\documentclass[12pt]{article}

\usepackage[margin=1in]{geometry}

\usepackage{amsthm,amsmath,amsfonts,amssymb}
\usepackage[authoryear,sort]{natbib}
\usepackage[colorlinks,citecolor=blue,urlcolor=blue]{hyperref}
\usepackage{graphicx}

\usepackage{tikz}
\usepackage{pgfplots}
\usepackage{pgfplots}
\pgfplotsset{compat=newest}
\usepgfplotslibrary{fillbetween}

\definecolor{pantoneRed}{RGB}{245,37,0}
\definecolor{pantoneBlue}{RGB}{0,84,166}
\definecolor{pantoneBlack}{RGB}{34,34,34}
\definecolor{pantoneLightGray}{RGB}{224,224,224}
\definecolor{highlightRed}{RGB}{255,0,0}
\definecolor{PantoneGreen}{RGB}{0,166,81}
\usetikzlibrary{patterns, arrows.meta,intersections,fillbetween,calc}

\usepackage{amsmath,bm}
\usepackage{mathptmx}
\usepackage{amsfonts}
\usepackage{amssymb}
\usepackage{graphicx}
\usepackage{rotating}
\usepackage{comment}

\theoremstyle{plain}

\newtheorem{theorem}{Theorem}[section]
\newtheorem{lemma}[theorem]{Lemma}

\theoremstyle{remark}

\newtheorem{example}{Example}

\newcommand{\yinchu}[1]{\textcolor{PantoneGreen}{Yinchu: #1}}
\newcommand{\jelena}[1]{\textcolor{pantoneRed}{Jelena: #1}}
\newcommand{\jelenax}[1]{\textcolor{black}{#1}}

\newtheorem{assumption}{Assumption}

\newtheorem{corollary}{Corollary}[theorem]

\newtheorem{remark}{Remark}[section]

\newcommand{\RR}[1]{\mathbb{R}}
\newcommand{\proj}{\mathrm{proj}}
\newcommand{\urltilde}{\kern -.15em\lower .7ex\hbox{~}\kern .04em}
\newcommand{\E}{\mathbb{E}}

\makeatletter
\def \@seccntformat#1{\csname the#1\endcsname.\quad}
\makeatother
\numberwithin{equation}{section}

\global\long\def\Mcal{\mathcal{M}}%

\global\long\def\RR{\mathbb{R}}%

\global\long\def\Acal{\mathcal{A}}%

\global\long\def\hphi{\hat{\phi}}%

\global\long\def\hpi{\hat{\pi}}%

\title{Minimax Semiparametric Learning With Approximate Sparsity}

\author{
Jelena Bradic\thanks{Department of Statistics and Data Science, Cornell University, \texttt{jelena.bradic@cornell.edu}} \and
Victor Chernozhukov\thanks{Department of Economics, MIT, \texttt{vchern@mit.edu}} \and
Whitney K. Newey\footnotemark[2]~\thanks{Department of Economics, MIT, \texttt{wnewey@mit.edu}} \and
Yinchu Zhu\thanks{Department of Economics, Brandeis University, \texttt{yinchuzhu@brandeis.edu}}
}

\date{July 31, 2025}  

\begin{document}
\maketitle

\begin{abstract}

Estimating linear, mean-square continuous functionals is a pivotal challenge in statistics. In high-dimensional contexts, this estimation is often performed under the assumption of exact model sparsity, meaning that only a small number of parameters are precisely non-zero. This excludes models where linear formulations only approximate the underlying data distribution, such as nonparametric regression methods that use basis expansion such as splines, kernel methods or polynomial regressions.
Many recent methods for root-$n$ estimation have been proposed, but the implications of exact model sparsity remain largely unexplored. In particular, minimax optimality for models that are not exactly sparse has not yet been developed.
This paper formalizes the concept of approximate sparsity through classical semi-parametric theory. We derive minimax rates under this formulation for a regression slope and an average derivative, finding these bounds to be substantially larger than those in low-dimensional, semi-parametric settings.
   We identify several new phenomena. We discover new regimes where rate double robustness does not hold, yet root-$n$ estimation is still possible. In these settings, we propose an estimator that achieves minimax optimal rates. Our findings further reveal distinct optimality boundaries for ordered versus unordered nonparametric regression estimation.

\end{abstract}

\section{Introduction}\label{sec intro}

\jelenax{High-dimensional statistical models have become increasingly pervasive in contemporary statistical research, driven by the natural availability of large-scale datasets across diverse application domains. Data collection now occurs on scales that were previously unimaginable, resulting in scenarios where the number of variables often exceeds the number of observations. This imbalance presents unique challenges for estimation and inference, necessitating the development of novel methodological techniques. The subject has been extensively studied, with significant advancements in regularization, sparsity, and ensemble methods \citep{hoerl1970ridge,tibshirani1996lasso,breiman1996bagging,freund1997boosting,fan2001scad}. } 

\jelenax{A prominent approach in high-dimensional settings is the assumption of underlying model's sparsity \citep{donoho2006sparsity}. }
However, in \jelenax{many} high-dimensional \jelenax{contexts}, strict \jelenax{model} sparsity might not be  realistic.  \jelenax{For instance, in genomics, natural language processing, finance, medical imaging, and climate science, the assumption that only a small number of parameters are non-zero often does not hold.} Allowing for many covariates or basis functions can potentially reduce approximation error of an unknown regression\jelenax{, as it provides the model with greater flexibility to capture complex relationships within the data.}
However, these covariates are not generally \jelenax{of equal importance}. \jelenax{For example, dummy variables created from categorical features or interaction terms may vary significantly in their relevance to the model.}  

\jelenax{Classical nonparametric methods, like splines and kernel regression, approximate non-linear data dependencies by selecting a specific cutoff for basis functions, assuming an ordered selection process where the initial basis functions are retained and the later ones are discarded; see for example,  \citep{wahba1990splines, silverman1984spline, schwarz2016spline, ruppert2003semiparametric}.}  \jelenax{ In principle, high-dimensional methods improve on this by allowing more basis functions and selecting an unordered subset that best predicts the outcome; rather than discarding higher-order terms outright, these methods might retain some of them as} the identity of the most important covariates  \jelenax{is not always apparent upfront.}  
 \jelenax{Open questions remain: Do high-dimensional approximately sparse models behave differently from low-dimensional, classical non-parametric methods? Is there a fundamental, perhaps irreconcilable, difference in the ability to estimate effectively in unordered versus ordered settings? Are minimax optimal rates the same in the two settings?  
 In this paper, we identify for the first time that these two settings are fundamentally different. } We find that allowing for unordered basis functions makes it harder to achieve the parametric rate for learning a linear functional. 
  
%

\subsection{Approximate sparsity}\label{sec11}

\jelenax{In the study of high-dimensional statistical models,  where the number of potential predictors \( p \) can be much larger than the sample size \( n \), a sparse regression model is one in which only a small subset of regressors have nonzero coefficients among a vast pool of potential regressors.}  

\jelenax{We focus here on approximately sparse models, a subject that naturally occurs as an extension of strict sparsity \citep{belloni2019valid, chernozhukov2022debiased, buhlmann2011statistics,koltchinskii2009sparse, zhang2012minimax, cai2021optimal}.} Approximate sparsity captures the idea that the regression is well approximated by a sparse specification without requiring that many true coefficients be exactly zero.
\jelenax{However, the existing formalizations are often inadequate, as they are limited to smaller than root-$n$ deviations from the model's exact sparsity or relaxations of $l_0$ to $l_q$, $q\leq 2$ constraints.} \jelenax{Our primary contribution is the establishment of a semi-parametric framework for approximate sparsity, which is grounded in approximation theory \citep{ward1927functions}}. 


\jelenax{Consider a nonparametric model $Y_i = f(X_i) + \varepsilon_i$ for $i=1,\dots,n$, where $\varepsilon_i$ are independent random variables, $\mathbb{E}(\varepsilon_i)=0$, and an outcome $Y_i\in \mathbb{R}$ with covariates $ X_i \in \mathbb{R}$ and a regression function $f: [0,1] \to \mathbb{R}$ such that $f(x)=\mathbb{E}[Y_i|X_i=x]$. The problem of nonparametric regression is to estimate $f$ given a priori that this function belongs to a nonparametric class of functions $\mathcal{F}$. For this discussion, we consider $\mathcal{F}$ to be $L_{2}[0,1]$, the space of square-integrable functions over the interval $[0,1]$, i.e., $E(f(X))^2<\infty$. 
For  a basis $\{\psi_j\}_{j=1}^\infty$ of $L_{2}[0,1]$, we will assume that $f$ can be expressed as an infinite series:
\[
f(x) = \sum_{j=1}^\infty \theta_j \psi_j(x)
\]
where the series converges for all $x \in [0,1]$}; see Sobolev classes in \cite{wasserman2006all} and \cite{Tsybakov2009}. \jelenax{A classical projection estimator of $f$ is based on a simple principle: approximate $f$ by its projection 
$
 \sum_{j=1}^p \theta_j \psi_j(x)
$
on the linear span of the first $p$ functions of the basis $\psi_1, \psi_2, \cdots,\psi_p $ and replace $\theta_j$ by their estimates.
A natural extension is to allow the order of the projection estimator, $p$, to be larger than the samples size $n$ as that would make the approximation error $\| f(\cdot) -  \sum_{j=1}^p \theta_j \psi_j(\cdot)\|_2$  smaller.  
Our work introduces two innovative sparse projection estimators. The {\it sparse ordered projection estimator} follows a traditional approach, using the first \(s\) terms to balance approximation and estimation errors. This method respects the natural order of approximation. In contrast, the {\it sparse non-ordered projection estimator} is better suited for high-dimensional settings, where \(p \gg n\). Instead of strictly following the order, it selects the best \(s\) terms from a larger set, potentially prioritizing higher-order terms when they provide a better balance between approximation and estimation.}

\jelenax{In this context, the concept of \textit{approximate sparsity} should emphasize  \textit{approximations} where the error scales with the sparsity level \( s \). } Our focus is particularly on sparse approximations characterized by a \textit{non-ordered} series. \jelenax{ With this perspective, we introduce a class of approximately sparse functions that disregard the natural order in selecting the optimal approximation, as follows:}
\begin{equation}
\mathcal{M}_{C,\xi}:=
\left\{
f: [0,1] \to \mathbb{R}:\ 
\min_{\Vert \theta \Vert_{0} \leq s} \ 
\bigl\|  f(\cdot)-\sum_{j=1}^{p}\theta_{j}\psi_{j}(\cdot)\bigl\|_{2}
\leq C s^{-\xi}, \forall s \in
\mathbb{N}\right\}, \label{eq:M}
\end{equation}
for some constant $C>0$. Here $s\jelenax{:=\|\theta\|_0 = \sum_{j=1}^p \mathbf{1} \{\theta_j \neq 0 \}}$ is the number of nonzero elements of $\theta$.   The decay rate index, denoted as \( \xi > 0 \), \jelenax{captures} the speed at which the approximation error decays to zero.  
In this context, functions \( f \) can be approximated in the mean square by sparse linear combinations, without \jelenax{restricting} the sparse elements \jelenax{to be within} the initial segments of the vector $\theta$.

\jelenax{Approximately} sparse specifications \jelenax{differ} from more \jelenax{conventional} nonparametric \jelenax{approaches} in ways that are \jelenax{advantageous} in high-dimensional settings. \jelenax{Approximate} sparsity \jelenax{permits} a \jelenax{very large} number of potential \jelenax{covariates}  (possibly \jelenax{exceeding} the sample size) when only a \jelenax{relatively small} subset of important \jelenax{ones} \jelenax{yields} a good approximation, but the \jelenax{identities} of these key \jelenax{features} are not known. In contrast, basis function \jelenax{approximations} are based on a \jelenax{comparatively limited} number of regressors, often \jelenax{significantly fewer} than the sample size. \jelenax{Both} approximately sparse and basis function \jelenax{approximations} are similar in that they both \jelenax{rely} on a few regressors to achieve a good approximation. They \jelenax{differ} in that basis function regression \jelenax{necessitates} knowledge of the \jelenax{specific} important \jelenax{features}, whereas approximate sparsity allows their identities to remain unknown. This \jelenax{distinction} is \jelenax{particularly beneficial} in high-dimensional settings, where there are potentially \jelenax{a vast} number of \jelenax{covariates} needed to approximate a function of many variables. Typically, there is \jelenax{minimal} guidance about which among the many regressors is important, such as interaction terms in a multivariate approximation. With approximate sparsity, such information is \jelenax{unnecessary}, since \jelenax{a large} number of terms can be included \jelenax{among} the potential \jelenax{features}.

\subsection{Main contributions} 
\jelenax{Our contributions fall into three areas: proposing a framework and establishing optimal estimation rates for approximately sparse models; presenting lower bound results indicating that double robustness may not be necessary for optimal inference; and achieving root‑$n$ estimation rates for linear functionals under some of the most general conditions.}

 We provide minimax lower bounds ftor high-dimensional linear models under approximate sparsity. 
\jelenax{The problem's complexity is governed by the regression model and its precision matrix, which mitigates the shrinkage bias introduced by regularization.
 }
We \jelenax{establish} that $\sqrt{n}$-rate for individual coefficients is impossible if $\max\{\xi_1,\xi_2\}  \leq 1/2$, where $\xi_1>0$ is the approximate sparsity index for the model and $\xi_2>0$ is that for the corresponding row in the precision matrix. \jelenax{To address this limitation, we introduce}  a novel estimator that \jelenax{attains} the $\sqrt{n}$-rate \jelenax{provided that } $\max\{\xi_1,\xi_2\} >1/2$\jelenax{, even under}  non-Gaussian \jelenax{covariates} and \jelenax{hteroscedastic} errors. \jelenax{Consequently}, the \jelenax{minimax optimal} rate for \jelenax{estimating
} individual parameters in approximately sparse models is $\sqrt{n}$ if and only if
\begin{equation}
\max\{\xi_{1},\xi_{2}\}>1/2. \label{box}%
\end{equation}

This is in stark contrast with the case in which the order of importance of the basis functions is known, \jelenax{i.e., the classical non-parametric setting}. \cite{donald1994series}  show that the $\sqrt{n}$-rate can be achieved \jelenax{therein} as long as 
\begin{equation}
\xi_{1}+\xi_{2}\geq1/2, \label{triangle}%
\end{equation}
which is \jelenax{inherently} weaker than $\max\{\xi_1,\xi_2\} >1/2$.  
\jelenax{Therefore, our results characterize the loss of having unordered covariates. We also derive \jelenax{an additional} result of independent interest: {\it for ordered basis functions, the condition of (\ref{triangle})  is sharp in that the $\sqrt{n}$-rate is impossible if $\xi_1+\xi_2 < 1/2$}. 
To the best of our knowledge, this is the first lower bound that  characterizes how  the approximation error rate affects the rate of learning a regression coefficient. Prior works such as \cite{robins2009semiparametric}  consider H\"older functions of low-dimensional covariates and do not consider the approximation error of basis functions.}


\jelenax{ Our results are also related to the estimation of other one-dimensional parameters. Given $(Y,X)$, many objects of interest can be written as linear functionals of $\mathbb{E}[Y\mid X]$, such as the average treatment effect, see Example \ref{example: ATE}. Given a $L_2$-continuous linear functional $L(\cdot)$, the Riesz representation theorem \citep{rudin1991functional} \jelenax{guarantees} the existence of a unique function $\alpha(\cdot)$ such that $L(f)=\mathbb{E}[\alpha(X)f(X)]$ for any $f$ with $\mathbb{E}[f(X)]^2<\infty$. We refer to this function $\alpha(\cdot)$ as the Riesz representer (of the functional under consideration). Now $\xi_1$ and $\xi_2$ describe the approximate sparsity of the conditional mean function $\mathbb{E}[Y\mid X=x]$ and the Riesz representer, respectively.   A popular condition for learning linear functionals at the $\sqrt{n}$-rate with high-dimensional models is the rate double robustness (RDR) condition of \cite{chernozhukov2018double},  \cite{smucler2019unifying}, \cite{balakrishnan2023fundamental}, and others. In low-dimensional settings, this condition has been referred to as the double robustness product-rate condition \citep{Seaman.2018, vansteelandt2022assumption}.
 }

 This condition requires that the product of the $\ell_2$-rates of estimating the two nuisance functions (one for the model and the other for the Riesz representer) be of the order $o_P(n^{-1/2})$. From our proof, one can see that the $\ell_2$-rate for estimation in 
$\mathcal{M}_{C,\xi}$ is $$(n^{-1}\log(p))^{\xi/(2\xi+1)}.$$ Thus, the \jelenax{RDR} condition translates into $(n^{-1}\log(p))^{\xi_1/(2\xi_1+1)+\xi_2/(2\xi_2+1)}\ll n^{-1/2}$, which implies $\xi_1/(2\xi_1+1)+\xi_2/(2\xi_2+1)>1/2$ or equivalently
\begin{equation}
\xi_1 \xi_2 >1/4. \label{hyperbola}%
\end{equation} 
\jelenax{We demonstrate that our estimator instead requires conditions weaker than RDR; see} Remark \ref{rem ATE}.  \jelenax{Specifically, it} achieves the $\sqrt{n}$-rate when either $\xi_1 \xi_2 > 1/4$ or $\xi_1 > 1/2$ (Areas $A \cup B$ in Figure \ref{fig:1}), \jelenax{which is less restrictive than the RDR condition} (Area $A$ in Figure \ref{fig:1}). 
\jelenax{Additionally, it permits an unbounded} Riesz representer, allowing practical applications such as propensity scores approaching zero or one in ATE estimation.

\begin{figure}
\caption{\label{fig:1}Minimal conditions for the root-n rate. This graph depicts three sets colored in blue, red and black.
The blue denotes points satisfying $ \xi_1+\xi_2=1/2$, the red denotes those satisfying $\max\{\xi_1,\xi_2\}=1/2$ and the black denotes those satisfying $\xi_1\xi_2=1/4$.  }
\centering
 \begin{tikzpicture}
    \begin{axis}[
        axis lines=middle,
        xlabel={$ \xi_1$},
        ylabel={$ \xi_2$},
        xmin=0, xmax=1,
        ymin=0, ymax=1,
        xtick={0,0.1,0.2,0.3,0.4,0.5,0.6,0.7,0.8,0.9,1},
        ytick={0,0.1,0.2,0.3,0.4,0.5,0.6,0.7,0.8,0.9,1},
        xlabel style={at={(axis description cs:0.5,-0.08)},anchor=north},
        ylabel style={at={(axis description cs:-0.12,0.5)},anchor=south},
        width=9cm, height=9cm,
        axis equal image
    ]
        \addplot [
            domain=0.25:1, 
            samples=100, 
            fill=pantoneLightGray,
            draw=none
        ] {0.25 / x} |- (1,1) -| (0.25, 4) -- cycle;

        \foreach \y in {0.505,0.51, 0.515, 0.52, 0.525, 0.53, 0.535, 0.54, 0.545, 0.55, 0.555, 0.56, 0.565, 0.57, 0.575, 0.58, 0.585, 0.59, 0.595, 0.6, 0.605, 0.61, 0.615, 0.62, 0.625, 0.63, 0.635, 0.64, 0.645, 0.65, 0.655, 0.66, 0.665, 0.67, 0.675, 0.68, 0.685, 0.69, 0.695, 0.7, 0.705, 0.71, 0.715, 0.72, 0.725, 0.73, 0.735, 0.74, 0.745, 0.75, 0.755, 0.76, 0.765, 0.77, 0.775, 0.78, 0.785, 0.79, 0.795, 0.8, 0.805, 0.81, 0.815, 0.82, 0.825, 0.83, 0.835, 0.84, 0.845, 0.85, 0.855, 0.86, 0.865, 0.87, 0.875, 0.88, 0.885, 0.89, 0.895, 0.9, 0.905, 0.91, 0.915, 0.92, 0.925, 0.93, 0.935, 0.94, 0.945, 0.95, 0.955, 0.96, 0.965, 0.97, 0.975, 0.98, 0.985, 0.99, 0.995, 1}
 {
            \addplot [
                ultra thick,  
                 fill=pantoneRed,
                color=pantoneRed!50
            ] coordinates {(0, \y) ({0.25/\y}, \y)};
        }

        \foreach \y in {0.005, 0.01, 0.015, 0.02, 0.025, 0.03, 0.035, 0.04, 0.045, 0.05, 0.055, 0.06, 0.065, 0.07, 0.075, 0.08, 0.085, 0.09, 0.095, 0.1, 0.105, 0.11, 0.115, 0.12, 0.125, 0.13, 0.135, 0.14, 0.145, 0.15, 0.155, 0.16, 0.165, 0.17, 0.175, 0.18, 0.185, 0.19, 0.195, 0.2, 0.205, 0.21, 0.215, 0.22, 0.225, 0.23, 0.235, 0.24, 0.245, 0.25, 0.255, 0.26, 0.265, 0.27, 0.275, 0.28, 0.285, 0.29, 0.295, 0.3, 0.305, 0.31, 0.315, 0.32, 0.325, 0.33, 0.335, 0.34, 0.345, 0.35, 0.355, 0.36, 0.365, 0.37, 0.375, 0.38, 0.385, 0.39, 0.395, 0.4, 0.405, 0.41, 0.415, 0.42, 0.425, 0.43, 0.435, 0.44, 0.445, 0.45, 0.455, 0.46, 0.465, 0.47, 0.475, 0.48, 0.485, 0.49, 0.495, 0.5}
{
            \addplot [
                ultra thick, 
                color=pantoneRed!50
            ] coordinates {(0.5, \y) ({0.25/\y}, \y)};
        }

        \addplot [
            domain=0:0.5,
            samples=100,
            fill=pantoneBlue!50, 
            draw=none
        ] coordinates {(0, 0.5) (0.5, 0.5) (0.5, 0) (0, 0.5)};
        
        \addplot [
            thick, 
            color=pantoneBlack, 
            domain=0.25:1, 
            samples=100
        ] {0.25 / x};
        
        \addplot[thick, color=pantoneRed!50] coordinates {(0, 0.5) (0.5, 0.5)};
        \addplot[thick, color=pantoneRed!50] coordinates {(0.5, 0.5) (0.5, 0)};
        
        \addplot[thick, color=pantoneBlue] coordinates {(0, 0.5) (0.5, 0)};
        
        \node at (axis cs:0.72,0.72) {\Huge \textbf{A}};
        \node at (axis cs:0.72,0.18) {\Huge \textbf{B}};
        \node at (axis cs:0.15,0.72) {\Huge \textbf{C}};
        \node at (axis cs:0.34,0.34) {\Huge \textbf{D}};
    \end{axis}
\end{tikzpicture}
 \end{figure}

\jelenax{To achieve \(\sqrt{n}\)-rate learning for individual coefficients, we have established the following minimal set of conditions, each representing a unique geometric region within Figure \ref{fig:1}:
}\begin{itemize}
    \item Ordered \jelenax{basis}: \jelenax{Represented by area} $A\bigcup B\bigcup C\bigcup D$ in Figure \ref{fig:1}, \jelenax{corresponding to } Eq (\ref{triangle}).
    \item Non-ordered \jelenax{basis} (approximate sparsity):  \jelenax{Represented by area}
$A\bigcup B\bigcup C$ in Figure \ref{fig:1}, \jelenax{corresponding to} Eq. (\ref{box}).
\end{itemize}
 
Hence, the RDR condition, \jelenax{widely adopted in sparse high-dimensional models (e.g., Assumption 11 in \cite{chernozhukov2022automatic})}-- is \jelenax{shown here} NOT \jelenax{to be}  necessary for \jelenax{achieving}  the parametric rate in approximately sparse  \jelenax{{\it non-ordered} high-dimensional settings}. \jelenax{This finding challenges a key assumption in the existing literature.}   
\jelenax{ Recently, \cite{balakrishnan2023fundamental} pointed out that the RDR condition is necessary for the parametric rate if estimation is done in a black-box manner. What we see here is that some structure can lead to big improvements over the RDR condition: {\it one can exploit approximate sparsity and achieve the parametric rate without RDR}. If the order of the dictionary is known, we can do even better; the cost of not knowing the dictionary order is that we lose the ability to perform root-n order inference in Area~D, Figure \ref{fig:1}.  \cite{Bonvini2405.08525} study a hybrid functional class by imposing  smoothness assumptions in addition to the rate conditions of black-box estimators. A key difference between their work and ours is that their structural assumption (i.e., smoothness condition) is conditional on the estimator, e.g., Proposition 3 therein; in contrast, our approximate sparsity class is not conditional on any estimator.  }

\subsection{Literature review}

\jelenax{For learning individual coefficients in linear models, \jelenax{numerous} prior works \jelenax{exist, including} \cite{belloni2014inference, zhang2014confidence, van2014asymptotically, javanmard2014confidence, javanmard2018debiasing, cai2017confidence, zhu2018linear, bradic2022testability, bellec2022debiasing}.
  Some of them investigates the minimal condition for the $\sqrt{n}$-rate under strictly sparse models or Riesz representers. \cite{cai2017confidence} show that if no sparsity is imposed on the Riesz representer (essentially the first column of the inverse of the covariance matrix), then the model parameter needs to be $s$-sparse with $s\ll \sqrt{n}/\log p$. \cite{bradic2022testability} show that if no sparsity is imposed on the model parameter, then the Riesz representer needs to be $s$-sparse with 
$s\ll \sqrt{n}/\log p$. \cite{javanmard2018debiasing} consider a balance between the sparsity condition of the model parameter and that of the Riesz representer, see Theorem 3.13 and Proposition 4.2 therein, but they need to assume bounded $\ell_1$-norm of the Riesz representer, see the quantity $\rho$ in their paper. This restriction is relaxed by \cite{bellec2022debiasing}, who also proposed an estimator based on degrees-of-freedom adjustments. All these papers work with strict sparsity and assume Gaussian designs and errors, whereas we consider approximate sparsity and allow for non-Gaussian covariates and errors and heteroscedasticity.   \jelenax{Lemma \ref{lem: strict vs approx sparse} confirms our results are not derived from exactly sparse methods.}
}

\jelenax{Our work is also related to the literature on learning other functionals of high-dimensional models, such as average derivative, ATE, instrumental variable regression and models defined by moment conditions; see \cite{farrell2015robust, chernozhukov2018double, chernozhukov2022debiased, chernozhukov2018learning, athey2018approximate, smucler2019unifying, bradic2019sparsity, tan2020regularized, ning2020robust} for related works. We contribute to this literature by proposing a new estimator for estimating a generic linear functional under approximate sparsity. 
 }

\jelenax{Some works in the literature consider models that are not strictly sparse, e.g, \cite{belloni2014inference} and Chapter 6 of \cite{buhlmann2011statistics}.
However, \jelenax{the optimality of their} assumptions on  approximation errors \jelenax{is unclear}. In fact, one can show that the assumptions in \cite{belloni2014inference} imply \jelenax{RDR condition}, which is not needed by our estimator; \jelenax{see Section \ref{sec22} for more details on exact vs. approximate sparsity. The Auto DML approach of \cite{chernozhukov2022automatic} also imposes conditions that imply  the RDR condition. The RDR condition is also needed by \cite{chernozhukov2022debiased}, who have weaker sparsity conditions under the assumption of bounded $\ell_1$-norm.  }
}

\jelenax{In this paper, we treat approximate sparsity as a function class and provide a complete theory (i.e., minimax lower and upper bounds) for learning individual coefficients and an upper bound for learning generic linear functionals. We show that the RDR condition is not needed to achieve the parametric rate even for generic linear functionals. We also do not impose any restrictions of the parameters in terms of their $\ell_1$-norm.}
\jelenax{Some papers, such as \cite{ye2010rate} and \cite{raskutti2011minimax}, use \(\ell_q\)-balls with \(q \in (0,1]\) as a way to approximate sparsity. However, their framework differs: while we approximate the true model with a sparse one up to some error, they assume exact linear models (with no approximation error) but allow these model parameters to deviate from strict sparsity.}
 
\jelenax{The proposed estimators allow for heavy-tailed error terms even for learning generic linear functionals. This is in sharp contrast with the majority of the literature on high-dimensional models, which typically relies on bounded or sub-Gaussian errors. We allow for unbounded Riesz representer, which means that the proposed estimators are valid under very weak assumptions. For example, in learning the average treatment effect, we only require the inverse propensity score to have bounded moments, rather than to be bounded by a constant.  }



\subsection{Notations and organization}

For any function $f$, we denote $\left\Vert f\right\Vert _{2}=\sqrt{\mathbb{E}[f(X)^{2}]}$. For any vector $v$,
 let $\Vert v\Vert_{0}$ denote the number of nonzero elements of $v$. For any vector $v=(v_1,...,v_p)^{\top}$ and $q \in [1,\infty)$, let $\|v\|_q=\left(\sum_{j=1}^{p} |v_j|^q \right)^{1/q}$ and $\|v\|_\infty=\max_{1\leq j\leq p} |v_j|$ . For $J \subset \{1,...,p\}$ and $v=(v_1,...,v_p)^{\top} \in\mathbb{R}^p$, $\left\vert J\right\vert $ denote
the number of elements of $J,$ $v_{J}$ be the vector in $\mathbb{R}^p$ consisting of
$(v_{J})_j=v_{j}$ for $j\in J$ and $(v_{J})_j=0$ otherwise, and
$v_{J^{c}}$ be the corresponding vector for $J^{c}$ (so that
$v=v_{J}+v_{J^{c}})$. We use $e_j$ to denote the $j$th column of the $p\times p$ identity matrix. A random variable $X$ is said to have a bounded sub-Gaussian norm if there exists a constant $c>0$ such that $\E \exp(tX)\leq \exp(ct^2)$ for any $t\in \RR$. A random vector $X$ is said to have a bounded sub-Gaussian norm if  there exists a constant $c>0$ such that $\E \exp(tv^{\top}X)\leq \exp(ct^2)$ for any $t\in \RR$ and any $v$ with $\|v\|_{2}=1$. We denote by $\mathbb{Z}$, $\mathbb{N}$, $\mathbb{C}$ and $\mathbb{R}$, the set of integers, the set of natural numbers, the set of complex numbers and the set of real numbers, respectively.

The rest of this paper is organized as follows. Section \ref{sec approx sparsity} discusses more details on the approximate sparsity. We then consider the problem of learning individual coefficients in high-dimensional linear models under approximate sparsity. Section \ref{sec lower bounds}  provides lower bounds and Section \ref{sec PLM upper bound} establishes upper bounds. Section \ref{sec general functionals} develops an estimator for generic linear functionals under approximate sparsity.

\section{Approximate Sparsity}\label{sec approx sparsity}  
 
 We note that \jelenax{the representation of }$\mathcal{M}_{C,\xi}$
 in \eqref{eq:M} \jelenax{implicitly depends on the choice} of the basis 
$b$.
 \jelenax{For a g}iven  basis $\{\psi_{j}(x)\}_{j=1}^{\infty}$, the set of functions that can be approximated \jelenax{is defined as $\mathcal{B}$.}
 \jelenax{This is the closure, with respect to the $\|\cdot\|_2$-norm, of the union of finite-dimensional subspaces:}
\[
\mathcal{B} = \operatorname{closure}\left(\bigcup_{p=1}^\infty \mathcal{B}_p\right), \quad \text{where} \quad \mathcal{B}_p = \left\{\sum_{j=1}^p a_j \psi_j(\cdot) : a_j \in \mathbb{R}\right\}.
\]
Clearly, $\mathcal{M}_{C,\xi} \subset \mathcal{B}$.
  Notice that $\mathcal{B}$ depends on the basis.  If  $\{\psi_{j}(x)\}_{j=1}^{\infty}$ does not span $L_2(X)$, then $\mathcal{B}\neq L_2(X)$. For example, one can often choose a basis that spans $L_2(X)$ and have $\mathcal{B}=L_2(X)$.

\subsection{{Comparison with classical semiparametrics}}
We can be precise about \jelenax{the} key difference between approximate sparsity and
basis approximations by defining the latter. \jelenax{For $\mathcal{M}_{C,\xi}$
the nonzero components of $a$ are allowed to be coefficients of any of the $p$
dictionary functions, where $p$ can be large, even larger than sample size$.$
In the definition of $\mathcal{M}_{C,\xi}$ it is unknown which dictionary
functions are used in the sparse approximate at rate $Ct^{-\xi}$. A series
approximation from the semiparametric literature would require that the
unknown function be well approximated by a linear combination of the first $t$
functions. The set of unknown $v$ allowed here would be
\[
\mathcal{S}_{C,\xi}=\left\{  f(\cdot):\ \min_{a=(a_{1}%
,...,a_{t},0,...,0)^{\top}}\Vert f(\cdot)-\psi(\cdot)^{\top}a\Vert_{2}\leq Ct^{-\xi
}\ \forall t\in\mathbb{N}\right\}
\]}

{For example, suppose that $X$ is continuously distributed with compact support
and that dictionary functions are products of all nonnegative powers of $x$
that are weakly increasing in order with $j$. If $x$ has dimension $d$ and $v$
is such that $\psi(x)^{\top}v$ has bounded derivatives of order $s$ then it is
well known that there is $C$ and an ordering of $\psi(x)$ such that the
inequality in the definition of $\mathcal{S}_{C,\xi}$ is satisfied with
$\xi=s/d.$ This $\xi$ is the well known rate for approximation of functions in
a H\"older class. }


{Notice that if $\xi\geq\tilde{\xi}$, then $\mathcal{M}_{C,\xi
}\subseteq\mathcal{M}_{C,\tilde{\xi}}$, similarly to $\mathcal{S}_{C,\xi}$
shrinking with $\xi$.}

\subsection{Comparison with strict sparsity} \label{sec22}
\jelenax{We note that the approximate sparsity is a generalization of ordered approximation, i.e., $\mathcal{S}_{C,\xi} \subset \mathcal{M}_{C,\xi}$. Moreover, approximate sparsity does not extend strict sparsity:   exact sparsity is not a strict subset of a set of approximately sparse models in that they are not nested sets.  Consider the set $\Theta(k)=\{v^\top \psi(\cdot) \in\mathbb{R}^p:\ \|v\|_0\leq k \} $ as in   \cite{cai2017confidence} and \cite{javanmard2018debiasing}. Clearly, since $\Theta(k)$ does not restrict the $\ell_2$-norm of the elements, it is impossible to say that $\Theta(k)\subset  \mathcal{M}_{C,\xi}$. However, even if we impose such a restriction, we still cannot fit it as a subset of  $\mathcal{M}_{C,\xi}$ due to the following result.}

\jelenax{\begin{lemma}\label{lem: strict vs approx sparse}
Suppose that $\mathbb{E}[\psi(X)\psi(X)^{\top}]=I_p$. Define the set \(\Theta(C, k)\) of sparse linear combinations of the basis functions with bounded $L_2$ norm:
$$
\Theta(C, k) = \Biggl\{ \theta^\top \psi(\cdot) \mid  \theta \in \mathbb{R}^p, \|\theta\|_0 \leq k, \|\theta\|_2 \leq C \Biggl\},
$$ for some constant $C>0$. \jelenax{As} the sparsity $k$ and the dimension $p$ grow to  $\infty$  there \jelenax{do}  not exist constants $D$ and $\xi>0$ such that 
$$\Theta(C,k) \subset \mathcal{M}_{D,\xi}.$$
\end{lemma}
 }

\jelenax{ Lemma \ref{lem: strict vs approx sparse} \jelenax{reveals that the}  approximate sparsity does not automatically generalize strict sparsity.   As discussed above, $  \mathcal{M}_{C,\xi}$ \jelenax{serves as} a natural extension of the ordered approximation, which is \jelenax{pivotal} in nonparametric model\jelenax{ing}. }

\subsection{\jelenax{Examples}}


\jelenax{We start with classical function classes and their approximations and develop examples  with  unordered basis functions.
Consider a convex $\Omega \subset [0,1]^D$ for $D \in \mathbb{N}$ (or a convex domain that can be mapped to $[1/2,3/4]^D$). 
For an integer $s \ge 0$, 
the Sobolev space $H^s(\Omega)$ is defined by
\[
H^s(\Omega) \;=\; \bigl\{ f : \Omega \to \mathbb{R} 
:\quad \|f\|_{H^s(\Omega)} < \infty \bigr\},
\]
where \jelenax{the norm $\|\cdot\|_{H^s(\Omega)}$ typically has the form}
\[
\|f\|_{H^s(\Omega)}^2 
\;=\; \sum_{|\alpha|\le s} \| D^\alpha f \|_{L^2(\Omega)}^2, \mbox{ with } 
\|g\|_{L^2(\Omega)}^2 
\;=\; \int_{\Omega} \bigl(g(x)\bigr)^2 \,dx.
\]
\jelenax{Above, $|\alpha| = \alpha_1 + \cdots + \alpha_D$ for a multi-index 
$\alpha = (\alpha_1, \ldots, \alpha_D)$, and $D^\alpha$ represents the 
weak partial derivatives of $f$ of order $\alpha$.
One may also define $H^s(\Omega)$ for non-integer $s$ (fractional Sobolev 
spaces); see \cite{adams2003sobolev} for more details.}
\jelenax{We can characterize membership in Sobolev 
spaces $H^s$ via the decay properties of expansion coefficients in an 
orthonormal eigenbasis.}
Consider a basis $\{\psi_{j}(\cdot)\}_{j=1}^{\infty}$
 of $L^{2}([0,1])$. For constants $Q,s>0$,
define  the \emph{Sobolev ellipsoid}
\begin{equation}\label{eq: sobolev ellipsoid}
\mathcal{F}_{s,Q}=\left\{ f(\cdot)=\sum_{j=1}^{\infty}\theta_{j}\psi_{j}(\cdot):\ \sum_{j=1}^{\infty}\theta_{j}^{2}j^{2s}\leq Q\right\} .
\end{equation}}

\jelenax{The factor $j^{2s}$ reflects how higher frequencies (indexed by $j$) 
are penalized to ensure smoothness. In many eigenfunction expansions 
related to the Laplacian or other differential operators, the eigenvalues 
grow on the order of $j^2$, so $j^{2s}$ is naturally tied to $s$ 
derivatives in a Sobolev sense.}

\jelenax{For $R>0$, consider the set
\[
B_{H^s}(R) 
\;=\;
\{\,f \in H^s(\Omega):\quad
\|f\|_{H^s(\Omega)} \le R
\}.
\]
 }


\jelenax{Below we show  that our condition for {\it approximate sparsity} is satisfied for the above Sobolev ellipsoids.}
\jelenax{We define for each multi-index 
$
  k \;=\; (k_1,k_2,\ldots,k_D)^\top \in \mathbb{Z}^D
$
the complex exponential function
\[
  \psi_k(x) 
  \;=\; \exp\bigl(2\pi \iota\, k^\top x\bigr),
  \quad 
  x=(x_1,\dots,x_D)\in \Omega,
\]
where $k^\top x = k_1x_1 + \cdots + k_D x_D$ is the usual dot product in 
$\mathbb{R}^D$, and $\iota = \sqrt{-1}$.
For $k, \ell \in \mathbb{Z}^D$,
\[
  \int_{\Omega} \psi_k(x)\,\overline{\psi_\ell(x)} \,dx
  \;=\;
  \int_{\Omega} 
    \exp\bigl(2\pi \iota\, (k-\ell)^\top \,x\bigr)\,dx
  \;=\;
  \begin{cases}
    1, & \text{if } k=\ell,\\
    0, & \text{otherwise}.
  \end{cases}
\]
For computational or theoretical reasons, one might consider a finite subset 
of these exponentials. Let $t \ge 1$ be a truncation parameter, and define
\[
  \mathcal{N}_t 
  \;=\;
  \bigl\{\, j \in \mathbb{Z} : |j| \le \lfloor t^{1/D}\rfloor \bigr\}, 
  \quad 
  \mathcal{N}_t^D 
  \;=\;
  \underbrace{\mathcal{N}_t \times \cdots \times \mathcal{N}_t}_{D \text{ times}}.
\]
Then we consider the set
\[
  \bigl\{\psi_k : k\in \mathcal{N}_t^D\bigr\}
\]
as a \emph{finite} collection of frequencies. (In fact, it has approximately $t$ elements.) This truncated basis approximates 
the full Fourier basis by including only those frequencies $\|k\|\lesssim t^{1/D}$.
The Sobolev space $H^s(\Omega)$ for an integer $s\ge 0$ can 
be characterized by
$
  \|f\|_{H^s([0,1]^D)}^2
  \;=\;
  \sum_{k\in \mathbb{Z}^D} 
    (1+\|k\|^2)^{s} \,\lvert \theta_k\rvert^2.
$
We show below that such collection approximates the unknown function $f$ in the sence of approximate sparsity as defined in Section \ref{sec11} and Eq.~\eqref{eq:M}. (The approximating function below will be real-valued but we use the notation in complex numbers to simplify the presentation.) }

\jelenax{\begin{lemma}
\label{lem: Sobolev approx}
For any $t\geq1$, consider the corresponding subset of basis functions $ \bigl\{\psi_k : k\in \mathcal{N}_t^D\bigr\}$.
Then, there exists a constant $C>0$, depending
only on $D,s,\Omega$, such that  
\begin{equation}\label{eq: basic ordered class}
\inf_{a_{k}\in\mathbb{C}} \Bigl\|  f-\sum_{k\in\mathcal{N}_{t}^{D}}a_{k}\psi_{k}\Bigl\|_{L^{2}(\Omega)}\leq CQ \ t^{-s/D}, \qquad \forall f\in B_{H^{s}(\Omega)}(Q). 
\end{equation}
\end{lemma}
\jelenax{This lemma bridges the gap between classical smoothness (measured by Sobolev norms) and our proposed sparse approximation theory by quantifying how smoothness translates into sparse representability.}
As a consequence of the above Lemma, we have
\[
B_{H^{s}(\Omega)}(Q)\subset \mathcal{M}_{C,\xi}, \qquad \mbox{with} \qquad \xi=s/D,
\]
}
\jelenax{which means that every function in the Sobolev ball \( B_{H^{s}(\Omega)}(Q) \) belongs to the model class \( \mathcal{M}_{C,\xi} \). 
  In other words, each function in the Sobolev ball is "approximately sparse" with respect to the given basis. This result is significant because it quantitatively links the smoothness of a function, measured in the Sobolev sense, with its sparse approximability using the index $\xi$.}

\jelenax{The above approximate sparsity is based on ordered basis functions, $\mathcal{N}_t^D$. Working with reproducing kernel Hilbert spaces (RKHS) is perhaps a more direct way of constructing a class of functions that can be approximated by ordered basis functions.
Classic Sobolev ellipsoids in nonparametric problems
would satisfy (\ref{eq: basic ordered class}),  
see e.g., \cite{Tsybakov2009} and many RKHS's are Sobolev ellipsoids, see e.g.,  \cite{zhang2022sieve}
and references therein. Consider a kernel function $\mathcal{K}(\cdot,\cdot)$
on $\Omega \times \Omega$ that admits the Mercer representation automatically
induces an ordered dictionary: 
\[
\mathcal{K}(s,t)=\sum_{j=1}^{\infty}\lambda_{j}\psi_{j}(s)\psi_{j}(t),
\]
where the eigenvalues $\{\lambda_{j}\}_{j=1}^{\infty}$ form a non-negative
decreasing sequence and $\{\psi_{j}\}_{j=1}^{\infty}$ form an orthonormal
basis. Then standard results (e.g., Corollary 12.26 in \cite{wainwright2019high})
imply that the RKHS induced by $\mathcal{K}$ consists functions of
the form $f=\sum_{j=1}^{\infty}a_{j}\psi_{j}$ such that $\sum_{j=1}^{\infty}a_{j}^{2}/\lambda_{j}<\infty$.
A ball of radius $C$ (in the Hilbert space norm) would be 
\[
\left\{ \sum_{j=1}^{\infty}a_{j}\psi_{j}:\sum_{j=1}^{\infty}\frac{a_{j}^{2}}{\lambda_{j}}\leq C\right\} .
\]
From a simple calculation, we can see that for eigenvalues decaying
at the rate $\lambda_{j}\lesssim j^{-2\xi}$, functions in the above
ball have the approximation error satisfying the rate in (\ref{eq: basic ordered class}). 
The Sobolev space on a bounded convex set $\Omega\subset \RR^d$ can be constructed this way. Consider the kernel $\mathcal{K}(s,t)=\Phi(s-t)$, where the Fourier transform of  $\Phi(x)$ is $(1+\|\cdot\|_2^2)^{-s}$ with $s>d/2$, see Theorem 6.13 of  \cite{wendland2004scattered} for the closed-end expression of $\Phi(\cdot)$. Then the corresponding RKHS is the Sobolev space $H^{s}(\Omega)$, see e.g., Corollary 10.48 of \cite{wendland2004scattered}.
}
\jelenax{
We now construct examples of approximate sparsity with unordered basis functions. }

 \subsubsection{Dimension-restricted Sobolev class}

\jelenax{Consider $x=(x_{1},...,x_{D})^\top \in I^{D}$, where $I=(1/4,3/4)$. Consider
the following class: 
\[
\mathcal{F}_{1}=\left\{ f(x)=f_{0}(x_{l_{1}},...,x_{l_{d}}):\|f_{0}\|_{H^{s}(I^{d})}\leq Q\right\}, 
\]
where $f_0: I^{d} \to \mathbb{R}$ is a lower-dimensional function $d \leq D$. Here, $f$ truly only depends on $d$ coordinates, but we do know know their order, i.e. not the first d coordinates. Since $f\in \mathcal{F}_{1}$ only depends on $d$ coordinates, Lemma \ref{lem: Sobolev approx} implies that we only need $t$ basis functions to get an approximation error of $O(t^{-s/d})$, which means that without knowing the $d$ ``true'' coordinates, it suffices to include  $t\cdot {D\choose d}$ \emph{unordered} basis functions. Thus, $\mathcal{F}_{1}\subset \mathcal{M}_{C,\xi}$ with $\xi=s/d$. Notice that $\xi$ depends only on $d$, the number of truly relevant variables, rather than $D$, the total number of variables under consideration. }

\jelenax{We notice that $D$ is allowed to go to infinity slowly here, in which case the it is impossible to solve the problem under the usual nonparametric methods due to the curse of dimensionality. We can certainly view $\mathcal{F}_{1}$ as a subset of $B_{H^{s}(I^D)}(Q)$. If $D\rightarrow \infty$, we can see that using an ordered dictionary (such as the Fourier basis) is infeasible. Even if $D$ is finite, it is still beneficial to  view $\mathcal{F}_{1}$ as an approximately sparse class by using an unordered dictionary. Suppose that we are interested in the regression coefficient in a partial linear model. Consider the setting of  (\ref{eq: PLM}) with $f$ and $g$  in $\mathcal{F}_{1}$ with $s=s_1$ and $s=s_2$, respectively. Then the lower bound in Theorem \ref{thm2} implies that the parametric rate is impossible if $s_1+s_2 <D/2$; in contrast, when we treat $f$ and $g$ as approximately sparse functions, Theorem \ref{thm: PLM main} shows that we only need $\max\{s_1,s_2\}>d/2$ to attain the parametric rate. Since $D$ can be much larger than $d$, the requirement of $\max\{s_1,s_2\}>d/2$ is weaker than $s_1+s_2 \geq D/2$ whenever $D>2d$. }

\subsubsection{Sliced Sobolev class}
\jelenax{Previous example can be seen as one-block of coordinates, i.e. one slice,  contributing to the function class. We can enlarge that by considering multiple blocks simultaneously, some of which can overlap. Hence, we consider 
\[
\mathcal{F}_{2}=\left\{ f(x)=\sum_{k\in\{1,...,D\}^{d}}c_{k}f_{k}(x_{k}):\|f_{k}\|_{H^{s}(I^{d})}\leq Q,c_{k}\in\mathbb{R}\right\} .
\]
In particular, if \(k = (k_{1},\ldots,k_{d})\), then \(x_{k}\) denotes the subvector of \(x\) indexed by \(k\), that is, \(x_{k} = \bigl(x_{k_{1}}, \ldots, x_{k_{d}}\bigr)\).
If each block \(x_k\) truly involves distinct coordinates (i.e., no coordinate overlaps between different \(k\)), then each block is a separate ``slice'' of dimension \(d\).  Functions in \(\mathcal{F}_2\) can be viewed as a sum of partial ``slices'' of dimension \(d\).  
Because \(k \in \{1,\dots,D\}^d\) can repeat coordinates or define overlapping subsets of \(\{1,\dots,D\}\), each \(f_k\) could share variables with another \(f_{k'}\).  The set \(\mathcal{F}_2\) might arise in a scenario where you suspect your true function is a linear superposition of local \(d\)-dimensional ``sub-features.''  
If $d=1$ we get back a classical additive model.
By Lemma
\ref{lem: Sobolev approx}, both of these classes satisfy the condition
that the approximation error decays with 
\[
{D \choose d}CQt^{-s/d}.
\]
Notice that this avoids the curse of dimensionality because the rate
essentially depends on $d$, not $D$. }

\subsubsection{Semi-parametric Sobolev class}
  \jelenax{Consider the class 
\[
\mathcal{F}_{3}=\left\{ f(x)=\sum_{j=1}^{d_{z}}c_{j}f_{j}(w)\cdot z_{j}:\|f_{j}\|_{H^{s}(I^{d_{w}})}\leq Q,\|c\|_{0}\leq q\right\} ,
\]
where $q$ is a bounded constant and $w\in\mathbb{R}^{d_{w}}$
and $z\in\mathbb{R}^{d_{z}}$.
Here,  the function \(f\) is \emph{linear} in the coordinates of \(z\), but each ``coefficient'' is replaced by a \emph{nonparametric} function \(f_j(w)\) lying in the Sobolev space \(H^s(I^{d_w})\). 
Because each \(z_j\) is multiplied by a function \(f_j(w)\), the effect of \(z_j\) on the output depends on \(w\), creating an \emph{interaction} between \((w,z)\). 
This setup also represents a \emph{varying-coefficient} perspective, since \(z_j\)'s contribution varies with \(w\). 
Additionally, the \(\ell^0\) constraint \(\|c\|_0 \le q\) enforces \emph{sparsity}, allowing at most \(q\) active interaction terms. 
Overall, \(\mathcal{F}_3\) can be viewed as a \emph{semi-parametric} class, combining a parametric (linear) form in \(z\) with a nonparametric Sobolev constraint in \(w\), ensuring \(s\)-smoothness in each coefficient function \(f_j\). 
This class satisfies the condition
that the approximation error decays with 
\[
{d_{z} \choose q}CQt^{-s/d_{w}}.
\]}


\section{Lower Bounds for Learning Individual Coefficients}\label{sec lower bounds}


In this section we give lower bounds for the convergence rates of estimators
of a slope coefficient in partial linear models.  We clarify the distinction between approximately sparse
models and those where identity of the important regressors is known and
explain the key conditions on which our results are based. We give lower
bounds under approximate sparsity and when the identity of the important
regressors is known. 


 Consider a partial linear model:
\begin{equation}\label{eq: PLM 3.1}
Y_i =D_i \theta +f(\tilde{Z}_i)+\varepsilon_i ,\qquad D_i =g(\tilde{Z}_i)+u_i,
\end{equation}
where $\E(\varepsilon_i\mid \tilde{Z}_i,D_i)=0$ and $\E (D_i \mid \tilde{Z}_i)=0$. Let
$\sigma_{u}^{2}=Eu_{i}^{2}$ and $\sigma_{\varepsilon}^{2}=\mathbb{E}\varepsilon_{i}^{2}$. In Sections \ref{sec:unordered} and \ref{sec:ordered}, we will derive lower bounds for learning $\theta$ based on different assumptions on the function class to which $f$ and $g$ belong.
Let $Z=\psi(\tilde{Z})\in\RR^p$. We observe a random sample of $W_1,W_2,\cdots, W_n \in \mathbb{R}^{p+2}$ with $W_i=(Y_i,X_i^\top)^\top$ and $X_i=(D_i,Z_i^\top)^\top \in  \mathbb{R}^{p+1}$. To derive the lower bounds, we consider the following special case of Gaussian data.
\begin{assumption} \label{assp 1}
Consider the model in (\ref{eq: PLM 3.1}). The random vector  $W_i=(Y_i,X_i^\top)^\top$ with $X_i=(D_i,Z_i^\top)^\top$ follows a
Gaussian distribution with mean zero and $\mathbb{E}[Z_{i}Z_{i}^{\top}]=I_{p}$. Furthermore, $f(\tilde{Z}_i)=Z_i^\top \beta$ and $g(\tilde{Z}_i)=Z_i^\top \phi$. 
\end{assumption}
The lower bound derived under this condition is a minimax lower bound for any model that nests this setting; likewise, any convergence‑rate lower bound for one model bounds every larger class containing it.

\subsection{Lower Bound Under Non-ordered Approximation} \label{sec:unordered}

Under the model in Assumption \ref{assp 1}, we first consider parameters that are approximately sparse, i.e., non-ordered basis functions.
We will assume that we are in a high-dimensional setting
where $p>n$ for large enough $n$ by imposing the condition that there exists a
constants, \jelenax{independent of the sample or covariate size}, $\kappa_1,\kappa_2>0$ such that $\kappa_1 \log p \leq n\leq \kappa_2 p\log p$ for large enough $n.$

\jelenax{Define
$$
\mathcal{M}^{*}_{C_{0},\xi}=\left\{ v\in\RR^{p}:\ \min_{\|a\|_0 \leq t} \|v-a\|_2 \leq C_{0} t^{-\xi}\quad \forall t\leq p \right\}.
$$ 
For $f(\tilde{Z})=Z^\top \beta$ and $g(\tilde{Z})=Z^\top \phi$  with  $Z=\psi(\tilde{Z})$, we can easily see that, under Assumption \ref{assp 1}, $f\in \mathcal{M}_{C_0,\xi_1}$ and $g\in \mathcal{M}_{C_0,\xi_2}$ if $\beta\in  \mathcal{M}^{*}_{C_0,\xi_1}$ and $\phi\in  \mathcal{M}^{*}_{C_0,\xi_2}$.}
The distribution of the data is now parameterized by $\lambda
=(\theta,\beta,\phi,\sigma_{u}^{2},\sigma_{\varepsilon}^{2})$. We use $\E_{\lambda}$ to denote the expectation under the distribution specified by $\lambda$. We define the following parameter space
\[
\Lambda_{\xi_{1},\xi_{2}}=\left\{  \lambda=(\theta,\beta,\phi,\sigma_{u}%
^{2},\sigma_{\varepsilon}^{2}):\theta\ \in\mathbb{R},\ \beta\in
\mathcal{M}^{*}_{C_{0},\xi_{1}},\ \phi\in\mathcal{M}^{*}_{C_{0},\xi_{2}},\ \ \{\sigma
_{u}^{2},\sigma_{\varepsilon}^{2}\}\subset\lbrack M^{-1},M]\right\}  ,
\]
where $C_{0},\xi_{1},\xi_{2}>0$ and $M\geq2$ are constants. 

Let 
\[
\mathcal{R}_{\xi_{1},\xi_{2}}=\operatorname{inf}_{\hat{\theta}}\sup
_{\lambda=(\theta,\beta,\phi,\sigma_{\varepsilon}^{2},\sigma_{u}^{2})\in \Lambda(\xi_{1}%
,\xi_{2})}\mathbb{E}_{\lambda}\left\vert \hat{\theta}(W)-\theta\right\vert ,
\]
where $\operatorname{inf}_{\hat{\theta}}$ is taken over all measurable
functions $\hat{\theta}$ of the data. The following result delivers a lower bound for
$\mathcal{R}_{\xi_{1},\xi_{2}}$.

\begin{theorem}\label{thm1}
Let $\tilde{\xi}=\max\{\xi_1,\xi_2 \}$. Suppose that Assumption  \ref{assp 1} is satisfied and there exist 
constants $\kappa_1,\kappa_2>0$  such that $ \kappa_{1}\log p\leq n\leq\kappa_{2}p\log p$  for
large enough $n$. If $\tilde{\xi}\leq 1/2$,   then 
\jelenax{\[
\mathcal{R}_{\xi_{1},\xi_{2}}\geq C(n^{-1}\log p)^{2\tilde{\xi}/(2\tilde{\xi}+1)},
\]}
where $C>0$  is a constant depending only on $\tilde{\xi}$\textit{, }$\kappa_1,\kappa_2,C_{0}$. 
\end{theorem}

\begin{figure}[h]
\caption{\label{fig:2} This graph depicts the exponent of the lower bound in Theorem \ref{thm1}. }
\begin{centering}
\includegraphics[scale=0.3]{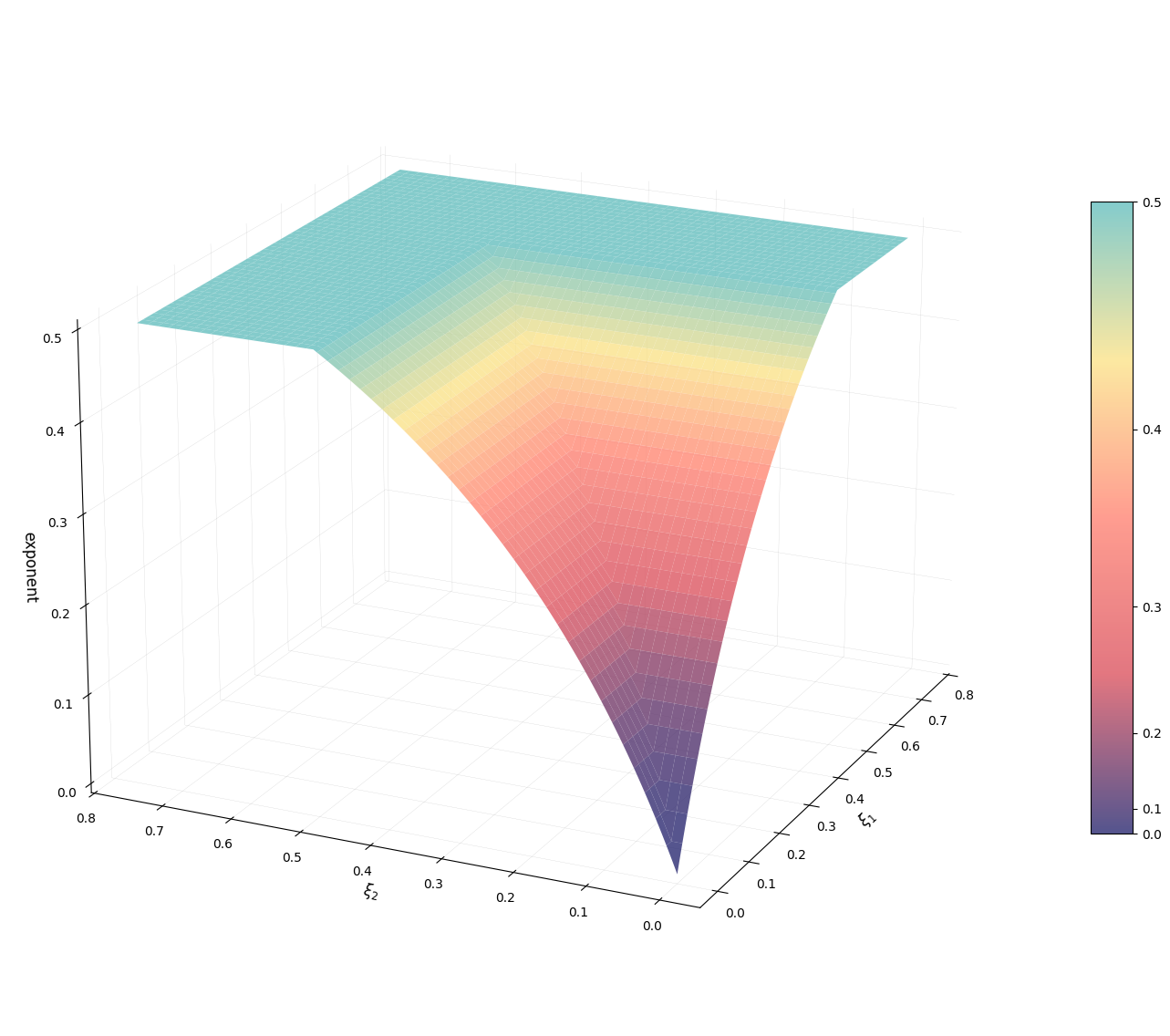}
\par\end{centering}
\end{figure}

 
 \jelenax{In this paper, our primary goal is to establish the most favorable and minimal conditions under which root‑\(n\) consistency and asymptotic normality can be achieved. To this end, we derive a lower bound that delineates when root‑\(n\) inference is possible and when it is not. Our results show that as long as at least one of the parameters \(\xi_1\) or \(\xi_2\) exceeds \(1/2\), root‑\(n\) inference can be attained. Specifically, when \(\max\{\xi_1, \xi_2\} = 1/2\), we achieve the optimal minimax rate of \(\sqrt{\log p/n}\); however, if the approximation constants \(\xi\) take on smaller values, the convergence rate deteriorates significantly, precluding root‑\(n\) inference.}

\jelenax{To illustrate our findings, Figure \ref{fig:2} plots the exponent of the lower bound from Theorem \ref{thm1}. Our result is conceptually similar to that of \cite{javanmard2018debiasing}, which also depends on the smaller of the two sparsity parameters—albeit under the much stronger assumption of a bounded \(\ell_1\) norm for the coefficients, Gaussian settings and within the framework of strictly sparse models only.} To the best of our knowledge, this is the only lower bound for approximate sparsity.

\jelenax{Our lower bounds do not follow from strictly sparse models (e.g., \cite{ren2015asymptotic} and \cite{cai2017confidence}), since, as Lemma \ref{lem: strict vs approx sparse} shows, strictly and approximately sparse models do not nest. Thus, a lower bound for strictly sparse models doesn't automatically apply to approximately sparse ones. While the exactly sparse minimax lower bound offers some intuition, the approximately sparse setting demands a fundamentally different approach. Our proof constructs a tailored prior—dependent on \(\xi_1\) and \(\xi_2\)— but goes well beyond a simple corollary.}

 In the proof of Theorem \ref{thm1}, we actually show a stronger result on the adaptivity of confidence intervals. In particular, suppose that $\max\{\tilde{\xi}_{1},\tilde{\xi}_2\}\leq 1/2<\max\{\xi_1,\xi_2\}$. By the rate in Theorem \ref{thm1}, it is impossible to construct a confidence interval that is valid on  $\Lambda_{\tilde{\xi}_1,\tilde{\xi}_2}$ and have length $O(n^{-1/2})$ on $ \Lambda_{\tilde{\xi}_1,\tilde{\xi}_2}$. However, we show later (Theorem \ref{thm: PLM main}) that it is possible to construct a confidence interval that is valid on  $\Lambda_{\xi_1,\xi_2}$ and have length $O(n^{-1/2})$ on $\Lambda_{\xi_1,\xi_2}$. Since $\Lambda_{\xi_1,\xi_2}\subset \Lambda_{\tilde{\xi}_1,\tilde{\xi}_2}$, a natural question  is whether one can have a confidence interval that is valid on  $\Lambda_{\tilde{\xi}_1,\tilde{\xi}_2}$ and have length $O(n^{-1/2})$ on $\Lambda_{\xi_1,\xi_2}$. In the proof of Theorem \ref{thm1}, we show that such a confidence interval does not exist. Any confidence interval that is valid on  $\Lambda_{\tilde{\xi}_1,\tilde{\xi}_2}$ must have length at least $O((n^{-1}\log p)^{2\tilde{\xi}/(2\tilde{\xi}+1)})$ even at points in $\Lambda_{\xi_1,\xi_2}$. In other words, no confidence interval can automatically adapt its length to the approximation indices. Due to the heavy notation needed to rigorously state this phenomenon, we leave the detailed discussions in the appendix.

\subsection{Lower Bound Under Ordered Approximation}\label{sec:ordered}

We now consider ordered basis functions for the model in  Assumption \ref{assp 1}.
\jelenax{For $f(\tilde{Z})=Z^\top \beta$ and $g(\tilde{Z})=Z^\top \phi$ with  $Z=\psi(\tilde{Z})$, we can observe that, under Assumption \ref{assp 1}, $f\in \mathcal{S}_{C_0,\xi_1}$ and $g\in \mathcal{S}_{C_0,\xi_2}$ if $\beta\in  \mathcal{S}^{*}_{C_0,\xi_1}$ and $\phi\in  \mathcal{S}^{*}_{C_0,\xi_2}$, where 
$$
\mathcal{S}^{*}_{C,\xi}=\left\{  v\in\RR^{p}:\ \min_{a=(a_{1}%
,...,a_{t},0,...,0)^{\top}} \|v-a\|_2 \leq C t^{-\xi}\quad \forall t\leq p \right\}.
$$
Therefore, under ordered basis functions, the distribution of the data is indexed by elements of the following set
\[
\Lambda^{\rm{Ordered}}(\xi_{1},\xi_{2})=\left\{
\begin{array}
[c]{c}%
\lambda=(\theta,\beta,\phi,\sigma_{\varepsilon}^{2},\sigma_{u}^{2}):\ \beta\in{\mathcal{S}
}_{M_1,\xi_{1}},\ \phi\in{\mathcal{S}}_{M_1,\xi_{2}},
 \sigma_{\varepsilon}^{2},\sigma_{u}^{2}\in\lbrack M_{2}^{-1},M_{2}]
\end{array}
\right\}  ,
\]
where $M_{1},\xi_1,\xi_2>0$ and $M_{2}>2$ are constants. }The minimax rate for estimating
$\beta$ is 
\[
\mathcal{R}^{\rm{Ordered}}_{\xi_{1},\xi_{2}}=\operatorname{inf}_{\hat{\theta}}\sup
_{\lambda=(\theta,\beta,\phi,\sigma_{\varepsilon}^{2},\sigma_{u}^{2})\in \Lambda^{\rm{Ordered}}(\xi_{1}%
,\xi_{2})}\mathbb{E}_{\lambda}\left\vert \hat{\theta}(W)-\theta\right\vert ,
\]
where $\operatorname{inf}_{\hat{\theta}}$ is taken over all measurable
functions $\hat{\theta}$ of the data. In other words, $\mathcal{R}^{\rm{Ordered}}_{\xi_{1},\xi_{2}}$ is the minimax risk—the smallest achievable worst‑case expected absolute error when the true parameter lies in $\Lambda^{\rm{Ordered}}(\xi_{1},\xi_{2})$.
The goal is to derive a lower bound for
$\mathcal{R}^{\rm{Ordered}}_{\xi_{1},\xi_{2}}$ when $\xi_{1}+\xi_{2}<1/2$.

\begin{theorem}\label{thm2}
\jelenax{If Assumption  \ref{assp 1} holds, and 
 $p\gtrsim n^{2}$, then, for all $\xi_{1}+\xi_{2}<1/2$, we have} $$\mathcal{R}^{\rm{Ordered}}_{\xi_{1},\xi_{2}}\gg n^{-1/2}.$$
\end{theorem}

Hence the bound is rate-optimal: \cite{donald1994series} construct a series estimator that attains the $\sqrt{n}$ convergence rate whenever the smoothness indices satisfy $\xi_{1}+\xi_{2}\ge 1/2$, exactly matching the bound and confirming its sharpness.
\jelenax{ Although this result might seem expected given minimax rates for H\"older classes in \cite{robins2009semiparametric}, existing results do not imply Theorem \ref{thm2}. The reason is that Theorem \ref{thm2} deals with models with increasing dimensions (similar to Gaussian sequence models) rather than \textit{fixed, low-dimensional} models.  As a result, the priors needed to obtain our bounds have a different structure. \jelenax{\cite{robins2009semiparametric} perturb kernel functions at various centers, whereas we perturb linear combinations with Gaussian deviations.}}

\jelenax{The proof of Theorem \ref{thm2} is based on a new construction of the priors. 
We sharpen the construction by contrasting the contiguous alternatives \((\beta,\phi)=(\zeta,\zeta)\) and \((\beta,\phi)=(0,\zeta)\), where \(\zeta\sim\mathcal N(0,\tau^{2})\) is tuned to balance bias and variance.  The proof then hinges on securing a sharp upper bound for the total‑variation distance between the associated mixture distributions, as established in Lemma \ref{lem: A5}.
Moreover, we  need to address the fact that these random perturbations can be outside $\Lambda^{\rm{Ordered}}(\xi_{1},\xi_{2})$ since Gaussian distributions are unbounded. Hence, another key argument establishes that enough of these perturbations are in $\Lambda^{\rm{Ordered}}(\xi_{1},\xi_{2})$. 
}

 \jelenax{Although Theorem \ref{thm1} addresses the optimal confidence interval length, it also implies that no estimator can achieve the parametric \(\sqrt{n}\)-rate when \(\max\{\xi_1,\xi_2\} \le 1/2\). If a \(\sqrt{n}\)-consistent estimator existed in this regime, one could construct a confidence interval with \(O(n^{-1/2})\) length. However, by Theorem \ref{thm1} (with \(\tilde{\xi}_1=\xi_1\) and \(\tilde{\xi}_2=\xi_2\)), such intervals cannot exist. Therefore, no estimator can be \(\sqrt{n}\)-consistent under \(\max\{\xi_1,\xi_2\} \le 1/2\).}

\begin{remark}
    \jelenax{C}omparing Theorems \ref{thm1} and \ref{thm2} \jelenax{highlights the crucial} role of ordering \jelenax{in the dictionary of basis functions.} \jelenax{In the {\bf{unordered}} case-- where} the most important terms might not be the first a few basis \jelenax{functions--} the parametric \jelenax{convergence} rate is only possible when $\max\{\xi_1,\xi_2\}>1/2 $. \jelenax{Conversely,} for {\bf ordered} models \jelenax{-- where the most significant terms appear first in the dictionary--}  the parametric rate \jelenax{can be achieved as long as} $\xi_1+\xi_2\geq 1/2$. As illustrated in Figure \ref{fig:1}, the \jelenax{condition} $\xi_1+\xi_2\geq 1/2$ is \jelenax{ considerably less restrictive} than $\max\{\xi_1,\xi_2\}>1/2 $ \jelenax{-- this difference is depicted in the figure as a triangle versus a square}. \jelenax{This disparity underscores the cost of} being unable to order the basis functions \jelenax{according to} their approximation power.   \end{remark}

\section{Upper bound for Learning Individual Coefficients} \label{sec PLM upper bound}
\jelenax{Through Section 3, we have determined that achieving the root-$n$ convergence rate for coefficient estimation in}  requires 
$\max \{ \xi_1,\xi_2 \}>1/2 $\jelenax{, when using an unordered dictionary of basis functions.} We now show that this condition is also sufficient by providing a \jelenax{new} estimator \jelenax{that achieves the lower bound of Section \ref{sec:unordered}}.
\jelenax{In this section, we will remove the Gaussian design assumption and work with the following models}
\begin{equation}
Y_{i}=D_{i}\theta+f(\tilde{Z}_{i})+\varepsilon_{i}\quad\text{and}\quad D_{i}=g(\tilde{Z}_{i})+u_{i}, \label{eq: PLM}
\end{equation}
where $ \mathbb{E}(\varepsilon_{i}\mid D_{i},\tilde{Z}_{i})=0$ and $\mathbb{E}(u_{i}\mid \tilde{Z}_{i})=0$. \jelenax{We consider $f\in \mathcal{M}_{C,\xi_1}$ and $g \in \mathcal{M}_{C,\xi_2}$ for some $C,\xi_1,\xi_2>0$. }For notational convenience, we denote the basis functions by $Z_{i}=(\psi_{1}(\tilde{Z}_{i}),...,\psi_{p}(\tilde{Z}_{i}))^{\top}\in\RR^{p}$ and $X_{i}=(D_{i},Z_{i}^{\top})^{\top}\in\RR^{p+1}$.

We approximate $f(\tilde{Z}_{i})$ and $g(\tilde{Z}_{i})$ by their
least-square projections on the basis functions: $Z_{i}^{\top}\beta$
and $Z_{i}^{\top}\phi$, where $Z_{i}=(\psi_{1}(\tilde{Z}_{i}),...,\psi_{p}(\tilde{Z}_{i}))^{\top}\in\RR^{p}$
is the basis functions, $\Sigma_{Z}=\mathbb{E}Z_{i}Z_{i}^{\top}$
and $\beta$ and $\phi$ are defined by 
\[
\mathbb{E}f(\tilde{Z}_{i})Z_{i}=\Sigma_{Z}\beta\quad\text{and}\quad\mathbb{E}g(\tilde{Z}_{i})Z_{i}=\Sigma_{Z}\phi.
\]
Let $\gamma=(\theta,\beta^{\top})^{\top}\in \RR^{p+1}$. For simplicity, we partition the \jelenax{set of} observation indices \(\{1, \ldots, n\}\) into two \jelenax{subsets of approximately} equal size, \jelenax{denoted by \(I_1\) and \(I_2\), with \(n_{\ell} = |I_{\ell}|\) representing the number of observations in subset \(I_{\ell}\) for \(\ell = 1, 2\). We define \(\mathbb{E}_{n,\ell}\)} as the sample average \jelenax{over} \(I_{\ell}\), given by \jelenax{\(\mathbb{E}_{n,\ell}[\cdot] = n_{\ell}^{-1} \sum_{i \in I_{\ell}} [\cdot]\). }For example, \(\mathbb{E}_{n,1}[Z_i Z_i^\top] = n_1^{-1} \sum_{i \in I_1} Z_i Z_i^\top\) and \(\mathbb{E}_{n,2}[Z_i D_i] = n_2^{-1} \sum_{i \in I_2} Z_i D_i\), and so on.

\jelenax{We introduce a novel multi-step estimation procedure. In the first stage, we estimate the parameter vector \(\beta  \) using a new Dantzig-type estimator that incorporates residual-balancing constraints with tailored cross-fitting for each moment condition. In the second stage, we obtain an estimator for \(\phi\), and finally, we construct a bias-corrected estimator for \(\theta\). }


\jelenax{Specifically, we define the estimator for 
 $\gamma$  as $\hat{\gamma} = (\hat{\theta}, \hat{\beta})$ and
\begin{align}
 (\hat{\theta}, \hat{\beta}) = \arg\min_{a \in \mathbb{R},\ b \in \mathbb{R}^{p}} \ & |a|+ \| b \|_{1} \label{eq: Dantzig_PLM} \\
\text{s.t.} \quad 
& \left| \mathbb{E}_{n,2}\left[D_{i} (Y_i - D_i a - Z_i^\top b ) \right]  \right| \leq \lambda_{1} \nonumber \\
& \left\| \mathbb{E}_{n,2}[Z_{i} (Y_i - D_{i} a)] -  \mathbb{E}_{n,1}[Z_{i} Z_{i}^\top] \, b  \right\|_{\infty} \leq \lambda_{1} \nonumber \\
& \left| \mathbb{E}_{n,2}[D_{i} (Y_i - D_i \, a)]  - \hat{\phi}^\top \mathbb{E}_{n,1}[Z_{i} Z_{i}^\top] \, b   \right| \leq \lambda_{2}, \nonumber
\end{align}}
\jelenax{where  
  $\lambda_{1}=\lambda_{0}\sqrt{n^{-1}\log p}$  and $\lambda_{2} = \lambda_{0} \sqrt{n^{-1/2}}$ for a large enough
constant $\lambda_{0}>0$. The nuisance parameter \(\phi\) is estimated via a  
 {\it cross-fitted} 
  Dantzig selector  defined by
  \begin{equation} \hat{\phi} = \arg\min_{v \in \mathbb{R}^{p}} \| v \|_{1} \quad \text{subject to} \quad \left\| \mathbb{E}_{n,1}[Z_i Z_i^\top] v - \mathbb{E}_{n,2}[Z_i D_i] \right\|_{\infty} \leq \lambda_{1}. \label{eq: Dantzig_PLM_part_0} 
  \end{equation}
 }
  \jelenax{The estimator is novel in several aspects. While the estimator for \(\phi\) (\(\hat{\phi}\)) is constructed using conventional cross-fitting methods, the estimator for \(\gamma\) (\(\hat{\gamma}\)) introduces a new approach. First,  it decouples the residual balancing in model \eqref{eq: PLM} into two distinct components: one for \(\theta\) and another for \(\beta\). The residual balancing for \(\theta\) is computed solely on the in-sample data from \(I_2\), whereas that for \(\beta\) employs a specialized cross-fitting technique that matches the right moments. Our approach is based on decomposing the moment condition
 $
 \mathbb{E}\bigl[Z\bigl(Y - D\theta - Z^\top\beta\bigr)\bigr]=0
$
into the equivalent equality
\[
 \mathbb{E}\bigl[Z\bigl(Y - D\theta\bigr)\bigr]= \mathbb{E}\bigl[ZZ^\top\beta\bigr].
\]
In practice, we approximate this by ensuring that the empirical moment computed on subsample \(I_1\), namely
$
\mathbb{E}_{n,1}\bigl[Z(Y-D\theta)\bigr],
$
is close to the moment computed on subsample \(I_2\),
$
\mathbb{E}_{n,2}\bigl[ZZ^\top\beta\bigr].
$
This splitting allows us to decouple the estimation of \(\theta\) and \(\beta\) effectively.}

 \jelenax{The final constraint in constructing \(\hat{\theta}\) is key to controlling second-order bias effects, which become problematic if only one of \(\phi\) or \(\beta\) is approximately sparse. We decompose
\[
\phi^\top \mathbb{E}[Z_i Z_i^\top]\beta - \hat{\phi}^\top \mathbb{E}[Z_i Z_i^\top]\hat{\beta} \approx \mathbb{E}\!\left[(D_i - Z_i^\top \hat{\phi})\, Z_i^\top\beta + Z_i^\top \hat{\phi}\,(Z_i^\top\beta - Z_i^\top\hat{\beta})\right],
\]
which further approximates to
\[
\mathbb{E}\!\left[D_i (Y_i - D_i\theta) - \hat{\phi}^\top Z_i Z_i^\top\beta\right] + \mathbb{E}\!\left[Z_i^\top \hat{\phi}\,(Z_i^\top\beta - Z_i^\top\hat{\beta})\right],
\]
using \(\mathbb{E}[D_i(Y_i-D_i\theta)] \approx \mathbb{E}[D_i Z_i^\top\beta]\). Our third constraint forces the first term to vanish by using cross-fitting to ensure that 
\[
\mathbb{E}_{n,1}\!\left[D_i (Y_i - D_i\theta)\right] \approx \mathbb{E}_{n,2}\!\left[\hat{\phi}^\top Z_i Z_i^\top\beta\right].
\]
The remaining term is straightforward to control. 
 }

\jelenax{We propose the following {\it in-sample} estimator of $\theta$}
\begin{equation}\label{eq:tildetheta}
\tilde{\theta}_{1} = \hat \theta - \hat{\pi}^\top \left\{  \hat{\Sigma}_{(1)} \hat{\gamma} - \mathbb{E}_{n,1}[ X_i Y_i ] \right\}.
\end{equation}
\jelenax{where $\hat \theta =e_{1}^\top \hat{\gamma}$,  from \eqref{eq: Dantzig_PLM},  \(\hat{\Sigma}_{(1)}\) is the sample covariance matrix computed from \(I_1\), $\hat{\Sigma}_{(1)} = \mathbb{E}_{n,1}[X_i X_i^\top] $, \(e_1\) is the first canonical basis vector in \(\mathbb{R}^{p+1}\) and  the estimate of  \(\pi\) is obtained using only the observations from the first sample, \(I_1\), as \begin{equation*}
\hat{\pi} = \arg\min_{v\in\mathbb{R}^{p+1}} \|v\|_1 \quad \text{subject to} \quad \|\hat{\Sigma}_{(1)}v - e_1\|_\infty \leq \lambda.
\label{eq: Riesz rep est PLM}
\end{equation*}
Cross-fitting, which dates back to \cite{schick1986asymptotically},  is a standard technique that uses sample splitting to reduce the bias. It is important to note that our estimator is not a standard cross-fitted estimator. Although the estimator \(\hat{\gamma}\) uses a (newly designed non-standard) cross-fitting strategy that combines information from both \(I_1\) and \(I_2\), both \(\hat{\Sigma}_{(1)}\) and \(\hat{\phi}\) are computed solely from \(I_1\). In a typical cross-fitting setup, one might use, say, \(I_2\) to construct both \(\hat{\gamma}\) and \(\hat{\phi}\) and then combine them using \(I_1\) to form the final estimator. Here, however, \(\hat{\gamma}\) is derived in a specially designed manner that utilizes both samples, while \(\hat{\phi}\) is based exclusively on \(I_1\). This justifies referring to \(\tilde{\theta}_1\) as an {\it in-sample} estimator.
}

Let $\tilde{\theta}_{2}$ be the same as $\tilde{\theta}_{1}$ except
that we swap the role of $I_{1}$ and $I_{2}$.  The final estimator
is \jelenax{ is then constructed as a weighted average of 
 $\tilde{\theta}_{1}$ and $\tilde{\theta}_{2}$ :}
\[
\tilde{\theta}=\frac{n_{1}}{n}\tilde{\theta}_{1}+\frac{n_{2}}{n}\tilde{\theta}_{2}.
\]

\jelenax{\begin{assumption}
\label{assu: PLM}Consider the model in (\ref{eq: PLM}). The
eigenvalues of $\Sigma=\mathbb{E}[X_{i}X_{i}^{\top}]$ and $\Sigma^{-1}$
bounded by constants. For some constant $\kappa>0$, $\mathbb{E}(\varepsilon_{i}^{2}\mid X_{i})$,
$\mathbb{E}|\varepsilon_{i}|^{2+\kappa}$, $\E (u_{i}^{2}\mid Z_{i})$, $(\E u_{i}^{2})^{-1}$  and $\mathbb{E}|Y_{i}|^{2+\kappa}$
are also bounded by constants. Moreover, either 
\[
\|X_{i}\|_{\infty}\ \text{and }\mathbb{E}|u_{i}|^{2+\kappa}\text{ are bounded}
\]
or 
\[
X_{i}\text{ has bounded sub-Gaussian norm}.
\]
Furthermore, $p\gg n^{1/(2\min\{\xi_{1},\xi_{2}\})}$ and $\log p \ll n^{\kappa/(2\kappa+4)}$. 
\end{assumption}
}

\jelenax{The requirement of $\log p \ll n^{\kappa/(2\kappa+4)}$ is due to the fact that the data can be heavy-tailed. For sub-Gaussian data, we can essentially take $\kappa=\infty$ and then the requirement becomes the usual $\log p \ll \sqrt{n}$. }

\jelenax{Notice that we do not maintain any sub-Gaussian assumption on $\varepsilon_{i}$ in Assumption \ref{assu: PLM} and $\varepsilon_{i}$ does not need to be independent of $X_i$. The error term
$\varepsilon_{i}$ is only assumed to have bounded conditional variance and bounded $(2+\kappa)$ moments for some $\kappa>0$,
which is much weaker than the usual condition that $\varepsilon_{i}$
is independent of $X_{i}$ and is sub-Gaussian.  Note that a large literature on robust estimation develops a variety of tools just to handle heavy tails in the error term $\varepsilon_i$, e.g., \cite{fan2014adaptive}, \cite{sun2020adaptive} and \cite{chen2020robust}. Theoretical analysis of models with heavy-tailed data can be much more complicated than the sub-Gaussian case even in the usual nonparametric settings,   see e.g., \cite{han2019convergence} and \cite{kuchibhotla2022least}.  Here, we exploit the fact that $\varepsilon_i$ only shows up in the product $X_i \varepsilon_i$ and derive a concentration bound (Lemma \ref{lem: C2}) that exploits this structure and the boundedness of $X_i$.   
}

\jelenax{The condition of $p\gg n^{1/(2\min\{\xi_{1},\xi_{2}\})}$ needs $p$ to be sufficiently large. Since the order of the basis functions is unknown, a large number of basis functions is needed to reduce the approximation error, leading to a high-dimensional setting. This condition is not needed if  the true function is a linear combination of the $p$ basis functions ($f(\tilde{Z}_i)=Z_i^\top \beta$ and $g(\tilde{Z}_i)=Z_i^\top \phi$ for some $\beta,\phi\in\RR^{p}$), which is an assumption in  most of the literature on high-dimensional models. We now state the main result on the upper bound for learning individual coefficients. }

\begin{theorem}
\label{thm: PLM main}\jelenax{Consider the model in (\ref{eq: PLM}). Let
Assumption \ref{assu: PLM} hold. }Suppose that 
$$\log p\ll n^{1/2-1/(4\max\{\xi_{1},\xi_{2}\})},$$
$f\in\Mcal_{C_{0},\xi_{1}}$ and $g\in\Mcal_{C_{0},\xi_{2}}$
for some constants $C_{0},\xi_{1},\xi_{2}>0$ such that $\max\{\xi_{1},\xi_{2}\}>1/2$.
Then 
\[
\sqrt{n}(\tilde{\theta}-\theta)=n^{-1/2}\sigma_{u}^{-2}\sum_{i=1}^{n}u_{i}\varepsilon_{i}+o_{P}(1)
\]
and $\tilde{\theta}-\theta=O_{P}(n^{-1/2})$. 
\end{theorem}

 \jelenax{This is one of the main contributions of our paper, as we establish the root‑\(n\) rate under the condition \(\max\{\xi_1,\xi_2\} > 1/2\).
 The condition \(\max\{\xi_1,\xi_2\} > 1/2\) is considerably more general than the rate double-robustness condition of \(\xi_1\xi_2 > 1/4\). While the latter requires both \(\xi_1\) and \(\xi_2\) to be moderately large—ensuring that neither parameter is too close to zero—the \(\max\) condition only necessitates that one of the two parameters exceeds \(1/2\). This means that even if one parameter is quite small, as long as the other is sufficiently large, the condition for root‑\(n\) inference is still satisfied. Consequently, our requirement is less restrictive and applies to a broader class of models than the traditional RDR; see Figure \ref{fig:1}.} \jelenax{Thus, we significantly broaden the scope of settings where root-\(n\) inference is attainable. } \jelenax{Moreover, our approach relaxes several common assumptions. First, we do not require the error term \(\varepsilon_i\) to be independent of the design; a zero conditional mean, \(\mathbb{E}(\varepsilon_i\mid X_i)=0\), suffices, thus allowing for heteroskedasticity. Second, the error term need not have sub-Gaussian tails; bounded \(2+\kappa\) moments are sufficient. This is an improvement over many existing works on linear models that assume  Gaussian or sub-Gaussian errors, e.g., \cite{van2014asymptotically}, 
 \cite{,javanmard2014confidence,javanmard2018debiasing}, \cite{cai2017confidence} and \cite{bellec2022debiasing}.}   This estimator achieves semiparametric efficiency under homoscedastic errors in the sense that it is efficient among all asymptotically linear estimators, see Theorem 6 of \cite{jankova2018semiparametric}.  

\section{Upper Bound for Linear Functionals of Approximately Sparse Models}\label{sec general functionals}

\jelenax{We use the Riesz representer formulation combined with Neyman's orthogonal loss function to achieve root‑n inference for linear functionals. While previous studies have explored Riesz representers within the framework of Neyman's orthogonality, our approach employs a novel non‑symmetric cross‑fitting procedure. We demonstrate that this method enables root‑n inference under more general conditions than those addressed in the existing literature.}

In this section, we aim to estimate
\[
\theta_0 = \mathbb{E}[m(W,\rho_0)],
\]
where \(\rho_0 = \mathbb{E}(Y \mid X)\) and \(m(\cdot, \cdot)\) is a given function. Assuming that the mapping \(\rho \mapsto \mathbb{E}[m(W,\rho)]\) is continuous in the \(\ell_{2}\) norm on a Hilbert space \(\mathcal{B}\), Theorem 5.25 in \cite{folland1999real} guarantees the existence of a unique element \(\alpha_0 \in \mathcal{B}\) (called the Riesz representer) satisfying \(\mathbb{E}[\alpha_0(X)^2] < \infty\) such that
\[
\mathbb{E}[m(W,\rho)] = \mathbb{E}[\alpha_0(X)\rho(X)]
\]
for all \(\rho \in \mathcal{B}\).
In this section, we assume the following model setting:
\jelenax{\[
\rho_0 \in \mathcal{M}_{C,\xi_1} \quad \text{and} \quad \alpha_0 \in \mathcal{M}_{C,\xi_2},
\]
for some positive constants \(C\), \(\xi_1\), and \(\xi_2\).}
\subsection{General methodology and theory}
 The newly proposed  estimator of $\theta_{0}$ ,  is given by $\hat{\theta}  = {n_{1}}\hat{\theta}_{1}/n + {n_{2}}\hat{\theta
}_{2}/n$  with an {\it asymmetric cross-fitted} estimator
\begin{align}
\hat{\theta}_{\ell} 
&= \E_{n,\ell} \Bigl[\hat{\alpha}_{\ell}(X_{i})Y_{i}   + m(W_{i},\hat{\rho}_{-\ell}) - \hat{\alpha}_{\ell}(X_{i}) \hat{\rho}_{-\ell}(X_{i})   \Bigl],
\label{ESTgen}
\end{align}
where $\ell \in \{1,2\}$ and $\hat{\alpha}_{\ell}$ is an estimator build on the sample $I_\ell$ whereas $\hat{\rho}_{-\ell}$ estimates $\rho_0$ and uses the complement set of $I_\ell$.  We utilize well known estimates of the  nuisance parameters $\alpha_0$ and $\rho_0$,
\begin{align}
\hat{\rho}_{\ell}(x)  &  :=\psi(x)^{\top}\hat{\gamma}_{\ell}, \quad  
\hat{\gamma}_{\ell}:=\arg\min_{\gamma}\left\{  \gamma^{\top}\hat{\Sigma
}_{\ell}\gamma-2\hat{\mu}_{\ell}^{\top}\gamma+2\lambda \left\Vert \gamma\right\Vert
_{1}\right\}  \label{pieces} \\ 
\hat{\alpha}_{\ell}(x)  &  =\psi(x)^{\top}\hat{\pi}_{\ell},
\quad   \hat{\pi
}_{\ell}=\arg\min_{\pi}\left\{  \pi^{\top}\hat{\Sigma}_{\ell}\pi-2\hat
{M}_{\ell}^{\top}\pi+2\lambda \left\Vert \pi\right\Vert _{1}\right\}  ,\text{
}\label{RRgen}
\end{align}
where 
$ \hat{\Sigma}_{\ell}     :=  \E_{n,\ell} [\psi(X_{i})\psi(X_{i})^{\top}]$, $ \hat{\mu}_{\ell}:=  \E_{n,\ell} [ \psi(X_{i})Y_{i}]$, as well as
$
\hat{M}_{\ell}    = \E_{n,\ell} [m(W_{i},b)]$ for $m(W,b)=(m(W,\psi_{1}),...,m(W,\psi_{p}))^{\top}.$
As, we are only cross-fitting $\hat{\alpha}_{\ell}$ in \eqref{ESTgen}, our estimator differs from \cite{chernozhukov2018double,chernozhukov2018learning}.
The above nuisance estimators themselves, are proxies for the population
parameters $\rho_n(X)=\psi(X)^\top\gamma$ and $\alpha_n(X)=\psi(X)^\top\pi$ with 
$\gamma$  and $\pi$ solving \[
\mu=\Sigma\gamma\text{, } \quad M=\Sigma\pi,
\]
where $\mu=\mathbb{E}[\psi(X)Y],$ $\Sigma=\mathbb{E}[\psi(X)\psi(X)^{\top}]$, and
$M=\mathbb{E}[m(W,b)]=\mathbb{E}[\psi(X)\alpha_0(X)]$.  For the above to be valid,  we  work with the following  assumption.

\jelenax{\begin{assumption}\label{assp 6}
The largest eigenvalues of $\Sigma$ and $\Sigma^{-1}$  are bounded by some constant $C>0$.  There exists a constant $\kappa>0$ such that $\E |\alpha_{0}(X_{i})|^{2+\kappa}$, $\E |Y|^{2+\kappa}$, $Var(Y\mid X)$ and the largest eigenvalues of  $\mathbb{E}[\alpha_0(X)^{2} \psi(X)\psi(X)^{\top}]$ and $\mathbb{E} [m(W,\psi)m(W,\psi)^{\top}]$ are  bounded. Moreover, either 
\[
\|\psi(X_{i})\|_{\infty}\ \text{and }\mathbb{E}|u_{i}|^{2+\kappa}\text{ are bounded}
\]
or 
\[
\psi(X_{i})\text{ has bounded sub-Gaussian norm}.
\]
The sub-Gaussian norm of $m(W,\psi_{j})$ is bounded uniformly in $j\in \{1,...,p \}$. Furthermore, $p\gg n^{1/(2\min\{\xi_{1},\xi_{2}\})}$, $\log p \ll n^{\kappa/(2\kappa+4)}$ and $\log p\ll n^{1/2-1/(4\max\{\xi_{1},\xi_{2}\})}$. 
\end{assumption}
}

\jelenax{Assumption \ref{assp 6} allows for unbounded Riesz presenters. This is particularly important in many problems. For example,  in the context of learning the average treatment effect, we allow the propensity score to be close to one or zero. This is a new contribution to the literature because the usual assumption states that the propensity score is in $[c,1-c]$ with probability one for some fixed constant $c>0$, e.g.,  \cite{farrell2015robust}, \cite{chernozhukov2018double}, \cite{ning2020robust},  \cite{wang2024debiased} and \cite{liu2023root} among others. See Remark \ref{remark ATE} for details. }

\jelenax{\begin{assumption}\label{assp 9}
Approximation parameters $\xi_1$ and $\xi_2$ satisfy
\begin{equation}
\xi_{1}>1/2, \quad n\gg (\log p)^{(2\xi_1)/(\xi_1-1/2)} \tag{i}
\end{equation}
or 
\begin{equation}
\xi_{1}\xi_{2}>1/4, \quad n^{\xi_1\xi_2-1/4}\gg (\log p)^{2\xi_1\xi_2+(\xi_1+\xi_2)/2}. \tag{ii}
\end{equation}
\end{assumption}
}

\jelenax{Only one of the above two conditions should hold at any given time. The former
  ensures that
$\sqrt{n}\left\Vert \rho_{0}-\rho_{n}\right\Vert _{2}\rightarrow0$, where as the latter is in line with RDR condition.  Note that most of the literature on high-dimensional models assumes that $\rho_0=\rho_n $, i.e., the true function is a linear combination of the $p$ basis functions.  }

 We now state the main result on the rate of convergence. 
\jelenax{\begin{theorem}\label{thm4}
Let Assumptions  \ref{assp 6} and \ref{assp 9} hold.  Then
$
\hat{\theta}-\theta_{0}=O_{P}(n^{-1/2}), 
$
and in particular, 
\[
\hat{\theta}-\theta_{0}=n^{-1}\sum_{i=1}^{n}\tilde{\varUpsilon}_{n}(W_{i})+o_{P}(n^{-1/2}),
\]
where $\tilde{\varUpsilon}_{n}(W_i)=m(W_i,\rho_{n})-\theta_{n}+\alpha_{n}(X_i)[Y_i-\rho_{0}(X_i)]$ such that  $\theta_n=\E m(W_i,\rho_n)$,   $\mathbb{E}\tilde{\varUpsilon}_{n}(W_i)=0$ and $\mathbb{E}\tilde{\varUpsilon}_{n}(W_i)^2$ is bounded by a constant.
\end{theorem}
}

By Theorem \ref{thm4}, under weak regularity conditions, the estimator for the generic problem achieves the parametric rate whenever $\xi_1>1/2$  or $\xi_1 \xi_2>1/4$. This is Area $A\bigcup B $ in Figure \ref{fig:1}. 

\begin{remark}\label{rem ATE}
    For a procedure that can be applied to a generic problem, this estimator achieves parametric rate under weaker conditions compared to some existing estimators that are designed for specific problems. In the problem of estimating ATE (Example \ref{example: ATE}), no existing estimator can guarantee the parametric rate in both $\xi_1>1/2$ and $\xi_1 \xi_2>1/4$. \cite{athey2018approximate} only impose the assumptions on $\rho_0$, which would roughly correspond to $\xi_1>1/2$ here; however, whether their estimator achieves the parameter rate under $\xi_1 \xi_2>1/4$ is unknown. Works such as  \cite{farrell2015robust}, \cite{ning2020robust} and \cite{chernozhukov2022debiased} among others impose assumptions that imply $\xi_1 \xi_2>1/4$; there is guarantee of the parametric rate under $\xi_1>1/2$. 

\end{remark}

\jelenax{Finally, we derive a procedure for statistical inference. For this, we consider a mild moment condition in order to establish asymptotic normality and consistent
estimation of the asymptotic variance. 
\begin{assumption}
\label{assu: asy normal}There exist a constant $c>0$ such that
$\E |m(W,\rho_{n})|^{2+c}$ is bounded.
\end{assumption}}

\jelenax{We estimate the asymptotic variance via cross-fitting: $\hat{V}=(n_{1}/n)\hat{V}_{1}+(n_{2}/n)\hat{V}_{2}$ 
with 
\[
\hat{V}_{\ell}=\frac{1}{n_{1}}\sum_{i\in I_{\ell}}\left[\left(m(W_{i},\hat{\rho}_{-\ell})-\hat{\theta}\right)^{2}+\hat{\alpha}_{\ell}(X_{i})^{2}[Y_{i}-\hat{\rho}_{-\ell}(X_{i})]^{2}\right].
\]}

\jelenax{This estimator is novel—it goes beyond the standard approach of simply computing the sample second moment of plug‑in estimates of the influence function. In many existing works (see, e.g., \cite{chernozhukov2022automatic}), plug‑in estimators require extra trimming procedures to achieve consistency. In the absence of guaranteed RDR, estimating the variance of the influence function becomes challenging because certain terms cannot be easily controlled.
 We utilize cross‑fitting to effectively eliminate bias terms that would otherwise compromise consistency.   
}

\jelenax{\begin{theorem}
\label{thm: asym norm gen}Let the Assumptions  \ref{assp 6}, \ref{assp 9} and  \ref{assu: asy normal} hold. Then 
\[
\frac{n^{1/2}(\hat{\theta}-\theta_0)}{\sqrt{\mathbb{E}\tilde{\varUpsilon}_{n}(W)^{2}}}\overset{d}{\rightarrow}N(0,1)
\]
and $\hat{V}=\mathbb{E}\tilde{\varUpsilon}_{n}(W)^{2}+o_{P}(1)$. 
\end{theorem}
 The proposed estimator achieves certain optimality properties. In Theorem 1 of \cite{hirshberg2021augmented}, semiparametric efficiency bound of linear functionals (including the average treatment effect) is derived under regularity conditions. Our asymptotic variance above attains this  bound.
}

\subsection{Discussion on Model Settings and Assumptions}\label{subsec: RR}
We now give more examples regarding the Riesz representer and give more detailed discussions on the assumptions.

\begin{example}[Average Treatment Effect]\label{example: ATE}
Suppose that we are interested in the average treatment effect of a binary $D\in\{0,1\}$ treatment variable: $\theta_0=\E[Y(1)-Y(0)]$, where $Y(1)$ and $Y(0)$ are potential outcomes. We observe $(Y,D,Z)$, where $Y=Y(1)\cdot D +Y(0)\cdot(1-D)$ and $Z$ is a vector of covariates. We work with the usual unconfoundedness condition:  $(Y(1),Y(0))$ and $D$ are independent conditional on $Z$. Let $X=(D,Z)$ and  $\psi(X)=\psi(D,Z)$. We define the basis functions as follows: $\psi_{2k-1}(d,z)=\tilde{b}_k(z) $ and $\psi_{2k}(d,z)=d\cdot \tilde{b}_k(z) $ for $k\geq 1$ where $\{\tilde{b}_j(z)\}_{j=1}^{\infty}$ 
is a dictionary of functions of $z$. Following the usual argument (e.g., \cite{rosenbaum1983central}), we can write the average treatment effect  as 
\[
\theta_{0}=\mathbb{E}[m(W,\rho_{0})],
\]
where $m(W,f)=f(1,Z)-f(0,Z)$. 
Let $\pi_{0}(Z)=P(D=1|Z)$ be the propensity score. Then by routine calculation the Riesz representer is 
$$
\alpha_0(X)=\frac{D}{\pi_{0}(Z)}-\frac{1-D}{1-\pi_{0}(Z)}
$$
assuming that $\alpha_0(X) \in \mathcal{B}$. 
This example was also considered in  \cite{chernozhukov2022locally} and \cite{hirshberg2021augmented}.
We include it here because of its fundamental importance and because the estimator we give is root-n consistent without the doubly robust rate condition.
\end{example}

\begin{remark}\label{remark ATE}

\jelenax{For Example \ref{example: ATE}, we notice that $m(W,\psi_{2k-1})=\psi_{2k-1}(1,Z)-\psi_{2k-1}(0,Z)=0$ and $m(W,\psi_{2k})=\psi_{2k}(1,Z)-\psi_{2k}(0,Z)=\tilde{b}_{k}(Z)$ for any $k\geq 1$. Most of the conditions in Assumption \ref{assp 6} are easy to verify. We focus on the boundedness of $\E|\alpha_0(X)|^{2+\kappa}$ and eigenvalues of $\mathbb{E}[\alpha_0(X)^{2} \psi(X)\psi(X)^{\top}]$. The condition of $\E |\alpha_{0}(X_{i})|^{2+\kappa}$ is needed because the parametric rate is impossible with $\E |\alpha_{0}(X_{i})|^{2}=\infty$ even if the propensity score is known, see e.g., \cite{hahn1998role}.  Moreover, the largest eigenvalues of $\mathbb{E}[\alpha_0(X)^{2} \psi(X)\psi(X)^{\top}]$ and $\mathbb{E} [m(W,\psi)m(W,\psi)^{\top}]$ are bounded whenever the largest eigenvalues of $\E(\tilde{b}(Z)\tilde{b}(Z)^{\top}/\pi_0(Z))$ and $\E(\tilde{b}(Z)\tilde{b}(Z)^{\top}/(1-\pi_0(Z)))$ are bounded, where $\tilde{b}(Z)=(\tilde{b}_1(Z),...,\tilde{b}_p(Z))^\top \in \RR^p$ is the basis function for $Z$ in the above example. This is significantly weaker than the commonly imposed strict overlap condition, which requires that $P(c \leq \pi_0(Z) \leq 1 - c) = 1$ for some constant $c > 0$. 
 As noted by \cite{d2021overlap}, 
  the strict overlap assumption can be unduly restrictive, particularly in high‑dimensional settings. 
  Our results 
   show that the parametric $n^{-1/2}$ rate 
    is achievable without strict overlap whenever $\max\{\xi_{1},\xi_{2}\}>1/2$.
    When  $\max\{\xi_{1},\xi_{2}\}\le 1/2 $, 
     the parametric rate remains out of reach even if strict overlap holds. 
}

\end{remark}

\begin{example}[Partially linear instrumental variable regression]\label{example: IV reg}

\jelenax{Consider the model $Y=\tau D+g(X)+U$ with $\E(U\mid X,Z)=0$, where
$D\in\RR$ is an endogenous variable and $Z\in\RR$ is an instrumental
variable. The goal is to learn the coefficient $\tau$ of the endogenous variable from observations of $W=(Y,D,Z,X)$. Notice that
$\E(Y\mid X,Z)-\E(Y\mid X)=\tau\cdot[\E(D\mid X,Z)-\E(D\mid X)]$.
Therefore, 
$$
\tau  =\frac{\E\left\{ Z\cdot[\E(Y\mid X,Z)-\E(Y\mid X)]\right\} }{\E\left\{ Z\cdot[\E(D\mid X,Z)-\E(D\mid X)]\right\} } =\frac{\E(YZ)-\E[Z\cdot\E(Y\mid X)]}{\E(DZ)-\E[Z\cdot\E(D\mid X)]}.
$$}

\jelenax{Notice that the term $\E[Z\cdot\E(Y\mid X)]$ is a linear functional
of the conditional mean function $\E(Y\mid X=\cdot)$, i.e., $\E [m(W,f)]$
with $m(W,f)=Z\cdot f(X)$. Moreover, the term $\E[Z\cdot\E(D\mid X)]$
is the same linear functional of the conditional mean function $\E(D\mid X=\cdot)$.
Since we can easily learn $\E(YZ)$ and $\E(DZ)$ at the parametric
rate, the problem of learning $\tau$ becomes solving the problem
of learning a linear functional of two conditional mean functions. 
}

\end{example}

\begin{example}[Continuous Treatment Effect in Panel Data]\label{example: PANEL}
\jelenax{Suppose that $Y_t$ is an outcome variable and $X_t=(D_t,Z_t)$ is a vector of regressors for two time periods $(t=1,2)$, where $D_t$ is a continuous treatment variable for period $t$ that varies with $t$. Also suppose that 
\[
Y_t=g_t(X_t)+A+\epsilon_t, \mathbb{E}[\epsilon_t|X_1,X_2]=0,
\]
where $A$ is a time invariant unobserved confounder, i.e. a fixed effect, such that $\mathbb{E}[A|X_1,X_2]$ is unrestricted. Define $Y=Y_2-Y_1$ and $X=(X_1,X_2)$ so that
\[
\rho_0(X)=\mathbb{E}[Y|X]=g_2(X_2)-g_1(X_1).
\]
A parameter of interest is 
\[
\theta_{0}=\mathbb{E}[m(W,\rho_{0})]=\mathbb{E}[\partial \rho_{0}(D_{1},Z_{1},D_{2},Z_{2})/\partial D_2],
\]
where $m(\cdot,\cdot)$ is defined by $m(w,f)=\partial f(x)/\partial d_{2}$. 
This parameter is a continuous average treatment effect for the second period as in for example \cite{klosin2023estimating}.
Consider a dictionary 
$\psi_{2k-1}(x)=\tilde{b}_k(x_2) $ and $\psi_{2k}(x)=-\tilde{b}_k(x_1) $ for $k\geq 1$ where $\{\tilde{b}_j(x)\}_{j=1}^{\infty}$ 
is a dictionary of functions of $x=(d,z)$ that are each differentiable in the first element of $x$. This dictionary imposes the same additive structure that $\rho_0(X)$ has while allowing the function $g_t(\cdot)$ to vary with $t$. This parameter was not considered as a linear functional of an approximately sparse specification in \cite{chernozhukov2022locally} or \cite{hirshberg2021augmented}.
The estimator given here will be root-n consistent without the doubly robust rate condition.
 }

\end{example}

\section{Discussion}

In this paper, we focus on the functional class of approximate sparsity, formulated as the error of sparse approximations decaying at a polynomial rate of the sparsity. This is motivated by classical nonparametric problems and can be viewed as series approximations with unordered basis functions, rather than the usual ordered basis functions (such as splines, Fourier basis). By providing the minimax lower and upper bounds, we show that the lack of order for the basis functions has a cost in terms of convergence rate for learning linear functionals. A proposed estimator  applies to the general problem of learning linear functionals of the conditional mean function in the approximately sparse class. Our work is related to several areas of future research. 

For example, how to construct estimators for general linear functionals under flexible conditions. Whereas our proposed estimator achieves the parametric rate under weaker conditions than the existing literature (such as the typical rate double-robustness condition), it is still an open question how to characterize the minimax optimality of estimators of general linear functionals. For a given linear functional $L$, we can apply the Riesz representer framework and build an estimator but it is unknown whether the difficulty of learning $L(\rho_0)$ ($\rho_0$ being the conditional mean function) is driven by the smoothness or approximate sparsity of the Riesz representer. 

  Moreover, one motivation for allowing approximate sparsity is to explicitly account for the approximation errors in high-dimensional sparse models, which connects to the literature on model misspecification. That literature often has broader aims—such as permitting one of the underlying models to be entirely wrong—and therefore demands new estimators for the nuisance parameters to address biases that standard loss functions fail to remove. In contrast, our estimators still rely on classical loss functions, leaving it an open question whether our approach can tolerate cases where one of the two function spaces is incorrectly specified (e.g., the wrong sparsity) or where the approximation error is unacceptably large.

\bigskip


\bibliographystyle{apalike}
\bibliography{SC_biblio}

\pagebreak

\pagenumbering{arabic}    
\setcounter{page}{1}      

\begin{center}
    {\Large \scshape Appendix to ``Minimax Semiparametric Learning With Approximate Sparsity''}\\[3ex] 
    {\large Jelena Bradic \quad Victor Chernozhukov \quad Whitney Newey \quad Yinchu Zhu}
\end{center}

\begin{appendix}


\section{Proof of Lemma \ref{lem: strict vs approx sparse}}

\begin{proof}[Proof of Lemma \ref{lem: strict vs approx sparse}]
We proceed by contradiction. Suppose that such $D,\xi$ exist. Consider $v\in\mathbb{R}^p$ with $v_1=v_2=\cdots=v_k=\sqrt{C/k} $ and $v_j=0$ for $j\geq k+1$. Clearly, $v\in\Theta(C,k)$. Consider $s=\left \lceil{k^{1/2}}\right \rceil$. Then $s\in \{1,...,k\}$ for large $k$. By  $\mathbb{E}[\psi(X)\psi(X)^{\top}]=I_p$ and the definition of  $  \mathcal{M}_{D,\xi}$, we have $(C/k)(k-s)\leq D^2 s^{-2\xi}$, which means that $1-s/k\leq (D^2/C) s^{-2\xi}$. This is impossible because the left-hand side tends to 1 and the right-hand side tends to zero as $k\rightarrow \infty$.
\end{proof}

\section{Proof of Lemma \ref{lem: Sobolev approx}}
\begin{proof}
\jelenax{We consider the periodic Sobolev space. Let $\Omega_{*}=[0,1]^{D}$
and 
\[
B_{H^{s}(\Omega_{*})}(Q_{*})=\left\{ f:\Omega_{*}\mapsto\mathbb{R}:\|f\|_{H^{s}(\Omega_{*})}\leq Q_{*}\right\} .
\]
}

\jelenax{Now we extend the Sobolev function $f\in H^{s}(\Omega)$ to a function
in $H^{s}(\mathbb{R}^{D})$. Notice that the boundary of $\Omega$
is minimally smooth as defined on page 189 of \citet{stein1970singular};
$\partial\Omega$ is minimally smooth since $\Omega$ is an open bounded
set. By Theorem 5 on page 181 of \citet{stein1970singular}, there
exists a function $f_{*}\in H^{s}(\mathbb{R}^{D})$ such that $f_{*}(x)=f(x)$
for any $x\in\Omega$ and 
\begin{equation}
\|f_{*}\|_{H^{s}(\mathbb{R}^{D})}\leq c_{1}\|f\|_{H^{s}(\Omega)},\label{eq: sobolev eq 2}
\end{equation}
where $c_{1}>0$ is a universal constant that depends only on $D,s,\Omega$. 
}

\jelenax{We now construct an extension of $f$ that is supported on $\Omega_{*}$.
To do so, let $\eta(\cdot)$ be an infinitely smooth function such
that $\eta$ is supported on ${\rm Int}(\Omega_{*})$ and $\eta(x)=1$
for $x\in\Omega$. This is possible by a smooth version of Urysohn's
lemma, e.g., page 38 of \citet{lieb2001analysis}. Then we define
$f_{**}(x)=f_{*}(x)\cdot\eta(x)$. Clearly, $f_{**}(x)=f(x)$ for
$x\in\Omega$ and $f_{**}$ is supported on ${\rm Int}(\Omega_{*})$. 
}

\jelenax{We now bound $\|f_{**}\|_{H^{s}(\Omega_{*})}$. By the product rule
of differentiation, $D^{\alpha}(f_{**})=D^{\alpha}(f_{*}\eta)=\sum_{q}(D^{q}f_{*})\cdot(D^{\alpha-q}\eta)$,
where the summation is taken over all $q\in\mathbb{N}_{0}^{D}$ such
that $0\leq q_{l}\leq a_{l}$ for all $l\in\{1,...,D\}$. Since $\eta(\cdot)$
is a fixed smooth function with compact support, $\sup_{x\in\mathbb{R}^{D}}|D^{\alpha-q}\eta(x)|$
is bounded by a constant that depends only on $\alpha-q$. Therefore,
there exist a constant $c_{2}>0$ depending only on $\alpha$ such
that $\|D^{\alpha}(f_{**})\|_{L^{2}(\Omega_{*})}\leq c_{2}\sum_{q}\|D^{q}f_{*}\|_{L^{2}(\Omega_{*})}$.
It follows that 
\[
\|f_{**}\|_{H^{s}(\Omega_{*})}=\sum_{|\alpha|\leq s}\|D^{\alpha}f_{**}\|_{L^{2}(\Omega_{*})}\leq c_{3}\sum_{|\alpha|\leq s}\|D^{\alpha}f_{*}\|_{L^{2}(\Omega_{*})}=c_{3}\|f_{*}\|_{H^{s}(\Omega_{*})}\leq c_{3}\|f_{*}\|_{H^{s}(\mathbb{R}^{D})}.
\]}

\jelenax{Thus, by (\ref{eq: sobolev eq 2}) and $\|f\|_{H^{s}(\Omega)}\leq Q$,
we have that 
\[
\|f_{**}\|_{H^{s}(\Omega_{*})}\leq c_{1}c_{3}Q.
\]}

\jelenax{We now study the approximation property of $f_{**}$. }

\jelenax{For any function $g\in L^{2}(\Omega_{*})$ and any $k\in\mathbb{Z}^{D}$,
define $s_{k}(g)=\int_{\Omega_{*}}g(x)\exp(2\pi\iota k'x)dx$, where
$\iota=\sqrt{-1}$. Let $\psi_{k}(x)=\exp(2\pi\iota k'x)$. We know
that 
\[
g=\sum_{k\in\mathbb{Z}^{D}}s_{k}(g)\cdot\psi_{k}\quad\text{and}\quad\|g\|_{L^{2}(\Omega_{*})}^{2}=\sum_{k\in\mathbb{Z}^{D}}|s_{k}(g)|^{2}.
\]}

\jelenax{We apply integration by parts and obtain that for $k_{l}\neq0$ and
$\alpha_{l}\geq1$,
\begin{align*}
 & \int_{0}^{1}D^{\alpha}f_{**}(x)\exp(2\pi\iota k_{l}x_{l})dx_{l}\\
 & =\int_{0}^{1}\exp(2\pi\iota k_{l}x_{l})dD^{\alpha-e_{l}}f_{**}(x)\\
 & =D^{\alpha-e_{l}}f_{**}(x)\cdot\exp(2\pi\iota k_{l}x_{l})\Bigg|_{0}^{1}-(2\pi\iota k_{l})\int_{0}^{1}\exp(2\pi\iota k_{l}x_{l})\cdot D^{\alpha-e_{l}}f(x)dx_{l}\\
 & \overset{\text{(i)}}{=}-(2\pi\iota k_{l})\int_{0}^{1}\exp(2\pi\iota k_{l}x_{l})\cdot D^{\alpha-e_{l}}f(x)dx_{l},
\end{align*}
where (i) follows by the fact that $f_{**}$ is supported on ${\rm Int}(\Omega_{*})$.
For $k_{l}=0$, $\int_{0}^{1}D^{\alpha}f_{**}(x)\exp(2\pi\iota k_{l}x_{l})dx_{l}=\int_{0}^{1}D^{\alpha}f_{**}(x)dx_{l}=0$
since $f_{**}$ is supported on ${\rm Int}(\Omega_{*})$. Thus, for
$a_{l}\geq1$, we have 
\[
\int_{0}^{1}D^{\alpha}f_{**}(x)\exp(2\pi\iota k_{l}x_{l})dx_{l}=-(2\pi\iota k_{l})\int_{0}^{1}\exp(2\pi\iota k_{l}x_{l})\cdot D^{\alpha-e_{l}}f(x)dx_{l}.
\]}

\jelenax{Therefore, for any $\alpha_{l}\geq1$, 
\begin{align*}
s_{k}\left(D^{\alpha}f_{**}\right) & =\int_{\Omega_{*}}D^{\alpha}f_{**}(x)\exp(2\pi\iota k'x)dx\\
 & =\int_{[0,1]^{D-1}}\left[\exp\left(2\pi\iota\sum_{j\neq l}k_{j}x_{j}\right)\left(\int_{0}^{1}D^{\alpha}f_{**}(x)\exp(2\pi\iota k_{l}x_{l})dx_{l}\right)\right]\prod_{j\neq l}dx_{j}\\
 & =-(2\pi\iota k_{l})\int_{[0,1]^{D-1}}\left[\exp\left(2\pi\iota\sum_{j\neq l}k_{j}x_{j}\right)\left(\int_{0}^{1}\exp(2\pi\iota k_{l}x_{l})\cdot D^{\alpha-e_{l}}f(x)dx_{l}\right)\right]\prod_{j\neq l}dx_{j}\\
 & =-(2\pi\iota k_{l})\int_{\Omega_{*}}D^{\alpha-e_{l}}f_{**}(x)\exp(2\pi\iota k'x)dx=-(2\pi\iota k_{l})s_{k}\left(D^{\alpha-e_{l}}f_{**}\right).
\end{align*}}

\jelenax{By induction, for any $\alpha_{l}\geq1$, 
\[
s_{k}\left(D^{\alpha}f_{**}\right)=(-2\pi\iota k_{l})^{a_{l}}s_{k}\left(D^{\alpha-\alpha_{l}e_{l}}f_{**}\right)
\]
and thus 
\[
s_{k}\left(D^{\alpha}f_{**}\right)=\prod_{l=1}^{D}(-2\pi\iota k_{l})^{a_{l}}s_{k}\left(f_{**}\right).
\]}

\jelenax{Therefore, 
\[
\|D^{\alpha}f_{**}\|_{L^{2}(\Omega_{*})}^{2}=\sum_{k\in\mathbb{Z}^{D}}\left|s_{k}\left(D^{\alpha}f_{**}\right)\right|^{2}=\sum_{k\in\mathbb{Z}^{D}}\prod_{l=1}^{D}|2\pi k_{l}|^{2a_{l}}\cdot\left|s_{k}\left(f_{**}\right)\right|^{2}.
\]}

\jelenax{Notice that $\|D^{\alpha}f_{**}\|_{L^{2}(\Omega_{*})}^{2}\leq c_{1}^{2}c_{3}^{2}Q^{2}$
for any $\alpha$ with $|\alpha|\leq s$. For any $b\in\{1,...,D\}$,
we choose $\alpha=se_{b}$. Thus, for any $b\in\{1,...,D\}$,
\[
c_{1}^{2}c_{3}^{2}Q^{2}\geq\|D^{se_{b}}f_{**}\|_{L^{2}(\Omega_{*})}^{2}=\sum_{k\in\mathbb{Z}^{D}}|2\pi k_{b}|^{2s}\cdot\left|s_{k}\left(f_{**}\right)\right|^{2}=(2\pi)^{2s}\cdot\sum_{k\in\mathbb{Z}^{D}}|k_{b}|^{2s}\cdot\left|s_{k}\left(f_{**}\right)\right|^{2}.
\]}

\jelenax{For any $b\in\{1,...,D\}$, let $k_{-b}\in\mathbb{Z}^{D-1}$ denote
the elements of $k=(k_{1},...,k_{D})'$ other than $k_{b}$. For $k$
such that $k_{b}=t$, define $m_{b}(t)=\sum_{k_{-b}\in\mathbb{Z}^{D-1}}|s_{k}(f_{**})|^{2}$.
Then the above means that for any $b\in\{1,...,D\}$, 
\[
c_{1}^{2}c_{3}^{2}Q^{2}\geq(2\pi)^{2s}\cdot\sum_{k_{b}\in\mathbb{Z}}|k_{b}|^{2s}\cdot\sum_{k_{-b}\in\mathbb{Z}^{D-1}}\left|s_{k}\left(f_{**}\right)\right|^{2}=(2\pi)^{2s}\cdot\sum_{k_{b}\in\mathbb{Z}}|k_{b}|^{2s}m_{b}(k_{b})=(2\pi)^{2s}\cdot\sum_{j\in\mathbb{Z}}|j|^{2s}m_{b}(j).
\]}

\jelenax{Hence, for any $r\geq1$ and for any $b\in\{1,...,D\}$,
\begin{multline}
\sum_{j\in\mathbb{Z},|j|\geq r+1}m_{b}(j)\leq\sum_{j\in\mathbb{Z},|j|\geq t+1}\left|\frac{j}{r+1}\right|^{2s}m_{b}(j)=(r+1)^{-2s}\sum_{j\in\mathbb{Z},|j|\geq r+1}|j|^{2s}m_{b}(j)\\
\leq(r+1)^{-2s}\sum_{j\in\mathbb{Z}}|j|^{2s}m_{b}(j)\leq(r+1)^{-2s}\cdot(2\pi)^{-2s}c_{1}^{2}c_{3}^{2}Q^{2}.\label{eq: sobolev eq 5}
\end{multline}}

\jelenax{For $t\geq1$, let $\mathcal{N}_{t}=\{j\in\mathbb{Z}:|j|\leq\left\lfloor t^{1/D}\right\rfloor \}$.
We observe that 
\begin{multline*}
\left(\mathbb{Z}^{D}\right)\backslash\left(\mathcal{N}_{t}^{D}\right)=\left\{ k\in\mathbb{Z}^{D}:|k_{b}|\geq\left\lfloor t^{1/D}\right\rfloor +1\ \text{for some }b\in\{1,...,D\}\right\} \\
=\bigcup_{b=1}^{D}\left\{ k\in\mathbb{Z}^{D}:|k_{b}|\geq\left\lfloor t^{1/D}\right\rfloor +1\right\} .
\end{multline*}}

\jelenax{It follows that 
\begin{align*}
\sum_{k\in\left(\mathbb{Z}^{D}\right)\backslash\left(\mathcal{N}_{t}^{D}\right)}\left|s_{k}\left(f_{**}\right)\right|^{2} & =\sum_{k\in\bigcup_{b=1}^{D}\left\{ k\in\mathbb{Z}^{D}:|k_{b}|\geq\left\lfloor t^{1/D}\right\rfloor +1\right\} }\left|s_{k}\left(f_{**}\right)\right|^{2}\\
 & \leq\sum_{b=1}^{D}\sum_{k\in\mathbb{Z}^{D}:|k_{b}|\geq\left\lfloor t^{1/D}\right\rfloor +1}\left|s_{k}\left(f_{**}\right)\right|^{2}\\
 & =\sum_{b=1}^{D}\sum_{|k_{b}|\geq\left\lfloor t^{1/D}\right\rfloor +1}\sum_{k_{-b}\in\mathbb{Z}^{D-1}}\left|s_{k}\left(f_{**}\right)\right|^{2}\\
 & =\sum_{b=1}^{D}\sum_{|k_{b}|\geq\left\lfloor t^{1/D}\right\rfloor +1}m_{b}(k_{b})\\
 & \overset{\text{(i)}}{\leq}D\left(\left\lfloor t^{1/D}\right\rfloor +1\right)^{-2s}\cdot(2\pi)^{-2s}c_{1}^{2}c_{3}^{2}Q^{2}\\
 & \overset{\text{(ii)}}{\leq}D\left(t^{1/D}\right)^{-2s}\cdot(2\pi)^{-2s}c_{1}^{2}c_{3}^{2}Q^{2}=(2\pi)^{-2s}c_{1}^{2}c_{3}^{2}Q^{2}D\cdot t^{-2s/D},
\end{align*}
where (i) follows by (\ref{eq: sobolev eq 5}) and (ii) follows by
$\left\lfloor t^{1/D}\right\rfloor +1\geq t^{1/D}$. Therefore, 
\begin{multline*}
\left\Vert f_{**}-\sum_{k\in\mathcal{N}_{t}^{D}}s_{k}(f_{**})\psi_{k}\right\Vert _{L^{2}(\Omega_{*})}^{2}=\left\Vert \sum_{k\in\left(\mathbb{Z}^{D}\right)\backslash\left(\mathcal{N}_{t}^{D}\right)}s_{k}(f_{**})\psi_{k}\right\Vert _{L^{2}(\Omega_{*})}^{2}\\
=\sum_{k\in\left(\mathbb{Z}^{D}\right)\backslash\left(\mathcal{N}_{t}^{D}\right)}\left|s_{k}\left(f_{**}\right)\right|^{2}\leq(2\pi)^{-2s}c_{1}^{2}c_{3}^{2}Q^{2}D\cdot t^{-2s/D}.
\end{multline*}}

\jelenax{Since $\Omega\subset\Omega_{*}$ and $f_{**}(x)=f(x)$ for $x\in\Omega$,
the desired result follows. }
\end{proof}

\section{Proofs for Section  \ref{sec lower bounds}}

In this section, we adopt the notation of $\tau(\lambda)=\theta$ for $ \lambda=(\theta,\beta,\phi,\sigma_{u}%
^{2},\sigma_{\varepsilon}^{2})$.

\subsection{Proof of Theorem  \ref{thm1}}
We prove a stronger result than Theorem  \ref{thm1} by deriving the optimal length of confidence intervals. 
Let $\mathcal{C}(\Lambda)$ be the set of $1-\alpha$ confidence intervals for
$\theta$ that are valid uniformly over $\lambda\in\Lambda$. A confidence interval $CI=[L,U]$ consists of  $L$ and $U$ measurable functions of the data $W_1,...,W_n$ and has length $|CI|=U-L$.  We are interested in how the minimax length of confidence interval in $\mathcal{C}(\Lambda_{\xi_1,\xi_2})$ depends on $(\xi_1,\xi_2)$. We also would like to investigate the adaptivity to $(\xi_1,\xi_2)$ in that whether knowledge of $(\xi_1,\xi_2)$ changes the efficiency of inference. To study these two questions simultaneous, we introduce
\[
\mathcal{L}(\Lambda_{\xi_1,\xi_2},\Lambda_{\tilde{\xi}_1,\tilde{\xi}_2})=\inf_{CI\in\mathcal{C}(\Lambda_{\tilde{\xi}_1,\tilde{\xi}_2}%
)}\sup_{\lambda\in \Lambda_{\xi_1,\xi_2}}\mathbb{E}_{\lambda}|CI|,
\] 
where $\xi_1\geq \tilde{\xi}_1$ and  $\xi_2\geq \tilde{\xi}_2$. Therefore, $\Lambda_{\xi_1,\xi_2}\subset \Lambda_{\tilde{\xi}_1,\tilde{\xi}_2}$. When $\xi_1= \tilde{\xi}_1$ and  $\xi_2 =\tilde{\xi}_2$, then the above measures the minimax length of confidence intervals in $\Lambda_{\xi_1,\xi_2}$. When  $\xi_1> \tilde{\xi}_1$ and  $\xi_2 >\tilde{\xi}_2$, the above
tells us the performance of confidence intervals that are valid in the enlarged class $\Lambda_{\tilde{\xi}_1,\tilde{\xi}_2}$. If a confidence interval that is valid in the enlarged class automatically achieves the same optimal length when the data distribution $\lambda$ is in the smaller $\Lambda_{\xi_1,\xi_2}$, then this confidence interval adapts to the true smoothness of the problem. If such a confidence interval does not exist, then we say that there is lack of adaptivity to $(\xi_1,\xi_2)$. Indeed, this is what we will show. It turns out that $\mathcal{L}(\Lambda_{\xi_1,\xi_2},\Lambda_{\tilde{\xi}_1,\tilde{\xi}_2})$ does not even depend on  $(\xi_1,\xi_2)$. Therefore,   if $\xi_1> \tilde{\xi}_1$ and  $\xi_2 >\tilde{\xi}_2$, then  $\Lambda_{\xi_1,\xi_2}$ is a very small subset of  $\Lambda_{\tilde{\xi}_1,\tilde{\xi}_2}$ but no confidence interval with validity on $\Lambda_{\tilde{\xi}_1,\tilde{\xi}_2}$  can automatically exploit the extra smoothness of the data distribution. Our result is in the following statement.

\begin{theorem}\label{old thm1}
If Assumption  \ref{assp 1} is satisfied and there exist 
constants $\kappa_1,\kappa_2>0$  such that $ \kappa_{1}\log p\leq n\leq\kappa_{2}p\log p$  for
large enough $n$,   then for any 
$$\xi_{1}\geq\tilde{\xi}_{1}\geq
0, \ \xi_{2}\geq\tilde{\xi}_{2}\geq0, \textit{ and } \tilde{\xi}%
=\max\{\tilde{\xi}_{1},\tilde{\xi}_{2}\}\leq1/2, $$%
\[
\mathcal{L}(\Lambda_{\xi_{1},\xi_{2}},\Lambda_{\tilde{\xi}_{1},\tilde{\xi}%
_{2}})
\geq C(n^{-1}\log p)^{2\tilde{\xi}/(2\tilde{\xi}+1)},
\]
where $C>0$  is a constant depending only on $\tilde{\xi}_{1}$\textit{, }$\tilde{\xi}_{2}$\textit{, }$\kappa_1,\kappa_2,\alpha,C_{0}$. 
\end{theorem}

 { Theorem \ref{old thm1} has two key implications. First, if  
$
\tilde{\xi} = \max\{\tilde{\xi}_{1},\tilde{\xi}_{2}\}
$
is significantly smaller than \(1/2\), then the rate  
$
\mathcal{L}\left(\Lambda_{\xi_{1},\xi_{2}},\Lambda_{\tilde{\xi}_{1},\tilde{\xi}_{2}}\right)
$
can be much slower than the standard parametric rate \(n^{-1/2}\). 
 Second, Theorem \ref{old thm1} shows that adaptivity to the sizes of \(\xi_{1}\) and \(\xi_{2}\) is impossible. Notice that the lower bound for  
\[
\mathcal{L}\left(\Lambda_{\xi_{1},\xi_{2}},\Lambda_{\tilde{\xi}_{1},\tilde{\xi}_{2}}\right)
\]
depends only on \(\max\{\tilde{\xi}_{1},\tilde{\xi}_{2}\}\) and not on the specific values of \((\xi_{1},\xi_{2})\). This implies that any confidence interval valid over \(\Lambda_{\tilde{\xi}_{1},\tilde{\xi}_{2}}\) with \(\max\{\tilde{\xi}_{1},\tilde{\xi}_{2}\}\le 1/2\) cannot achieve an expected width of \(n^{-1/2}\) even when restricted to a smaller subspace \(\Lambda_{\xi_{1},\xi_{2}}\), no matter how narrow that space is. In other words, no confidence interval can simultaneously (1) be valid over \(\Lambda_{\tilde{\xi}_{1},\tilde{\xi}_{2}}\) with \(\max\{\tilde{\xi}_{1},\tilde{\xi}_{2}\} < 1/2\) and (2) have an expected width of \(O(n^{-1/2})\) on a (possibly much smaller) subspace \(\Lambda_{\xi_{1},\xi_{2}}\). Consequently, it is impossible to distinguish between the cases \(\max\{\xi_{1},\xi_{2}\} < 1/2\) and \(\max\{\xi_{1},\xi_{2}\} > 1/2\) from the data. In particular, no test can yield a confidence interval of length \(1/\sqrt{n}\) when \(\max\{\xi_{1},\xi_{2}\} < 1/2\).
}
 {We emphasize that although \(\mathcal{L}(\Lambda_{\xi_{1},\xi_{2}},\Lambda_{\tilde{\xi}_{1},\tilde{\xi}_{2}})\) is defined as the expected confidence interval length, our rates are not driven by rare extreme values. The interval lengths are  bounded, so the expected length is not skewed by outliers.  }

\begin{lemma} 
\label{lem: hypergeometric}Let $k\in\mathbb{N}$ and define $\mathcal{Q}_{k}=\{v\in\{0,1\}^{p}:\ \|v\|_{0}=k\}$.
Let $v$ and $u$ be two independent vectors that have a uniform distribution
on $\mathcal{Q}_{k}$. Then for any $D\geq0$, 
\[
\mathbb{E}\exp\left(Du'v\right)<\exp\left(\exp\left(D+\log(k^{2}/p)\right)\right).
\]
\end{lemma}

\begin{proof}[Proof of Lemma \ref{lem: hypergeometric}]
Let $N=|\mathcal{Q}_{k}|$. We list elements in $\mathcal{Q}_{k}$, i.e., $\mathcal{Q}_{k}=\{x_{1},...,x_{N}\}$.
Then 
\[
\mathbb{E}\exp\left(Du'v\right)=N^{-2}\sum_{j_{2}=1}^{N}\sum_{j_{1}=1}^{N}\exp\left(Dx_{j_{1}}'x_{j_{2}}\right)\overset{\text{(i)}}{=}N^{-1}\sum_{j=1}^{N}\exp\left(Dx_{1}'x_{j}\right)=\mathbb{E}\exp(Dx_{1}'v),
\]
(i) follows by the observation that $\sum_{j_{1}=1}^{N}\exp(Dx_{j_{1}}'x_{j_{2}})$
does not depend on $j_{2}$. Without loss of generality, we take $x_{1}=(1,...,1,0,...,0)^{\top}$,
i.e., the vector whose first $k$ entries are nonzero. 

Let $\mathcal{C}_{n,k}$ be the population that consists of $n$ elements
with $n-k$ elements being 0 and the remaining $k$ being 1. Let $\{\xi_{i}\}_{i=1}^{k}$
be a random sample without replacement from the population of $\mathcal{C}_{n,k}$.
We observe that $x_{1}'v$ has the same distribution as $\sum_{i=1}^{k}\xi_{i}$.
Then 
\[
\mathbb{E}\exp(Dx_{1}'v)=\mathbb{E}\exp\left(D\sum_{i=1}^{k}\xi_{i}\right).
\]

Let $\{\zeta_{i}\}_{i=1}^{k}$ be a random sample with replacement
from $\mathcal{C}_{n,k}$. In other words, $\{\zeta_{i}\}_{i=1}^{k}$ is
i.i.d Bernoulli with $\mathbb{E}(\zeta_{i})=k/p$. Since $x\mapsto\exp(Dx)$
is a convex function, we can use Theorem 4 of \citet{hoeffding1963probability}
and obtain that 
\begin{multline*}
\mathbb{E}\exp\left(D\sum_{i=1}^{k}\xi_{i}\right)
\leq 
\mathbb{E}\exp\left(D\sum_{i=1}^{k}\zeta_{i}\right)
=\left[\mathbb{E}\exp(D\zeta_{1})\right]^{k}
\\
\overset{\text{(i)}}{=}
\left(1-\frac{k}{p}+\frac{k}{p}\exp(D)\right)^{k}\overset{\text{(ii)}}
{\leq}\exp\left(\frac{k^{2}}{p}[\exp(D)-1]\right)\\
<\exp\left(\frac{k^{2}}{p}\exp(D)\right)=\exp\left(\exp\left(D+\log(k^{2}/p)\right)\right) \hskip 90pt
\end{multline*}
where (i) follows by the moment generating function of Bernoulli distributions,
(ii) follows by the elementary inequality $1+x\leq\exp(x)$ for $x\geq0$.
The proof is complete. 
\end{proof}

\begin{lemma}
\label{lem: minimax part 1}Let $k$ be a positive integer and define
$\mathcal{Q}_{k}=\{v\in\{0,1\}^{p}:\ \|v\|_{0}=k\}$. Let $N=|\mathcal{Q}_{k}|$.
We list elements in $\mathcal{Q}_{k}$, i.e., $\mathcal{Q}_{k}=\{\delta_{1},...,\delta_{N}\}$.
Let 
\[
\Sigma_{j}=\begin{pmatrix}1 & 0 & q_{1}\delta_{j}'\\
0 & 1 & q_{2}\delta_{j}'\\
q_{1}\delta_{j} & q_{2}\delta_{j} & I_{p}
\end{pmatrix}
\]
and $\Sigma_{*}=I_{p+2}$. Let $P_{j}$ denote the distribution of
$N(0,\Sigma_{j})$ and $P_{*}$ denote the distribution of $N(0,I_{p+2})$.
If $n^{-1}k\log(p/k^{2})\leq3$ and $6n(q_{1}^{2}+q_{2}^{2})\leq\log(p/k^{2})$,
then 
\[
{\rm TV}\left(N^{-1}\sum_{j=1}^{N}P_{j},\ P_{*}\right)\leq\sqrt{\exp\left(\sqrt{k^{2}/p}\right)-1}.
\]
\end{lemma}

\begin{proof}[Proof of Lemma \ref{lem: minimax part 1}]
Let $\mathbb{E}_{*}$ denotes the expectation under $P_{*}$. By Lemma 3 in
\citet{cai2017confidence}, we have that 
\begin{equation}
\mathbb{E}_{*}\frac{dP_{j_{1}}}{dP_{*}}\frac{dP_{j_{2}}}{dP_{*}}=\left[\det\left(I_{p+2}-(\Sigma_{*}^{-1}\Sigma_{j_{1}}-I_{p+2})(\Sigma_{*}^{-1}\Sigma_{j_{2}}-I_{p+2})\right)\right]^{-n/2}.\label{eq: lower bnd 2-1}
\end{equation}

By the definitions of $\Sigma_{j}$ and $\Sigma_{*}$, we have that
\[
\Sigma_{*}^{-1}\Sigma_{j}-I_{p+2}=\begin{pmatrix}0 & 0 & q_{1}\delta_{j}'\\
0 & 0 & q_{2}\delta_{j}'\\
q_{1}\delta_{j} & q_{2}\delta_{j} & 0
\end{pmatrix}.
\]

Then we have 
\begin{align*}
 & I_{p+2}-(\Sigma_{*}^{-1}\Sigma_{j_{1}}-I_{p+2})(\Sigma_{*}^{-1}\Sigma_{j_{2}}-I_{p+2})\\
 & =I_{p+2}-\begin{pmatrix}0 & 0 & q_{1}\delta_{j_{1}}'\\
0 & 0 & q_{2}\delta_{j_{1}}'\\
q_{1}\delta_{j_{1}} & q_{2}\delta_{j_{1}} & 0
\end{pmatrix}\begin{pmatrix}0 & 0 & q_{1}\delta_{j_{2}}'\\
0 & 0 & q_{2}\delta_{j_{2}}'\\
q_{1}\delta_{j_{2}} & q_{2}\delta_{j_{2}} & 0
\end{pmatrix}\\
 & =I_{p+2}-\begin{pmatrix}q_{1}^{2}\delta_{j_{1}}'\delta_{j_{2}} & q_{1}q_{2}\delta_{j_{1}}'\delta_{j_{2}} & 0\\
q_{1}q_{2}\delta_{j_{1}}'\delta_{j_{2}} & q_{2}^{2}\delta_{j_{1}}'\delta_{j_{2}} & 0\\
0 & 0 & (q_{1}^{2}+q_{2}^{2})\delta_{j_{1}}\delta_{j_{2}}'
\end{pmatrix}\\
 & =\begin{pmatrix}1-q_{1}^{2}\delta_{j_{1}}'\delta_{j_{2}} & -q_{1}q_{2}\delta_{j_{1}}'\delta_{j_{2}} & 0\\
-q_{1}q_{2}\delta_{j_{1}}'\delta_{j_{2}} & 1-q_{2}^{2}\delta_{j_{1}}'\delta_{j_{2}} & 0\\
0 & 0 & I_{p}-(q_{1}^{2}+q_{2}^{2})\delta_{j_{1}}\delta_{j_{2}}'
\end{pmatrix}.
\end{align*}

Therefore, 
\begin{align*}
 & \det\left[I_{p+2}-(\Sigma_{*}^{-1}\Sigma_{j_{1}}-I_{p+2})(\Sigma_{*}^{-1}\Sigma_{j_{2}}-I_{p+2})\right]\\
 & =\det\left(I_{p}-(q_{1}^{2}+q_{2}^{2})\delta_{j_{1}}\delta_{j_{2}}'\right)\times\det\begin{pmatrix}1-q_{1}^{2}\delta_{j_{1}}'\delta_{j_{2}} & -q_{1}q_{2}\delta_{j_{1}}'\delta_{j_{2}}\\
-q_{1}q_{2}\delta_{j_{1}}'\delta_{j_{2}} & 1-q_{2}^{2}\delta_{j_{1}}'\delta_{j_{2}}
\end{pmatrix}\\
 & \overset{\text{(i)}}{=}\left(1-(q_{1}^{2}+q_{2}^{2})\delta_{j_{1}}'\delta_{j_{2}}\right)\times\det\begin{pmatrix}1-q_{1}^{2}\delta_{j_{1}}'\delta_{j_{2}} & -q_{1}q_{2}\delta_{j_{1}}'\delta_{j_{2}}\\
-q_{1}q_{2}\delta_{j_{1}}'\delta_{j_{2}} & 1-q_{2}^{2}\delta_{j_{1}}'\delta_{j_{2}}
\end{pmatrix}\\
 & =\left(1-(q_{1}^{2}+q_{2}^{2})\delta_{j_{1}}'\delta_{j_{2}}\right)\times\left\{ \left(1-q_{1}^{2}\delta_{j_{1}}'\delta_{j_{2}}\right)\left(1-q_{2}^{2}\delta_{j_{1}}'\delta_{j_{2}}\right)-\left(-q_{1}q_{2}\delta_{j_{1}}'\delta_{j_{2}}\right)^{2}\right\} \\
 & =\left(1-(q_{1}^{2}+q_{2}^{2})\delta_{j_{1}}'\delta_{j_{2}}\right)^{2},
\end{align*}
where (i) follows by Sylvester's determinant identity. By (\ref{eq: lower bnd 2-1}),
we have that 
\begin{equation}
\mathbb{E}_{*}\frac{dP_{j_{1}}}{dP_{*}}\frac{dP_{j_{2}}}{dP_{*}}=\left(1-(q_{1}^{2}+q_{2}^{2})\delta_{j_{1}}'\delta_{j_{2}}\right)^{-n}\overset{\text{(i)}}{\leq}\exp\left(3n(q_{1}^{2}+q_{2}^{2})\delta_{j_{1}}'\delta_{j_{2}}\right),\label{eq: lower bnd 3-1}
\end{equation}
where (i) follows by $(q_{1}^{2}+q_{2}^{2})\delta_{j_{1}}'\delta_{j_{2}}j\leq(q_{1}^{2}+q_{2}^{2})k\leq1/2$
and the fact that $(1-x)^{-n}<\exp(3xn)$ for any $x\in[0,1/2]$.

To see the above, define $f(x)=-3x-\log(1-x)$. Notice that $f(\cdot)$
is convex on $[0,1/2]$ by checking $f''(\cdot)$. Also notice that
$f(0)<0$ and $f(1/2)<0$. Hence, $f(x)<0$ on $[0,1/2]$. This means
$-\log(1-x)\leq3x$. Multiplying both sides by $n$ and taking exponential,
we obtain $(1-x)^{-n}\leq\exp(3xn)$.

By (\ref{eq: lower bnd 3-1}), we have
\begin{align}
\mathbb{E}_{*}\left(N^{-1}\sum_{j=1}^{N}\frac{dP_{j}}{dP_{*}}-1\right)^{2} & =N^{-2}\sum_{j_{2}=1}^{N}\sum_{j_{1}=1}^{N}\mathbb{E}_{*}\frac{dP_{j_{1}}}{dP_{*}}\frac{dP_{j_{2}}}{dP_{*}}-1
\\
&\leq N^{-2}\sum_{j_{2}=1}^{N}\sum_{j_{1}=1}^{N}\exp\left(3n(q_{1}^{2}+q_{2}^{2})\delta_{j_{1}}'\delta_{j_{2}}\right)-1\nonumber \\
 & \overset{\text{(i)}}{<}\exp\left(\exp\left(3n(q_{1}^{2}+q_{2}^{2})+\log(k^{2}/p)\right)\right)-1 \nonumber \\
 & \overset{\text{(ii)}}{=}\exp\left(\exp\left((1/2)\log\frac{k^{2}}{p}\right)\right)-1=\exp\left(\sqrt{k^{2}/p}\right)-1,\label{eq: lower bnd 4-1}
\end{align}
where (i) follows by Lemma \ref{lem: hypergeometric} and (ii) follows
by $q_{1}^{2}+q_{2}^{2}\leq(6n)^{-1}\log(p/k^{2})$. The desired result
follows by Equation (2.27) of \citet{Tsybakov2009}.
\end{proof}

\begin{lemma}
\label{lem: minimax part 2}Suppose that $q\neq0$ and $\delta\in\{0,1\}^{p}$
with $\|\delta\|_{0}=k$. Let $\xi>0$. If $|q|k^{\xi+1/2}\leq\sqrt{C}$,
then $(k-t)q^{2}\leq Ct^{-2\xi}$ for $1\leq t\leq k$. 
\end{lemma}

\begin{proof}[Proof of Lemma \ref{lem: minimax part 2}]
We shall show that for any $0\leq t\leq k$, we have
\[
kt^{2\xi}-t^{2\xi+1}\leq Cq^{-2}.
\]

The maximum of this mapping on $[0,k]$ is either at the end points
or at an interior point. Taking the first-order condition of the mapping
$t\mapsto kt^{2\xi}-t^{2\xi+1}$ and setting it zero, we have that
$t=\alpha k$ with $\alpha=2\xi/(2\xi+1)$, which corresponds to a
function value of $k^{2\xi+1}(\alpha^{2\xi}-\alpha^{2\xi+1})$. Notice
that the requirement of $k^{2\xi+1}(\alpha^{2\xi}-\alpha^{2\xi+1})\leq Cq^{-2}$
is satisfied by $k^{2\xi+1}\leq Cq^{-2}$ since $\alpha\in(0,1)$.
Clearly, the value of the mapping $t\mapsto kt^{2\xi}-t^{2\xi+1}$
at the end points ($t=0$ and $t=k$) is less than $Cq^{-2}$. The
desired result follows.
\end{proof}

\begin{proof}[Proof of Theorems  \ref{thm1} and  \ref{old thm1}]Let 
\[
k=\left\lfloor c_{0}\left(n/\log p\right)^{1/(2\tilde{\xi}+1)}\right\rfloor ,
\]
where $c_{0}>0$ is a constant to be chosen. Let $q_{n}=c_{1}\sqrt{n^{-1}\log p}$
and $c_{1}>0$ is a constant to be chosen.

Define $\mathcal{Q}_{k}=\{v\in\{0,1\}^{p}:\ \|v\|_{0}=k\}$. Let $N=|\mathcal{Q}_{k}|$.
Clearly, $N={p \choose k}$. We list elements in $\mathcal{Q}_{k}$, i.e.,
$\mathcal{Q}_{k}=\{\delta_{1},...,\delta_{N}\}$. For $1\leq j\leq N$,
define $\beta_{j}=-2q_{n}\delta_{j}$ and $\phi_{j}=q_{n}\delta_{j}$.
Clearly, we can choose $c_{0},c_{1}$ such that 
\[
2q_{n}k^{\tilde{\xi}+1/2}\leq2c_{1}c_{0}^{\tilde{\xi}+1/2}\leq C_{0}\qquad{\rm and}\qquad\exp(c_{0}^{2}\kappa_{2})-1\leq\alpha^{2}.
\]

By Lemma \ref{lem: minimax part 2}, this implies that $-2q_{n}\delta_{j},q_{n}\delta_{j}\in\mathcal{M}_{C_{0},\tilde{\xi}}$.
We define $\lambda_{*}=(0,0,0,1,1)$ and for $1\leq j\leq N$, $\lambda_{j}=(2q_{n}^{2}k,-2q_{n}\delta_{j},q_{n}\delta_{j},1-q_{n}^{2}k,\sigma_{\varepsilon,1}^{2})$
with $\sigma_{\varepsilon,1}^{2}=1-4q_{n}^{2}k(1-q_{n}^{2}k)^{2}-4q_{n}^{4}k^{2}(1-q_{n}^{2}k)$.
Notice that 
\[
q_{n}^{2}k\leq c_{1}^{2}c_{0}\left(n/\log p\right)^{-1+1/(2\tilde{\xi}+1)}=c_{1}^{2}c_{0}\left(n/\log p\right)^{-2\tilde{\xi}/(2\tilde{\xi}+1)}\leq c_{1}^{2}c_{0}\kappa_{1}^{-2\tilde{\xi}/(2\tilde{\xi}+1)}.
\]

Then we can choose $c_{1},c_{0}$ small enough such that $1-q_{n}^{2}k,\sigma_{\varepsilon,1}^{2}\in[1/2,2]\subseteq[M^{-1},M]$
(due to $M\geq2$). Hence, $\lambda_{*},\lambda_{j}\in\Lambda_{\tilde{\xi}_{1},\tilde{\xi}_{2}}$. 

Then under $P_{\lambda_{*}}$, the distribution of $(Y_{i},Z_{i},X_{i})$
is $N(0,I_{p+2})$. Under $P_{\lambda_{j}}$, $Y_{i}=Z_{i}\cdot(2q_{n}^{2}k)-2q_{n}X_{i}'\delta_{j}+\varepsilon_{i}$
and $Z_{i}=q_{n}X_{i}'\delta_{j}+u_{i}$, where $\mathbb{E}\varepsilon_{i}^{2}=\sigma_{\varepsilon,1}^{2}$
and $Eu_{i}^{2}=1-q_{n}^{2}k$. After simple calculations, we have
that under $P_{\lambda_{j}}$, the distribution of $(Y_{i},Z_{i},X_{i})$
is $N(0,\Sigma_{j})$ with
\[
\Sigma_{j}=\begin{pmatrix}1 & 0 & -q_{n}\delta_{j}'\\
0 & 1 & q_{n}\delta_{j}'\\
-q_{n}\delta_{j} & q_{n}\delta_{j} & I_{p}
\end{pmatrix}.
\]

Clearly, we can choose $c_{0},c_{1}$ small enough such that
\begin{eqnarray*}
n^{-1}k\log(p/k^{2})
&\leq n^{-1}k\log p\leq c_{0}n^{-1}\left(n/\log p\right)^{1/(2\tilde{\xi}+1)}\log p
\\
&=c_{0}\left(n/\log p\right)^{-2\tilde{\xi}/(2\tilde{\xi}+1)}\leq c_{0}\kappa_{1}^{-2\tilde{\xi}/(2\tilde{\xi}+1)}\leq3
\end{eqnarray*}
and 
\[
\frac{6n(q_{n}^{2}+q_{n}^{2})}{\log(p/k^{2})}\leq\frac{6n(q_{n}^{2}+q_{n}^{2})}{\log(p)}=12c_{1}^{2}\leq1.
\]

Hence, the assumptions of Lemma \ref{lem: minimax part 1} hold. Therefore,
\begin{equation}
{\rm TV}\left(N^{-1}\sum_{j=1}^{N}P_{\lambda_{j}},\ P_{\lambda_{*}}\right)\leq\sqrt{\exp\left(\sqrt{k^{2}/p}\right)-1}.\label{eq: lower bnd 4}
\end{equation}

Let $CI_{n}=[l_{n},u_{n}]$ be an arbitrary confidence interval for
$\theta=\tau(\lambda)$ with nomimal coverage probability $1-\alpha$
on $\Lambda_{\tilde{\xi}_{1},\tilde{\xi}_{2}}$. In other words, 
\begin{equation}
\inf_{\lambda\in\Lambda_{\tilde{\xi}_{1},\tilde{\xi}_{2}}}P_{\lambda}\left(l_{n}\leq\tau(\lambda)\leq u_{n}\right)=\inf_{\lambda\in\Lambda_{\tilde{\xi}_{1},\tilde{\xi}_{2}}}P_{\lambda}\left(\tau(\lambda)\in CI_{n}\right)\geq1-\alpha.\label{eq: lower bnd 4.5}
\end{equation}

We now define the random variable
\[
\Psi_{n}=\mathbf{1}\left\{ 2q_{n}^{2}k\in CI_{n}\right\} .
\]

Thus, by (\ref{eq: lower bnd 4}), $n\leq\kappa_{2}p\log p$ and $\exp(c_{0}^{2}\kappa_{2})-1\leq\alpha^{2}$,
we have
\begin{eqnarray*}
\left|N^{-1}\sum_{j=1}^{N}\mathbb{E}_{\lambda_{j}}\Psi_{n}-\mathbb{E}_{\lambda_{*}}\Psi_{n}\right|
&\leq&{\rm TV}\left(N^{-1}\sum_{j=1}^{N}P_{\lambda_{j}},\ P_{\lambda_{*}}\right)\leq\sqrt{\exp\left(\sqrt{k^{2}/p}\right)-1}\\
&\leq&\sqrt{\exp\left(\frac{c_{0}^{2}n}{p\log p}\right)-1}\leq\sqrt{\exp(c_{0}^{2}\kappa_{2})-1}\leq\alpha.
\end{eqnarray*}

Therefore, 
\begin{eqnarray*}
P_{\lambda_{*}}\left(2q_{n}^{2}k\in CI_{n}\right)
&=&\mathbb{E}_{\lambda_{*}}\Psi_{n}\geq N^{-1}\sum_{j=1}^{N}\mathbb{E}_{\lambda_{j}}\Psi_{n}-\alpha=N^{-1}\sum_{j=1}^{N}P_{\lambda_{j}}\left(2q_{n}^{2}k\in CI_{n}\right)-\alpha\\
&\overset{\text{(i)}}{=}&N^{-1}\sum_{j=1}^{N}P_{\lambda_{j}}\left(\tau(\lambda_{j})\in CI_{n}\right)-\alpha\overset{\text{(ii)}}{\geq}1-2\alpha,
\end{eqnarray*}
where (i) follows by $\tau(\lambda_{j})=2q_{n}^{2}k$ and (ii) follows
by (\ref{eq: lower bnd 4.5}). On the other hand, by $\tau(\lambda_{*})=0$
and (\ref{eq: lower bnd 4.5}), we have 
\[
P_{\lambda_{*}}\left(0\in CI_{n}\right)\geq1-\alpha.
\]

The above two displays imply that $P_{\lambda_{*}}\left(\{0,2q_{n}^{2}k\}\subset CI_{n}\right)\geq1-3\alpha$.
Since $CI_{n}$ is an interval, the event of $\{0,2q_{n}^{2}k\}\subset CI_{n}$
is the same as the event $[0,2q_{n}^{2}k]\subset CI_{n}$. Hence,
we have $P_{\lambda_{*}}\left([0,2q_{n}^{2}k]\subset CI_{n}\right)\geq1-3\alpha$,
which means
\[
\frac{\mathbb{E}_{\lambda_{*}}|CI_{n}|}{2q_{n}^{2}k}\geq P_{\lambda_{*}}\left(|CI_{n}|\geq2q_{n}^{2}k\right)\geq1-3\alpha.
\]

Since $2q_{n}^{2}k\asymp(n/\log p)^{-2\tilde{\xi}/(2\tilde{\xi}+1)}$ (by the definitions
of $q_{n}$ and $k$), we have 
$$
\mathbb{E}_{\lambda_{*}}|CI_{n}| \gtrsim (n^{-1} \log p )^{2\tilde{\xi}/(2\tilde{\xi}+1)}.
$$

This proves Theorem \ref{old thm1}. 

Let $a=\mathcal{R}_{\xi_{1},\xi_{2}}$. Then there exists an estimator $\hat{\theta}$ such that $\mathbb{E}_{\lambda} | \hat{\theta}(W)-\theta|\leq 2a$ for any $\lambda=(\theta,\beta,\phi,\sigma_{\varepsilon}^{2},\sigma_{u}^{2})\in \Lambda(\xi_{1},\xi_{2})$. Consider the confidence interval $CI=[\hat{\theta}-2a/\alpha,\hat{\theta}+2a/\alpha]$. Then for any $\lambda=(\theta,\beta,\phi,\sigma_{\varepsilon}^{2},\sigma_{u}^{2})\in \Lambda(\xi_{1},\xi_{2})$,
$$
P_{\lambda} \left( \theta\in CI\right)=P_{\lambda} \left( |\hat{\theta}-\theta|>2a/\alpha \right)\leq \mathbb{E}_{\lambda} | \hat{\theta}(W)-\theta|/(2a/\alpha)\leq 2a/(2a/\alpha)=\alpha.
$$
Therefore, $CI$ is a confidence interval that has coverage probability $1-\alpha$. Thus, $|CI|\gtrsim  (n^{-1} \log p )^{2\tilde{\xi}/(2\tilde{\xi}+1)}$. On the other hand, $|CI|=4a/\alpha$, which means that $a \gtrsim (n^{-1} \log p )^{2\tilde{\xi}/(2\tilde{\xi}+1)}$. This proves Theorem  \ref{thm1}. 
\end{proof}

\subsection{Proof of Theorem  \ref{thm2}}

We next give several Lemmas that are useful in the proof of Theorem \ref{thm2}.

\bigskip


\begin{lemma}\label{lem: A3}
 Let $C_{0},C_{1},C_{2},\xi_{1},\xi_{2}>0$ 
be constants. Define $k=\left\lfloor C_{0}n^{1/(\xi_{1}+\xi_{2}%
+1/2)}\right\rfloor $, $c_{1}=C_{1}k^{-\xi_{1}-1/2}$  and
 $c_{2}=C_{2}k^{-\xi_{2}-1/2}$, where $\left\lfloor \cdot
\right\rfloor $  denotes the integer part. Let $\omega\sim
N(0,I_{k})$ . There exist a constant $D>0$  depending only on
$(\xi_{1},\xi_{2},M_{1})$  such that for any $C_{1},C_{2}\in
(0,D)$ , we have $P(c_{1}\omega\in \mathcal{S}^{*}_{M_1,\xi_{1}})=1-o(1)$  and
$P(c_{2}\omega\in \mathcal{S}^{*}_{M_1,\xi_{2}})=1-o(1)$.
\end{lemma}

\begin{proof}[Proof of Lemma \ref{lem: A3}]

 We now show $P(c_{1}\omega\in\mathcal{S}^{*}_{M_1,\xi_{1}}%
)=1-o(1)$. Let $\omega=(\omega_{1},...,\omega_{k})^{\top}\in{\mathbb{R}%
}^{k}$ follow $N(0,I_{k})$. The goal is to show that the following event
occurs with probability $1-o(1)$,
\[
\bigcap_{t=1}^{k-1}\left\{  c_{1}\sqrt{\sum_{j=t+1}^{k}\omega_{j}^{2}}\leq
M_{1}t^{-\xi_{1}}\right\}  .
\]

We can rewrite this event as
\[
\bigcap_{t=1}^{k-1}\left\{  \sum_{j=t+1}^{k}(\omega_{j}^{2}-1)\leq\left(
c_{1}^{-2}M_{1}^{2}t^{-2\xi_{1}}-(k-t)\right)  \right\}  .
\]

Notice that $\mathbb{E}(\omega_{j}^{2}-1)^{2}=3$. By Kolmogorov's maximal inequality,
we have that for $x>0$,
\[
P\left(  \max_{1\leq t\leq k-1}\left|  \sum_{j=t+1}^{k}(\omega_{j}%
^{2}-1)\right|  >x\right)  \leq\frac{3k}{x^{2}}.
\]

Hence,
\[
P\left(  \max_{1\leq t\leq k-1}\left|  \sum_{j=t+1}^{k}(\omega_{j}%
^{2}-1)\right|  >\sqrt{k\log k}\right)  =o(1).
\]

Therefore, we only need to show
\[
\min_{1\leq t\leq k}c_{1}^{-2}M_{1}^{2}t^{-2\xi_{1}}-(k-t)\geq\sqrt{k\log k}.
\]

Recall $c_{1}=C_{1}k^{-\xi_{1}-1/2}$ for a constant $C_{1}>0$. We need to
show
\[
\min_{1\leq t\leq k}M_{1}^{2}C_{1}^{-2}k^{2\xi_{1}+1}t^{-2\xi_{1}}+t\geq
k+\sqrt{k\log k}.
\]

Notice that the left-hand side is the minimum of a convex function of $t$.
Therefore, the minimum occurs at $t=t_{*}$ with $t_{*}$ being the solution of
the first-order condition
\[
-M_{1}^{2}C_{1}^{-2}k^{2\xi_{1}+1}t_{*}^{-2\xi_{1}-1}/(2\xi_{1})+1=0.
\]

This means $t_{*}=(2\xi_{1}M_{1}^{-2}C_{1}^{2})^{-1/(2\xi_{1}+1)}k$, which
means
\begin{eqnarray*}
&&\min_{t\in{\mathbb{R}}} \left\{ M_{1}^{2}C_{1}^{-2}k^{2\xi_{1}+1}t^{-2\xi_{1}%
}+t \right\}
\\
&=&\left[  (2\xi_{1}M_{1}^{-2})^{-1/(2\xi_{1}+1)}+M_{1}^{2}(2\xi_{1}%
M_{1}^{-2})^{2\xi_{1}/(2\xi_{1}+1)}\right]  k\times C_{1}^{-2/(2\xi_{1}+1)}.
\end{eqnarray*}

Hence, the right-hand side is larger than $k+\sqrt{k\log k}$ for small enough
$C_{1}$. Therefore, we have proved that for $c_{1}=C_{1}k^{-\xi_{1}-1/2}$ with
a small constant $C_{1}>0$, $P(c_{1}\omega\in\mathcal{S}^{*}_{M_1,\xi_{1}}%
)=1-o(1)$. The result for $P(c_{2}\lambda\in\mathcal{S}^{*}_{M_1,\xi_{2}})=1-o(1)$
is analogous.
\end{proof}

\bigskip


\begin{lemma}\label{lem: A4}
 
Let $W\sim N(0,I_{k})$. Suppose that
$A\in \RR^{k\times k}$ is positive semi-definite and $b\in R^{k}$. Then
\[
\mathbb{E}\operatorname{exp}\left(  -\frac{1}{2}W^{\top}AW+b^{\top}W\right)
=\left(  \det(A+I_{k})\right)  ^{-1/2}\operatorname{exp}\left(  b^{\top
}(A+I_{k})^{-1}b/2\right).
\]

\end{lemma}

\begin{proof} [Proof of Lemma \ref{lem: A4}]

Let $A=UBU^{\top}$, where $B$ is diagonal with $B_{jj}\geq0$
and $U^{\top}U=I_{k}$. This is possible because $A$ is positive
semi-definite. Define $T=U^{\top}W$ and $c=U^{\top}b$. Clearly, $T$ is
Gaussian with mean zero and variance $U^{\top}I_{k}U=I_{k}$. Then we need to
compute
\[
\mathbb{E}\operatorname{exp}\left(  -\frac{1}{2}T^{\top}BT+c^{\top}T\right)
=\mathbb{E}\operatorname{exp}\left(  \sum_{j=1}^{k}(-T_{j}^{2}B_{jj}/2+c_{j}%
T_{j})\right)  .
\]

Since $T_{j}\sim N(0,1)$ is independent across $j$, it follows that
\begin{align*}
\E \operatorname{exp}\left(  -\frac{1}{2}T^{\top}BT+c^{\top}T\right)   &
=\E \operatorname{exp}\left(  \sum_{j=1}^{k}(-T_{j}^{2}B_{jj}/2+c_{j}%
T_{j})\right) \\
&  =\prod_{j=1}^{k}\E \operatorname{exp}\left(  -T_{j}^{2}B_{jj}/2+c_{j}%
T_{j}\right) \\
&  =\prod_{j=1}^{k}\left[  \int_{-\infty}^{\infty} \operatorname{exp}\left(
-x^{2}B_{jj}/2+c_{j}x\right)  \cdot\frac{1}{\sqrt{2\pi}} \operatorname{exp}%
(-x^{2}/2)dx\right] \\
&  =\prod_{j=1}^{k}\left[  \frac{1}{\sqrt{B_{jj}+1}} \operatorname{exp}\left(
\frac{c_{j}^{2}}{2(B_{jj}+1)}\right)  \right] \\
&  =\left[  \prod_{j=1}^{k}(B_{jj}+1)\right]  ^{-1/2} \operatorname{exp}%
\left(  \frac{1}{2}\cdot\sum_{j=1}^{k}\frac{c_{j}^{2}}{B_{jj}+1}\right)  .
\end{align*}

Since $U^{\top}U=I_{k}$, we have
\[
\prod_{j=1}^{k}(B_{jj}+1)=\det(B+I_{k})=\det(UBU^{\top}+I_{k})=\det
(A+I_{k}).
\]

Moreover, since $B$ is diagonal, we have
\[
\sum_{j=1}^{k}\frac{c_{j}^{2}}{B_{jj}+1}=c^{\top}(B+I_{k})^{-1}c=(U^{\top
}b)^{\top}\left(  U^{\top}AU+I_{k}\right)  ^{-1}(U^{\top}b)=b^{\top
}(A+I_{k})^{-1}b.
\]

The above three displays imply
\[
\E \operatorname{exp}\left(  -\frac{1}{2}T^{\top}BT+c^{\top}T\right)
=\left(  \det(A+I_{k})\right)  ^{-1/2} \operatorname{exp}\left(  b^{\top
}(A+I_{k})^{-1}b/2\right)  .
\]

The proof is complete. 
\end{proof}

\bigskip

\begin{lemma}\label{lem: A5}
Let $C_{0},C_{1},C_{2},\xi_{1},\xi_{2}>0$ 
be constants. Define $k=\left\lfloor C_{0}n^{1/(\xi_{1}+\xi_{2}%
+1/2)}\right\rfloor $\textit{, }$c_{1}=C_{1}k^{-\xi_{1}-1/2}$ and
 $c_{2}=C_{2}k^{-\xi_{2}-1/2}$, where $\left\lfloor \cdot
\right\rfloor $ denotes the integer part. For any $\omega\in \RR^{k}$, let $\Phi_{1,\omega}$ denote the distribution of i.i.d
$\{(Y_{i},D_{i},Z_{i})\}_{i=1}^{n}$  given by $Z_{i}\sim N(0,I_{k}%
)$, 
\[%
\begin{pmatrix}
Y_{i}\\
D_{i}%
\end{pmatrix}
\mid Z_i \sim N\left(
\begin{pmatrix}
c_{1}Z_{i}^{\top}\omega\\
0
\end{pmatrix}
,\ \Sigma\right)  \qquad{\text{and}}\qquad\Sigma=%
\begin{pmatrix}
1 & c_{1}c_{2}k\\
c_{1}c_{2}k & 1
\end{pmatrix}.
\]
Let $\Phi_{2,\omega}$  denote the distribution of i.i.d
$\{(Y_{i},D_{i},Z_{i})\}_{i=1}^{n}$  given by $Z_{i}\sim N(0,I_{k}%
)$,
\[%
\begin{pmatrix}
Y_{i}\\
D_{i}%
\end{pmatrix}
\mid Z_i \sim N\left(
\begin{pmatrix}
c_{1}Z_{i}^{\top}\omega\\
c_{2}Z_{i}^{\top}\omega
\end{pmatrix}
,\
\begin{pmatrix}
1 \quad & 0\\
0 \quad & 1
\end{pmatrix}
\right)  .
\]
Assume that $\xi_{2}\geq1/4>\xi_{1}$  and $\xi_{1}+\xi
_{2}<1/2$. Let $Z_{i}$  be i.i.d from $N(0,I_{k})$.
Then 
\[
{\text{TV}}(\bar{\phi}_{1},\bar{\phi}_{2})\leq2C_{1}C_{2}C_{0}^{-(\xi_{1}%
+\xi_{2}+1/2)}+o(1),
\]
where $TV$  denotes the total variation distance, $\bar{\phi
}_{1}=\int_{{\mathbb{R}}^{k}}\Phi_{1,\omega}\mu(\omega)d\omega$, $\bar{\phi}_{2}=\int_{{\mathbb{R}}^{k}}\Phi_{2,\omega}\mu(\omega)d\omega
$  and $\mu$ is the density function of $N(0,I_{k})$,
i.e., $\mu(\omega)=(2\pi)^{-k/2}\operatorname{exp}(-\Vert\omega\Vert
_{2}^{2}/2)$.
\end{lemma}


\begin{proof} [Proof of Lemma \ref{lem: A5}]
For simplicity, we shall use $|\cdot|$ to denote the determinant of a
square matrix. Since $\xi_{1}+\xi_{2}<1/2$, we have $k\gg n$. Let $\otimes$ denote the Kronecker product of matrices.

Let $Y=(Y_1,...,Y_n)^{\top}\in\RR^{n}$, $D=(D_1,...,D_n)^{\top}\in\RR^{n}$, $\varepsilon=(\varepsilon_1,...,\varepsilon_n)^{\top}\in\RR^{n}$, $u=(u_1,...,u_n)^{\top}\in\RR^{n}$  and $Z=(Z_1,...,Z_n)^{\top}\in\RR^{n\times k}$. For $j\in\{1,2\}$, let $f_{j,\omega}$ denote the conditional distribution of
$(Y,D)$ given $Z$ under $\Phi_{j,\omega}$ and define $\bar{f}_{j}%
=\int_{{\mathbb{R}}^{k}}f_{j,\omega}\mu(\omega)d\omega$. Let ${\text{KL}}$
denote the Kullback-Leibler divergence. We will proceed in four steps. In the
first three steps, we bound ${\text{KL}}(\bar{f}_{1},\bar{f}_{2})$; in the
last step, we prove the final result.

\textbf{Step 1:} Deriving the mixture $\bar{f}_{1}$

Notice that
\[
R=
\begin{pmatrix}
Y\\
D
\end{pmatrix}
=
\begin{pmatrix}
c_{1}Z\\
0
\end{pmatrix}
\omega+
\begin{pmatrix}
\varepsilon\\
u
\end{pmatrix}
=G\omega+
\begin{pmatrix}
\varepsilon\\
u
\end{pmatrix}
,
\]
where $G=
\begin{pmatrix}
c_{1}Z\\
0
\end{pmatrix}
=g\otimes Z$ with $g=
\begin{pmatrix}
c_{1}\\
0
\end{pmatrix}
$ and
\[
\E
\begin{pmatrix}
\varepsilon\\
u
\end{pmatrix}
\begin{pmatrix}
\varepsilon\\
u
\end{pmatrix}
^{\top}=
\begin{pmatrix}
\mathbb{E}\varepsilon\varepsilon^{\top} & \mathbb{E}\varepsilon u^{\top}\\
\E u\varepsilon^{\top} & \E uu^{\top}%
\end{pmatrix}
=
\begin{pmatrix}
\Sigma_{1}I_{n} & \Sigma_{3}I_{n}\\
\Sigma_{3}I_{n} & \Sigma_{2}I_{n}%
\end{pmatrix}
=\Sigma\otimes I_{n}
\]
with%

\[
\Sigma=
\begin{pmatrix}
1 & c_{1}c_{2}k\\
c_{1}c_{2}k & 1
\end{pmatrix}
.
\]

Let $S=\Sigma^{-1}$. Notice that
\[
S=
\begin{pmatrix}
S_{1} & S_{3}\\
S_{3} & S_{2}%
\end{pmatrix}
=
\begin{pmatrix}
1 & -c_{1}c_{2}k\\
-c_{1}c_{2}k & 1
\end{pmatrix}
\frac{1}{1-c_{1}^{2}c_{2}^{2}k^{2}}
\]
and $|\Sigma|=1-c_{1}^{2}c_{2}^{2}k^{2}$.

Thus, under $f_{1,\omega}$, $R$ conditional on $Z$ is Gaussian mean $G\omega$ and variance
$\Sigma\otimes I_{n}$. Hence, the density is
\begin{align*}
f_{1,\omega}&=f_{1,\omega}(R,Z)    =(2\pi)^{-n}|\Sigma\otimes I_{n}|^{-1/2}
\operatorname{exp}\left(  -\frac{1}{2}(R-G\omega)^{\top}(\Sigma\otimes
I_{n})^{-1}(R-G\omega)\right) \\
&  =(2\pi)^{-n}|S\otimes I_{n}|^{1/2} \operatorname{exp}\left(  -\frac{1}%
{2}(R-G\omega)^{\top}(S\otimes I_{n})(R-G\omega)\right) \\
&  =(2\pi)^{-n}|S\otimes I_{n}|^{1/2} \operatorname{exp}\left(  -\frac{1}%
{2}R^{\top}(S\otimes I_{n})R-\frac{1}{2}\omega^{\top}G^{\top}(S\otimes
I_{n})G\omega+R^{\top}(S\otimes I_{n})G\omega\right)
\end{align*}

Hence, we can apply Lemma \ref{lem: A4} and obtain
\begin{align*}
\bar{f}_{1}  &  =\int f_{1,\omega}\mu(\omega)d\omega\\
&  =\mathbb{E}_{\omega\sim N(0,I_{k})}f_{1,\omega}\\
&  =(2\pi)^{-n}|S\otimes I_{n}|^{1/2}\operatorname{exp}\left(  -\frac{1}%
{2}R^{\top}(S\otimes I_{n})R\right) \\
& \qquad \times \mathbb{E}_{\omega\sim N(0,I_{k}%
)}\operatorname{exp}\left(  -\frac{1}{2}\omega^{\top}G^{\top}(S\otimes
I_{n})G\omega+R^{\top}(S\otimes I_{n})G\omega\right) \\
&  =(2\pi)^{-n}|S\otimes I_{n}|^{1/2}\operatorname{exp}\left(  -\frac{1}%
{2}R^{\top}(S\otimes I_{n})R\right)  \times|I+G^{\top}(S\otimes
I_{n})G|^{-1/2}\\
&  \qquad \times \operatorname{exp}\left(  \frac{1}{2}R^{\top}(S\otimes
I_{n})G(I+G^{\top}(S\otimes I_{n})G)^{-1}G^{\top}(S\otimes I_{n})R\right)
\\
&  =(2\pi)^{-n}|S\otimes I_{n}|^{1/2}\cdot|I+G^{\top}(S\otimes
I_{n})G|^{-1/2}\operatorname{exp}\left(  -\frac{1}{2}R^{\top}MR\right)  ,
\end{align*}
where
\begin{align*}
M  &  =S\otimes I_{n}-(S\otimes I_{n})G\left(  I+G^{\top}(S\otimes
I_{n})G\right)  ^{-1}G^{\top}(S\otimes I_{n})\\
&  =S\otimes I_{n}-(Sg\otimes Z)\left(  I+g^{\top}SgZ^{\top}Z\right)
^{-1}(g^{\top}S\otimes Z^{\top}).
\end{align*}

We notice that
\begin{align*}
|M|  &  =\left|  S\otimes I_{n}-(S\otimes I_{n})G\left(  I+G^{\top}(S\otimes
I_{n})G\right)  ^{-1}G^{\top}(S\otimes I_{n})\right| \\
&  \overset{{\text{(i)}}}{=}|S\otimes I_{n}|\cdot\left|  I-G\left(
I+G^{\top}(S\otimes I_{n})G\right)  ^{-1}G^{\top}(S\otimes I_{n})\right|
\\
&  \overset{{\text{(ii)}}}{=}|S\otimes I_{n}|\cdot\left|  I-\left(
I+G^{\top}(S\otimes I_{n})G\right)  ^{-1}G^{\top}(S\otimes I_{n})G\right|
\\
&  =|S\otimes I_{n}|\cdot\left|  I-\left(  I+G^{\top}(S\otimes
I_{n})G\right)  ^{-1}\left(  I+G^{\top}(S\otimes I_{n})G-I\right)  \right|
\\
&  =|S\otimes I_{n}|\cdot\left|  \left(  I+G^{\top}(S\otimes I_{n})G\right)
^{-1}\right| \\
&  =|S\otimes I_{n}|\cdot\left|  I+G^{\top}(S\otimes I_{n})G\right|  ^{-1},
\end{align*}
where (i) follows by the identity $|H-AB|=|H|\cdot|I-BH^{-1}A|$ with
$H=A=S\otimes I_{n}$ and $B=G[ I+G^{\top}(S\otimes I_{n})G]
^{-1}G^{\top}(S\otimes I_{n})$ and (ii) follows by the same identity with
$H=I$ and $A=G$. Therefore,
\[
\bar{f}_{1}=(2\pi)^{-n}|M|^{1/2} \operatorname{exp}\left(  -\frac{1}%
{2}R^{\top}MR\right)  .
\]

\textbf{Step 2:} Deriving the mixture $\bar{f}_{2}$

We can write
\[
R=
\begin{pmatrix}
Y\\
D
\end{pmatrix}
=
\begin{pmatrix}
c_{1}Z\\
c_{2}Z
\end{pmatrix}
\omega+
\begin{pmatrix}
\varepsilon\\
u
\end{pmatrix}
=\tilde{G}\omega+
\begin{pmatrix}
\varepsilon\\
u
\end{pmatrix}
,
\]
where $\tilde{G}=
\begin{pmatrix}
c_{1}Z\\
c_{2}Z
\end{pmatrix}
=\tilde{g}\otimes Z$ with $\tilde{g}=
\begin{pmatrix}
c_{1}\\
c_{2}%
\end{pmatrix}
$. Observe that $E
\begin{pmatrix}
\varepsilon\\
u
\end{pmatrix}
\begin{pmatrix}
\varepsilon\\
u
\end{pmatrix}
^{\top}=I_{2n}$.

Therefore, under $f_{2,\omega}$, $R$ conditional on $Z$ is Gaussian with mean $\tilde{G}\omega$
and variance $I_{2n}$. This means that the density is
\begin{align*}
f_{2,\omega}&=f_{2,\omega}(R,Z)=(2\pi)^{-n} \operatorname{exp}\left(  -\frac
{1}{2}(R-\tilde{G}\omega)^{\top}(R-\tilde{G}\omega)\right) \\
& =(2\pi)^{-n}
\operatorname{exp}\left(  -\frac{1}{2}R^{\top}R-\frac{1}{2}\omega^{\top
}\tilde{G}^{\top}\tilde{G}\omega+R^{\top}\tilde{G}\omega\right)  .
\end{align*}

By a similar calculation as in step 1, we have
\begin{align*}
\bar{f}_{2}  &  =\int f_{2,\omega}\mu(\omega)d\omega\\
&  =\mathbb{E}_{\omega\sim N(0,I_{k})}f_{2,\omega}\\
&  =(2\pi)^{-n} \operatorname{exp}\left(  -\frac{1}{2}R^{\top}R\right)
\mathbb{E}_{\omega\sim N(0,I_{k})} \operatorname{exp}\left(  -\frac{1}{2}%
\omega^{\top}\tilde{G}^{\top}\tilde{G}\omega+R^{\top}\tilde{G}%
\omega\right) \\
&  =(2\pi)^{-n} \operatorname{exp}\left(  -\frac{1}{2}R^{\top}R\right)
\times|I+\tilde{G}^{\top}\tilde{G}|^{-1/2} \operatorname{exp}\left(  \frac
{1}{2}R^{\top}\tilde{G}\left(  I+\tilde{G}^{\top}\tilde{G}\right)
^{-1}\tilde{G}^{\top}R\right) \\
&  =(2\pi)^{-n}|\tilde{M}|^{1/2} \operatorname{exp}\left(  -\frac{1}%
{2}R^{\top}\tilde{M}R\right)  ,
\end{align*}
where
\[
\tilde{M}=I_{2n}-\tilde{G}\left(  I+\tilde{G}^{\top}\tilde{G}\right)
^{-1}\tilde{G}^{\top}=I-(\tilde{g}\otimes Z)\left(  I+\|\tilde{g}\|_{2}%
^{2}Z^{\top}Z\right)  ^{-1}(\tilde{g}^{\top}\otimes Z^{\top}).
\]

\textbf{Step 3:} Computeing KL divergence between the two mixtures

Hence, $\bar{f}_{1}$ and $\bar{f}_{2}$ represents $N(0,\Sigma_{1})$ and
$N(0,\Sigma_{2})$, respectively, where $\Sigma_{1}=M^{-1}$ and $\Sigma
_{2}=\tilde{M}^{-1}$. Now we consider the KL divergence conditional on $Z$
\[
{\text{KL}}(\bar{f}_{1}(\cdot,Z),\bar{f}_{2}(\cdot,Z))=\frac{1}{2}\left(  \log\frac{|\Sigma
_{2}|}{|\Sigma_{1}|}-2n+\mathrm{trace}(\Sigma_{2}^{-1}\Sigma_{1})\right)
=\frac{1}{2}\left(  \log\frac{|M|}{|\tilde{M}|}-2n+\mathrm{trace}(M^{-1}%
\tilde{M})\right)  .
\]

Using the Woodbury's identity, we have
\[
M^{-1}=\Sigma\otimes I_{n}+(gg^{\top})\otimes(ZZ^{\top}).
\]

Therefore,
\begin{align*}
&  \mathrm{trace}(M^{-1}\tilde{M})\\
&  =\mathrm{trace}\left[  \left(  \Sigma\otimes I_{n}+(gg^{\top}%
)\otimes(ZZ^{\top})\right)  \left(  I-(\tilde{g}\otimes Z)\left(
I+\|\tilde{g}\|_{2}^{2}Z^{\top}Z\right)  ^{-1}(\tilde{g}^{\top}\otimes
Z^{\top})\right)  \right] \\
&  =\mathrm{trace}\left[  \Sigma\otimes I_{n}+(gg^{\top})\otimes(ZZ^{\top
})-\left(  \Sigma\otimes I_{n}+(gg^{\top})\otimes(ZZ^{\top})\right)
(\tilde{g}\otimes Z)\left(  I+\|\tilde{g}\|_{2}^{2}Z^{\top}Z\right)
^{-1}(\tilde{g}^{\top}\otimes Z^{\top})\right] \\
&  =n\cdot\mathrm{trace}(\Sigma)+\|g\|_{2}^{2}\mathrm{trace}(Z^{\top
}Z)
\\
& \qquad -\mathrm{trace}\left[  (\tilde{g}^{\top}\otimes Z^{\top})\left(
\Sigma\otimes I_{n}+(gg^{\top})\otimes(ZZ^{\top})\right)  (\tilde
{g}\otimes Z)\left(  I+\|\tilde{g}\|_{2}^{2}Z^{\top}Z\right)  ^{-1}\right]
\\
&  =n\cdot\mathrm{trace}(\Sigma)+\|g\|_{2}^{2}\mathrm{trace}(Z^{\top
}Z)-\mathrm{trace}\left[  \left(  \tilde{g}^{\top}\Sigma\tilde{g}Z^{\top
}Z+(g^{\top}\tilde{g})^{2}Z^{\top}ZZ^{\top}Z\right)  \left(
I+\|\tilde{g}\|_{2}^{2}Z^{\top}Z\right)  ^{-1}\right]  .
\end{align*}

By Sylvester's theorem, we have
\begin{align*}
|\tilde{M}|  &  =\left|  I-(\tilde{g}\otimes Z)\left(  I+\|\tilde{g}\|_{2}%
^{2}Z^{\top}Z\right)  ^{-1}(\tilde{g}^{\top}\otimes Z^{\top})\right| \\
&  =\left|  I-(\tilde{g}^{\top}\otimes Z^{\top})(\tilde{g}\otimes
Z)\left(  I+\|\tilde{g}\|_{2}^{2}Z^{\top}Z\right)  ^{-1}\right| \\
&  =\left|  I-\|\tilde{g}\|_{2}^{2}Z^{\top}Z\left(  I+\|\tilde{g}\|_{2}%
^{2}Z^{\top}Z\right)  ^{-1}\right|  =\left|  \left(  I+\|\tilde{g}\|_{2}%
^{2}Z^{\top}Z\right)  ^{-1}\right|  =\left|  I+(c_{1}^{2}+c_{2}%
^{2})Z^{\top}Z\right|  ^{-1}%
\end{align*}
and
\begin{align*}
|M|  &  =\left|  S\otimes I_{n}-(Sg\otimes Z)\left(  I+g^{\top}SgZ^{\top
}Z\right)  ^{-1}(g^{\top}S\otimes Z^{\top})\right| \\
&  =\left|  S\otimes I_{n}\right|  \cdot\left|  I-\left(  I+g^{\top
}SgZ^{\top}Z\right)  ^{-1}(g^{\top}S\otimes Z^{\top-1}\otimes
I_{n})(Sg\otimes Z)\right| \\
&  =\left|  S\otimes I_{n}\right|  \cdot\left|  I-\left(  I+g^{\top
}SgZ^{\top}Z\right)  ^{-1}g^{\top}SgZ^{\top}Z\right| \\
&  =\left|  S\right|  ^{n}\cdot\left|  \left(  I+g^{\top}SgZ^{\top
}Z\right)  ^{-1}\right|  =\left|  S\right|  ^{n}\cdot\left|  I+c_{1}^{2}%
S_{1}Z^{\top}Z\right|  ^{-1}.
\end{align*}

Let $\sigma_{1},...,\sigma_{n}$ be the eigenvalues of $Z^{\top}Z$. By
Corollary 5.35 of Vershynin (2012) (with $t=\sqrt{n}$), we have that
\[
P\left(  \max_{1\leq i\leq n}|\sqrt{\sigma}_{i}-\sqrt{k}|>2\sqrt{n}\right)
\leq2 \operatorname{exp}(-n/2).
\]

Thus, with probability at least $1-2 \operatorname{exp}(-n/2)$, $\sqrt
{\sigma_{i}}+\sqrt{k}\leq2\sqrt{k}+2\sqrt{n}$ for all $i$. By $k\gg n$, it
follows that $P({\mathcal{A}})\geq1-o(1)$, where
\[
{\mathcal{A}}=\left\{  \max_{1\le i\leq n}|\sigma_{i}-k|\leq\min\{k,4\sqrt
{nk}\}\right\}  .
\]

Now we notice
\begin{align*}
\log\frac{|M|}{|\tilde{M}|}  &  =n\log|S|+\log\frac{\left|  I+(c_{1}^{2}%
+c_{2}^{2})Z^{\top}Z\right|  }{\left|  I+c_{1}^{2}S_{1}Z^{\top}Z\right|
}\\
&  =n\log|S|+\sum_{i=1}^{n}\log\frac{1+(c_{1}^{2}+c_{2}^{2})\sigma_{i}%
}{1+c_{1}^{2}S_{1}\sigma_{i}}\\
&  =n\log|S|+\sum_{i=1}^{n}\log\left(  1+\frac{(c_{1}^{2}(1-S_{1})+c_{2}%
^{2})\sigma_{i}}{1+c_{1}^{2}S_{1}\sigma_{i}}\right) \\
&  \overset{{\text{(i)}}}{\leq}n\log|S|+\sum_{i=1}^{n}\frac{(c_{1}^{2}%
(1-S_{1})+c_{2}^{2})\sigma_{i}}{1+c_{1}^{2}S_{1}\sigma_{i}}\\
&  \overset{{\text{(ii)}}}{=}n\log\left(  \frac{1}{1-c_{1}^{2}c_{2}^{2}k^{2}%
}\right)  +\sum_{i=1}^{n}\frac{(c_{1}^{2}(1-S_{1})+c_{2}^{2})\sigma_{i}%
}{1+c_{1}^{2}S_{1}\sigma_{i}}\\
&  \overset{{\text{(iii)}}}{\leq}nc_{1}^{2}c_{2}^{2}k^{2}+nc_{1}^{4}c_{2}%
^{4}k^{4}+\sum_{i=1}^{n}\frac{(c_{1}^{2}(1-S_{1})+c_{2}^{2})\sigma_{i}%
}{1+c_{1}^{2}S_{1}\sigma_{i}},
\end{align*}
where (i) follows by $\log(1+x)\leq x$ for $x\geq-1$ and (ii) follows by
$|S|=(1-c_{1}^{2}c_{2}^{2}k^{2})^{-1}$ and (iii) follows by the elementary
inequality $\log(\frac{1}{1-x})\leq x+x^{2}$ for $x\in[0,1/2]$ and $c_{1}%
^{2}c_{2}^{2}k^{2}\asymp k^{-2\xi_{1}-2\xi_{2}}=o(1)$. Observe that on the
event ${\mathcal{A}}$,%

\begin{align*}
&  2\cdot{\text{KL}}(\bar{f}_{1}(\cdot,Z),\bar{f}_{2}(\cdot,Z))\\
&  =\log\frac{|M|}{|\tilde{M}|}-2n+\mathrm{trace}(M^{-1}\tilde{M})\\
&  \leq nc_{1}^{2}c_{2}^{2}k^{2}+nc_{1}^{4}c_{2}^{4}k^{4}+\sum_{i=1}^{n}%
\frac{(c_{1}^{2}(1-S_{1})+c_{2}^{2})\sigma_{i}}{1+c_{1}^{2}S_{1}\sigma_{i}%
}-2n+n\mathrm{trace}(\Sigma)\\
&  \qquad+\Vert g\Vert_{2}^{2}\mathrm{trace}(Z^{\top}Z)-\mathrm{trace}%
\left[  \left(  \tilde{g}^{\top}\Sigma\tilde{g}Z^{\top}Z+c_{1}%
^{4}Z^{\top}ZZ^{\top}Z\right)  \left(  I+\Vert\tilde{g}\Vert_{2}%
^{2}Z^{\top}Z\right)  ^{-1}\right] \\
&  =nc_{1}^{2}c_{2}^{2}k^{2}+nc_{1}^{4}c_{2}^{4}k^{4}+\sum_{i=1}^{n}%
\frac{(c_{1}^{2}(1-S_{1})+c_{2}^{2})\sigma_{i}}{1+c_{1}^{2}S_{1}\sigma_{i}}\\
&  \qquad+c_{1}^{2}\sum_{i=1}^{n}\sigma_{i}-\sum_{i=1}^{n}\frac{(c_{1}%
^{2}+c_{2}^{2}+2c_{1}^{2}c_{2}^{2}k)\sigma_{i}+c_{1}^{4}\sigma_{i}^{2}%
}{1+(c_{1}^{2}+c_{2}^{2})\sigma_{i}}\\
&  =nc_{1}^{2}c_{2}^{2}k^{2}+nc_{1}^{4}c_{2}^{4}k^{4}+\sum_{i=1}^{n}%
\frac{(c_{1}^{2}(1-S_{1})+c_{2}^{2})\sigma_{i}}{1+c_{1}^{2}S_{1}\sigma_{i}}\\
&  \qquad+\sum_{i=1}^{n}\left(  \frac{c_{1}^{2}\sigma_{i}[1+(c_{1}^{2}%
+c_{2}^{2})\sigma_{i}]-(c_{1}^{2}+c_{2}^{2}+2c_{1}^{2}c_{2}^{2}k)\sigma
_{i}-c_{1}^{4}\sigma_{i}^{2}}{1+(c_{1}^{2}+c_{2}^{2})\sigma_{i}}\right) \\
&  =nc_{1}^{2}c_{2}^{2}k^{2}+nc_{1}^{4}c_{2}^{4}k^{4}+\sum_{i=1}^{n}\left(
\frac{(c_{1}^{2}(1-S_{1})+c_{2}^{2})\sigma_{i}}{1+c_{1}^{2}S_{1}\sigma_{i}%
}+\frac{c_{2}^{2}c_{1}^{2}\sigma_{i}^{2}-(c_{2}^{2}+2c_{1}^{2}c_{2}%
^{2}k)\sigma_{i}}{1+(c_{1}^{2}+c_{2}^{2})\sigma_{i}}\right) \\
&  =nc_{1}^{4}c_{2}^{4}k^{4}+\sum_{i=1}^{n}\left(  \frac{(c_{1}^{2}%
(1-S_{1})+c_{2}^{2})\sigma_{i}}{1+c_{1}^{2}S_{1}\sigma_{i}}+\frac{c_{2}%
^{2}c_{1}^{2}\sigma_{i}^{2}-(c_{2}^{2}+2c_{1}^{2}c_{2}^{2}k)\sigma_{i}%
}{1+(c_{1}^{2}+c_{2}^{2})\sigma_{i}}+c_{1}^{2}c_{2}^{2}k^{2}\right) \\
&  =nc_{1}^{4}c_{2}^{4}k^{4}+\sum_{i=1}^{n}\left(  \frac{(c_{1}^{2}%
(1-S_{1})+c_{2}^{2})\sigma_{i}}{1+c_{1}^{2}S_{1}\sigma_{i}}+\frac{c_{1}%
^{2}c_{2}^{2}(\sigma_{i}-k)^{2}-c_{2}^{2}\sigma_{i}+c_{1}^{2}c_{2}^{2}%
k^{2}(c_{1}^{2}+c_{2}^{2})\sigma_{i}}{1+(c_{1}^{2}+c_{2}^{2})\sigma_{i}%
}\right) \\
&  =nc_{1}^{4}c_{2}^{4}k^{4}+\sum_{i=1}^{n}\left(  \frac{c_{1}^{2}c_{2}%
^{2}(\sigma_{i}-k)^{2}}{1+(c_{1}^{2}+c_{2}^{2})\sigma_{i}}\right)  +\sum
_{i=1}^{n}\left(  \frac{c_{1}^{2}(1-S_{1})+c_{2}^{2}}{1+c_{1}^{2}S_{1}%
\sigma_{i}}+\frac{-c_{2}^{2}+c_{1}^{2}c_{2}^{2}k^{2}(c_{1}^{2}+c_{2}^{2}%
)}{1+(c_{1}^{2}+c_{2}^{2})\sigma_{i}}\right)  \sigma_{i}\\
&  \overset{{\text{(i)}}}{=}nc_{1}^{4}c_{2}^{4}k^{4}+\sum_{i=1}^{n}\left(
\frac{c_{1}^{2}c_{2}^{2}(\sigma_{i}-k)^{2}}{1+(c_{1}^{2}+c_{2}^{2})\sigma_{i}%
}\right)  +\sum_{i=1}^{n}\left(  \frac{-c_{1}^{4}c_{2}^{2}k^{2}S_{1}+c_{2}%
^{2}}{1+c_{1}^{2}S_{1}\sigma_{i}}+\frac{-c_{2}^{2}+c_{1}^{2}c_{2}^{2}%
k^{2}(c_{1}^{2}+c_{2}^{2})}{1+(c_{1}^{2}+c_{2}^{2})\sigma_{i}}\right)
\sigma_{i}\\
&  =nc_{1}^{4}c_{2}^{4}k^{4}+\sum_{i=1}^{n}\left(  \frac{c_{1}^{2}c_{2}%
^{2}(\sigma_{i}-k)^{2}}{1+(c_{1}^{2}+c_{2}^{2})\sigma_{i}}\right)  +\sum
_{i=1}^{n}\frac{-c_{1}^{4}c_{2}^{2}k^{2}S_{1}+c_{1}^{2}c_{2}^{2}k^{2}%
(c_{1}^{2}+c_{2}^{2})}{\left(  1+c_{1}^{2}S_{1}\sigma_{i}\right)  \left(
1+(c_{1}^{2}+c_{2}^{2})\sigma_{i}\right)  }\sigma_{i}\\
&  \qquad+\sum_{i=1}^{n}\sigma_{i}^{2}\frac{(-c_{1}^{4}c_{2}^{2}k^{2}%
S_{1}+c_{2}^{2})(c_{1}^{2}+c_{2}^{2})-c_{2}^{2}c_{1}^{2}S_{1}+c_{1}^{4}%
c_{2}^{2}k^{2}(c_{1}^{2}+c_{2}^{2})S_{1}}{\left(  1+c_{1}^{2}S_{1}\sigma
_{i}\right)  \left(  1+(c_{1}^{2}+c_{2}^{2})\sigma_{i}\right)  }\\
&  =nc_{1}^{4}c_{2}^{4}k^{4}+\sum_{i=1}^{n}\left(  \frac{c_{1}^{2}c_{2}%
^{2}(\sigma_{i}-k)^{2}}{1+(c_{1}^{2}+c_{2}^{2})\sigma_{i}}\right)  +\sum
_{i=1}^{n}\frac{c_{1}^{4}c_{2}^{2}k^{2}(1-S_{1})+c_{1}^{2}c_{2}^{4}k^{2}%
}{\left(  1+c_{1}^{2}S_{1}\sigma_{i}\right)  \left(  1+(c_{1}^{2}+c_{2}%
^{2})\sigma_{i}\right)  }\sigma_{i}\\
&  \qquad+c_{2}^{2}\sum_{i=1}^{n}\sigma_{i}^{2}\frac{c_{2}^{2}+c_{1}%
^{2}(1-S_{1})}{\left(  1+c_{1}^{2}S_{1}\sigma_{i}\right)  \left(  1+(c_{1}%
^{2}+c_{2}^{2})\sigma_{i}\right)  }%
\end{align*}

\begin{align*}
&  \overset{{\text{(ii)}}}{=}nc_{1}^{4}c_{2}^{4}k^{4}+\sum_{i=1}^{n}\left(
\frac{c_{1}^{2}c_{2}^{2}(\sigma_{i}-k)^{2}}{1+(c_{1}^{2}+c_{2}^{2})\sigma_{i}%
}\right)  +\sum_{i=1}^{n}\frac{-c_{1}^{6}c_{2}^{4}k^{4}S_{1}+c_{1}^{2}%
c_{2}^{4}k^{2}}{\left(  1+c_{1}^{2}S_{1}\sigma_{i}\right)  \left(
1+(c_{1}^{2}+c_{2}^{2})\sigma_{i}\right)  }\sigma_{i}\\
&  \qquad+c_{2}^{4}\sum_{i=1}^{n}\sigma_{i}^{2}\frac{1-c_{1}^{4}k^{2}S_{1}%
}{\left(  1+c_{1}^{2}S_{1}\sigma_{i}\right)  \left(  1+(c_{1}^{2}+c_{2}%
^{2})\sigma_{i}\right)  }\\
&  \leq nc_{1}^{4}c_{2}^{4}k^{4}+\sum_{i=1}^{n}c_{1}^{2}c_{2}^{2}(\sigma
_{i}-k)^{2}+c_{1}^{2}c_{2}^{4}k^{2}\sum_{i=1}^{n}\sigma_{i}+c_{2}^{4}%
\sum_{i=1}^{n}\sigma_{i}^{2}\\
&  \overset{{\text{(iii)}}}{\leq}nc_{1}^{4}c_{2}^{4}k^{4}+16c_{1}^{2}c_{2}%
^{2}n^{2}k+c_{1}^{2}c_{2}^{4}k^{2}\cdot(n\cdot2k)+c_{2}^{4}\cdot(n\cdot
4k^{2})\\
&  \overset{{\text{(iv)}}}{=}O\left(  n^{[1/2-3(\xi_{1}+\xi_{2})]/(\xi_{1}%
+\xi_{2}+1/2)}\right)  +16C_{1}^{2}C_{2}^{2}C_{0}^{-(2\xi_{1}+2\xi_{2}+1)}\\
&  \qquad+O\left(  n^{(1/2-\xi_{1}-3\xi_{2})/(\xi_{1}+\xi_{2}+1/2)}\right)
+O\left(  n^{(1/2+\xi_{1}-3\xi_{2})/(\xi_{1}+\xi_{2}+1/2)}\right) \\
&  \overset{{\text{(v)}}}{=}16C_{1}^{2}C_{2}^{2}C_{0}^{-(2\xi_{1}+2\xi_{2}%
+1)}+o(1),
\end{align*}
where (i) and (ii) follow by $1-S_{1}=-c_{1}^{2}c_{2}^{2}k^{2}S_{1}$ (due to
$S_{1}=(1-c_{1}^{2}c_{2}^{2}k^{2})^{-1}$), (iii) follows by the definition of
${\mathcal{A}}$, (iv) follows by $c_{1}=C_{1}k^{-\xi_{1}-1/2}$, $c_{2}%
=C_{2}k^{-\xi_{2}-1/2}$ and $k=C_{0}n^{1/(\xi_{1}+\xi_{2}+1/2)}$ and (vii)
follows by $\xi_{2}\geq1/4>\xi_{1}$. Thus, we have proved that on the event
${\mathcal{A}}$,
\[
{\text{KL}}(\bar{f}_{1}(\cdot,Z),\bar{f}_{2}(\cdot,Z))\leq8C_{1}^{2}C_{2}^{2}C_{0}^{-(2\xi
_{1}+2\xi_{2}+1)}+o(1).
\]

\textbf{Step 4:} Deriving the final result.

By first Pinsker's inequality (Lemma 2.5 of \cite{Tsybakov2009}), we have that on
the event ${\mathcal{A}}$,
\begin{equation}
{\text{TV}}(\bar{f}_{1}(\cdot,Z),\bar{f}_{2}(\cdot,Z))\leq\sqrt{{\text{KL}}(\bar{f}_{1}(\cdot,Z),\bar
{f}_{2}(\cdot,Z))/2}\leq2C_{1}C_{2}C_{0}^{-(\xi_{1}+\xi_{2}+1/2)}+o(1).
\label{eq:TVcalculationeq5}%
\end{equation}

By Scheff\'e's theorem (Lemma 2.1 of \cite{Tsybakov2009}),
\begin{equation}
{\text{TV}}(\bar{f}_{1}(\cdot,Z),\bar{f}_{2}(\cdot,Z))=\frac{1}{2}\int|\bar{f}_{1}(r,Z)-\bar
{f}_{2}(r,Z)|dr. \label{eq:TVcaleq6}%
\end{equation}

Since $Z_{i}$ is i.i.d $N(0,I_{k})$, we have that for $j\in\{1,2\}$,
$\Phi_{j,\omega}(r,z)=f_{j,\omega}(r,z)\mu_{n}(z)$ with $r=(y,d)$, where
$\mu_{n}(z)=\prod_{i=1}^{n}\mu(z_{i})$ with $z=(z_{1},...,z_{n})$. This means
that
\[
\bar{\phi}_{j}(r,z)=\int_{{\mathbb{R}}^{k}}\Phi_{j,\omega}(r,z)\mu
(\omega)d\omega=\int_{{\mathbb{R}}^{k}}f_{j,\omega}(r,z)\mu_{n}%
(z)\mu(\omega)d\omega=\bar{f}_{j}(r,z)\mu_{n}(z).
\]

In the following argument, $E$ denotes expectation with respect to the
randomness in $Z$ given by $Z_{i}\sim i.i.d.\ N(0,I_{k})$. Then
\begin{align*}
{\text{TV}}(\bar{\phi}_{1},\bar{\phi}_{2})  &  =\frac{1}{2}\int\mu_{n}%
(z)|\bar{f}_{1}(r,z)-\bar{f}_{2}(r,z)|drdz\\
&  =\frac{1}{2}\mathbb{E}\int|\bar{f}_{1}(r,Z)-\bar{f}_{2}(r,Z)|dr\\
&  =\frac{1}{2}\mathbb{E}\int|\bar{f}_{1}(r,Z)-\bar{f}_{2}(r,Z)|dr\cdot{\mathbf{1}%
}\{Z\in{\mathcal{A}}\}+\frac{1}{2}\mathbb{E}\int|\bar{f}_{1}(r,Z)-\bar{f}%
_{2}(r,Z)|dr\cdot{\mathbf{1}}\{Z\notin{\mathcal{A}}\}\\
&  \overset{{\text{(i)}}}{\leq}2C_{1}C_{2}C_{0}^{-(\xi_{1}+\xi_{2}%
+1/2)}+o(1)+\frac{1}{2}\mathbb{E}\int|\bar{f}_{1}(r,Z)-\bar{f}_{2}(r,Z)|dr\cdot
{\mathbf{1}}\{Z\notin{\mathcal{A}}\}\\
&  \overset{{\text{(ii)}}}{\leq}2C_{1}C_{2}C_{0}^{-(\xi_{1}+\xi_{2}%
+1/2)}+o(1)+P(Z\notin{\mathcal{A}}),
\end{align*}
where (i) follows by (\ref{eq:TVcalculationeq5}) and (\ref{eq:TVcaleq6}) and
(ii) follows by the fact that $\int\bar{f}_{j}(r,Z)dr=1$ almost surely. Since
$P({\mathcal{A}})=1-o(1)$, we have ${\text{TV}}(\bar{\phi}_{1},\bar{\phi}%
_{2})\leq2C_{1}C_{2}C_{0}^{-(\xi_{1}+\xi_{2}+1/2)}+o(1)$. 
\end{proof}

\bigskip


\begin{lemma}\label{lem: A6}
Assume that $p\gtrsim n^{2}$. If $\xi
_{1}+\xi_{2}<1/2$ and $\xi_{2}\geq1/4>\xi_{1}$, then
$$\mathcal{R}^{\rm{Ordered}}_{\xi_{1},\xi_{2}}\gtrsim n^{-(\xi_{1}+\xi_{2})/(\xi_{1}+\xi
_{2}+1/2)}.$$
\end{lemma}

\begin{proof} [Proof of Lemma \ref{lem: A6}]

 Define 
 $$k=\left\lfloor C_{0}n^{1/(\xi_{1}+\xi_{2}%
+1/2)}\right\rfloor, \ c_{1}=C_{1}k^{-\xi_{1}-1/2} \mbox{ and } c_{2}=C_{2}%
k^{-\xi_{2}-1/2},$$ where $C_{0},C_{1},C_{2}>0$ are constants to be determined.
Since $\xi_{1}+\xi_{2}<1/2$, we have $k\gg n$. Since $p\gtrsim n^{2}$ and
$k\asymp n^{1/(\xi_{1}+\xi_{2}+1/2)}\ll n^{2}$, we have $p\gg k$.

In the following analysis, we set $\beta_{j}=\phi_{j}=0$ for $j>k$. This
effectively reduces the dimension from $p$ to $k$ in our analysis.

For any $\omega\in{\mathbb{R}}^{k}$, let 
$$\lambda_{\omega}=(c_{1}c_{2}%
k,c_{1}\omega,0,1-c_{1}^{2}c_{2}^{2}k^{2},1), \mbox{ and } \tilde{\lambda}_{\omega
}=(0,c_{1}\omega,c_{2}\omega,1,1).$$ Notice that $1-c_{1}^{2}c_{2}^{2}k^{2}$
is well defined as the value of a variance (i.e., $1-c_{1}^{2}c_{2}^{2}%
k^{2}>0$) if we choose $C_{0},C_{1},C_{2}>0$ small enough; to see this notice
that $c_{1}c_{2}k\asymp k^{-(\xi_{1}+\xi_{2})}=o(1)$.

We consider two probability measures: $P_{1}=\int P_{\lambda_{\omega}}%
\mu(\omega)d\omega$ and $P_{2}=\int P_{\tilde{\lambda}_{\omega}}\mu
(\omega)d\omega$, where $\mu$ is the density of $N(0,I_{k})$, i.e.,
$\mu(\omega)=(2\pi)^{-k/2} \operatorname{exp}(-\|\omega\|_{2}^{2}/2)$.

Under $P_{\lambda_{\omega}}$, $Z_{i}\sim N(0,I_{k})$ and
\[%
\begin{pmatrix}
Y_{i}\\
D_{i}%
\end{pmatrix}
=
\begin{pmatrix}
-D_{i}c_{1}c_{2}k+c_{1}Z_{i}^{\top}\omega+\varepsilon_{i}\\
D_{i}%
\end{pmatrix}
=
\begin{pmatrix}
c_{1}Z_{i}^{\top}\omega+c_{1}c_{2}ku_{i}+\varepsilon_{i}\\
u_{i}%
\end{pmatrix}
\]
with $\E u_{i}^{2}=1$ and $\E v_{i}^{2}=1-c_{1}^{2}c_{2}^{2}k^{2}$. Therefore,
under $P_{\lambda_{\omega}}$,
\[%
\begin{pmatrix}
Y_{i}\\
D_{i}%
\end{pmatrix}
\mid Z\sim N\left(
\begin{pmatrix}
c_{1}Z_{i}^{\top}\omega\\
0
\end{pmatrix}
,\ \Sigma\right)  \qquad{\text{and}}\qquad\Sigma=
\begin{pmatrix}
1 & c_{1}c_{2}k\\
c_{1}c_{2}k & 1
\end{pmatrix}
.
\]

Similarly, we observe that under $P_{\tilde{\lambda}_{\omega}}$, $Z_{i}\sim
N(0,I_{k})$ and
\[%
\begin{pmatrix}
Y_{i}\\
D_{i}%
\end{pmatrix}
=
\begin{pmatrix}
c_{1}Z_{i}^{\top}\omega+\varepsilon_{i}\\
c_{2}Z_{i}^{\top}\omega+u_{i}%
\end{pmatrix}
,
\]
which means that
\[%
\begin{pmatrix}
Y_{i}\\
D_{i}%
\end{pmatrix}
\mid Z\sim N\left(
\begin{pmatrix}
c_{1}Z_{i}^{\top}\omega\\
c_{2}Z_{i}^{\top}\omega
\end{pmatrix}
,\ I_{2}\right)  .
\]

By Lemma  \ref{lem: A5}, we have ${\text{TV}}(P_{1},P_{2})\leq2C_{1}C_{2}C_{0}^{-(\xi
_{1}+\xi_{2}+1/2)}+o(1)$. We can choose constants $C_{1},C_{2},C_{0}>0$ small
enough such that $C_{1}C_{2}C_{0}^{-(\xi_{1}+\xi_{2}+1/2)}\leq\alpha/4$.
Thus,
\begin{equation}
{\text{TV}}(P_{1},P_{2})\leq\alpha/2+o(1). \label{eq:lemcase1eq5}%
\end{equation}

Let $\psi$ be a test for
\[
H_{0}:\ \lambda\in\Lambda_{0}\qquad{\text{versus}}\qquad H_{1}:\ \lambda\in
\Lambda_{1},
\]
where 
$$\Lambda_{1}=\{\lambda=(\beta,\phi,\pi,\sigma_{\varepsilon}^{2},\sigma_{u}^{2}%
)\in\Lambda^{\rm{Ordered}}(\xi_{1},\xi_{2}):\ |\beta|\geq C_{3}n^{-(\xi_{1}+\xi_{2})/(\xi
_{1}+\xi_{2}+1/2)}\}$$ and $$\Lambda_{0}=\{\lambda=(\theta,\beta,\phi,\sigma_{\varepsilon}%
^{2},\sigma_{u}^{2})\in\Lambda^{\rm{Ordered}}(\xi_{1},\xi_{2}):\ \theta=0\}$$ with $C_{3}%
=C_{1}C_{2}C_{0}^{-(\xi_{1}+\xi_{2})}/2$. Let $\eta(\lambda)=\mathbb{E}_{\lambda}(\psi)$
be the power function. By definition of $\psi$, we have $\sup_{\lambda\in
\Lambda_{0}}\eta(\lambda)\leq\alpha$.

By Lemma  \ref{lem: A3}, we can choose small enough constants $C_{1},C_{2}>0$ such that
\begin{equation}
\int_{{\mathbb{R}}^{k}}{\mathbf{1}}\{\lambda_{\omega}\in\Lambda_{0}%
\}\mu(\omega)d\omega=1-o(1). \label{eq:lemcase1eq6}%
\end{equation}
and%

\begin{equation}
\int_{{\mathbb{R}}^{k}}{\mathbf{1}}\{\tilde{\lambda}_{\omega}\in\Lambda
_{1}\}\mu(\omega)d\omega=1-o(1). \label{eq:lemcase1eq7}%
\end{equation}

Then
\begin{align*}
\left\vert \int_{{\mathbb{R}}^{k}}\eta(\tilde{\lambda}_{\omega})\mu
(\omega)d\omega-\int_{{\mathbb{R}}^{k}}\eta(\lambda_{\omega})\mu
(\omega)d\omega\right\vert  &  =\left\vert \int\left(  \int\psi
dP_{\lambda_{\omega}}\right)  \mu(\omega)d\omega-\int\left(  \int\psi
dP_{\tilde{\lambda}_{\omega}}\right)  \mu(\omega)d\omega\right\vert \\
&  =\left\vert \int\psi dP_{1}-\int\psi dP_{2}\right\vert \overset{{\text{(i)}%
}}{\leq}2\cdot{\text{TV}}(P_{1},P_{2})\overset{{\text{(ii)}}}{\leq}%
\alpha+o_{P}(1),
\end{align*}
where (i) follows by $|\psi|\leq1$ and Scheff\'e's theorem (Lemma 2.1 of
\cite{Tsybakov2009}) and (ii) follows by (\ref{eq:lemcase1eq5}).

Since the power function is bounded by 1, (\ref{eq:lemcase1eq6}) and
(\ref{eq:lemcase1eq7}) imply
\begin{align*}
&  \left\vert \int_{{\mathbb{R}}^{k}}\eta(\lambda_{\omega})\mu(\omega
)d\omega-\int_{{\mathbb{R}}^{k}}{\mathbf{1}}\{\lambda_{\omega}\in\Lambda
_{0}\}\eta(\lambda_{\omega})\mu(\omega)d\omega\right\vert \\
&  =\int_{{\mathbb{R}}^{k}}{\mathbf{1}}\{\lambda_{\omega}\notin\Lambda
_{0}\}\eta(\lambda_{\omega})\mu(\omega)d\omega\leq\int_{{\mathbb{R}}^{k}%
}{\mathbf{1}}\{\lambda_{\omega}\notin\Lambda_{0}\}\mu(\omega)d\omega=o(1)
\end{align*}
and similarly
\[
\left\vert \int_{{\mathbb{R}}^{k}}\eta(\tilde{\lambda}_{\omega})\mu
(\omega)d\omega-\int_{{\mathbb{R}}^{k}}{\mathbf{1}}\{\tilde{\lambda}%
_{\omega}\in\Lambda_{1}\}\eta(\tilde{\lambda}_{\omega})\mu(\omega
)d\omega\right\vert =o(1).
\]

By the above three displays, we have
\[
\int_{{\mathbb{R}}^{k}}{\mathbf{1}}\{\tilde{\lambda}_{\omega}\in\Lambda
_{1}\}\eta(\tilde{\lambda}_{\omega})\mu(\omega)d\omega\leq\alpha
+\int_{{\mathbb{R}}^{k}}{\mathbf{1}}\{\lambda_{\omega}\in\Lambda_{0}%
\}\eta(\lambda_{\omega})\mu(\omega)d\omega+o(1)\overset{{\text{(i)}}}{\leq
}2\alpha+o(1),
\]
where (i) follows by $\sup_{\lambda\in\Lambda_{0}}\eta(\lambda)\leq\alpha$. On
the other hand, by (\ref{eq:lemcase1eq7}),
\begin{align*}
& \int_{{\mathbb{R}}^{k}}{\mathbf{1}}\{\tilde{\lambda}_{\omega}\in\Lambda
_{1}\}\eta(\tilde{\lambda}_{\omega})\mu(\omega)d\omega
\\
&\geq\left(
\operatorname{inf}_{\lambda\in\Lambda_{1}}\eta(\lambda)\right)  \times
\int_{{\mathbb{R}}^{k}}{\mathbf{1}}\{\tilde{\lambda}_{\omega}\in\Lambda
_{1}\}\mu(\omega)d\omega=\left(  \operatorname{inf}_{\lambda\in\Lambda_{1}%
}\eta(\lambda)\right)  \times(1-o(1)).
\end{align*}

Combining the above two displays, we have that
\[
\operatorname{inf}_{\lambda\in\Lambda_{1}}\eta(\lambda)\leq\frac{2\alpha
+o(1)}{1-o(1)}=2\alpha+o(1).
\]

Therefore, for testing $H_{0}:\ \tau(\lambda)=0$ versus $H_{1}:\ |\tau
(\lambda)|\geq C_{3}n^{-(\xi_{1}+\xi_{2})/(\xi_{1}+\xi_{2}+1/2)}$, any test of
size $\alpha$ can has the worst-case power bounded by $2\alpha+o(1)$. The
desired result follows. 
\end{proof}

\bigskip

\begin{lemma}\label{lem: A7}
\textit{ Let }$C_{0},C_{1},C_{2},\xi_{1},\xi_{2}>0$\textit{
be constants. Define }$k=\left\lfloor C_{0}n^{1/(\xi_{1}+\xi_{2}%
+1/2)}\right\rfloor $\textit{, }$c_{1}=C_{1}k^{-\xi_{1}-1/2}$\textit{ and
}$c_{2}=C_{2}k^{-\xi_{2}-1/2}$\textit{, where }$\left\lfloor \cdot
\right\rfloor $\textit{ denotes the integer part. For any }$\omega\in \RR^{k}%
$\textit{, let }$\Phi_{1,\omega}$\textit{ denote the distribution of i.i.d
}$\{(Y_{i},D_{i},Z_{i})\}_{i=1}^{n}$\textit{ given by }$Z_{i}\sim N(0,I_{k})$
\[%
\begin{pmatrix}
D_{i}\\
Y_{i}%
\end{pmatrix}
\mid Z\sim N\left(
\begin{pmatrix}
c_{2}Z_{i}^{\top}\omega\\
0
\end{pmatrix}
,\ \Sigma\right)  \qquad{\text{and}}\qquad\Sigma=%
\begin{pmatrix}
1 & c_{1}c_{2}k\\
c_{1}c_{2}k & 1
\end{pmatrix}
.
\]
\textit{Let }$\Phi_{2,\omega}$\textit{ denote the distribution of i.i.d
}$\{(Y_{i},D_{i},Z_{i})\}_{i=1}^{n}$\textit{ given by }$Z_{i}\sim N(0,I_{k}%
)$,
\[%
\begin{pmatrix}
D_{i}\\
Y_{i}%
\end{pmatrix}
\mid Z\sim N\left(
\begin{pmatrix}
c_{2}Z_{i}^{\top}\omega\\
c_{1}Z_{i}^{\top}\omega
\end{pmatrix}
,\
\begin{pmatrix}
1 & 0\\
0 & 1
\end{pmatrix}
\right)  .
\]
\textit{Assume that }$\xi_{1}\geq1/4>\xi_{2}$\textit{ and }$\xi_{1}+\xi
_{2}<1/2$\textit{. Let }$Z_{i}$\textit{ be i.i.d from }$N(0,I_{k})$\textit{.
Then}
\[
{\text{TV}}(\bar{\phi}_{1},\bar{\phi}_{2})\leq2C_{1}C_{2}C_{0}^{-(\xi_{1}%
+\xi_{2}+1/2)}+o(1),
\]
\textit{where }$TV$\textit{ denotes the total variation distance, }$\bar{\phi
}_{1}=\int_{{\mathbb{R}}^{k}}\Phi_{1,\omega}\mu(\omega)d\omega$\textit{,
}$\bar{\phi}_{2}=\int_{{\mathbb{R}}^{k}}\Phi_{2,\omega}\mu(\omega)d\omega
$\textit{ and }$\mu$\textit{ is the density function of }$N(0,I_{k})$\textit{,
i.e., }$\mu(\omega)=(2\pi)^{-k/2}\operatorname{exp}(-\Vert\omega\Vert
_{2}^{2}/2)$.
\end{lemma}

\begin{proof} [Proof of Lemma \ref{lem: A7}]

 The proof is analogous to that of Lemma  \ref{lem: A5}. We simply swapped
the role of $c_{1}$ and $c_{2}$ as well as the role of $Y_{i}$ and $D_{i}$.
\end{proof}

\bigskip

\begin{lemma}\label{lem: A8}
\textit{Assume that }$p\gtrsim n^{2}$\textit{. If }$\xi
_{1}+\xi_{2}<1/2$\textit{ and }$\xi_{1}\geq1/4>\xi_{2}$\textit{, then }%
$$\mathcal{R}^{\rm{Ordered}}_{\xi_{1},\xi_{2}}\gtrsim n^{-(\xi_{1}+\xi_{2})/(\xi_{1}+\xi_{2}+1/2)}.$$
\end{lemma}

\begin{proof} [Proof of Lemma \ref{lem: A8}]

The argument is analogous to that of Lemma  \ref{lem: A6} with $\xi_{1}$ and
$\xi_{2}$ swapped. Define 
$$k=C_{0}n^{1/(\xi_{1}+\xi_{2}+1/2)}, \qquad c_{1}%
=C_{1}k^{-\xi_{1}-1/2}, \mbox{ and } c_{2}=C_{2}k^{-\xi_{2}-1/2},$$ where $C_{0}%
,C_{1},C_{2}>0$ are constants to be determined. Since $\xi_{1}+\xi_{2}<1/2$,
we have $k\gg n$. Since $p\gtrsim n^{2}$ and $k\asymp n^{1/(\xi_{1}+\xi
_{2}+1/2)}\ll n^{2}$, we have $p\gg k$.

In the following analysis, we set $\beta_{j}=\phi_{j}=0$ for $j>k$. This
effectively reduces the dimension from $p$ to $k$ in our analysis.

For any $\omega\in{\mathbb{R}}^{k}$, let 
$$\lambda_{\omega}=(c_{1}%
c_{2}k,-c_{1}c_{2}^{2}k\omega,c_{2}\omega,1-c_{1}^{2}c_{2}^{2}k^{2},1) \mbox{  and }
 \tilde{\lambda}_{\omega}=(0,c_{1}\omega,c_{2}\omega,1,1).$$  Notice that
$1-c_{1}^{2}c_{2}^{2}k^{2}$ is well defined as the value of a variance (i.e.,
$1-c_{1}^{2}c_{2}^{2}k^{2}>0$) if we choose $C_{0},C_{1},C_{2}>0$ small
enough; to see this notice that $c_{1}c_{2}k\asymp k^{-(\xi_{1}+\xi_{2}%
)}=o(1)$.

We consider two probability measures: 
$$P_{1}=\int P_{\lambda_{\omega}}%
\mu(\omega)d\omega  \mbox{ and  } P_{2}=\int P_{\tilde{\lambda}_{\omega}}\mu
(\omega)d\omega,$$ where $\mu$ is the density of $N(0,I_{k})$, i.e.,
$\mu(\omega)=(2\pi)^{-k/2} \operatorname{exp}(-\|\omega\|_{2}^{2}/2)$.

Under $P_{\lambda_{\omega}}$, $Z_{i}\sim N(0,I_{k})$ and
\[%
\begin{pmatrix}
D_{i}\\
Y_{i}%
\end{pmatrix}
=
\begin{pmatrix}
D_{i}\\
D_{i}c_{1}c_{2}k-c_{1}c_{2}^{2}kZ_{i}^{\top}\omega+\varepsilon_{i}%
\end{pmatrix}
=
\begin{pmatrix}
c_{2}Z_{i}^{\top}\omega+u_{i}\\
c_{1}c_{2}ku_{i}+\varepsilon_{i}%
\end{pmatrix}
\]
with $Eu_{i}^{2}=1$ and $Ev_{i}^{2}=1-c_{1}^{2}c_{2}^{2}k^{2}$. Therefore,
under $P_{\lambda_{\omega}}$,
\[%
\begin{pmatrix}
D_{i}\\
Y_{i}%
\end{pmatrix}
\mid Z\sim N\left(
\begin{pmatrix}
c_{2}Z_{i}^{\top}\omega\\
0
\end{pmatrix}
,\ \Sigma\right)  \qquad{\text{and}}\qquad\Sigma=
\begin{pmatrix}
1 & c_{1}c_{2}k\\
c_{1}c_{2}k & 1
\end{pmatrix}
.
\]

Similarly, we observe that under $P_{\tilde{\lambda}_{\omega}}$, $Z_{i}\sim
N(0,I_{k})$ and
\[%
\begin{pmatrix}
D_{i}\\
Y_{i}%
\end{pmatrix}
=
\begin{pmatrix}
c_{2}Z_{i}^{\top}\omega+u_{i}\\
c_{1}Z_{i}^{\top}\omega+\varepsilon_{i}%
\end{pmatrix}
,
\]
which means that
\[%
\begin{pmatrix}
D_{i}\\
Y_{i}%
\end{pmatrix}
\mid Z\sim N\left(
\begin{pmatrix}
c_{2}Z_{i}^{\top}\omega\\
c_{1}Z_{i}^{\top}\omega
\end{pmatrix}
,\ I_{2}\right)  .
\]

By Lemma  \ref{lem: A7}, we have ${\text{TV}}(P_{1},P_{2})\leq2C_{1}C_{2}C_{0}^{-(\xi
_{1}+\xi_{2}+1/2)}+o(1)$. We can choose constants $C_{1},C_{2},C_{0}>0$ small
enough such that $C_{1}C_{2}C_{0}^{-(\xi_{1}+\xi_{2}+1/2)}\leq\alpha/4$.
Thus,
\begin{equation}
{\text{TV}}(P_{1},P_{2})\leq\alpha/2+o(1). \label{eq:lemcase1eq5-1}%
\end{equation}

Let $\psi$ be a test for
\[
H_{0}:\ \lambda\in\Lambda_{0}\qquad{\text{versus}}\qquad H_{1}:\ \lambda\in
\Lambda_{1},
\]
where 
$$\Lambda_{1}=\{\lambda=(\theta,\beta,\phi,\sigma_{\varepsilon}^{2},\sigma_{u}^{2}%
)\in\Lambda^{\rm{Ordered}}(\xi_{1},\xi_{2}):\ |\theta|\geq C_{3}n^{-(\xi_{1}+\xi_{2})/(\xi
_{1}+\xi_{2}+1/2)}\}$$
 and 
 $$\Lambda_{0}=\{\lambda=(\theta,\beta,\phi,\sigma_{\varepsilon}%
^{2},\sigma_{u}^{2})\in\Lambda^{\rm{Ordered}}(\xi_{1},\xi_{2}):\ \theta=0\}$$ with $C_{3}%
=C_{1}C_{2}C_{0}^{-(\xi_{1}+\xi_{2})}/2$. Let $\eta(\lambda)=\mathbb{E}_{\lambda}(\psi)$
be the power function. By definition of $\psi$, we have $\sup_{\lambda\in
\Lambda_{0}}\eta(\lambda)\leq\alpha$.

By Lemma \ref{lem: A3}, we can choose small enough constants $C_{1},C_{2}>0$ such that
\begin{equation}
\int_{{\mathbb{R}}^{k}}{\mathbf{1}}\{\lambda_{\omega}\in\Lambda_{0}%
\}\mu(\omega)d\omega=1-o(1). \label{eq:lemcase1eq6-1}%
\end{equation}
and%

\begin{equation}
\int_{{\mathbb{R}}^{k}}{\mathbf{1}}\{\tilde{\lambda}_{\omega}\in\Lambda
_{1}\}\mu(\omega)d\omega=1-o(1). \label{eq:lemcase1eq7-1}%
\end{equation}

Then
\begin{align*}
\left\vert \int_{{\mathbb{R}}^{k}}\eta(\tilde{\lambda}_{\omega})\mu
(\omega)d\omega-\int_{{\mathbb{R}}^{k}}\eta(\lambda_{\omega})\mu
(\omega)d\omega\right\vert  &  =\left\vert \int\left(  \int\psi
dP_{\lambda_{\omega}}\right)  \mu(\omega)d\omega-\int\left(  \int\psi
dP_{\tilde{\lambda}_{\omega}}\right)  \mu(\omega)d\omega\right\vert \\
&  =\left\vert \int\psi dP_{1}-\int\psi dP_{2}\right\vert \overset{{\text{(i)}%
}}{\leq}2\cdot{\text{TV}}(P_{1},P_{2})\overset{{\text{(ii)}}}{\leq}%
\alpha+o_{P}(1),
\end{align*}
where (i) follows by $|\psi|\leq1$ and Scheff\'e's theorem (Lemma 2.1 of
 \cite{Tsybakov2009}) and (ii) follows by (\ref{eq:lemcase1eq5-1}).

Since the power function is bounded by 1, (\ref{eq:lemcase1eq6-1}) and
(\ref{eq:lemcase1eq7-1}) imply
\begin{align*}
&  \left\vert \int_{{\mathbb{R}}^{k}}\eta(\lambda_{\omega})\mu(\omega
)d\omega-\int_{{\mathbb{R}}^{k}}{\mathbf{1}}\{\lambda_{\omega}\in\Lambda
_{0}\}\eta(\lambda_{\omega})\mu(\omega)d\omega\right\vert \\
&  =\int_{{\mathbb{R}}^{k}}{\mathbf{1}}\{\lambda_{\omega}\notin\Lambda
_{0}\}\eta(\lambda_{\omega})\mu(\omega)d\omega\leq\int_{{\mathbb{R}}^{k}%
}{\mathbf{1}}\{\lambda_{\omega}\notin\Lambda_{0}\}\mu(\omega)d\omega=o(1)
\end{align*}
and similarly
\[
\left\vert \int_{{\mathbb{R}}^{k}}\eta(\tilde{\lambda}_{\omega})\mu
(\omega)d\omega-\int_{{\mathbb{R}}^{k}}{\mathbf{1}}\{\tilde{\lambda}%
_{\omega}\in\Lambda_{1}\}\eta(\tilde{\lambda}_{\omega})\mu(\omega
)d\omega\right\vert =o(1).
\]

By the above three displays, we have
\[
\int_{{\mathbb{R}}^{k}}{\mathbf{1}}\{\tilde{\lambda}_{\omega}\in\Lambda
_{1}\}\eta(\tilde{\lambda}_{\omega})\mu(\omega)d\omega\leq\alpha
+\int_{{\mathbb{R}}^{k}}{\mathbf{1}}\{\lambda_{\omega}\in\Lambda_{0}%
\}\eta(\lambda_{\omega})\mu(\omega)d\omega+o(1)\overset{{\text{(i)}}}{\leq
}2\alpha+o(1),
\]
where (i) follows by $\sup_{\lambda\in\Lambda_{0}}\eta(\lambda)\leq\alpha$. On
the other hand, by (\ref{eq:lemcase1eq7-1}),
\begin{align*}
&\int_{{\mathbb{R}}^{k}}{\mathbf{1}}\{\tilde{\lambda}_{\omega}\in\Lambda
_{1}\}\eta(\tilde{\lambda}_{\omega})\mu(\omega)d\omega
\\
&\geq\left(
\operatorname{inf}_{\lambda\in\Lambda_{1}}\eta(\lambda)\right)  \times
\int_{{\mathbb{R}}^{k}}{\mathbf{1}}\{\tilde{\lambda}_{\omega}\in\Lambda
_{1}\}\mu(\omega)d\omega=\left(  \operatorname{inf}_{\lambda\in\Lambda_{1}%
}\eta(\lambda)\right)  \times(1-o(1)).
\end{align*}

Combining the above two displays, we have that
\[
\operatorname{inf}_{\lambda\in\Lambda_{1}}\eta(\lambda)\leq\frac{2\alpha
+o(1)}{1-o(1)}=2\alpha+o(1).
\]

Therefore, for testing $H_{0}:\ \tau(\lambda)=0$ versus $H_{1}:\ |\tau
(\lambda)|\geq C_{3}n^{-(\xi_{1}+\xi_{2})/(\xi_{1}+\xi_{2}+1/2)}$, any test of
size $\alpha$ can has the worst-case power bounded by $2\alpha+o(1)$. The
desired result follows. 
\end{proof}

\bigskip

\begin{proof}[Proof of Theorem \ref{thm2}]

We consider two cases: (1) $\max\{\xi_{1},\xi
_{2}\}\geq1/4$ and (2) $\max\{\xi_{1},\xi_{2}\}<1/4$. We first consider the
case of $\max\{\xi_{1},\xi_{2}\}\geq1/4$, which means that either $\xi_1\geq 1/4>\xi_2$ or $\xi_2\geq 1/4>\xi_1$ (due to $\xi_1+\xi_2<1/2$). Then by Lemmas  \ref{lem: A6} and  \ref{lem: A8}, we have proved
that  $\mathcal{R}^{\rm{Ordered}}_{\xi_{1},\xi_{2}}\gtrsim
n^{-(\xi_{1}+\xi_{2})/(\xi_{1}+\xi_{2}+1/2)} \gg n^{-1/2}$.
Thus the result follows for the case of $\max\{\xi_{1},\xi_{2}\}\geq1/4$.

Now we consider the case of $\max\{\xi_{1},\xi_{2}\}<1/4$. Without loss of
generality, we assume $\xi_{1}\geq\xi_{2}$. Then we have $1/4>\xi_{1}\geq
\xi_{2}$. Now we set $\tilde{\xi}_{1}=1/4$ and $\tilde{\xi}_{2}=\xi_{2}$.
Thus, we have $\tilde{\xi}_{1}+\tilde{\xi}_{2}=1/4+\xi_{2}<1/4+1/4=1/2$ and
$\max\{\tilde{\xi}_{1},\tilde{\xi}_{2}\}=1/4$. The results for the case of
$\max\{\xi_{1},\xi_{2}\}\geq1/4$ implies $\mathcal{R}^{\rm{Ordered}}_{\tilde{\xi}_{1},\tilde{\xi}_{2}}\gtrsim n^{-(\tilde{\xi}_{1}+\tilde{\xi}_{2})/(\tilde{\xi
}_{1}+\tilde{\xi}_{2}+1/2)}\gg n^{-1/2}$. By $\tilde{\xi}_{1}>\xi_{1}$,
$\Lambda^{\rm{Ordered}}(\tilde{\xi}_{1},\tilde{\xi}_{2})\subset\Lambda^{\rm{Ordered}}(\xi_{1},\xi_{2})$, which
means that $\mathcal{R}^{\rm{Ordered}}_{\xi_{1},\xi_{2}}\geq \mathcal{R}^{\rm{Ordered}}_{\tilde{\xi}_{1},\tilde{\xi}_{2}}$. Hence, $\mathcal{R}^{\rm{Ordered}}_{\xi_{1},\xi_{2}}\gg n^{-1/2}$. 
\end{proof}

\section{Proof of Theorem \ref{thm: PLM main}}


For the proof of Theorem \ref{thm: PLM main}, recall $Z_{i}=(\psi_{1}(\tilde{Z}_{i}),...,\psi_{p}(\tilde{Z}_{i}))^{\top}\in\RR^{p}$
and $X_{i}=(D_{i},Z_{i}^{\top})^{\top}\in\RR^{p+1}$. Let $\tilde{X}_{i}=(D_{i},\tilde{Z}_{i}^{\top})^{\top}$.
Define $\alpha_{0}(\tilde{X}_{i})=\sigma_{u}^{-2}D_{i}-\sigma_{u}^{-2}g(\tilde{Z}_{i})$,
where $\sigma_{u}^{2}=\mathbb{E}(u_{i}^{2})$. Then 
\[
\mathbb{E}\alpha_{0}(\tilde{X}_{i})X_{i}=\begin{pmatrix}\mathbb{E}\alpha_{0}(\tilde{X}_{i})D_{i}\\
\mathbb{E}\alpha_{0}(\tilde{X}_{i})u_{i}
\end{pmatrix}=e_{1}.
\]

Since $g\in\mathcal{M}_{C_{0},\xi_{2}}$, we have that $\alpha_{0}\in\mathcal{M}_{C,\xi_{2}}$
for some $C>0$. 

Define
\begin{equation}\label{eq:epsilon}
\varepsilon_{n}=\sqrt{\log(p)/n} \qquad \mbox{and} \qquad  
s_{0}\geq C\varepsilon_{n}^{-2/(2\xi_{2}+1)}, 
\end{equation} and $\pi$ be coefficients
of the least squares projection of $\alpha_{0}(\tilde{X}_{i})$ on
$X_{i}$, satisfying
\[
e_{1}-\Sigma\pi=\mathbb{E}[X_{i}\{\alpha_{0}(\tilde{X_{i}})-X_{i}^{\top}\pi\}]=0.
\]

Since $\alpha_{0}(\cdot)\in\mathcal{M}_{C,\xi_{2}}$ for some $C,\xi_{2}>0$,
we have $\mathbb{E}(\alpha_{0}(\tilde{X}_{i})-X_{i}^{\top}\pi)^{2}\leq(Cp^{-\xi_{2}})^{2}$.
Moreover, $\alpha_{0}(\cdot)\in\mathcal{M}_{C,\xi_{2}}$ also implies
that $\mathbb{E}(\alpha_{0}(\tilde{X}_{i})-X_{i}^{\top}\tilde{\pi})^{2}\leq(Cs_{0}^{-\xi_{2}})^{2}$
for some $\tilde{\pi}\in\RR^{p}$ such that for some $J_{0}\subset\{1,...,p\}$
with $|J_{0}|=s_{0}$ and $\tilde{\pi}_{j}=0$ for $j\notin J_{0}$.
Therefore, 
$$\mathbb{E}(X_{i}^{\top}\pi-X_{i}^{\top}\tilde{\pi})^{2}\leq2(Cp^{-\xi_{2}})^{2}+2(Cs^{-\xi_{2}})^{2}\leq4(Cs^{-\xi_{2}})^{2}\lesssim\varepsilon_{n}^{2\xi_{2}/(2\xi_{2}+1)}.$$
Thus, 
\begin{equation}
(\pi-\tilde{\pi})^{\top}\Sigma(\pi-\tilde{\pi})\leq C\varepsilon_{n}^{4\xi_{2}/(2\xi_{2}+1)}.\label{pistar2D}
\end{equation}
Since the eigenvalues of $\Sigma$ are bounded, we have 
\[
\left\Vert \pi-\tilde{\pi}\right\Vert _{2}\leq C\varepsilon_{n}^{2\xi_{2}/(2\xi_{2}+1)}.
\]
We define $\pi_{\ast}$ as 
\begin{equation}
\pi_{\ast}\in\arg\min_{v}\;(\pi-v)^{\top}\Sigma(\pi-v)+2\varepsilon_{n}\sum_{j\in J_{0}^{c}}|v_{j}|\text{.}\label{pistarD}
\end{equation}

\jelenax{Let $J$ be the vector of indices of nonzero elements of $\pi_{\ast}.$
Define $\gamma=(\theta,\beta^{\top})^{\top}$, where $\beta$ is defined
as in the main text. We notice that 
\[
\mathbb{E}Y_{i}X_{i}=\mathbb{E}(X_{i}^{\top}\gamma+\eta_{f,i}+\varepsilon_{i})X_{i}=\Sigma\gamma+\mathbb{E}\eta_{f,i}X_{i}=\Sigma\gamma+\begin{pmatrix}\mathbb{E}\eta_{f,i}D_{i}\\
\mathbb{E}\eta_{f,i}Z_{i}
\end{pmatrix}=\Sigma\gamma+\begin{pmatrix}\mathbb{E}\eta_{f,i}\eta_{g,i}\\
0
\end{pmatrix},
\]
where (i) follows by $\mathbb{E}\eta_{f,i}Z_{i}=0$ (due to the definition
of $\beta$) and $\mathbb{E}\eta_{f,i}D_{i}=\mathbb{E}\eta_{f,i}u_{i}+\mathbb{E}\eta_{f,i}Z_{i}^{\top}\phi+\mathbb{E}\eta_{f,i}\eta_{g,i}=\mathbb{E}\eta_{f,i}\eta_{g,i}$
(due to the definition of $\beta$, which implies $\mathbb{E}\eta_{f,i}Z_{i}=0$).
Then 
\begin{equation}
\mathbb{E}Y_{i}X_{i}=\Sigma\gamma+e_{1}\varsigma,\label{eq: plm gamma}
\end{equation}
where $\varsigma=\mathbb{E}\eta_{f,i}\eta_{g,i}$ and $e_{1}\in\RR^{p+1}$
is the first column of the $(p+1)\times(p+1)$ identity matrix. Since
$p\gg\max\{n^{1/(2\xi_{1})},n^{1/(2\xi_{2})}\}$, we have 
\begin{equation}
|\varsigma|\leq\sqrt{\mathbb{E}\eta_{f,i}^{2}}\cdot\sqrt{\mathbb{E}\eta_{g,i}^{2}}=\sqrt{O(p^{-2\xi_{1}})}\cdot\sqrt{O(p^{-2\xi_{2}})}=o(n^{-1/2})=o(\varepsilon_{n}).\label{eq: plm gamma 2}
\end{equation}
We will use the approximation $g(\tilde{Z}_{i})=Z_{i}^{\top}\phi+\eta_{g,i}$. }

We first prove some preliminary results (Lemmas \ref{lem: PLM tool 1}-\ref{lem: PLM tool 5})
and then present the proof of Theorem \ref{thm: PLM main}.
\begin{lemma}
\label{lem: PLM tool 1}Let Assumption \ref{assu: PLM} hold. Then
the following hold:\\
 (1) $\|\Sigma(\pi_{*}-\pi)\|_{\infty}\leq\varepsilon_{n}$, $(\pi-\pi_{*})^{\top}\Sigma(\pi-\pi_{*})\lesssim\varepsilon_{n}^{4\xi_{2}/(2\xi_{2}+1)}$
and $|J|\lesssim\varepsilon_{n}^{-2/(2\xi_{2}+1)}$.\\
 (2) If $\xi_{2}>1/2$, then $\|\pi_{*}-\pi\|_{1}\lesssim\varepsilon_{n}^{(2\xi_{2}-1)/(2\xi_{2}+1)}$. 
\end{lemma}

\begin{proof}[Proof of Lemma \ref{lem: PLM tool 1}]
This follows by the same argument as in Lemmas \ref{lem: B1}, \ref{lem: B2},
\ref{lem: B3} and \ref{lem: B4}. 
\end{proof}
\begin{lemma}
\label{lem: PLM tool 2}Let Assumption \ref{assu: PLM} hold. Then
with probability at least $1-o(1)$, we have $\|\hat{\Sigma}_{(1)}\pi_{*}-\Sigma\pi_{*}\|_{\infty}\leq C_{1}\varepsilon_{n}$,
$\|\mathbb{E}_{n,2}Z_{i}Y_{i}-\mathbb{E}Z_{i}Y_{i}\|_{\infty}\leq C_{1}\varepsilon_{n}$,
$\|\mathbb{E}_{n,2}Z_{i}D_{i}-\mathbb{E}Z_{i}D_{i}\|_{\infty}\leq C_{1}\varepsilon_{n}$
and $\|\mathbb{E}_{n,2}Z_{i}D_{i}-\mathbb{E}_{n,1}Z_{i}D_{i}\|_{\infty}\leq C_{1}\varepsilon_{n}$
for some constant $C_{1}>0$. 
\end{lemma}

\begin{proof}[Proof of Lemma \ref{lem: PLM tool 2}]
Notice that 
\[
\sqrt{\mathbb{E}(X_{i}^{\top}\pi_{*})^{2}}=\|\Sigma^{1/2}\pi_{*}\|_{2}\leq\|\Sigma^{1/2}\pi\|_{2}+\|\Sigma^{1/2}(\pi-\pi_{*})\|_{2}.
\]
Since 
\begin{equation}
(\pi-\pi_{*})^{\top}\Sigma(\pi-\pi_{*})\lesssim\varepsilon_{n}^{4\xi_{2}/(2\xi_{2}+1)}\label{eq: lem plm tool 2 eq 3}
\end{equation}
(due to Lemma \ref{lem: PLM tool 1}), we have that $\|\Sigma^{1/2}(\pi-\pi_{*})\|_{2}\lesssim\varepsilon_{n}^{2\xi_{2}/(2\xi_{2}+1)}=o(1)$
(by $\xi_{2}>0$). 

\jelenax{By assumption, $\mathbb{E}|Y_{i}|^{2+\kappa}$ is bounded and
entries of $Z_{i}$ are sub-Gaussian. Then Lemma
\ref{lem: C2} yields the bound for $\|\mathbb{E}_{n,2}Z_{i}Y_{i}-\mathbb{E}Z_{i}Y_{i}\|_{\infty}$.}

Similarly, since entries of $Z_{i}$ are sub-Gaussian and $D_{i}$
is sub-Gaussian, the same argument yields the bound for $\|\mathbb{E}_{n,2}Z_{i}D_{i}-\mathbb{E}Z_{i}D_{i}\|_{\infty}$
and $\|\mathbb{E}_{n,1}Z_{i}D_{i}-\mathbb{E}_{n,1}Z_{i}D_{i}\|_{\infty}$.
Hence, $\|\mathbb{E}_{n,2}Z_{i}D_{i}-\mathbb{E}_{n,1}Z_{i}D_{i}\|_{\infty}\lesssim\varepsilon_{n}$
with probability $1-o(1)$. 

\jelenax{In the case of sub-Gaussian $X_{i}$, $X_{i}^{\top}\pi_{*}$ is sub-Gaussian
and thus has bounded moments. Therefore, Lemma \ref{lem: C2}
implies that
\[
\left\Vert \mathbb{E}_{n,1}X_{i}X_{i}^{\top}\pi_{*}-\mathbb{E} X_{i}X_{i}^{\top} \pi_{*}\right\Vert _{\infty}\lesssim\varepsilon_{n}
\]
with probability $1-o(1)$. It remains to show the above for the case
of bounded $X_{i}$. }

\jelenax{Recall that $\alpha_{0}(\tilde{X}_{i})=\sigma_{u}^{-2}D_{i}-\sigma_{u}^{-2}g(\tilde{Z}_{i})=\sigma_{u}^{-2}u_{i}$. Let $\eta_{\alpha,i}=\alpha_{0}(\tilde{X}_{i})-X_{i}^{\top}\pi$. 
Then 
\[
X_{i}^{\top}\pi_{*}-\alpha_{0}(\tilde{X}_{i})=X_{i}^{\top}\pi-\alpha_{0}(\tilde{X}_{i})+X_{i}^{\top}(\pi_{*}-\pi)=-\eta_{\alpha,i}+X_{i}^{\top}(\pi_{*}-\pi).
\]}

\jelenax{By (\ref{eq: lem plm tool 2 eq 3}), we have that $\mathbb{E}(X_{i}^{\top}(\pi-\pi_{*}))^{2}=(\pi-\pi_{*})^{\top}\Sigma(\pi-\pi_{*})\lesssim\varepsilon_{n}^{4\xi_{2}/(2\xi_{2}+1)}\ll1/\log(pn)$.
We also know that $\mathbb{E}\eta_{\alpha,i}^{2} \lesssim p^{-\xi_{2}} \ll1/\log(pn)$ (due to $\alpha_{0}\in\mathcal{M}_{C,\xi_{2}}$ and $p\gg\max\{n^{1/(2\xi_{1})},n^{1/(2\xi_{2})}\}$). Then $\mathbb{E}(X_{i}^{\top}\pi_{*}-\alpha_{0}(\tilde{X}_{i}))^{2}\ll1/\log(pn)$. }

By the boundedness of $\|X_{i}\|_{\infty}$ and Lemma \ref{lem: C3},
\[
\left\Vert \mathbb{E}_{n,1}X_{i}(X_{i}^{\top}\pi_{*}-\alpha_{0}(\tilde{X}_{i}))-\mathbb{E}[X_{i}(X_{i}^{\top}\pi_{*}-\alpha_{0}(\tilde{X}_{i}))\right\Vert _{\infty}\lesssim\varepsilon_{n}
\]
with probability $1-o(1)$. By $\alpha_{0}(\tilde{X}_{i})=\sigma_{u}^{-2}u_{i}$ and the assumption of $\mathbb{E}|u_{i}|^{2+\kappa}$
being bounded, the boundedness of $\|X_{i}\|_{\infty}$,  Lemma
\ref{lem: C2} implies 
\[
\left\Vert \mathbb{E}_{n,1}X_{i}\alpha_{0}(\tilde{X}_{i})-\mathbb{E}X_{i}\alpha_{0}(\tilde{X}_{i})\right\Vert _{\infty}\lesssim\varepsilon_{n}
\]
with probability $1-o(1)$. Therefore,
\[
\left\Vert \mathbb{E}_{n,1}X_{i}X_{i}^{\top}\pi_{*}-\mathbb{E}[X_{i}X_{i}^{\top}\pi_{*}\right\Vert _{\infty}\lesssim\varepsilon_{n}
\]
with probability $1-o(1)$. 
\end{proof}

\begin{lemma}
\label{lem: PLM tool 3}Let Assumption \ref{assu: PLM} hold. Suppose
that $\log p\ll n^{1/2-1/(4\max\{\xi_{1},\xi_{2}\})}$. Then $\hat{\pi}$
is well defined and 
\[
\|\hat{\pi}-\pi\|_{2}\lesssim\varepsilon_{n}^{2\xi_{2}/(2\xi_{2}+1)}
\]
with probability at least $1-o(1)$ and 
\[
(\hat{\pi}-\pi)^{\top}\hat{\Sigma}_{(1)}(\hat{\pi}-\pi)=O_{P}(\varepsilon_{n}^{4\xi_{2}/(2\xi_{2}+1)}).
\]
Moreover, if $\xi_{2}>1/2$, then with probability at least $1-o(1)$,
$\|\hat{\pi}-\pi\|_{1}\lesssim\varepsilon_{n}^{(2\xi_{2}-1)/(2\xi_{2}+1)}$. 
\end{lemma}

\begin{proof}[Proof of Lemma \ref{lem: PLM tool 3}]
Recall the estimator in (\ref{eq: Riesz rep est PLM}): 
\begin{align*}
\hat{\pi}=\arg\min_{v\in\RR^{p+1}} & \|v\|_{1}\qquad s.t.\qquad\|\hat{\Sigma}_{(1)}v-e_{1}\|_{\infty}\leq\lambda_{1}
\end{align*}
where $\lambda_{1}=\lambda_{0}\varepsilon_{n}$ for a large enough
constant $\lambda_{0}>0$. Let $J$ be the vector of indices of nonzero elements of $\pi_{\ast}$
defined in (\ref{pistarD}). By Lemma \ref{lem: PLM tool 1}, $|J|\lesssim\varepsilon_{n}^{-2/(2\xi_{2}+1)}\lesssim(n/\log p)^{1/(2\xi_{2}+1)}$. 
For constants $\kappa_{1},\kappa_{2}>0$ (to be chosen), define the event 
\[
\Acal=\left\{ \|\hat{\Sigma}_{(1)}\pi_{*}-\Sigma\pi_{*}\|_{\infty}\leq C_{1}\varepsilon_{n}\right\} \bigcap\left\{ \min_{|A|\leq \kappa_{1} (n/\log p)^{1/(2\xi_{2}+1)} } \min_{\|v_{A^{c}}\|_{1}\leq\|v_{A}\|_{1}}\frac{v^{\top}\hat{\Sigma}_{(1)}v}{\|v_{A}\|_{2}^{2}}\geq\kappa_{2}\right\} ,
\]
where $C_{1}>0$ is a constant that makes the first part of $\Acal$ holds with probability $1-o(1)$ (this is possible due to  Lemma \ref{lem: PLM tool 2}). 
The rest of two proof consists of two steps. We first show $P(\Acal)=1-o(1)$
and then derive the final result.

\textbf{Step 1:} show that for any $\kappa_{1}>0$, there exists $\kappa_{2}>0$ such that $P(\Acal)\geq1-o(1)$.

We first show this in the case of bounded $X_{i}$. We apply Theorem
22 of \citet{rudelson2013reconstruction}, which requires $n\gg|A|(\log p)\log^{3}(|A|\log p)$.
By  $|A| \lesssim(n/\log p)^{1/(2\xi_{2}+1)}$, we notice that $|A|(\log p)\log^{3}(|A|\log p)\lesssim(n/\log p)^{1/(2\xi_{2}+1)}(\log p)(\log n)^{3}$.
Hence, we only need $\log p\ll n(\log n)^{-3(2\xi_{2}+1)/(2\xi_{2})}$,
which holds due to $\log p\ll n^{1/2-1/(4\max\{\xi_{1},\xi_{2}\})}\ll n^{1/2}$.
Therefore, for bounded $X_{i}$, Theorem 22 of \citet{rudelson2013reconstruction}
implies the second part of $\Acal$ has probability approaching one for some $\kappa_{2}>0$.

For the case of sub-Gaussian $X_{i}$, we apply Theorem 16 of \citet{rudelson2013reconstruction},
which requires $n\gg|A|\log(p/|A|)$. This  holds due to  $|A|\lesssim(n/\log p)^{1/(2\xi_{2}+1)}$ and $\log p \ll n^{1/2}$. Thus, $P(\Acal)\geq1-o(1)$ for some $\kappa_{2}>0$.

\textbf{Step 2:} show the final result.

The rest of the argument applies on the event $\Acal$ with a large $\kappa_{1}$ and a corresponding $\kappa_{2}>$ such that $P(\Acal)=1-o(1)$. 
Notice that $\Sigma\pi=e_{1}$. By Lemma \ref{lem: PLM tool 1}, we
have that for a large constant $\lambda_{0}>0$
\[
\|\hat{\Sigma}_{(1)}\pi_{*}-e_{1}\|_{\infty}=\|\hat{\Sigma}_{(1)}\pi_{*}-\Sigma\pi\|_{\infty}\leq\|\hat{\Sigma}_{(1)}\pi_{*}-\Sigma\pi_{*}\|_{\infty}+\|\Sigma\pi_{*}-\Sigma\pi\|_{\infty}\leq\varepsilon_{n}+C_{1}\varepsilon_{n}\leq\lambda_{0}\varepsilon_{n}=\lambda_{1}.
\]

Hence, the constraint $\|\hat{\Sigma}_{(1)}v-e_{1}\|_{\infty}\leq\lambda_{1}$
is satisfied by $v=\pi_{*}$ so $\hpi$ is well defined and $\|\hpi\|_{1}\leq\|\pi_{*}\|_{1}$.
Let $\Delta=\hat{\pi}-\pi_{*}$. Since $J={\rm supp}(\pi_{*})$, we
have that $\|\pi_{*}+\Delta_{J}\|_{1}+\|\Delta_{J^{c}}\|_{1}\leq\|\pi_{*}\|_{1}$
and thus 
\begin{equation}
\|\Delta_{J^{c}}\|_{1}\leq\|\Delta_{J}\|_{1}.\label{eq: lem PLM tool 3 eq 5}
\end{equation}

Since $\|\hat{\Sigma}_{(1)}\pi_{*}-e_{1}\|_{\infty}\leq\lambda_{1}$
and $\|\hat{\Sigma}_{(1)}\hphi-e_{1}\|_{\infty}\leq\lambda_{1}$,
we have that $\|\hat{\Sigma}_{(1)}\Delta\|_{\infty}\leq2\lambda_{1}$.
Then 
\begin{equation}
\Delta^{\top}\hat{\Sigma}_{(1)}\Delta\leq\|\Delta\|_{1}\|\hat{\Sigma}_{(1)}\Delta\|_{\infty}\overset{\text{(i)}}{\leq}2\sqrt{|J|}\|\Delta_{J}\|_{2}\cdot2\lambda_{1},\label{eq: lem PLM tool 3 eq 6}
\end{equation}
where (i) follows by $\|\Delta\|_{1}=\|\Delta_{J}\|_{1}+\|\Delta_{J^{c}}\|_{1}\leq2\|\Delta_{J}\|_{1}\leq2\sqrt{|J|}\|\Delta_{J}\|_{2}$
due to (\ref{eq: lem PLM tool 3 eq 5}). From the definition of $\Acal$,
we have 
\[
\Delta^{\top}\hat{\Sigma}_{(1)}\Delta\geq\kappa_{2}\|\Delta_{J}\|_{2}^{2}.
\]

The above two displays imply 
\begin{equation}
\|\Delta_{J}\|_{2}\leq4\kappa_{2}^{-1}\lambda_{1}\sqrt{|J|}.\label{eq: lem PLM tool 3 eq 7}
\end{equation}

Let $B$ denote the indices corresponding to the largest (in magnitude)
$|J|$ entries in $\Delta_{J^{c}}$. Then $|B\bigcup J|=2\cdot|J|$
and $\|\Delta_{(B\bigcup J)^{c}}\|_{1}\leq\|\Delta_{J^{c}}\|_{1}\leq\|\Delta_{J}\|_{1}\leq\|\Delta_{B\bigcup J}\|_{1}$.
Thus, we can repeat the same argument as above and obtain 
\[
\|\Delta_{B\bigcup J}\|_{2}\leq4\kappa_{2}^{-1}\lambda_{1}\sqrt{|B\bigcup J|}=4\kappa_{2}^{-1}\lambda_{1}\sqrt{2|J|}.
\]

By Lemma 6.9 of \citet{buhlmann2011statistics}, we have that 
\begin{equation}
\|\Delta_{(B\bigcup J)^{c}}\|_{2}\leq|B\bigcup J|^{-1/2}\|\Delta\|_{1}.\label{eq: lem PLM tool 3 eq 8}
\end{equation}

Since $\|\Delta\|_{1}=\|\Delta_{J}\|_{1}+\|\Delta_{J^{c}}\|_{1}\leq2\|\Delta_{J}\|_{1}$,
(\ref{eq: lem PLM tool 3 eq 7}) implies that 
\begin{equation}
\|\Delta\|_{1}\leq2\|\Delta_{J}\|_{1}\leq2|J|^{1/2}\|\Delta_{J}\|_{2}\leq8\kappa_{2}^{-1}\lambda_{1}|J|.\label{eq: lem PLM tool 3 eq 9}
\end{equation}

The above two displays and $|B\bigcup J|=2\cdot|J|$ imply 
\[
\|\Delta_{(B\bigcup J)^{c}}\|_{2}\leq2|B\bigcup J|^{-1/2}\|\Delta_{J}\|_{1}\leq16|B\bigcup J|^{-1/2}\kappa_{2}^{-1}\lambda_{1}|J|=8\sqrt{2}\kappa_{2}^{-1}\lambda_{1}\sqrt{|J|}.
\]

Hence, (\ref{eq: lem PLM tool 3 eq 8}) and the above display imply
\begin{equation}
\|\Delta\|_{2}\leq\|\Delta_{B\bigcup J}\|_{2}+\|\Delta_{(B\bigcup J)^{c}}\|_{2}\leq4\kappa_{2}^{-1}\lambda_{1}\sqrt{2|J|}+8\sqrt{2}\kappa_{2}^{-1}\lambda_{1}\sqrt{|J|}=12\sqrt{2}\kappa_{2}^{-1}\lambda_{1}\sqrt{|J|}.\label{eq: lem PLM tool 3 eq 10}
\end{equation}

Since $|J|\lesssim\varepsilon_{n}^{-2/(2\xi_{2}+1)}$ and $\lambda_{1}=\lambda_{0}\varepsilon_{n}$,
the above display implies 
\[
\|\Delta\|_{2}\lesssim\varepsilon_{n}^{2\xi_{2}/(2\xi_{2}+1)}
\]
and (\ref{eq: lem PLM tool 3 eq 9}) implies 
\[
\|\Delta\|_{1}\lesssim\varepsilon_{n}^{(2\xi_{2}-1)/(2\xi_{2}+1)}.
\]
If $\xi_{2}>1/2$, then $\|\hpi-\pi\|_{1}\leq\|\hpi-\pi_{*}\|_{1}+\|\pi_{*}-\pi\|_{1}=\|\Delta\|_{1}+\|\pi_{*}-\pi\|_{1}\lesssim\varepsilon_{n}^{(2\xi_{2}-1)/(2\xi_{2}+1)}$
due to $\|\pi_{*}-\pi\|_{1}\lesssim\varepsilon_{n}^{(2\xi_{2}-1)/(2\xi_{2}+1)}$
(Lemma \ref{lem: PLM tool 1}).

We now bound $\|\hpi-\pi\|_{2}$. Since we already showed that $\|\hpi-\pi_{*}\|_{2}=\|\Delta\|_{2}\lesssim\varepsilon_{n}^{2\xi_{2}/(2\xi_{2}+1)}$,
it remains to show that $\|\pi_{*}-\pi\|_{2}\lesssim\varepsilon_{n}^{2\xi_{2}/(2\xi_{2}+1)}$.
This is true because $(\pi-\pi_{*})^{\top}\Sigma(\pi-\pi_{*})\lesssim\varepsilon_{n}^{4\xi_{2}/(2\xi_{2}+1)}$
(due to Lemma \ref{lem: PLM tool 1}) and the eigenvalues of $\Sigma$
are bounded below by $\kappa_{3}>0$.

Finally, we bound $(\hat{\pi}-\pi)^{\top}\hat{\Sigma}_{(1)}(\hat{\pi}-\pi)$.
By (\ref{eq: lem PLM tool 3 eq 6}) and (\ref{eq: lem PLM tool 3 eq 7}),
we have that 
\[
\Delta^{\top}\hat{\Sigma}_{(1)}\Delta\leq4\sqrt{|J|}\|\Delta_{J}\|_{2}\lambda\leq16\kappa_{2}^{-1}\lambda^{2}|J|\lesssim\varepsilon_{n}^{2}\cdot\varepsilon_{n}^{-2/(2\xi_{2}+1)}\lesssim\varepsilon_{n}^{4\xi_{2}/(2\xi_{2}+1)}
\]
and thus 
\begin{align*}
 & (\hat{\pi}-\pi)^{\top}\hat{\Sigma}_{(1)}(\hat{\pi}-\pi)\\
 & =[\Delta+(\pi_{*}-\pi)]^{\top}\hat{\Sigma}_{(1)}[\Delta+(\pi_{*}-\pi)]\\
 & \leq2\Delta^{\top}\hat{\Sigma}_{(1)}\Delta+2(\pi_{*}-\pi)^{\top}\hat{\Sigma}_{(1)}(\pi_{*}-\pi)\\
 & =O_{P}(\varepsilon_{n}^{4\xi_{2}/(2\xi_{2}+1)})+2(\pi_{*}-\pi)^{\top}\hat{\Sigma}_{(1)}(\pi_{*}-\pi)\overset{\text{(i)}}{=}O_{P}(\varepsilon_{n}^{4\xi_{2}/(2\xi_{2}+1)}),
\end{align*}
where (i) follows by 
\[
\mathbb{E}(\pi_{*}-\pi)^{\top}\hat{\Sigma}_{(1)}(\pi_{*}-\pi)=(\pi_{*}-\pi)^{\top}\Sigma(\pi_{*}-\pi)\lesssim\varepsilon_{n}^{4\xi_{2}/(2\xi_{2}+1)}.
\]

The proof is complete. 
\end{proof}
\begin{lemma}
\label{lem: PLM tool 4}Recall $\hat{\phi}$ from (\ref{eq: Dantzig_PLM_part_0}):
\[
\hat{\phi}=\arg\min_{v\in\RR^{p}}\|v\|_{1}\qquad s.t.\qquad\|\mathbb{E}_{n,1}Z_{i}Z_{i}^{\top}v-\mathbb{E}_{n,2}Z_{i}D_{i}\|_{\infty}\leq\lambda_{1},
\]
where $\lambda_{1}=\lambda_{0}\varepsilon_{n}$. Let Assumption \ref{assu: PLM}
hold. Suppose that $\log p\ll n^{1/2-1/(4\max\{\xi_{1},\xi_{2}\})}$.
Then for large enough $\lambda_{0}>0$, $\hat{\phi}$ is well defined
and $\|\mathbb{E}_{n,1}Z_{i}Z_{i}^{\top}(\hat{\phi}-\phi)\|_{\infty}\leq2\lambda_{1}$
with probability at least $1-o(1)$ and 
\[
(\hat{\phi}-\phi)^{\top}\mathbb{E}_{n,1}[Z_{i}Z_{i}^{\top}](\hat{\phi}-\phi)=O_{P}(\varepsilon_{n}^{4\xi_{2}/(2\xi_{2}+1)}).
\]
Moreover, if $\xi_{2}>1/2$, then $\|\hat{\phi}-\phi\|_{1}\lesssim\varepsilon_{n}^{(2\xi_{2}-1)/(2\xi_{2}+1)}$
and $\|\hat{\phi}-\phi\|_{2}\lesssim\varepsilon_{n}^{2\xi_{2}/(2\xi_{2}+1)}$
with probability at least $1-o(1)$. 
\end{lemma}

\begin{proof}[Proof of Lemma \ref{lem: PLM tool 4}]
This follows by the same argument for Dantzig selector as the proof
of Lemma \ref{lem: PLM tool 3}. We only need to verify 
\begin{equation}
P\left(\|\mathbb{E}_{n,1}Z_{i}Z_{i}^{\top}\phi-\mathbb{E}_{n,2}Z_{i}D_{i}\|_{\infty}\leq\lambda\right)\geq1-o(1)\label{eq: lem PLM tool 4 eq 2}
\end{equation}
and $\mathbb{E}_{n,1}Z_{i}Z_{i}^{\top}$ satisfies the restricted
eigenvalue condition. Notice that 
\[
\hat{\Sigma}_{(1)}=\begin{pmatrix}\mathbb{E}_{n,1}D_{i}^{2} & \mathbb{E}_{n,1}D_{i}Z_{i}^{\top}\\
\mathbb{E}_{n,1}D_{i}Z_{i} & \mathbb{E}_{n,1}Z_{i}Z_{i}^{\top}
\end{pmatrix}.
\]

In the proof of Lemma \ref{lem: PLM tool 3}, we show that $\hat{\Sigma}_{(1)}$
satisfies the restricted eigenvalue condition. This means that $\mathbb{E}_{n,1}Z_{i}Z_{i}^{\top}$
also satisfies the restricted eigenvalue condition. It remains to
show (\ref{eq: lem PLM tool 4 eq 2}). We do so in two cases: bounded
$X_{i}$ and sub-Gaussian $X_{i}$.

\textbf{Case 1:} show (\ref{eq: lem PLM tool 4 eq 2}) for bounded
$X_{i}$.

By Lemma \ref{lem: PLM tool 2}, $P(\|\mathbb{E}_{n,2}Z_{i}D_{i}-\mathbb{E}Z_{i}D_{i}\|_{\infty}\leq C_{1}\varepsilon_{n})\geq1-o(1)$.
By the same argument as Lemma \ref{lem: PLM tool 2}, we have
that $\|\mathbb{E}_{n,1}Z_{i}Z_{i}^{\top}\phi-\mathbb{E}Z_{i}Z_{i}^{\top}\phi\|_{\infty}\leq C_{2}\varepsilon_{n}$
with probability approaching one for some constant $C_{2}>0$. Hence,
with probability approaching one, 
\begin{multline*}
\|\mathbb{E}_{n,1}Z_{i}Z_{i}^{\top}\phi-\mathbb{E}_{n,2}Z_{i}D_{i}\|_{\infty}\leq\|\mathbb{E}_{n,1}Z_{i}Z_{i}^{\top}\phi-\mathbb{E}Z_{i}Z_{i}^{\top}\phi\|_{\infty}+\|\mathbb{E}_{n,2}Z_{i}D_{i}-\mathbb{E}Z_{i}D_{i}\|_{\infty}\\
+\|\mathbb{E}Z_{i}Z_{i}^{\top}\phi-\mathbb{E}Z_{i}D_{i}\|_{\infty}\overset{\text{(i)}}{\leq}C_{1}\varepsilon_{n}+C_{2}\varepsilon_{n},
\end{multline*}
where (i) follows by $\mathbb{E}Z_{i}Z_{i}^{\top}\phi-\mathbb{E}Z_{i}D_{i}=\mathbb{E}Z_{i}[Z_{i}^{\top}\phi-g(\tilde{Z}_{i})]-\mathbb{E}Z_{i}u_{i}=0$
(due to $\mathbb{E}Z_{i}[Z_{i}^{\top}\phi-g(\tilde{Z}_{i})]$ by the
definition of $\phi$). Since $\lambda_{1}=\lambda_{0}\varepsilon_{n}$
with a large constant $\lambda_{0}$, (\ref{eq: lem PLM tool 4 eq 2})
follows. 

\jelenax{\textbf{Case 2:} show (\ref{eq: lem PLM tool 4 eq 2}) for sub-Gaussian
$X_{i}$.}

\jelenax{Again, by $\mathbb{E}Z_{i}Z_{i}^{\top}\phi-\mathbb{E}Z_{i}D_{i}=0$,
we have 
\[
\|\mathbb{E}_{n,1}Z_{i}Z_{i}^{\top}\phi-\mathbb{E}_{n,2}Z_{i}D_{i}\|_{\infty}\leq\|\mathbb{E}_{n,1}Z_{i}Z_{i}^{\top}\phi-\mathbb{E}Z_{i}Z_{i}^{\top}\phi\|_{\infty}+\|\mathbb{E}_{n,2}Z_{i}D_{i}-\mathbb{E}Z_{i}D_{i}\|_{\infty}.
\]}

\jelenax{By Lemma \ref{lem: PLM tool 2}, $P(\|\mathbb{E}_{n,2}Z_{i}D_{i}-\mathbb{E}Z_{i}D_{i}\|_{\infty}\leq C_{1}\varepsilon_{n})\geq1-o(1)$.
Since $Z_{i}$ is a sub-vector of $X_{i}$, which is sub-Gaussian,
$Z_{i}^{\top}\phi$ is sub-Gaussian and hence has bounded third moments.
By Lemma \ref{lem: C2}, $\|\mathbb{E}_{n,1}Z_{i}Z_{i}^{\top}\phi-\mathbb{E}Z_{i}Z_{i}^{\top}\phi\|_{\infty}\leq C\varepsilon_{n}$
with probability $1-o(1)$ for a large enough $C>0$. Hence, (\ref{eq: lem PLM tool 4 eq 2})
follows for a large enough constant $\lambda_{0}$. }
\end{proof}
\begin{lemma}
\label{lem: PLM tool 5}
\jelenax{Recall  $\hat{\gamma}=(\hat{\theta},\hat{\beta}^{\top})^{\top}$
from (\ref{eq: Dantzig_PLM}): 
\begin{align*}
\hat{\gamma}=(\hat{\theta},\hat{\beta}^{\top})^{\top}=\arg\min_{b\in\RR,g\in\RR^{p}} & \|(b,g^{\top})^{\top}\|_{1}\\
s.t. & \|\mathbb{E}_{n,2}Z_{i}D_{i}b+\mathbb{E}_{n,1}Z_{i}Z_{i}^{\top}g-\mathbb{E}_{n,2}Z_{i}Y_{i}\|_{\infty}\leq\lambda_{1}\\
 & |\mathbb{E}_{n,2}D_{i}^{2}b+\mathbb{E}_{n,2}D_{i}Z_{i}^{\top}g-\mathbb{E}_{n,2}D_{i}Y_{i}|\leq\lambda_{1}\\
 & |\mathbb{E}_{n,2}D_{i}^{2}b+\hat{\phi}^{\top}\mathbb{E}_{n,1}Z_{i}Z_{i}^{\top}g-\mathbb{E}_{n,2}D_{i}Y_{i}|\leq\lambda_{2},
\end{align*}
where $\lambda_{1}=\lambda_{0}\varepsilon_{n}$ and $\lambda_{2}=\lambda_{n}n^{-1/4}$
for a large $\lambda_{0}>0$. Let Assumption \ref{assu: PLM} hold.
Suppose that that $\max\{\xi_{1},\xi_{2}\}>1/2$ and $\log p\ll n^{1/2-1/(4\max\{\xi_{1},\xi_{2}\})}$.
Then $\hat{\gamma}$ is well defined with probability at least $1-o(1)$,
\[
\|\hat{\gamma}-\gamma\|_{2}=O_{P}(\varepsilon_{n}^{2\xi_{1}/(2\xi_{1}+1)}),
\]
\[
(\hat{\beta}-\beta)^{\top}\mathbb{E}_{n,1}[Z_{i}Z_{i}^{\top}](\hat{\beta}-\beta)=O_{P}(\varepsilon_{n}^{4\xi_{1}/(2\xi_{1}+1)})
\]
and $\|\mathbb{E}_{n,2}Z_{i}D_{i}(\hat{\theta}-\theta)+\mathbb{E}_{n,1}Z_{i}Z_{i}^{\top}(\hat{\beta}-\beta)\|_{\infty}=O_{P}(\varepsilon_{n})$.
Moreover, if $\xi_{1}>1/2$, then 
\[
\|\hat{\gamma}-\gamma\|_{1}=O_{P}\left(\varepsilon_{n}^{(2\xi_{1}-1)/(2\xi_{1}+1)}\right).
\]
Furthermore, if $\xi_{2}>1/2$, then $|\mathbb{E}_{n,2}D_{i}^{2}(\hat{\theta}-\theta)+\hat{\phi}^{\top}\mathbb{E}_{n,1}Z_{i}Z_{i}^{\top}(\hat{\beta}-\beta)|=O_{P}(n^{-1/4})$. }
\end{lemma}

\begin{proof}[Proof of Lemma \ref{lem: PLM tool 5}]
By the same argument as in Lemma \ref{lem: PLM tool 1}, we can show
that there exists $\gamma_{*}$ with $J_{*}={\rm supp}(\gamma_{*})$
such that (1) $\|\Sigma(\gamma_{*}-\gamma)\|_{\infty}\leq2\varepsilon_{n}$,
$(\gamma-\gamma_{*})^{\top}\Sigma(\gamma-\gamma_{*})\lesssim\varepsilon_{n}^{4\xi_{1}/(2\xi_{1}+1)}$
and $|J_{*}|\lesssim\varepsilon_{n}^{-2/(2\xi_{1}+1)}$ and (2) if
$\xi_{1}>1/2$, $\|\gamma_{*}-\gamma\|_{1}\lesssim\varepsilon_{n}^{(2\xi_{1}-1)/(2\xi_{1}+1)}$.
Notice that the only difference is that now we have $\mathbb{E}Y_{i}X_{i}=\Sigma\gamma+e_{1}\varsigma$
from (\ref{eq: plm gamma}), rather than $\mathbb{E}Y_{i}X_{i}=\Sigma\gamma$;
however, since $|\varsigma|\ll\varepsilon_{n}$ (due to (\ref{eq: plm gamma 2})),
all the calculations involved in Lemma \ref{lem: PLM tool 1} hold,
potentially at the cost of larger constants. Partition $\gamma_{*}=(\theta_{*},\beta_{*}^{\top})^{\top}$.
Notice that $|J_{*}|\lesssim\varepsilon_{n}^{-2/(2\xi_{1}+1)}\ll\varepsilon_{n}^{-2}$
(due to $\xi_{1}>0$). For constants $\kappa_{1},\kappa_{2}>0$ (to be chosen), we define the event $\Acal=\Acal_{1}\bigcap\Acal_{2}$,
where 
\[
\Acal_{1}=\left\{ \text{Constraints of }\hat{\gamma}\ \text{are satisfied with}\ \gamma_{*}\right\} 
\]
and 
\[
\Acal_{2}=\left\{ \min_{|A|\leq \kappa_{1} (n/\log p)^{1/(2\xi_{2}+1)} } \min_{\|v_{A^{c}}\|_{1}\leq\|v_{A}\|_{1}} \frac{v^{\top}\tilde{\Sigma}v}{\|v_{A}\|_{2}^{2}}\geq\kappa_{2} \right\} \quad{\rm with}\quad\tilde{\Sigma}=\begin{pmatrix}\mathbb{E}_{n,2}D_{i}^{2} & \mathbb{E}_{n,2}D_{i}Z_{i}^{\top}\\
\mathbb{E}_{n,2}Z_{i}D_{i} & \mathbb{E}_{n,1}Z_{i}Z_{i}^{\top}
\end{pmatrix}.
\]

Notice that 
\begin{align*}
\hat{\gamma}=(\hat{\theta},\hat{\beta}^{\top})^{\top}=\arg\min_{b\in\RR,v\in\RR^{p}} & \|(b,v^{\top})^{\top}\|_{1}\\
s.t. & \left\Vert \tilde{\Sigma}\begin{pmatrix}b\\
v
\end{pmatrix}\right\Vert _{\infty}\leq\lambda_{1}\\
 & |\mathbb{E}_{n,2}D_{i}^{2}b+\hat{\phi}^{\top}\mathbb{E}_{n,1}Z_{i}Z_{i}^{\top}v-\mathbb{E}_{n,2}D_{i}Y_{i}|\leq\lambda_{2}.
\end{align*}

The rest of the proof proceeds in four steps. We first show two preliminary
results in two steps and then present the proof of the final bounds
in additional steps.

\jelenax{\textbf{Step 1:} show that $P(\Acal_{1})\geq1-o(1)$.}

\jelenax{We can follow the same argument as in Lemma \ref{lem: PLM tool 2} and show $P(\|\mathbb{E}_{n,1}Z_{i}X_{i}^{\top}\gamma_{*}-\mathbb{E}Z_{i}X_{i}^{\top}\gamma_{*}\|_{\infty}\leq C_{2}\varepsilon_{n})=1-o(1)$ for some $C_{2}>0$.
Then with probability approaching one, 
\begin{align*}
 & \|\mathbb{E}_{n,2}Z_{i}D_{i}\theta_{*}+\mathbb{E}_{n,1}Z_{i}Z_{i}^{\top}\beta_{*}-\mathbb{E}_{n,2}Z_{i}Y_{i}\|_{\infty}\\
 & =\|\mathbb{E}_{n,2}Z_{i}D_{i}\theta_{*}-\mathbb{E}_{n,1}Z_{i}D_{i}\theta_{*}+\mathbb{E}_{n,1}Z_{i}D_{i}\theta_{*}+\mathbb{E}_{n,1}Z_{i}Z_{i}^{\top}\beta_{*}-\mathbb{E}_{n,2}Z_{i}Y_{i}\|_{\infty}\\
 & =\|\mathbb{E}_{n,2}Z_{i}D_{i}\theta_{*}-\mathbb{E}_{n,1}Z_{i}D_{i}\theta_{*}+\mathbb{E}_{n,1}Z_{i}X_{i}^{\top}\gamma_{*}-\mathbb{E}_{n,2}Z_{i}Y_{i}\|_{\infty}\\
 & \leq|\theta_{*}|\cdot\|\mathbb{E}_{n,2}Z_{i}D_{i}-\mathbb{E}_{n,1}Z_{i}D_{i}\|_{\infty}+\|\mathbb{E}_{n,1}Z_{i}X_{i}^{\top}\gamma_{*}-\mathbb{E}Z_{i}X_{i}^{\top}\gamma_{*}\|_{\infty}\\
 & \quad+\|\mathbb{E}_{n,2}Z_{i}Y_{i}-\mathbb{E}Z_{i}Y_{i}\|_{\infty}+\|\mathbb{E}Z_{i}X_{i}^{\top}\gamma_{*}-\mathbb{E}Z_{i}Y_{i}\|_{\infty}\\
 & \overset{\text{(i)}}{\leq}O(1)\cdot C_{1}\varepsilon_{n}+C_{2}\varepsilon_{n}+C_{1}\varepsilon_{n}+\|\mathbb{E}Z_{i}X_{i}^{\top}\gamma_{*}-\mathbb{E}Z_{i}Y_{i}\|_{\infty}\\
 &  \overset{\text{(ii)}}{=}O(1)\cdot C_{1}\varepsilon_{n}+C_{2}\varepsilon_{n}+C_{1}\varepsilon_{n}+\|\mathbb{E}Z_{i}X_{i}^{\top}\gamma_{*}-\mathbb{E}Z_{i}X_{i}^{\top}\gamma\|_{\infty}\\
 & =O(1)\cdot C_{1}\varepsilon_{n}+C_{2}\varepsilon_{n}+C_{1}\varepsilon_{n}+\max_{2\leq j\leq p+1}\left|\mathbb{E}X_{i,j}X_{i}^{\top}\gamma_{*}-\mathbb{E}X_{i,j}X_{i}^{\top}\gamma\right|\\
 & \leq O(1)\cdot C_{1}\varepsilon_{n}+C_{2}\varepsilon_{n}+C_{1}\varepsilon_{n}+\|\Sigma(\gamma_{*}-\gamma)\|_{\infty}\\
 & \overset{\text{(iii)}}{\leq}O(1)\cdot C_{1}\varepsilon_{n}+C_{2}\varepsilon_{n}+C_{1}\varepsilon_{n}+2\varepsilon_{n},
\end{align*}
where (i) follows by $\|\mathbb{E}_{n,1}Z_{i}X_{i}^{\top}\gamma_{*}-\mathbb{E}Z_{i}X_{i}^{\top}\gamma_{*}\|_{\infty}\leq C_{2}\varepsilon_{n}$,
Lemma \ref{lem: PLM tool 2} and $|\theta_{*}|\leq\|\gamma_{*}\|_{2}\lesssim1$, (ii) follows by $Y_{i}=\varepsilon_{i}+X_{i}^{\top}\gamma+\eta_{f,i}$ and $\E Z_{i} \eta_{f,i}=0$ (due to the definition of $\beta$) and  (iii) follows by $\|\Sigma(\gamma_{*}-\gamma)\|_{\infty}\leq2\varepsilon_{n}$.
Since $\lambda_{1}=\lambda_{0}\varepsilon_{n}$ with a large constant
$\lambda_{0}>0$, we have that 
\[
P\left(\|\mathbb{E}_{n,2}Z_{i}D_{i}\theta_{*}+\mathbb{E}_{n,1}Z_{i}Z_{i}^{\top}\beta_{*}-\mathbb{E}_{n,2}Z_{i}Y_{i}\|_{\infty}\leq\lambda_{1}\right)\geq1-o(1).
\]}

\jelenax{Moreover, we have 
\begin{align*}
 & \left|\mathbb{E}_{n,2}D_{i}^{2}\theta_{*}+\mathbb{E}_{n,2}D_{i}Z_{i}^{\top}\beta_{*}-\mathbb{E}_{n,2}D_{i}Y_{i}\right|\\
 & =\left|\mathbb{E}_{n,2}D_{i}(D_{i}\theta_{*}+Z_{i}^{\top}\beta_{*}-Y_{i})\right|\\
 & =\left|\mathbb{E}_{n,2}D_{i}(X_{i}^{\top}\gamma_{*}-X_{i}^{\top}\gamma-\varepsilon_{i}-\eta_{f,i})\right|\\
 & \overset{\text{(i)}}{\leq}\left|\mathbb{E}_{n,2}D_{i}X_{i}^{\top}(\gamma_{*}-\gamma)\right|+O_{P}(n^{-1/2})\\
 & \leq\left|\mathbb{E}D_{i}X_{i}^{\top}(\gamma_{*}-\gamma)\right|+\left|\left(\mathbb{E}_{n,2}D_{i}X_{i}-\mathbb{E}D_{i}X_{i}\right)^{\top}(\gamma_{*}-\gamma)\right|+O_{P}(n^{-1/2})\\
 & \overset{\text{(ii)}}{\leq}\left|\mathbb{E}D_{i}X_{i}^{\top}(\gamma_{*}-\gamma)\right|+o_{P}(n^{-1/2})+O_{P}(n^{-1/2})\\
 & \overset{\text{(iii)}}{\leq}\|\Sigma(\gamma_{*}-\gamma)\|_{\infty}+o_{P}(n^{-1/2})+O_{P}(n^{-1/2})\overset{\text{(iv)}}{\leq}2\varepsilon_{n}+o_{P}(n^{-1/2})+O_{P}(n^{-1/2}),
\end{align*}
where (i) follows by 
$$\mathbb{E}\left(\mathbb{E}_{n,2}D_{i}\varepsilon_{i}\right)^{2}=n_{2}^{-2}\sum_{i\in I_{2}}\mathbb{E}D_{i}^{2}\varepsilon_{i}^{2}\lesssim n_{2}^{-2}\sum_{i\in I_{2}}\mathbb{E}\varepsilon_{i}^{2}\lesssim n^{-1}$$
and $\left(\mathbb{E}_{n,2}D_{i}\eta_{f,i}\right)^{2}\leq[\mathbb{E}_{n,2}D_{i}^{2}]\cdot[\mathbb{E}_{n,2}\eta_{f,i}^{2}]=O_{P}(1)\cdot O_{P}(p^{-2\xi_{1}})$
together with $p\gg n^{1/(2\xi_{1})}$, (ii) follows by 
\begin{eqnarray*}
&&\mathbb{E}[\left(\mathbb{E}_{n,2}D_{i}X_{i}-\mathbb{E}D_{i}X_{i}\right)^{\top}(\gamma_{*}-\gamma)]^{2}
\\
&=&n_{2}^{-1}\mathbb{E}[\left(D_{i}X_{i}-\mathbb{E}D_{i}X_{i}\right)^{\top}(\gamma_{*}-\gamma)]^{2}
\\
&\leq &n_{2}^{-1}\mathbb{E}[D_{i}X_{i}^{\top}(\gamma_{*}-\gamma)]^{2}\leq n_{2}^{-1} \mathbb{E} (D_{i}^2) \mathbb{E}[X_{i}^{\top}(\gamma_{*}-\gamma)]^{2}\ll n^{-1}
\end{eqnarray*}
due to $\mathbb{E}[X_{i}^{\top}(\gamma_{*}-\gamma)]^{2}=(\gamma-\gamma_{*})^{\top}\Sigma(\gamma-\gamma_{*})\lesssim\varepsilon_{n}^{4\xi_{1}/(2\xi_{1}+1)}=o(1)$,
(iii) follows by the fact that $\mathbb{E}D_{i}X_{i}^{\top}(\gamma_{*}-\gamma)$
is the first entry of $\Sigma(\gamma_{*}-\gamma)$ and (iv) follows
by $\|\Sigma(\gamma_{*}-\gamma)\|_{\infty}\leq2\varepsilon_{n}$.
Notice that $\lambda_{1}=\lambda_{0}\varepsilon_{n}$ with a large
enough $\lambda_{0}$ and $\varepsilon_{n}\gg n^{-1/2}$. Thus, 
\begin{equation}
P\left(\left|\mathbb{E}_{n,2}D_{i}^{2}\theta_{*}+\mathbb{E}_{n,2}D_{i}Z_{i}^{\top}\beta_{*}-\mathbb{E}_{n,2}D_{i}Y_{i}\right|\leq\lambda_{1}\right)\geq1-o(1).\label{eq: lem PLM tool 5 eq 4}
\end{equation}}

\jelenax{Finally, we need to verify that if $\max\{\xi_{1},\xi_{2}\}>1/2$,
then 
\begin{equation}
P\left(\left|\mathbb{E}_{n,2}D_{i}^{2}\theta_{*}+\hat{\phi}^{\top}\mathbb{E}_{n,1}Z_{i}Z_{i}^{\top}\beta_{*}-\mathbb{E}_{n,2}D_{i}Y_{i}\right|\leq\lambda_{2}\right)\geq1-o(1).\label{eq: lem PLM tool 5 eq 5}
\end{equation}}

\jelenax{Notice that $\mathbb{E}_{n,2}D_{i}^{2}\theta_{*}=\mathbb{E}D_{i}^{2}\theta_{*}+O_{P}(n^{-1/2})$
and $\mathbb{E}_{n,2}D_{i}Y_{i}=\mathbb{E}D_{i}Y_{i}+O_{P}(n^{-1/2})$
are bounded since $D_{i}$ is sub-Gaussian and $\mathbb{E}Y_{i}^{2}$
is bounded by assumption. Moreover, notice that 
\begin{align*}
\mathbb{E}(D_{i}Z_{i}^{\top}\beta_{*})^{2}&=\mathbb{E}(u_{i}Z_{i}^{\top}\beta_{*}+\phi^{\top}Z_{i}Z_{i}^{\top}\beta_{*}+\eta_{g,i}Z_{i}^{\top}\beta_{*})^{2}
\\
&=\mathbb{E}(u_{i}Z_{i}^{\top}\beta_{*})^{2}+\mathbb{E}(\phi^{\top}Z_{i}Z_{i}^{\top}\beta_{*})^{2}+\mathbb{E}(\eta_{g,i}Z_{i}^{\top}\beta_{*})^{2}
\end{align*}
due to $\mathbb{E}(u_{i}\mid Z_{i})=0$ and $\mathbb{E}\eta_{g,i}Z_{i}=0$
(from the definition of $\eta_{g,i}$). For bounded $X_{i}$, since
$\mathbb{E}(X_{i}^{\top}\gamma_{*})^{2}$ has been shown to be bounded
and $|D_{i}|\leq\|X_{i}\|_{\infty}$ is bounded, so is $\mathbb{E}(D_{i}Z_{i}^{\top}\beta_{*})^{2}\lesssim\mathbb{E}(Z_{i}^{\top}\beta_{*})^{2}=\mathbb{E}(X_{i}^{\top}\gamma_{*}-D_{i}\theta_{*})^{2}\leq2\mathbb{E}(X_{i}^{\top}\gamma_{*})^{2}+2\theta_{*}^{2}\mathbb{E}D_{i}^{2}$.
Thus, $\mathbb{E}(\phi^{\top}Z_{i}Z_{i}^{\top}\beta_{*})^{2}$ is
bounded. For sub-Gaussian $X_{i}$, both $Z_{i}^{\top}\phi$ and $Z_{i}^{\top}\beta_{*}$
are sub-Gaussian, which means that $\mathbb{E}(\phi^{\top}Z_{i}Z_{i}^{\top}\beta_{*})^{2}$
is bounded. Hence,
\begin{align*}
\mathbb{E}\left(\phi^{\top}\mathbb{E}_{n,1}Z_{i}Z_{i}^{\top}\beta_{*}-\phi^{\top}\mathbb{E}Z_{i}Z_{i}^{\top}\beta_{*}\right)^{2} & =n_{1}^{-2}\sum_{i\in I_{1}}\mathbb{E}[\phi^{\top}Z_{i}Z_{i}^{\top}\beta_{*}-\mathbb{E}\phi^{\top}Z_{i}Z_{i}^{\top}\beta_{*}]^{2}\\
 & \leq n_{1}^{-2}\sum_{i\in I_{1}}\mathbb{E}[\phi^{\top}Z_{i}Z_{i}^{\top}\beta_{*}]^{2}\lesssim n_{1}^{-1}\mathbb{E}(Z_{i}^{\top}\beta_{*})^{2}\lesssim n^{-1}.
\end{align*}}

\jelenax{Therefore, 
\begin{align*}
 & \left|\mathbb{E}_{n,2}D_{i}^{2}\theta_{*}+\phi^{\top}\mathbb{E}_{n,1}Z_{i}Z_{i}^{\top}\beta_{*}-\mathbb{E}_{n,2}D_{i}Y_{i}\right|\\
 & \leq\left|\mathbb{E}_{n,2}D_{i}^{2}-\mathbb{E}D_{i}^{2}\right|\cdot|\theta_{*}|+\left|\phi^{\top}\mathbb{E}_{n,1}Z_{i}Z_{i}^{\top}\beta_{*}-\phi^{\top}\mathbb{E}Z_{i}Z_{i}^{\top}\beta_{*}\right|\\
 & \quad+|\mathbb{E}_{n,2}D_{i}Y_{i}-\mathbb{E}D_{i}Y_{i}|+\left|\mathbb{E}D_{i}^{2}\theta_{*}+\phi^{\top}\mathbb{E}Z_{i}Z_{i}^{\top}\beta_{*}-\mathbb{E}D_{i}Y_{i}\right|\\
 & =O_{P}(n^{-1/2})+\left|\mathbb{E} D_{i}^{2}\theta_{*}+\phi^{\top}\mathbb{E}Z_{i}Z_{i}^{\top}\beta_{*}-\mathbb{E}D_{i}Y_{i}\right|\\
 & =O_{P}(n^{-1/2})+\left|\mathbb{E}D_{i}^{2}\theta_{*}+\phi^{\top}\mathbb{E}Z_{i}Z_{i}^{\top}\beta_{*}-\mathbb{E}D_{i}(D_{i}\theta+Z_{i}^{\top}\beta+\eta_{f,i}+\varepsilon_{i})\right|\\
 & \leq  O_{P}(n^{-1/2})+\left|\mathbb{E}D_{i}^{2}(\theta_{*}-\theta)+\phi^{\top}\mathbb{E}Z_{i}Z_{i}^{\top}(\beta_{*}-\beta)\right| +\left|\mathbb{E} D_i \eta_{f,i} \right| \\
 & \overset{\text{(i)}}{=} O_{P}(n^{-1/2})+\left|\mathbb{E}D_{i}^{2}(\theta_{*}-\theta)+\phi^{\top}\mathbb{E}Z_{i}Z_{i}^{\top}(\beta_{*}-\beta)\right| +\left|\mathbb{E} \eta_{g,i} \eta_{f,i} \right| \\
 & \overset{\text{(ii)}}{=}O_{P}(n^{-1/2})+\left|\mathbb{E}D_{i}X_{i}^{\top}(\gamma_{*}-\gamma)\right|  +\left|\mathbb{E} \eta_{g,i} \eta_{f,i} \right|\\
 & \overset{\text{(iii)}}{=}O_{P}(n^{-1/2})+O\left(\sqrt{n^{-1}\log p}\right)  +\left|\mathbb{E} \eta_{g,i} \eta_{f,i} \right| \overset{\text{(iv)}}{=}o_{P}(n^{-1/4}),
\end{align*}
where (i) follows by $D_{i}=u_{i}+Z_{i}^{\top}\phi+\eta_{g,i}$ and $\mathbb{E} Z_{i}\eta_{f,i}=0$ (due to the definition of $\beta$),
 (ii) follows by $\phi^{\top}\mathbb{E}Z_{i}Z_{i}^{\top}=\mathbb{E}D_{i}Z_{i}^{\top}$,
(iii) follows by the fact that $\mathbb{E}D_{i}X_{i}^{\top}(\gamma_{*}-\gamma)$
is the first entry of $\Sigma(\gamma_{*}-\gamma)$ and $\|\Sigma(\gamma_{*}-\gamma)\|_{\infty}\leq2\varepsilon_{n}\asymp\sqrt{n^{-1}\log p}$
and (iv) follows by $\log p\ll\sqrt{n}$, $|\mathbb{E} \eta_{g,i} \eta_{f,i}|^2 \leq \mathbb{E}  \eta_{f,i}^2 \cdot \mathbb{E}  \eta_{g,i}^2 \lesssim p^{-2\xi_{1}} p^{-2\xi_{2}}$ and  $p\gg n^{1/(2\xi_{1})}$. Thus, to show (\ref{eq: lem PLM tool 5 eq 5}),
it suffices to verify that if $\max\{\xi_{1},\xi_{2}\}>1/2$, then
\begin{equation}
\left|(\hat{\phi}-\phi)^{\top}\mathbb{E}_{n,1}Z_{i}Z_{i}^{\top}\beta_{*}\right|=o_{P}(n^{-1/4}).\label{eq: lem PLM tool 5 eq 6}
\end{equation}}

\jelenax{We verify (\ref{eq: lem PLM tool 5 eq 6}) in two cases: (a) $\xi_{1}>1/2$
and (b) $\xi_{2}>1/2$.}

\jelenax{In Case (a), we have that $\log p\ll n^{1/2-1/(4\xi_{1})}$ and 
\begin{multline*}
\left|(\hat{\phi}-\phi)^{\top}\mathbb{E}_{n,1}Z_{i}Z_{i}^{\top}\beta_{*}\right|\leq\|\beta_{*}\|_{1}\|\mathbb{E}_{n,1}Z_{i}Z_{i}^{\top}(\hat{\phi}-\phi)\|_{\infty}\\
\leq\|\gamma_{*}\|_{1}\cdot\|\mathbb{E}_{n,1}Z_{i}Z_{i}^{\top}(\hat{\phi}-\phi)\|_{\infty}\overset{\text{(i)}}{=}O_{P}\left(\varepsilon_{n}^{-1/(2\xi_{1}+1)}\cdot\varepsilon_{n}\right)\overset{\text{(ii)}}{=}o_{P}(n^{-1/4})
\end{multline*}
where (i) follows by $\|\mathbb{E}_{n,1}Z_{i}Z_{i}^{\top}(\hat{\phi}-\phi)\|_{\infty}\leq2\lambda_{1}\lesssim\varepsilon_{n}$
(Lemma \ref{lem: PLM tool 4}) and $\|\gamma_{*}\|_{1}\leq\sqrt{|J_{*}|}\cdot\|\gamma_{*}\|_{2}\lesssim\sqrt{|J_{*}|}\lesssim\varepsilon_{n}^{-1/(2\xi_{1}+1)}$
due to $\|\gamma_{*}\|_{2}\lesssim1$, $J_{*}={\rm supp}(\gamma_{*})$
and $|J_{*}|\lesssim\varepsilon_{n}^{-2/(2\xi_{1}+1)}$ (shown above)
and (ii) follows by $\log p\ll n^{1/2-1/(4\xi_{1})}$.}

\jelenax{In Case (b), we have that $\log p\ll n^{1/2-1/(4\xi_{2})}$, $\|\hat{\phi}-\phi\|_{1}\lesssim\varepsilon_{n}^{(2\xi_{2}-1)/(2\xi_{2}+1)}$
and $\|\hat{\phi}-\phi\|_{2}\lesssim\varepsilon_{n}^{2\xi_{2}/(2\xi_{2}+1)}$
with probability $1-o(1)$ (due to Lemma \ref{lem: PLM tool 4}).
Then with probability $1-o(1)$, 
\begin{align}
 & \left|(\hat{\phi}-\phi)^{\top}\mathbb{E}_{n,1}Z_{i}Z_{i}^{\top}\beta_{*}\right|\nonumber \\
 & \leq\left|(\hat{\phi}-\phi)^{\top}\mathbb{E}Z_{i}Z_{i}^{\top}\beta_{*}\right|+\left|(\hat{\phi}-\phi)^{\top}(\mathbb{E}_{n,1}Z_{i}Z_{i}^{\top}-\mathbb{E}Z_{i}Z_{i}^{\top})\beta_{*}\right|\nonumber \\
 & \leq\|\hat{\phi}-\phi\|_{2}\cdot\|\mathbb{E}Z_{i}Z_{i}^{\top}\beta_{*}\|_{2}+\|\hat{\phi}-\phi\|_{1}\cdot\|\mathbb{E}_{n,1}Z_{i}Z_{i}^{\top}\beta_{*}-\mathbb{E}Z_{i}Z_{i}^{\top}\beta_{*}\|_{\infty}\nonumber \\
 & \lesssim\varepsilon_{n}^{2\xi_{2}/(2\xi_{2}+1)}\cdot\|\mathbb{E}Z_{i}Z_{i}^{\top}\beta_{*}\|_{2}+\varepsilon_{n}^{(2\xi_{2}-1)/(2\xi_{2}+1)}\cdot\|\mathbb{E}_{n,1}Z_{i}Z_{i}^{\top}\beta_{*}-\mathbb{E}Z_{i}Z_{i}^{\top}\beta_{*}\|_{\infty}\nonumber \\
 & \overset{\text{(i)}}{\lesssim}\varepsilon_{n}^{2\xi_{2}/(2\xi_{2}+1)}\cdot1+\varepsilon_{n}^{(2\xi_{2}-1)/(2\xi_{2}+1)}\cdot\|\mathbb{E}_{n,1}Z_{i}Z_{i}^{\top}\beta_{*}-\mathbb{E}Z_{i}Z_{i}^{\top}\beta_{*}\|_{\infty}\nonumber \\
 & \overset{\text{(ii)}}{\lesssim}\varepsilon_{n}^{2\xi_{2}/(2\xi_{2}+1)}\cdot1+\varepsilon_{n}^{(2\xi_{2}-1)/(2\xi_{2}+1)}\cdot\varepsilon_{n}\overset{\text{(iii)}}{\ll}n^{-1/4},\label{eq: lem PLM tool 5 eq 16}
\end{align}
where (i) follows by $\|\mathbb{E}Z_{i}Z_{i}^{\top}\beta_{*}\|_{2}\leq\lambda_{\max}(\mathbb{E}Z_{i}Z_{i}^{\top})\cdot\|\beta_{*}\|_{2}\leq\lambda_{\max}(\Sigma)\cdot\|\gamma_{*}\|_{2}\lesssim1$
(due to $\|\gamma_{*}\|_{2}\lesssim1$ as shown above) and (ii) follows
by $\|\mathbb{E}_{n,1}Z_{i}Z_{i}^{\top}\beta_{*}-\mathbb{E}Z_{i}Z_{i}^{\top}\beta_{*}\|_{\infty}\lesssim\varepsilon_{n}$;
to see this, follow the same argument in Lemma \ref{lem: PLM tool 2}. Moreover, (iii) in the above display follows
by $\log p\ll n^{1/2-1/(4\xi_{2})}$. Therefore, we have verified (\ref{eq: lem PLM tool 5 eq 6}). Hence,
we have proved that $\gamma_{*}$ satisfies all the constraints with
probability at least $1-o(1)$.}

\jelenax{\textbf{Step 2:} show that for any $\kappa_{1}>0$, there exists $\kappa_{2}>0$ such that $P(\Acal_{2})\geq1-o(1)$.}

\jelenax{Notice that 
\[
\tilde{\Sigma}-\hat{\Sigma}_{(1)}=\begin{pmatrix}\mathbb{E}_{n,2}D_{i}^{2}-\mathbb{E}_{n,1}D_{i}^{2} & \mathbb{E}_{n,2}D_{i}Z_{i}^{\top}-\mathbb{E}_{n,1}D_{i}Z_{i}^{\top}\\
\mathbb{E}_{n,2}Z_{i}D_{i}-\mathbb{E}_{n,1}Z_{i}D_{i} & 0
\end{pmatrix}.
\]
In Lemma \ref{lem: PLM tool 3}, we have proved that for any $\kappa_{1}$, there exists $\kappa_{2}>0$ such that 
$$
P \left( \min_{|A|\leq \kappa_{1} (n/\log p)^{1/(2\xi_{2}+1)} } \min_{\|v_{A^{c}}\|_{1}\leq\|v_{A}\|_{1}}\frac{v^{\top}\hat{\Sigma}_{(1)}v}{\|v_{A}\|_{2}^{2}}\geq 2 \kappa_{2} \right)\geq 1-o(1).
$$
Therefore, we only need to show that 
\begin{equation}
P\left(\max_{|A|\leq \kappa_{1} (n/\log p)^{1/(2\xi_{2}+1)} } \max_{\|v_{A^{c}}\|_{1}\leq\|v_{A}\|_{1}}\frac{\left|v^{\top}\left(\tilde{\Sigma}-\hat{\Sigma}_{(1)}\right)v\right|}{\|v_{A}\|_{2}^{2}}\leq\kappa_{2} \right)\geq1-o(1).\label{eq: lem PLM tool 5 eq 17}
\end{equation}}

\jelenax{Fix a set $A$. Notice that if $v=(v_{1},v_{-1}^{\top})^{\top}$ with $v_{1}=0$,
then $v^{\top}\left(\tilde{\Sigma}-\hat{\Sigma}_{(1)}\right)v=0$.
Hence, we only need to consider $v=(v_{1},v_{-1}^{\top})^{\top}$
with $v_{1}\neq0$. Without loss of generality, we can normalize $v=(1,v_{-1}^{\top})^{\top}$
since the ratio $\left|v^{\top}\left(\tilde{\Sigma}-\hat{\Sigma}_{(1)}\right)v\right|/\|v_{A}\|_{2}^{2}$
depends only on $v_{-1}/v_{1}$. Notice that for $v=(1,v_{-1}^{\top})^{\top}$,
\[
\left|v^{\top}\left(\tilde{\Sigma}-\hat{\Sigma}_{(1)}\right)v\right|=\left|\left(\mathbb{E}_{n,2}D_{i}^{2}-\mathbb{E}_{n,1}D_{i}^{2}\right)+2v_{-1}^{\top}\left(\mathbb{E}_{n,2}Z_{i}D_{i}-\mathbb{E}_{n,1}Z_{i}D_{i}\right)\right|\leq\tau_{1}+\|v_{-1}\|_{1}\tau_{2},
\]
where $\tau_{1}=|\mathbb{E}_{n,2}D_{i}^{2}-\mathbb{E}_{n,1}D_{i}^{2}|$
and $\tau_{2}=\|\mathbb{E}_{n,2}Z_{i}D_{i}-\mathbb{E}_{n,1}Z_{i}D_{i}\|_{\infty}$.}

\jelenax{We discuss two cases: (a) $1\in A$ and (b) $1\notin A$.}

\jelenax{In Case (a), we have 
\begin{multline*}
\frac{\left|v^{\top}\left(\tilde{\Sigma}-\hat{\Sigma}_{(1)}\right)v\right|}{\|v_{A}\|_{2}^{2}}\leq\frac{\tau_{1}+\|v_{-1}\|_{1}\tau_{2}}{\|v_{A}\|_{2}^{2}}\overset{\text{(i)}}{\leq}\frac{\tau_{1}+2\|v_{A}\|_{1}\tau_{2}}{\|v_{A}\|_{2}^{2}}\\
\leq\frac{\tau_{1}+2\sqrt{|A|}\cdot\|v_{A}\|_{2}\tau_{2}}{\|v_{A}\|_{2}^{2}}=\frac{\tau_{1}}{\|v_{A}\|_{2}^{2}}+\frac{2\sqrt{|A|}\tau_{2}}{\|v_{A}\|_{2}}\overset{\text{(ii)}}{\leq}\tau_{1}+2\sqrt{|A|}\tau_{2},
\end{multline*}
where (i) follows by $\|v_{-1}\|_{1}\leq\|v\|_{1}=\|v_{A}\|_{1}+\|v_{A^{c}}\|_{1}\leq2\|v_{A}\|_{1}$
due to $\|v_{A^{c}}\|_{1}\leq\|v_{A}\|_{1}$ and (ii) follows
by $\|v_{A}\|_{2}\geq|v_{1}|=1$ due to $v_{1}=1$ and $1\in A$.}

\jelenax{In Case (b), we have that $1\in A^{c}$ and thus $\|v_{A^{c}}\|_{1}\geq|v_{1}|=1$,
which means that $\|v_{A}\|_{1}\geq\|v_{A^{c}}\|_{1}=1$.
Thus, 
\begin{multline*}
\frac{\left|v^{\top}\left(\tilde{\Sigma}-\hat{\Sigma}_{(1)}\right)v\right|}{\|v_{A}\|_{2}^{2}}\leq\frac{\tau_{1}+\|v_{-1}\|_{1}\tau_{2}}{\|v_{A}\|_{2}^{2}}\leq\frac{\tau_{1}+2\|v_{A}\|_{1}\tau_{2}}{\|v_{A}\|_{2}^{2}}\\
\leq\frac{\tau_{1}+2\sqrt{|A|}\cdot\|v_{A}\|_{2}\tau_{2}}{\|v_{A}\|_{2}^{2}}=\frac{\tau_{1}}{\|v_{A}\|_{2}^{2}}+\frac{2\sqrt{|A|}\tau_{2}}{\|v_{A}\|_{2}}\leq\tau_{1}+2\sqrt{|A|}\tau_{2}.
\end{multline*}}

\jelenax{Therefore, we have proved that 
\begin{equation}
\max_{\|v_{A^{c}}\|_{1}\leq\|v_{A}\|_{1}}\frac{\left|v^{\top}\left(\tilde{\Sigma}-\hat{\Sigma}_{(1)}\right)v\right|}{\|v_{A}\|_{2}^{2}}\leq\tau_{1}+2\sqrt{|A|}\tau_{2}.\label{eq: lem PLM tool 5 eq 18}
\end{equation}}

\jelenax{Notice that $\tau_{1}=o_{P}(1)$,  $\tau_{2}=O_{P}(\sqrt{n^{-1}\log p})=O_{P}(\varepsilon_{n})$
(due to Lemma \ref{lem: PLM tool 2}) and $|A|\lesssim \varepsilon_{n}^{-2/(2\xi_{2}+1)}$ The above display is $o_{P}(1)+o_{P}(\varepsilon_{n}^{-1/(2\xi_{2}+1)}\cdot\varepsilon_{n})=o_{P}(1)$.
Thus, (\ref{eq: lem PLM tool 5 eq 17}) follows by (\ref{eq: lem PLM tool 5 eq 18}).}

\jelenax{\textbf{Step 3:} show the bounds for $\|\hat{\gamma}-\gamma\|_{2}$,
$\|\hat{\gamma}-\gamma\|_{1}$ and $(\hat{\beta}-\beta)^{\top}\mathbb{E}_{n,1}Z_{i}Z_{i}^{\top}(\hat{\beta}-\beta)$.}

\jelenax{In the previous two steps, we have verified that for any $\kappa_{1}>0$, there exists $\kappa_{2}>0$ such that $P(\Acal)\geq1-o(1)$.
On the event $\Acal$, $\hat{\gamma}$ is well defined. The same argument
as in Lemma \ref{lem: PLM tool 3} yields the bound for $\|\hat{\gamma}-\gamma\|_{2}$
and $\|\hat{\gamma}-\gamma\|_{1}$.}

\jelenax{We now show the bound for $(\hat{\beta}-\beta)^{\top}\mathbb{E}_{n,1}Z_{i}Z_{i}^{\top}(\hat{\beta}-\beta)$.
Let $\Delta:=\hat{\gamma}-\gamma_{*}$ be partitioned as $\Delta=(\Delta_{\theta},\Delta_{\beta}^{\top})^{\top}$.
From the same argument as (\ref{eq: lem PLM tool 3 eq 9}) and (\ref{eq: lem PLM tool 3 eq 10})
in Lemma \ref{lem: PLM tool 3}, we can derive that on the event $\Acal$,
$\|\Delta\|_{1}\leq8\kappa_{2}^{-1}\lambda_{1}|J|$ and $\|\Delta\|_{2}\leq12\sqrt{2}\kappa_{2}^{-1}\lambda_{1}\sqrt{|J|}$
with $|J|\lesssim\varepsilon_{n}^{-2/(2\xi_{1}+1)}$ and $\lambda_{1}\asymp\varepsilon_{n}$.
Hence, $\|\Delta\|_{1}=O_{P}(\varepsilon_{n}^{(2\xi_{1}-1)/(2\xi_{1}+1)})$
and $\|\Delta\|_{2}=O_{P}(\varepsilon_{n}^{4\xi_{1}/(2\xi_{1}+1)})$.}

\jelenax{On the event $\Acal$, both $\gamma_{*}$ and $\hat{\gamma}$ satisfy
the constraints of $\hat{\gamma}$, which means that 
\[
\begin{cases}
\|\mathbb{E}_{n,2}Z_{i}D_{i}\theta_{*}+\mathbb{E}_{n,1}Z_{i}Z_{i}^{\top}\beta_{*}-\mathbb{E}_{n,2}Z_{i}Y_{i}\|_{\infty}\leq\lambda_{1}\\
\|\mathbb{E}_{n,2}Z_{i}D_{i}\hat{\theta}+\mathbb{E}_{n,1}Z_{i}Z_{i}^{\top}\hat{\beta}-\mathbb{E}_{n,2}Z_{i}Y_{i}\|_{\infty}\leq\lambda_{1}
\end{cases}
\]
and 
\[
\begin{cases}
|\mathbb{E}_{n,2}D_{i}^{2}\theta_{*}+\mathbb{E}_{n,2}D_{i}Z_{i}^{\top}\beta_{*}-\mathbb{E}_{n,2}D_{i}Y_{i}|\leq\lambda_{1}\\
|\mathbb{E}_{n,2}D_{i}^{2}\hat{\theta}+\mathbb{E}_{n,2}D_{i}Z_{i}^{\top}\hat{\beta}-\mathbb{E}_{n,2}D_{i}Y_{i}|\leq\lambda_{1}.
\end{cases}
\]}

\jelenax{Therefore, by the triangular inequality, we have that on the event
$\Acal$, 
\begin{equation}
\|\mathbb{E}_{n,2}Z_{i}D_{i}\Delta_{\theta}+\mathbb{E}_{n,1}Z_{i}Z_{i}^{\top}\Delta_{\beta}\|_{\infty}\leq2\lambda_{1}\label{eq: lem PLM tool 5 eq 19}
\end{equation}
and 
\begin{equation}
|\mathbb{E}_{n,2}D_{i}^{2}\Delta_{\theta}+\mathbb{E}_{n,2}D_{i}Z_{i}^{\top}\Delta_{\beta}|\leq2\lambda_{1}.\label{eq: lem PLM tool 5 eq 20}
\end{equation}}

\jelenax{By (\ref{eq: lem PLM tool 5 eq 19}) and Hölder's inequality, we have
that $|\Delta_{\beta}^{\top}\mathbb{E}_{n,2}Z_{i}D_{i}\Delta_{\theta}+\Delta_{\beta}^{\top}\mathbb{E}_{n,1}Z_{i}Z_{i}^{\top}\Delta_{\beta}|\leq2\lambda_{1}\|\Delta_{\beta}\|_{1}$,
which means that 
\[
\Delta_{\beta}^{\top}\mathbb{E}_{n,1}Z_{i}Z_{i}^{\top}\Delta_{\beta}\leq2\lambda_{1}\|\Delta_{\beta}\|_{1}+|\Delta_{\beta}^{\top}\mathbb{E}_{n,2}Z_{i}D_{i}|\cdot|\Delta_{\theta}|.
\]}

\jelenax{By (\ref{eq: lem PLM tool 5 eq 20}), we have $|\mathbb{E}_{n,2}D_{i}Z_{i}^{\top}\Delta_{\beta}|\leq\mathbb{E}_{n,2}D_{i}^{2}\cdot|\Delta_{\theta}|+2\lambda_{1}$.
Hence, 
\begin{align*}
\Delta_{\beta}^{\top}\mathbb{E}_{n,1}Z_{i}Z_{i}^{\top}\Delta_{\beta} & \leq2\lambda_{1}\|\Delta_{\beta}\|_{1}+2\lambda_{1}\cdot|\Delta_{\theta}|+\mathbb{E}_{n,2}D_{i}^{2}\cdot\Delta_{\theta}^{2}\\
 & \overset{\text{(i)}}{\leq}4\lambda_{1}\|\Delta\|_{1}+\mathbb{E}_{n,2}D_{i}^{2}\cdot\|\Delta\|_{2}^{2}\overset{\text{(ii)}}{=}O_{P}(\varepsilon_{n}^{4\xi_{1}/(2\xi_{1}+1)}),
\end{align*}
where (i) follows by $\|\Delta_{\beta}\|_{1}\leq\|\Delta\|_{1}$,
$|\Delta_{\theta}|\leq\|\Delta\|_{1}$ and $\Delta_{\theta}^{2}\leq\|\Delta\|_{2}^{2}$
and (ii) follows by the sub-Gaussianity of $D_{i}$, $\lambda_{1}=O(\varepsilon_{n})$,
$\|\Delta\|_{1}=O_{P}(\varepsilon_{n}^{(2\xi_{1}-1)/(2\xi_{1}+1)})$
and $\|\Delta\|_{2}=O_{P}(\varepsilon_{n}^{4\xi_{1}/(2\xi_{1}+1)})$.
Now we notice that 
\begin{align*}
 & (\hat{\beta}-\beta)^{\top}\mathbb{E}_{n,1}Z_{i}Z_{i}^{\top}(\hat{\beta}-\beta)\\
 & =[\Delta_{\beta}+(\beta_{*}-\beta)]^{\top}\mathbb{E}_{n,1}Z_{i}Z_{i}^{\top}[\Delta_{\beta}+(\beta_{*}-\beta)]\\
 & \leq2\Delta_{\beta}^{\top}\mathbb{E}_{n,1}Z_{i}Z_{i}^{\top}\Delta_{\beta}+2(\beta_{*}-\beta)^{\top}\mathbb{E}_{n,1}Z_{i}Z_{i}^{\top}(\beta_{*}-\beta)\\
 & \overset{\text{(i)}}{=}O_{P}(\varepsilon_{n}^{4\xi_{1}/(2\xi_{1}+1)})+2(\beta_{*}-\beta)^{\top}\mathbb{E}_{n,1}Z_{i}Z_{i}^{\top}(\beta_{*}-\beta)\\
 & \overset{\text{(ii)}}{=}O_{P}(\varepsilon_{n}^{4\xi_{1}/(2\xi_{1}+1)}),
\end{align*}
where (i) follows by $\Delta_{\beta}^{\top}\mathbb{E}_{n,1}Z_{i}Z_{i}^{\top}\Delta_{\beta}=O_{P}(\varepsilon_{n}^{4\xi_{1}/(2\xi_{1}+1)})$
and (ii) follows by 
\begin{multline*}
\mathbb{E}(\beta_{*}-\beta)^{\top}\mathbb{E}_{n,1}Z_{i}Z_{i}^{\top}(\beta_{*}-\beta)=(\beta_{*}-\beta)^{\top}(\mathbb{E}Z_{i}Z_{i}^{\top})(\beta_{*}-\beta)\leq\lambda_{\max}(\mathbb{E}Z_{i}Z_{i}^{\top})\cdot\|\beta_{*}-\beta\|_{2}^{2}\\
\leq\lambda_{\max}(\Sigma)\cdot\|\gamma_{*}-\gamma\|_{2}^{2}\leq\frac{\lambda_{\max}(\Sigma)}{\lambda_{\min}(\Sigma)}(\gamma_{*}-\gamma)^{\top}\Sigma(\gamma_{*}-\gamma)\lesssim\varepsilon_{n}^{4\xi_{1}/(2\xi_{1}+1)}.
\end{multline*}}

\jelenax{\textbf{Step 4:} show the bounds for $\|\mathbb{E}_{n,2}Z_{i}D_{i}(\hat{\theta}-\theta)+\mathbb{E}_{n,1}Z_{i}Z_{i}^{\top}(\hat{\beta}-\beta)\|_{\infty}$
and $|\mathbb{E}_{n,2}D_{i}^{2}\delta_{\theta}+\hat{\phi}^{\top}\mathbb{E}_{n,1}Z_{i}Z_{i}^{\top}\delta_{\beta}|$.}

\jelenax{To bound $\|\mathbb{E}_{n,2}Z_{i}D_{i}(\hat{\theta}-\theta)+\mathbb{E}_{n,1}Z_{i}Z_{i}^{\top}(\hat{\beta}-\beta)\|_{\infty}$,
we observe that on the event $\Acal$, 
\begin{multline}
\|\mathbb{E}_{n,2}Z_{i}D_{i}(\hat{\theta}-\theta_{*})+\mathbb{E}_{n,1}Z_{i}Z_{i}^{\top}(\hat{\beta}-\beta_{*})\|_{\infty}\leq\|\mathbb{E}_{n,2}Z_{i}D_{i}\hat{\theta}+\mathbb{E}_{n,1}Z_{i}Z_{i}^{\top}\hat{\beta}-\mathbb{E}_{n,2}Z_{i}Y_{i}\|_{\infty}\\
+\|\mathbb{E}_{n,2}Z_{i}D_{i}\theta_{*}+\mathbb{E}_{n,1}Z_{i}Z_{i}^{\top}\beta_{*}-\mathbb{E}_{n,2}Z_{i}Y_{i}\|_{\infty}\leq2\lambda_{1}.\label{eq: lem PLM tool 5 eq 21}
\end{multline}}

\jelenax{We notice that $\mathbb{E}Z_{i}D_{i}(\theta_{*}-\theta)+\mathbb{E}Z_{i}Z_{i}^{\top}(\beta_{*}-\beta)=\mathbb{E}Z_{i}X_{i}^{\top}(\gamma_{*}-\gamma)$
are the last $p$ elements of $\Sigma(\gamma_{*}-\gamma)$, which
means that 
\[
\|\mathbb{E}Z_{i}D_{i}(\theta_{*}-\theta)+\mathbb{E}Z_{i}Z_{i}^{\top}(\beta_{*}-\beta)\|_{\infty}\leq\|\Sigma(\gamma_{*}-\gamma)\|_{\infty}\leq\varepsilon_{n}.
\]}

\jelenax{By Lemma \ref{lem: PLM tool 2}, $\|\mathbb{E}_{n,2}Z_{i}D_{i}-\mathbb{E}Z_{i}D_{i}\|_{\infty}=O_{P}(\varepsilon_{n})$.
Notice that $\|\beta_{*}-\beta\|_{2}\leq\|\gamma_{*}-\gamma\|_{2}\leq\sqrt{(\gamma-\gamma_{*})^{\top}\Sigma(\gamma-\gamma_{*})/\lambda_{\min}(\Sigma)}\lesssim1$.
Thus, the same argument as in Lemma \ref{lem: PLM tool 2} yields
$\|\mathbb{E}_{n,1}Z_{i}Z_{i}^{\top}(\beta_{*}-\beta)-\mathbb{E}Z_{i}Z_{i}^{\top}(\beta_{*}-\beta)\|_{\infty}=O_{P}(\varepsilon_{n})$.
Therefore, 
\begin{multline*}
\|\mathbb{E}_{n,2}Z_{i}D_{i}(\theta_{*}-\theta)+\mathbb{E}_{n,1}Z_{i}Z_{i}^{\top}(\beta_{*}-\beta)\|_{\infty}\\
\leq\|\mathbb{E}Z_{i}D_{i}(\theta_{*}-\theta)+\mathbb{E}Z_{i}Z_{i}^{\top}(\beta_{*}-\beta)\|_{\infty}+O_{P}(\varepsilon_{n})=O_{P}(\varepsilon_{n}).
\end{multline*}}

\jelenax{The above display and (\ref{eq: lem PLM tool 5 eq 21}) imply $\|\mathbb{E}_{n,2}Z_{i}D_{i}(\hat{\theta}-\theta)+\mathbb{E}_{n,1}Z_{i}Z_{i}^{\top}(\hat{\beta}-\beta)\|_{\infty}=O_{P}(\varepsilon_{n})$.}

\jelenax{To bound $|\mathbb{E}_{n,2}D_{i}^{2}\delta_{\theta}+\hat{\phi}^{\top}\mathbb{E}_{n,1}Z_{i}Z_{i}^{\top}\delta_{\beta}|$
for $\xi_{2}>1/2$, we similarly observe that on the event $\Acal$,
the triangular inequality yields 
\begin{multline}
|\mathbb{E}_{n,2}D_{i}^{2}(\hat{\theta}-\theta_{*})+\hat{\phi}^{\top}\mathbb{E}_{n,1}Z_{i}Z_{i}^{\top}(\hat{\beta}-\beta_{*})|\leq|\mathbb{E}_{n,2}D_{i}^{2}\hat{\theta}+\hat{\phi}^{\top}\mathbb{E}_{n,1}Z_{i}Z_{i}^{\top}\hat{\beta}-\mathbb{E}_{n,2}D_{i}Y_{i}|\\
+|\mathbb{E}_{n,2}D_{i}^{2}\theta_{*}+\hat{\phi}^{\top}\mathbb{E}_{n,1}Z_{i}Z_{i}^{\top}\beta_{*}-\mathbb{E}_{n,2}D_{i}Y_{i}|\leq2n^{-1/4}.\label{eq: lem PLM tool 5 eq 22}
\end{multline}}

\jelenax{We notice that $\mathbb{E}D_{i}^{2}(\theta_{*}-\theta)+\phi^{\top}\mathbb{E}Z_{i}Z_{i}^{\top}(\beta_{*}-\beta)=\mathbb{E}D_{i}^{2}(\theta_{*}-\theta)+\mathbb{E}D_{i}Z_{i}^{\top}(\beta_{*}-\beta)=\mathbb{E}D_{i}X_{i}^{\top}(\gamma_{*}-\gamma)$
is the first element of $\Sigma(\gamma_{*}-\gamma)$, which means
that 
\[
|\mathbb{E}D_{i}^{2}(\theta_{*}-\theta)+\phi^{\top}\mathbb{E}Z_{i}Z_{i}^{\top}(\beta_{*}-\beta)|\leq\|\Sigma(\gamma_{*}-\gamma)\|_{\infty}\leq\varepsilon_{n}.
\]}

\jelenax{We notice that 
\begin{align*}
 & \left|\hat{\phi}^{\top}\mathbb{E}_{n,1}Z_{i}Z_{i}^{\top}(\beta_{*}-\beta)-\phi^{\top}\mathbb{E}Z_{i}Z_{i}^{\top}(\beta_{*}-\beta)\right|\\
 & \leq\left|(\hat{\phi}-\phi)^{\top}\mathbb{E}_{n,1}Z_{i}Z_{i}^{\top}(\beta_{*}-\beta)\right|+\left|\phi^{\top}\mathbb{E}_{n,1}Z_{i}Z_{i}^{\top}(\beta_{*}-\beta)-\phi^{\top}\mathbb{E}Z_{i}Z_{i}^{\top}(\beta_{*}-\beta)\right|\\
 & \overset{\text{(i)}}{=}o_{P}(n^{-1/4})+\left|\phi^{\top}\mathbb{E}_{n,1}Z_{i}Z_{i}^{\top}(\beta_{*}-\beta)-\phi^{\top}\mathbb{E}Z_{i}Z_{i}^{\top}(\beta_{*}-\beta)\right|\\
 & \overset{\text{(ii)}}{=}o_{P}(n^{-1/4})+\left|\mathbb{E}_{n,1}(D_{i}-u_{i})Z_{i}^{\top}(\beta_{*}-\beta)-\mathbb{E}D_{i}Z_{i}^{\top}(\beta_{*}-\beta)\right|\\
 & \leq o_{P}(n^{-1/4})+\left|\mathbb{E}_{n,1}u_{i}Z_{i}^{\top}(\beta_{*}-\beta)\right|+\left|\mathbb{E}_{n,1}D_{i}Z_{i}^{\top}(\beta_{*}-\beta)-\mathbb{E}D_{i}Z_{i}^{\top}(\beta_{*}-\beta)\right|\\
 & \overset{\text{(iii)}}{\leq}o_{P}(n^{-1/4})+O_{P}(n^{-1/2}),
\end{align*}
where (i) follows by the same argument as in (\ref{eq: lem PLM tool 5 eq 16})
with $\beta_{*}$ replaced by $\beta_{*}-\beta$ for $\xi_{2}>1/2$,
(ii) follows by $D_{i}=\phi^{\top}Z_{i}+\eta_{g,i}+u_{i}$ and $\mathbb{E}Z_{i}\eta_{g,i}=0$, 
and (iii) follows by Markov's inequality and 
$$\mathbb{E}[D_{i}Z_{i}^{\top}(\beta_{*}-\beta)]^{2}\leq \mathbb{E}[D_{i}^2] \cdot \mathbb{E}[Z_{i}^{\top}(\beta_{*}-\beta)]^{2}\lesssim1$$
and 
$$\mathbb{E}[u_{i}Z_{i}^{\top}(\beta_{*}-\beta)]^{2}\lesssim\mathbb{E}[Z_{i}^{\top}(\beta_{*}-\beta)]^{2}\lesssim1$$
(due to $\mathbb{E}(u_{i}^{2}\mid Z_{i})\lesssim1$);
to see the last step, we simply notice that $n_{1}^{-1}\sum_{i\in I_{1}}Q_{i}=O_{P}(n^{-1/2})$
for independent variables $Q_{i}$ if $\mathbb{E}(Q_{i})=0$ and $\mathbb{E}(Q_{i}^{2})\lesssim1$.}

\jelenax{Also notice that $|\mathbb{E}_{n,2}D_{i}^{2}(\theta_{*}-\theta)-\mathbb{E}D_{i}^{2}(\theta_{*}-\theta)|=O_{P}(n^{-1/2})$
(due to the law of large numbers). This and the above two displays imply
\[
|\mathbb{E}_{n,2}D_{i}^{2}(\theta_{*}-\theta)+\hat{\phi}^{\top}\mathbb{E}_{n,1}Z_{i}Z_{i}^{\top}(\beta_{*}-\beta)|=O_{P}(\varepsilon_{n})+o_{P}(n^{-1/4}).
\]
}

\jelenax{The above display and (\ref{eq: lem PLM tool 5 eq 22}) imply $|\mathbb{E}_{n,2}D_{i}^{2}(\hat{\theta}-\theta)+\hat{\phi}^{\top}\mathbb{E}_{n,1}Z_{i}Z_{i}^{\top}(\hat{\beta}-\beta)|=O_{P}(\varepsilon_{n})+O_{P}(n^{-1/4})=O_{P}(n^{-1/4})$.
The proof is complete. }
\end{proof}
\begin{proof}[Proof of Theorem \ref{thm: PLM main}]
Let $\delta_{\gamma}=\hat{\gamma}-\gamma$ and $\delta_{\pi}=\hat{\pi}-\pi$.
We partition $\delta_{\gamma}=(\delta_{\theta},\delta_{\beta}^{\top})^{\top}$
and $\delta_{\pi}=(\delta_{\pi_{1}},\delta_{\pi_{-1}}^{\top})^{\top}$
such that the first component is a scalar and the second component
is an element in $\RR^{p}$. \jelenax{ Let $\eta_{f,i}=f(\tilde{Z}_{i})-Z_{i}^{\top}\beta$
and $\eta_{\alpha,i}=\alpha_{0}(\tilde{X}_{i})-X_{i}^{\top}\pi$.
We show that if $\max\{\xi_{1},\xi_{2}\}>1/2$, then 
\begin{equation}
\tilde{\theta}_{1}-\theta=\pi^{\top}\mathbb{E}_{n,1}X_{i}\varepsilon_{i}+o_{P}(n^{-1/2}).\label{eq: thm PLM eq 2}
\end{equation}
Notice that this would mean $\tilde{\theta}_{1}-\theta=\sigma_{u}^{-2}\mathbb{E}_{n,1}u_{i}\varepsilon_{i}+o_{P}(n^{-1/2})$ because $\E (\pi^{\top}X_{i}-\alpha_{0}(\tilde{X}_{i}))^{2} \lesssim p^{-2\xi_{2}} \ll n^{-1}$ (due to $\alpha_{0}\in \mathcal{M}_{C,\xi_{2}}$ and $p \gg n^{1/(2\xi_{2}}$).}

An analogous argument (with $I_{1}$ and $I_{2}$ swapped) yields
$$\tilde{\theta}_{2}-\theta=\pi^{\top}\mathbb{E}_{n,2}X_{i}\varepsilon_{i}+o_{P}(n^{-1/2})=\sigma_{u}^{-2}\mathbb{E}_{n,2}u_{i}\varepsilon_{i}+o_{P}(n^{-1/2}).$$
Since $\hat{\theta}=\frac{n_{1}}{n}\tilde{\theta}_{1}+\frac{n_{2}}{n}\tilde{\theta}_{2}$,
we have that $\sqrt{n}(\tilde{\theta}-\theta)=n^{-1/2}\sum_{i=1}^{n}\pi^{\top}X_{i}\varepsilon_{i}+o_{P}(1)$.
Notice that 
$$\mathbb{E}(\pi^{\top}X_{i}\varepsilon_{i})^{2}=\mathbb{E}[\mathbb{E}(\varepsilon_{i}^{2}\mid X_{i})\cdot(X_{i}^{\top}\pi)^{2}]\lesssim\mathbb{E}[(X_{i}^{\top}\pi)^{2}]=\pi^{\top}\Sigma\pi \lesssim \| \pi\|_{2}^{2}\lesssim1.$$
We have $\mathbb{E}(n^{-1/2}\sum_{i=1}^{n}\pi^{\top}X_{i}\varepsilon_{i})^{2}\lesssim1$,
which means that $\sqrt{n}(\tilde{\theta}-\theta)=O_{P}(1)$.

Since either $\xi_{1}\geq\xi_{2}$ or $\xi_{1}<\xi_{2}$, we show
(\ref{eq: thm PLM eq 2}) in these two cases in two steps.

\textbf{Step 1:} show (\ref{eq: thm PLM eq 2}) assuming $\xi_{1}\geq\xi_{2}$.

Since $\max\{\xi_{1},\xi_{2}\}>1/2$ , $\xi_{1}\geq\xi_{2}$ implies
$\xi_{1}>1/2$. The assumption of $\log p\ll n^{1/2-1/(4\max\{\xi_{1},\xi_{2}\})}$
and $\xi_{1}\ge\xi_{2}$ imply $\log p\ll n^{1/2-1/(4\xi_{1})}$.
Notice that 
\begin{align}
\tilde{\theta}_{1}-\theta & =e_{1}^{\top}\hat{\gamma}-e_{1}^{\top}\gamma+\hat{\pi}^{\top}\mathbb{E}_{n,1}X_{i}Y_{i}-\hat{\pi}^{\top}\hat{\Sigma}_{(1)}\hat{\gamma}\nonumber \\
 & =e_{1}^{\top}\delta_{\gamma}+\hat{\pi}^{\top}\mathbb{E}_{n,1}X_{i}Y_{i}-\hat{\pi}^{\top}\hat{\Sigma}\gamma-\hat{\pi}^{\top}\hat{\Sigma}_{(1)}\delta_{\gamma}\nonumber \\
 & =(e_{1}^{\top}-\hat{\pi}^{\top}\hat{\Sigma}_{(1)})\delta_{\gamma}+\hat{\pi}^{\top}\mathbb{E}_{n,1}X_{i}(X_{i}^{\top}\gamma+\varepsilon_{i}+\eta_{f,i})-\hat{\pi}^{\top}\hat{\Sigma}_{(1)}\gamma\nonumber \\
 & =(e_{1}^{\top}-\hat{\pi}^{\top}\hat{\Sigma}_{(1)})\delta_{\gamma}+\pi^{\top}\mathbb{E}_{n,1}X_{i}\varepsilon_{i}+(\hat{\pi}-\pi)^{\top}\mathbb{E}_{n,1}X_{i}\varepsilon_{i}+\hat{\pi}^{\top}\mathbb{E}_{n,1}X_{i}\eta_{f,i}.\label{eq: thm PLM eq 3}
\end{align}

We notice that 
\begin{align}
\mathbb{E}\left([(\hat{\pi}-\pi)^{\top}\mathbb{E}_{n,1}X_{i}\varepsilon_{i}]^{2}\mid\{X_{i}\}_{i\in I_{1}}\right) & \overset{\text{(i)}}{=}n_{1}^{-2}\sum_{i\in I_{1}}[X_{i}^{\top}(\hat{\pi}-\pi)]^{2}\mathbb{E}(\varepsilon_{i}^{2}\mid\{X_{i}\}_{i\in I_{1}})\nonumber \\
 & \overset{\text{(ii)}}{\lesssim }n_{1}^{-2}\sum_{i\in I_{1}}[X_{i}^{\top}(\hat{\pi}-\pi)]^{2}\nonumber \\
 & =n_{1}^{-1}(\hat{\pi}-\pi)^{\top}\hat{\Sigma}_{(1)}(\hat{\pi}-\pi)\nonumber \\
 & \overset{\text{(iii)}}{=}n_{1}^{-1}O_{P}\left(\varepsilon_{n}^{4\xi_{2}/(2\xi_{2}+1)}\right)=o_{P}(n^{-1}),\label{eq: thm PLM eq 4}
\end{align}
where (i) follows by the fact that $\hat{\pi}$ only depends on $\{X_{i}\}_{i\in I_{1}}$
and $\mathbb{E}(\varepsilon_{i}\mid\{X_{i}\}_{i\in I_{1}})=0$, (ii)
follows by the assumption of $\mathbb{E}(\varepsilon_{i}^{2}\mid X_{i})$ being bounded
and (iii) follows by $(\hat{\pi}-\pi)^{\top}\hat{\Sigma}_{(1)}(\hat{\pi}-\pi)=O_{P}\left(\varepsilon_{n}^{4\xi_{2}/(2\xi_{2}+1)}\right)$
from Lemma \ref{lem: PLM tool 3}.

By Lemma \ref{lem: PLM tool 3}, with probability approaching one,
$\hat{\pi}$ is well defined and thus $\|\hat{\Sigma}_{(1)}\hat{\pi}-e_{1}\|_{\infty}\leq\lambda_{1}\lesssim\varepsilon_{n}$.
By Lemma \ref{lem: PLM tool 5} and $\xi_{1}>1/2$, $\|\delta_{\gamma}\|_{1}=O_{P}\left(\varepsilon_{n}^{(2\xi_{1}-1)/(2\xi_{1}+1)}\right)$.
Hence, 
\[
\left|(e_{1}^{\top}-\hat{\pi}^{\top}\hat{\Sigma}_{(1)})\delta_{\gamma}\right|\leq\|\hat{\Sigma}_{(1)}\hat{\pi}-e_{1}\|_{\infty}\|\delta_{\gamma}\|_{1}=O_{P}(\varepsilon_{n})\cdot O_{P}\left(\varepsilon_{n}^{(2\xi_{1}-1)/(2\xi_{1}+1)}\right)\overset{\text{(i)}}{=}o_{P}(n^{-1/2}),
\]
where (i) follows by $\xi_{1}>1/2$ and $\log p\ll n^{1/2-1/(4\xi_{1})}$.
Moreover, by $\mathbb{E}\eta_{f,i}^{2}\lesssim p^{-2\xi_{1}}$, $\mathbb{E}_{n,1}(X_{i}^{\top}\hat{\pi})^{2}=O_{P}(1)$
and $p\gg n^{1/(2\min\{\xi_{1},\xi_{2}\})}\geq n^{1/(2\xi_{1})}$, we have 
\[
\left|\hat{\pi}^{\top}\mathbb{E}_{n,1}X_{i}\eta_{f,i}\right|\leq\sqrt{\mathbb{E}_{n,1}(X_{i}^{\top}\hat{\pi})^{2}}\cdot\sqrt{\mathbb{E}_{n,1}\eta_{f,i}^{2}}=\sqrt{O_{P}(1)}\cdot\sqrt{O_{P}(p^{-2\xi_{1}})}=o_{P}(n^{-1/2}).
\]

The above four displays imply (\ref{eq: thm PLM eq 2}).

\jelenax{\textbf{Step 2:} show (\ref{eq: thm PLM eq 2}) assuming $\xi_{1}<\xi_{2}$.}

\jelenax{Since $\max\{\xi_{1},\xi_{2}\}>1/2$ , $\xi_{2}>\xi_{1}$ implies
$\xi_{2}>1/2$. The assumption of $\log p\ll n^{1/2-1/(4\max\{\xi_{1},\xi_{2}\})}$
and $\xi_{2}>\xi_{1}$ imply $\log p\ll n^{1/2-1/(4\xi_{2})}$. Notice
that (\ref{eq: thm PLM eq 3}) and (\ref{eq: thm PLM eq 4}) still
hold and thus 
\begin{align}
\tilde{\theta}_{1}-\theta & =(e_{1}^{\top}-\hat{\pi}^{\top}\hat{\Sigma}_{(1)})\delta_{\gamma}+\pi^{\top}\mathbb{E}_{n,1}X_{i}\varepsilon_{i}+o_{P}(n^{-1/2})\nonumber \\
 & =(e_{1}^{\top}-\pi^{\top}\hat{\Sigma}_{(1)})\delta_{\gamma}-\delta_{\pi}^{\top}\hat{\Sigma}_{(1)}\delta_{\gamma}+\pi^{\top}\mathbb{E}_{n,1}X_{i}\varepsilon_{i}+o_{P}(n^{-1/2}).\label{eq: thm PLM eq 6}
\end{align}}

\jelenax{Then we observe
that 
\begin{align}
(e_{1}^{\top}-\pi^{\top}\hat{\Sigma}_{(1)})\delta_{\gamma} & =\left[1-(\pi_{1}\mathbb{E}_{n,1}D_{i}^{2}+\mathbb{E}_{n,1}\pi_{-1}^{\top}Z_{i}D_{i})\right]\delta_{\theta}-\pi_{1}\mathbb{E}_{n,1}u_{i}Z_{i}^{\top}\delta_{\beta}\nonumber \\
 & \overset{\text{(i)}}{=}\left[1-\mathbb{E}_{n,1}\pi^{\top}X_{i}D_{i}\right]\delta_{\theta}+O_{P}\left(n^{-1/2}\mathbb{E}_{n,1}(Z_{i}^{\top}\delta_{\beta})^{2}\right)\nonumber \\
 & \overset{\text{(ii)}}{=}\left[1-\mathbb{E}_{n,1}\pi^{\top}X_{i}D_{i}\right]\delta_{\theta}+o_{P}(n^{-1/2})\nonumber \\
 & \overset{\text{(iii)}}{=}O_{P}(n^{-1/2})|\delta_{\theta}|+o_{P}(n^{-1/2}),\label{eq: thm PLM eq 6.5}
\end{align}
where (i) follows by the fact that $\delta_{\beta}$ only depends
on $\{Z_{i}\}_{i\in I_{1}},\{(Z_{i},Y_{i},D_{i})\}_{i\in I_{2}}$
and 
\begin{eqnarray*}
&& \mathbb{E}[(\mathbb{E}_{n,1}u_{i}Z_{i}^{\top}\delta_{\beta})^{2}\mid\{Z_{i}\}_{i\in I_{1}},\{(Z_{i},Y_{i},D_{i})\}_{i\in I_{2}}]
\\
&=&n_{1}^{-2}\sum_{i\in I_{1}}\mathbb{E}(u_{i}^{2}\mid Z_{i})(Z_{i}^{\top}\delta_{\beta})^{2}\lesssim n_{1}^{-2}\sum_{i\in I_{1}}(Z_{i}^{\top}\delta_{\beta})^{2}\lesssim n^{-1}\mathbb{E}_{n,1}(Z_{i}^{\top}\delta_{\beta})^{2}
\end{eqnarray*}
(due to $\mathbb{E}(u_{i}\mid Z_{i})=0$), (ii) follows by $\mathbb{E}_{n,1}(Z_{i}^{\top}\delta_{\beta})^{2}=O_{P}(\varepsilon_{n}^{4\xi_{1}/(2\xi_{1}+1)})=o_{P}(1)$
(due to Lemma \ref{lem: PLM tool 5}) and (iii) follows by $1-\mathbb{E}_{n,1}\pi^{\top}X_{i}D_{i}=O_{P}(n^{-1/2})$;
to see the last step, notice that $\mathbb{E}\pi^{\top}X_{i}D_{i}=1$ (due to $e_{1}=\Sigma \pi$) and  we will show $\mathbb{E}\left(\pi^{\top}X_{i}D_{i}\right)^{2} \lesssim 1$. For sub-Gaussian $X_{i}$, both $\pi^{\top}X_{i}$ and $D_{i}$ are sub-Gaussian, which means that  $\mathbb{E}\left(\pi^{\top}X_{i}D_{i}\right)^{2} \lesssim 1$; for bounded $X_{i}$, we have that $\mathbb{E}\left(\pi^{\top}X_{i}D_{i}\right)^{2} \lesssim \mathbb{E}\left(\pi^{\top}X_{i}\right)^{2} =\pi^{\top}\Sigma \pi \lesssim 1 $.}


\jelenax{Since $|\delta_{\theta}|\leq\|\delta_{\gamma}\|_{2}=O_{P}(\varepsilon_{n}^{2\xi_{2}/(2\xi_{2}+1)})=o_{P}(1)$
due to Lemma \ref{lem: PLM tool 5}, (\ref{eq: thm PLM eq 6.5}) implies
that 
\begin{equation}
(e_{1}^{\top}-\pi^{\top}\hat{\Sigma}_{(1)})\delta_{\gamma}=O_{P}(n^{-1/2})|\delta_{\theta}|+o_{P}(n^{-1/2})=o_{P}(n^{-1/2}).\label{eq: thm PLM eq 7}
\end{equation}}

\jelenax{We decompose 
\begin{equation}
\delta_{\pi}^{\top}\hat{\Sigma}_{(1)}\delta_{\gamma}=\underset{T_{1}}{\underbrace{\delta_{\pi_{1}}\left(\mathbb{E}_{n,1}D_{i}^{2}\delta_{\theta}+\mathbb{E}_{n,1}D_{i}Z_{i}^{\top}\delta_{\beta}\right)}}+\underset{T_{2}}{\underbrace{\delta_{\pi_{-1}}^{\top}\left(\mathbb{E}_{n,1}D_{i}Z_{i}\delta_{\theta}+\mathbb{E}_{n,1}Z_{i}Z_{i}^{\top}\delta_{\beta}\right)}}.\label{eq: thm PLM eq 8}
\end{equation}}

\jelenax{By Lemma \ref{lem: PLM tool 5}, $\|\mathbb{E}_{n,2}Z_{i}D_{i}\delta_{\theta}+\mathbb{E}_{n,1}Z_{i}Z_{i}^{\top}\delta_{\beta}\|_{\infty}=O_{P}(\varepsilon_{n})$,
which means that 
\begin{align}
T_{2} & =\delta_{\pi_{-1}}^{\top}\left(\mathbb{E}_{n,2}Z_{i}D_{i}\delta_{\theta}+\mathbb{E}_{n,1}Z_{i}Z_{i}^{\top}\delta_{\beta}\right)+\delta_{\pi_{-1}}^{\top}\left(\mathbb{E}_{n,1}Z_{i}D_{i}-\mathbb{E}_{n,2}Z_{i}D_{i}\right)\delta_{\theta}\nonumber \\
 & \overset{\text{(i)}}{=}O_{P}\left(\varepsilon_{n}\|\delta_{\pi_{-1}}\|_{1}\right)+\delta_{\pi_{-1}}^{\top}\left(\mathbb{E}_{n,1}Z_{i}D_{i}-\mathbb{E}_{n,2}Z_{i}D_{i}\right)\delta_{\theta}\nonumber \\
 & \overset{\text{(ii)}}{=}O_{P}\left(\varepsilon_{n}\|\delta_{\pi_{-1}}\|_{1}\right)+O_{P}\left(\|\delta_{\pi_{-1}}\|_{1}\cdot\varepsilon_{n}\cdot|\delta_{\theta}|\right)\overset{\text{(iii)}}{=}o_{P}(n^{-1/2}),\label{eq: thm PLM eq 9}
\end{align}
where (i) follows by Hölder's inequality, (ii) follows by Hölder's
inequality and $\|\mathbb{E}_{n,1}Z_{i}D_{i}-\mathbb{E}_{n,2}Z_{i}D_{i}\|_{\infty}=O(\varepsilon_{n})$
(due to Lemma \ref{lem: PLM tool 2}) and (iii) follows by 
$$|\delta_{\theta}|\leq\|\delta_{\gamma}\|_{2}=O_{P}(\varepsilon_{n}^{2\xi_{1}/(2\xi_{1}+1)})=o_{P}(1)$$
(due to $\xi_{1}>0$), $\xi_{2}>1/2$, $\|\delta_{\pi_{-1}}\|_{1}\leq\|\delta_{\pi}\|_{1}=O_{P}(\varepsilon_{n}^{(2\xi_{2}-1)/(2\xi_{2}+1)})$
(due to Lemma \ref{lem: PLM tool 3}) and $\log p\ll n^{1/2-1/(4\xi_{2})}$.}

\jelenax{For $T_{1}$, we observe
\begin{align*}
T_{1} & =\delta_{\pi_{1}}\left(\mathbb{E}_{n,1}D_{i}^{2}\delta_{\theta}+\mathbb{E}_{n,1}D_{i}Z_{i}^{\top}\delta_{\beta}\right)\\
 & =\delta_{\pi_{1}}\left(\mathbb{E}_{n,1}D_{i}^{2}\delta_{\theta}+\phi^{\top}\mathbb{E}_{n,1}Z_{i}Z_{i}^{\top}\delta_{\beta}\right)+\delta_{\pi_{1}}\mathbb{E}_{n,1}u_{i}Z_{i}^{\top}\delta_{\beta} + \delta_{\pi_{1}}\mathbb{E}_{n,1}\eta_{g,i}Z_{i}^{\top}\delta_{\beta}  \\
 & \overset{\text{(i)}}{=}\delta_{\pi_{1}}\left(\mathbb{E}_{n,1}D_{i}^{2}\delta_{\theta}+\phi^{\top}\mathbb{E}_{n,1}Z_{i}Z_{i}^{\top}\delta_{\beta}\right)+o_{P}(n^{-1/2})\\
 & =\delta_{\pi_{1}}\left(\mathbb{E}_{n,2}D_{i}^{2}\delta_{\theta}+\hat{\phi}^{\top}\mathbb{E}_{n,1}Z_{i}Z_{i}^{\top}\delta_{\beta}\right)+\delta_{\pi_{1}}\left(\mathbb{E}_{n,1}D_{i}^{2}-\mathbb{E}_{n,2}D_{i}^{2}\right)\delta_{\theta}\\
 & \qquad+\delta_{\pi_{1}}(\phi-\hat{\phi})^{\top}\mathbb{E}_{n,1}Z_{i}Z_{i}^{\top}\delta_{\beta}+o_{P}(n^{-1/2})\\
 & \overset{\text{(ii)}}{=}\delta_{\pi_{1}}\left(\mathbb{E}_{n,2}D_{i}^{2}\delta_{\theta}+\hat{\phi}^{\top}\mathbb{E}_{n,1}Z_{i}Z_{i}^{\top}\delta_{\beta}\right)+o_{P}(n^{-1/2})+\delta_{\pi_{1}}(\phi-\hat{\phi})^{\top}\mathbb{E}_{n,1}Z_{i}Z_{i}^{\top}\delta_{\beta}+o_{P}(n^{-1/2})
\end{align*}
where (i) follows by $\mathbb{E}_{n,1}u_{i}Z_{i}^{\top}\delta_{\beta}=o_{P}(n^{-1/2})$
as shown above in (\ref{eq: thm PLM eq 6.5}) and $\mathbb{E}_{n,1}\eta_{g,i}^2 =O_{P}(p^{-2\xi_2}) =o_{P}(n^{-1}) $ (due to $p\gg n^{1/(2\xi_2)}$) and (ii) follows by
$\mathbb{E}_{n,1}D_{i}^{2}-\mathbb{E}_{n,2}D_{i}^{2}=O_{P}(n^{-1/2})$
(due to the sub-Gaussianity of $D_{i}$). Notice that 
\begin{multline*}
\left|\delta_{\pi_{1}}(\phi-\hat{\phi})^{\top}\mathbb{E}_{n,1}Z_{i}Z_{i}^{\top}\delta_{\beta}\right|\leq|\delta_{\pi_{1}}|\cdot\sqrt{\mathbb{E}_{n,1}(Z_{i}^{\top}(\hat{\phi}-\phi))^{2}}\times\sqrt{\mathbb{E}_{n,1}(Z_{i}^{\top}\delta_{\beta})^{2}}\\
\overset{\text{(i)}}{\leq}O_{P}(\varepsilon_{n}^{2\xi_{2}/(2\xi_{2}+1)})\cdot O_{P}\left(\varepsilon_{n}^{2\xi_{2}/(2\xi_{2}+1)}\right)\cdot o_{P}(1)\overset{\text{(ii)}}{=}o_{P}(n^{-1/2}),
\end{multline*}
where (i) follows by $|\delta_{\pi_{1}}|\leq\|\delta_{\pi}\|_{2}=O_{P}(\varepsilon_{n}^{2\xi_{2}/(2\xi_{2}+1)})$
(due to Lemma \ref{lem: PLM tool 3}), $(\hat{\phi}-\phi)^{\top}\mathbb{E}_{n,1}Z_{i}Z_{i}^{\top}(\hat{\phi}-\phi)=O_{P}(\varepsilon_{n}^{4\xi_{2}/(2\xi_{2}+1)})$
(due to Lemma \ref{lem: PLM tool 4}) and $\mathbb{E}_{n,1}(Z_{i}^{\top}\delta_{\beta})^{2}=O_{P}(\varepsilon_{n}^{4\xi_{1}/(2\xi_{1}+1)})=o_{P}(1)$
(due to Lemma \ref{lem: PLM tool 5}) and (ii) follows by $\log p\ll n^{1/2-1/(4\xi_{2})}$.}

\jelenax{We also notice that since $\xi_{2}>1/2$, we have 
\begin{align*}
\left|\delta_{\pi_{1}}\left(\mathbb{E}_{n,2}D_{i}^{2}\delta_{\theta}+\hat{\phi}^{\top}\mathbb{E}_{n,1}Z_{i}Z_{i}^{\top}\delta_{\beta}\right)\right| & \leq|\delta_{\pi}|\cdot|\mathbb{E}_{n,2}D_{i}^{2}\delta_{\theta}+\hat{\phi}^{\top}\mathbb{E}_{n,1}Z_{i}Z_{i}^{\top}\delta_{\beta}|\\
 & \overset{\text{(i)}}{=}O_{P}(\varepsilon_{n}^{2\xi_{2}/(2\xi_{2}+1)})\cdot O_{P}(n^{-1/4})\overset{\text{(ii)}}{=}o_{P}(n^{-1/2}),
\end{align*}
where (i) follows by $|\mathbb{E}_{n,2}D_{i}^{2}(\hat{\theta}-\theta)+\hat{\phi}^{\top}\mathbb{E}_{n,1}Z_{i}Z_{i}^{\top}(\hat{\beta}-\beta)|=O_{P}(n^{-1/4})$
(due to Lemma \ref{lem: PLM tool 5}) and $|\delta_{\pi_{1}}|\leq\|\delta_{\pi}\|_{2}=O_{P}(\varepsilon_{n}^{2\xi_{2}/(2\xi_{2}+1)})$
(due to Lemma \ref{lem: PLM tool 3}) and (ii) follows by $\log p\ll n^{1/2-1/(4\xi_{2})}$.}

\jelenax{The above three displays imply $T_{1}=o_{P}(n^{-1/2})$, which together
with (\ref{eq: thm PLM eq 8}) and (\ref{eq: thm PLM eq 9}) yields
$\delta_{\pi}^{\top}\hat{\Sigma}_{(1)}\delta_{\gamma}=o_{P}(n^{-1/2})$.
Then the desired result follows by (\ref{eq: thm PLM eq 6}) and (\ref{eq: thm PLM eq 7}). }
\end{proof}

\section{Proofs for Section \ref{sec general functionals}}

For the proof of results in Section  \ref{sec general functionals}, again define $\varepsilon_{n}=\sqrt{\log(p)/n}$ and let 
$s_{0}\geq C\varepsilon_{n}^{-2/(2\xi_{2}+1)}$, and  $\pi$ be coefficients of
the least squares projection of $\alpha_0(X)$ on $\psi(X)$, satisfying%
\[
M-\Sigma\pi=\mathbb{E}[\psi(X)\{\alpha_0(X)-\psi(X)^{\top}\pi\}]=0.
\]

Suppose that $\alpha_0(\cdot)\in\mathcal{M}_{C,\xi_2}$ for some $C,\xi_2>0$. Then $\mathbb{E}(\alpha_0(X)-\psi(X)^{\top}\pi)^2\leq (C p^{-\xi_2})^2$. Moreover, $\alpha_0(\cdot)\in\mathcal{M}_{C,\xi_2}$ also implies that   $\mathbb{E}(\alpha_0(X)-\psi(X)^{\top}\tilde{\pi})^2\leq (C s_{0}^{-\xi_2})^2$ for some $\tilde{\pi}\in\RR^{p}$ such that for some $J_{0}\subset\{1,...,p\}$  with $|J_0|=s_0$ and $\tilde{\pi}_j=0$ for $j\notin J_0$. Therefore, 
$ \mathbb{E}(\psi(X)^{\top}\pi-\psi(X)^{\top}\tilde{\pi})^2\leq 2  (C p^{-\xi_2})^2+2 (C s^{-\xi_2})^2\leq  4(C s^{-\xi_2})^2 \lesssim \varepsilon_n^{2\xi_2/(2\xi_2+1)} $. Thus, 
\begin{equation}
(\pi-\tilde{\pi})^{\top}\Sigma (\pi-\tilde{\pi}) \leq C \varepsilon_n^{4\xi_2/(2\xi_2+1)}. \label{pistar2}%
\end{equation}
Since the eigenvalues of $\Sigma$ are bounded, we have 
\[
\left\Vert \pi-\tilde{\pi}\right\Vert _{2} \leq 
C\varepsilon_{n}^{2\xi_{2}/(2\xi_{2}+1)}.
\]
We define $\pi_{\ast}$ as
\begin{equation}
\pi_{\ast}\in\arg\min_{v}\;(\pi-v)^{\top}\Sigma(\pi-v)+2\varepsilon_{n}%
\sum_{j\in J_{0}^{c}}|v_{j}|\text{.} \label{pistar}%
\end{equation}

Let $J$ be the vector of indices of nonzero elements of $\pi_{\ast}.$


\begin{lemma}\label{lem: B1} If $\alpha_0(\cdot)\in \mathcal{M}_{C,\xi_2}$, then 
 $\Vert\Sigma(\pi_{\ast}-\pi)\Vert_{\infty}\leq
C\varepsilon_{n}.$
\end{lemma}

\begin{proof} [Proof of Lemma \ref{lem: B1}]

 Let $e_{j}\in\mathbb{R}^{p}$ denote the $j$-th column of
$I_{p}$. The first-order condition for $\pi^{\ast}$ implies that for $j\in
J_{0}$, we have $e_{j}^{\top}\Sigma(\pi_{\ast}-\pi)=0$; for $j\in J_{0}^{c}%
$, we have that $e_{j}^{\top}\Sigma(\pi_{\ast}-\pi)+\varepsilon_{n}z_{j}=0$,
where $z_{j}=\text{sign}(\pi_{\ast,j})$ if $\pi_{\ast,j}\neq0$ and $z_{j}%
\in\lbrack-1,1]$ if $\pi_{\ast,j}=0$. Therefore, for any $j$, we have that
$|e_{j}^{\top}\Sigma(\pi_{\ast}-\pi)|\leq\varepsilon_{n}$. Hence,
$\Vert\Sigma(\pi_{\ast}-\pi)\Vert_{\infty}\leq\varepsilon_{n}$. 
\end{proof}

\begin{lemma}\label{lem: B2}Let Assumption \ref{assp 6} hold.
If $\alpha_0(\cdot)\in \mathcal{M}_{C,\xi_2}$, then  
$$(\pi-\pi_{\ast})^{\top}\Sigma(\pi-\pi_{\ast})\leq
C\varepsilon_{n}^{4\xi_{2}/(2\xi_{2}+1)}.$$

\end{lemma}

\begin{proof}[Proof of Lemma \ref{lem: B2}]

By the definition of $\pi_{\ast}$, we have that by the bounded
eigenvalues of $\Sigma$,
\begin{align*}
(\pi-\pi_{\ast})^{\top}\Sigma(\pi-\pi_{\ast})+\varepsilon_{n}\sum_{j\in
J_{0}^{c}}|\pi_{\ast,j}|  &  \leq(\pi-\tilde{\pi})^{\top}\Sigma(\pi
-\tilde{\pi})+\varepsilon_{n}\sum_{j\in J_{0}^{c}}|\tilde{\pi}_{j}%
|=(\pi-\tilde{\pi})^{\top}\Sigma(\pi-\tilde{\pi})\\
&  \leq C\left\Vert \pi-\tilde{\pi}\right\Vert _{2}^{2}\leq C\varepsilon
_{n}^{4\xi_{2}/(2\xi_{2}+1)}.
\end{align*}
\end{proof}


\begin{lemma}\label{lem: B3} Let Assumption \ref{assp 6} hold.
If $\alpha_0(\cdot)\in \mathcal{M}_{C,\xi_2}$, then 
$\left\vert J\right\vert \leq C\varepsilon_{n}%
^{-2/(2\xi_{2}+1)}$.
\end{lemma}

\begin{proof}[Proof of Lemma \ref{lem: B3}]

 For all $j\in J\backslash J_{0}$ the first order conditions to
Equation (\ref{pistar}) imply $|e_{j}^{\top}\Sigma(\pi_{\ast}-\pi
)|=\varepsilon_{n}$. Therefore, It follows that
\[
\sum_{j\in J\backslash J_{0}}\left(  e_{j}^{\top}\Sigma(\pi_{\ast}%
-\pi)\right)  ^{2}=\frac{1}{4}\varepsilon_{n}^{2}|J\backslash J_{0}|\text{.}%
\]
In addition,
\begin{align*}
\sum_{j\in J\backslash J_{0}}\left(  e_{j}^{\top}\Sigma(\pi_{\ast}%
-\pi)\right)  ^{2}  &  \leq\sum_{j=1}^{p}\left(  e_{j}^{\top}\Sigma
(\pi_{\ast}-\pi)\right)  ^{2}=(\pi_{\ast}-\pi)^{\top}\Sigma\left(
\sum_{j=1}^{p}e_{j}e_{j}^{\top}\right)  \Sigma(\pi_{\ast}-\pi)\\
&  =(\pi_{\ast}-\pi)^{\top}\Sigma^{2}(\pi_{\ast}-\pi)\leq\lambda_{\max
}(\Sigma)\{(\pi-\pi_{\ast})^{\top}\Sigma(\pi-\pi_{\ast})\}\leq
C\varepsilon_{n}^{4\xi_{2}/(2\xi_{2}+1)}\text{,}%
\end{align*}
where the last inequality follows by Lemma \ref{lem: B2} and $\lambda_{\max}(\Sigma)\leq
C.$ Combining the above two displays, we obtain
\[
\frac{1}{4}\varepsilon_{n}^{2}|J\backslash J_{0}|\leq C\varepsilon_{n}%
^{4\xi_{2}/(2\xi_{2}+1)}\text{.}%
\]
Dividing through by $\varepsilon_{n}^{2}$ gives $|J\backslash J_{0}|\leq
C\varepsilon_{n}^{-2/(2\xi_{2}+1)}$. Thus by $s_{0}\leq C\varepsilon
_{n}^{-2/(2\xi_{2}+1)},$%
\[
|J|=|J_{0}|+|J\backslash J_{0}|=s_{0}+|J\backslash J_{0}|\leq s_{0}%
+C\varepsilon_{n}^{-2/(2\xi_{2}+1)}\leq C\varepsilon_{n}^{-2/(2\xi_{2}%
+1)}.
\]
\end{proof}


\begin{lemma}\label{lem: B4}
Let Assumption \ref{assp 6} hold.
If $\alpha_0(\cdot)\in \mathcal{M}_{C,\xi_2}$ with $\xi_{2}>1/2$, then $\Vert\pi_{\ast}%
-\pi\Vert_{1}\lesssim\varepsilon_{n}^{(2\xi_{2}-1)/(2\xi_{2}+1)}$.
\end{lemma}

\begin{proof} [Proof of Lemma \ref{lem: B4}]

By Lemma \ref{lem: C1} and $\xi_{2}>1/2$, we have that
\[
\left\Vert \pi_{J_{0}^{c}}\right\Vert _{1}\leq C s_{0}^{1/2-\xi_{2}}\leq
C\varepsilon_{n}^{(2\xi_{2}-1)/(2\xi_{2}+1)}.
\]
Let $J_{1}=J\cup J_{0}$ note that $J\subset J_{1}$ and $J_{0}\subset J_{1}$
imply $J_{1}^{c}\subset J^{c}$ and $J_{1}^{c}\subset J_{0}^{c},$ so that
\[
\left\Vert (\pi_{\ast})_{J_{1}^{c}}-\pi_{J_{1}^{c}}\right\Vert _{1}=\left\Vert
\pi_{J_{1}^{c}}\right\Vert _{1}\leq\left\Vert \pi_{J_{0}^{c}}\right\Vert
_{1}.
\]
Also, by Lemma \ref{lem: B3},
\[
\left\vert J_{1}\right\vert \leq\left\vert J\right\vert +\left\vert
J_{0}\right\vert \leq C\varepsilon_{n}^{-2/(2\xi_{2}+1)}+s_{0}\leq
C\varepsilon_{n}^{-2/(2\xi_{2}+1)}%
\]
Therefore we have
\begin{align*}
\Vert\pi_{\ast}-\pi\Vert_{1}  &  =\left\Vert (\pi_{\ast})_{J_{1}}-\pi_{J_{1}%
}\right\Vert _{1}+\left\Vert (\pi_{\ast})_{J_{1}^{c}}-\pi_{J_{1}^{c}%
}\right\Vert _{1}\leq\left\Vert (\pi_{\ast})_{J_{1}}-\pi_{J_{1}}\right\Vert
_{1}+\left\Vert \pi_{J_{0}^{c}}\right\Vert _{1}\\
&  \leq\sqrt{|J_{1}|}\left\Vert (\pi_{\ast})_{J_{1}}-\pi_{J_{1}}\right\Vert
_{2}+C\varepsilon_{n}^{(2\xi_{2}-1)/(2\xi_{2}+1)}\\
&  \leq C\varepsilon_{n}^{-1/(2\xi_{2}+1)}\Vert\pi_{\ast}-\pi\Vert
_{2}+C\varepsilon_{n}^{(2\xi_{2}-1)/(2\xi_{2}+1)}\leq C\varepsilon_{n}%
^{(2\xi_{2}-1)/(2\xi_{2}+1)}.
\end{align*}
\end{proof}


\begin{lemma}\label{lem: B5} Let Assumption \ref{assp 6} hold.
If $\alpha_0(\cdot)\in \mathcal{M}_{C,\xi_2}$, then 
 $P\left(\Vert\hat{\Sigma}\pi_{\ast}-\Sigma\pi_{\ast}\Vert_{\infty
} > C\varepsilon_n\right)=o_{P}(1)$ and $P\left(\Vert\hat{\Sigma}\pi-\Sigma\pi\Vert_{\infty} >C\varepsilon_n \right)
=o_{P}(1)$, where $C>0$ is a constant. 
\end{lemma}

\begin{proof}[Proof of Lemma \ref{lem: B5}]
\jelenax{The  results follow by the same argument as in the proof of Lemma \ref{lem: PLM tool 2}. For sub-Gaussian $\psi(X_{i})$, we apply Lemma  \ref{lem: C2}. For bounded $\psi(X_{i})$, we  apply Lemma  \ref{lem: C2} for $\psi(X_{i})\alpha_{0}(X_{i})$ and Lemma  \ref{lem: C3} for $\psi(X_{i})(\alpha_{0}(X_{i})-\psi(X_{i})^{\top}\pi_{\ast})$ and $\psi(X_{i})(\alpha_{0}(X_{i})-\psi(X_{i})^{\top}\pi)$.  We omit the details for brevity. }
\end{proof}

\bigskip

Next let%
\[
\hat{\pi}=\arg\min_{\pi}\{-2\hat{M}^{\top}\pi+\pi^{\top}\hat{\Sigma}%
\pi+2r\left\Vert \pi\right\Vert _{1}\},
\]
for $\hat{M}$ to be specified later in this appendix. Recall that $r=\kappa \varepsilon_n$. Let  $\Delta=\hat{\pi}-\pi^{\ast}$.


\begin{lemma}\label{lem: B6}Let Assumption \ref{assp 6} hold. Suppose that $\alpha_0(\cdot)\in \mathcal{M}_{C,\xi_2}$. Suppose that $\kappa>0$ is large enough. If $P\left(\left\Vert \hat{M}-M\right\Vert _{\infty
}\leq C_0 \varepsilon_{n} \right)\geq 1-o(1)$, then for
 any $\hat{J}$  such that
$(\pi^{\ast})_{\hat{J}^{c}}=0,$  with probability approaching one,
\[
\Delta^{\top}\hat{\Sigma}\Delta\leq3r\Vert\Delta\Vert_{1}\text{\ and \  }%
\Vert\Delta_{\hat{J}^{c}}\Vert_{1}\leq3\Vert\Delta_{\hat{J}}\Vert_{1}.
\]
\end{lemma}

\begin{proof}[Proof of Lemma \ref{lem: B6}]

By the definition of the estimator, we have
\[
\hat{\pi}^{\top}\hat{\Sigma}\hat{\pi}-2\hat{M}^{\top}\hat{\pi}+2r\Vert
\hat{\pi}\Vert_{1}\leq\pi_{\ast}^{\top}\hat{\Sigma}\pi_{\ast}-2\hat
{M}^{\top}\pi_{\ast}+4r\Vert\pi_{\ast}\Vert_{1}\text{.}%
\]
Plugging $\hat{\pi}=\pi_{\ast}+\Delta$ into the above equation and rearranging
the terms gives%
\begin{equation}
\Delta^{\top}\hat{\Sigma}\Delta+2r\Vert\pi_{\ast}+\Delta\Vert_{1}\leq
2r\Vert\pi_{\ast}\Vert_{1}+2(\hat{M}-\hat{\Sigma}\pi_{\ast})^{\top}\Delta.
\label{key ineq}%
\end{equation}
By assumption, $P(\mathcal{A}_1)\geq 1-o(1)$, where  $\mathcal{A}_1=\{\Vert\hat{M}-M\Vert_{\infty}\leq C_0 \varepsilon_n\}$. By Lemma
\ref{lem: B5},  $P(\mathcal{A}_2)\geq 1-o(1)$, where 
$\mathcal{A}_2=\{\Vert\hat{\Sigma}\pi_{\ast}-\Sigma\pi_{\ast}\Vert_{\infty}
\leq C_1\varepsilon_{n}\}$. Then by Lemma \ref{lem: B1}, $M=\Sigma\pi,$ and the triangle
inequality, on the event $\mathcal{A}_1\bigcap\mathcal{A}_2$,
\begin{align*}
\Vert\hat{M}-\hat{\Sigma}\pi_{\ast}\Vert_{\infty}  &  \leq\Vert\hat{\Sigma}%
\pi_{\ast}-\Sigma\pi_{\ast}\Vert_{\infty}+\Vert\hat{M}-M\Vert_{\infty
}+\left\Vert M-\Sigma\pi_{\ast}\right\Vert _{\infty}\\
&  \leq (C_0+C_1)\varepsilon_{n}+\Vert M-\Sigma\pi\Vert_{\infty}+\Vert\Sigma
(\pi_{\ast}-\pi)\Vert_{\infty}\leq (C_0+C_1+C)\varepsilon_n ,
\end{align*}
Therefore, by H\"older's inequality we have $\left\vert (\hat{M}-\hat{\Sigma
}\pi_{\ast})^{\top}\Delta\right\vert \leq\Vert\hat{M}-\hat{\Sigma}\pi_{\ast
}\Vert_{\infty}\Vert\Delta\Vert_{1},$ so that on the event $\mathcal{A}_1\bigcap\mathcal{A}_2$,
\begin{equation}
\Delta^{\top}\hat{\Sigma}\Delta+2r\Vert\pi_{\ast}+\Delta\Vert_{1}\leq
2r\Vert\pi_{\ast}\Vert_{1}+C_2\varepsilon_{n}\Vert\Delta\Vert_{1},
\nonumber\label{eq: thm bnd pi 2}%
\end{equation}
where $C_2=2(C_0+C_1+C)$. By $r=\kappa\varepsilon_n$ and $\kappa>C_2$, it follows that with probability approaching one
\begin{equation}
\Delta^{\top}\hat{\Sigma}\Delta+2r\Vert\pi_{\ast}+\Delta\Vert_{1}\leq
2r\Vert\pi_{\ast}\Vert_{1}+r\Vert\Delta\Vert_{1}.\label{eq thm bnd pi 2.5}
\end{equation}
The triangle inequality implies $\Vert\pi_{\ast}\Vert_{1}=\Vert\pi_{\ast
}+\Delta-\Delta\Vert_{1}\leq\Vert\pi_{\ast}+\Delta\Vert_{1}+\Vert\Delta
\Vert_{1}$ so subtracting $2r\Vert\pi_{\ast}+\Delta\Vert_{1}$ from both sides
gives the first conclusion.

Next, since $\Delta^{\top}\hat{\Sigma}\Delta\geq0$ it also follows from
Equation (\ref{eq thm bnd pi 2.5}) that $2r\Vert\pi_{\ast}+\Delta\Vert_{1}\leq
2r\Vert\pi_{\ast}\Vert_{1}+r\Vert\Delta\Vert_{1}$ with probability approaching
one, so dividing through by $r$ gives%
\[
2\Vert\pi_{\ast}+\Delta\Vert_{1}\leq2\Vert\pi_{\ast}\Vert_{1}+\Vert\Delta
\Vert_{1}\text{.}%
\]
It follows by $(\pi_{\ast})_{\hat{J}^{c}}=0$ that $\Vert\pi_{\ast}+\Delta
\Vert_{1}=\Vert(\pi_{\ast})_{\hat{J}}+\Delta_{\hat{J}}\Vert_{1}+\Vert
\Delta_{\hat{J}^{c}}\Vert_{1}\geq \|(\pi_{*})_{\hat{J}}\|_{1}-\|\Delta_{\hat{J}}\|_{1}+\|\Delta_{\hat{J}^{c}}\|_{1} $ and $\Vert\pi_{\ast}\Vert_{1}=\Vert(\pi_{\ast
})_{\hat{J}}\Vert_{1}$. Substituting in the previous display then gives%
$$
2\|(\pi_{*})_{\hat{J}}\|_{1}-2\|\Delta_{\hat{J}}\|_{1}+2\|\Delta_{\hat{J}^{c}}\|_{1}\leq 2\|(\pi_{*})_{\hat{J}}\|_{1}+ \|\Delta_{\hat{J}}\|_{1}+\|\Delta_{\hat{J}^{c}}\|_{1}.
$$
The second conclusion follows.
\end{proof}


\begin{lemma}\label{lem: B7}Let Assumption  \ref{assp 6} hold. Suppose that $\kappa>0$ is large enough and $n\gg (\log p)^{3+3/(2\xi_2)}$. 
If $P\left(\left\Vert \hat{M}-M\right\Vert _{\infty
}\leq C_0 \varepsilon_{n} \right)\geq 1-o(1)$, then
$$\Delta^{\top}\hat{\Sigma}\Delta=O_{P}(
\varepsilon_{n}^{4\xi_{2}/(2\xi_{2}+1)}), \quad \left\Vert \Delta\right\Vert
_{1}=O_{P}(\varepsilon_{n}^{(2\xi_{2}-1)/(2\xi_{2}+1)}) \mbox{ and  }
 \left\Vert \Delta\right\Vert _{2}=O_{P}(\varepsilon
_{n}^{2\xi_{2}/(2\xi_{2}+1)}).$$

\end{lemma}

\begin{proof}[Proof of Lemma \ref{lem: B7}]

We apply the argument in Step 1 of the proof of Lemma \ref{lem: PLM tool 3} and obtain that for any $C_{0}>0$, there exists $C_{1}>0$ such that
\begin{equation}
    P\left(\min_{|S|\leq C_{0} (n/\log p)^{1/(2\xi_{2}+1)}}\min_{\|v_{S^{c}}\|_{1}\leq3\|v_{S}\|_{1}}\frac{v'\hat{\Sigma}v}{\|v_{S}\|_{2}^{2}}\geq C_{1}\right)\geq1-o(1),\label{eq RE cond}
\end{equation}
By Lemma \ref{lem: B5}, $|J| \lesssim \varepsilon_n^{-2/(2\xi_2+1)} \lesssim  (n/\log p)^{1/(2\xi_{2}+1)}$.

For $\hat{J}=J$ it follows from the above display and Lemma \ref{lem: B6} that with high probability%
\begin{align*}
\left\Vert \Delta_{J}\right\Vert _{2}^{2}  &  \leq C\Delta^{\top}\hat
{\Sigma}\Delta\leq Cr\left\Vert \Delta\right\Vert _{1}\leq Cr\left\Vert
\Delta\right\Vert _{1}=Cr(\left\Vert \Delta_{J}\right\Vert _{1}+\left\Vert
\Delta_{J}^{c}\right\Vert _{1})\leq Cr\left\Vert \Delta_{J}\right\Vert _{1}\\
&  \leq Cr\sqrt{\left\vert J\right\vert }\left\Vert \Delta_{J}\right\Vert
_{2}\leq Cr\varepsilon_{n}^{-1/(2\xi_{2}+1)}\left\Vert \Delta_{J}\right\Vert
_{2}=C\varepsilon_{n}^{2\xi_{2}/(2\xi_{2}+1)}\left\Vert
\Delta_{J}\right\Vert _{2}.
\end{align*}
Dividing through by $\left\Vert \Delta_{J}\right\Vert _{2}$ then gives%
\[
\left\Vert \Delta_{J}\right\Vert _{2}\leq C\varepsilon
_{n}^{2\xi_{2}/(2\xi_{2}+1)}.
\]
Plugging this back in the final expression in the previous inequality gives
the first conclusion.

For the second conclusion note that by Lemma \ref{lem: B6},
\begin{align*}
\left\Vert \Delta\right\Vert _{1}  &  =\left\Vert \Delta_{J^{c}}\right\Vert
_{1}+\left\Vert \Delta_{J}\right\Vert _{1}\leq4\left\Vert \Delta
_{J}\right\Vert _{1}\leq4\sqrt{\left\vert J\right\vert }\left\Vert \Delta
_{J}\right\Vert _{2}\\
&  \leq C\varepsilon_{n}^{-1/(2\xi_{2}+1)}\varepsilon
_{n}^{2\xi_{2}/(2\xi_{2}+1)}=C\varepsilon_{n}^{(2\xi
_{2}-1)/(2\xi_{2}+1)}.
\end{align*}
For the third conclusion let $N$ denote the indices corresponding to the
largest $\left\vert J\right\vert $ entries in $\Delta_{J^{c}}$, so that
$N\subset J^{c}$, $|N|=\left\vert J\right\vert $ and $|\Delta_{j}|\geq
|\Delta_{k}|$ for any $j\in J^{c}\cap N$ and $k\in J^{c}\backslash N$. By
Lemma \ref{lem: B6} for $\hat{J}=J\cup N$ it follows exactly as in second previous
display that
\[
\left\Vert \Delta_{\hat{J}}\right\Vert _{2}\leq C\varepsilon_{n}^{2\xi_{2}/(2\xi_{2}+1)}.
\]
By Lemma 6.9 of \cite{buhlmann2011statistics} and Lemma \ref{lem: B6},
\[
\Vert\Delta_{\hat{J}^{c}}\Vert_{2}\leq(\left\vert J\right\vert )^{-1/2}%
\Vert\Delta_{\hat{J}^{c}}\Vert_{1}\leq(\left\vert J\right\vert )^{-1/2}%
3\Vert\Delta_{\hat{J}}\Vert_{1}\leq3(\left\vert J\right\vert )^{-1/2}%
\sqrt{\left\vert J\right\vert }\Vert\Delta_{J}\Vert_{2}\leq C \varepsilon_{n}^{2\xi_{2}/(2\xi_{2}+1)}.
\]
Therefore, by the triangle inequality,
\[
\Vert\Delta\Vert_{2}\leq\Vert\Delta_{\hat{J}}\Vert_{2}+\Vert\Delta_{\hat
{J}^{c}}\Vert_{2}\leq C \varepsilon_{n}^{2\xi_{2}/(2\xi
_{2}+1)},
\]
giving the third conclusion. 
\end{proof}


\begin{lemma}\label{lem: B8} Let Assumption \ref{assp 6} hold. Suppose that $\alpha_0(\cdot)\in \mathcal{M}_{C,\xi_2}$. 
If $P\left(\left\Vert \hat{M}-M\right\Vert _{\infty
}\leq C_0 \varepsilon_{n} \right)\geq 1-o(1)$, then 
$P\left(\left\Vert \hat{\Sigma}(\hat{\pi}-\pi)\right\Vert _{\infty}\leq C_1 \varepsilon_n\right)\geq 1-o(1)$ for some $C_1>0$.
\end{lemma}

\begin{proof}[Proof of Lemma \ref{lem: B8}]

The Lasso first order conditions imply $\left\Vert \hat
{\Sigma}\hat{\pi}-\hat{M}\right\Vert _{\infty}\leq r.$ By Lemma \ref{lem: B5}, $\left\Vert
\hat{\Sigma}\pi-\Sigma\pi\right\Vert _{\infty}\leq C\varepsilon_{n}$ with high probability. Then
by the triangle inequality, with high probability
\begin{align*}
\left\Vert \hat{\Sigma}(\hat{\pi}-\pi)\right\Vert _{\infty}  &  \leq\left\Vert
\hat{\Sigma}\hat{\pi}-\hat{M}\right\Vert _{\infty}+\left\Vert \hat
{M}-M\right\Vert _{\infty}+\left\Vert M-\Sigma\pi\right\Vert _{\infty
}+\left\Vert \left(  \Sigma-\hat{\Sigma}\right)  ^{\top}\pi\right\Vert
_{\infty}\\
&  \leq r + C_0\varepsilon_n +0 +C\varepsilon_n .
\end{align*}
\end{proof}


\begin{lemma}\label{lem: B9} Let Assumption  \ref{assp 6} hold. 
If $P\left(\left\Vert \hat{M}_1-M\right\Vert _{\infty
}\leq C_0 \varepsilon_{n} \right)\geq 1-o(1)$, then 
$$(\hat{\pi
}_{1}-\pi)^{\top}\hat{\Sigma}_1(\hat{\pi}_{1}-\pi)=O_{P}(\varepsilon_{n}^{4\xi_{2}/(2\xi_{2}+1)})=o_{P}(1).$$
\end{lemma}

\begin{proof}[Proof of Lemma \ref{lem: B9}]

By Lemma \ref{lem: B2}, 
$(\pi-\pi_{*})^{\top}\Sigma(\pi-\pi_{*})\lesssim\varepsilon_{n}^{4\xi_{2}/(2\xi_{2}+1)}
$. Since $\mathbb{E}(\pi-\pi_{*})^{\top}\hat{\Sigma}_{1}(\pi-\pi_{*})=(\pi-\pi_{*})^{\top}\Sigma(\pi-\pi_{*})$, Markov's inequality implies that 
$$
(\pi-\pi_{*})^{\top}\hat{\Sigma}_{1}(\pi-\pi_{*})=O_{P}(\varepsilon_{n}^{4\xi_{2}/(2\xi_{2}+1)}).
$$

By Lemma \ref{lem: B7} (applied to $\hat{\Sigma}_1$ and $\hat{\pi}_1$), $$(\hat{\pi}_{1}-\pi_{*})^{\top}\hat{\Sigma}_{1}(\hat{\pi}_{1}-\pi_{*})=O_{P}(\varepsilon_{n}^{4\xi_{2}/(2\xi_{2}+1)}).$$

The result follows by the above two displays. 
\end{proof} 

Before proving Theorem \ref{thm4} we first give an intermediate result that takes care
of important remainders in a stochastic expansion of $\hat{\theta}$ when
$\xi_{1}>1/2.$ Recall that $n_{\ell}$ denotes the number of elements of
$I_{\ell}$ and define
\begin{align*}
\theta_{n} &  :=\gamma^{\top}\Sigma\pi=\mu^{\top}\pi,\text{ }\rho
_{n}(x):=\psi(x)^{\top}\gamma,\text{ }\alpha_{n}(x):=\psi(x)^{\top}\pi.\\
\hat{U} &  :=\hat{\mu}_{1}-\hat{\Sigma}_{1}\gamma,\text{ }\hat{R}:=\hat{M}%
_{1}-\hat{\Sigma}_{1}\pi,\\
\varUpsilon_{n}(w) &  :=m(w,\rho_{n})-\theta_{n}+\alpha_{n}(x)[y-\rho_{n}(x)].
\end{align*}


\begin{lemma}\label{lem: B11}
Let Assumptions  \ref{assp 6} and \ref{assp 9}(i) hold. Then 
\[
\hat{\theta}_{1}=\theta_{n}+\frac{1}{n_{1}}\sum_{i\in I_{1}}\varUpsilon_{n}%
(W_{i})+o_{P}(n^{-1/2})=\theta_{n}+O_{P}(n^{-1/2}).
\]
\end{lemma}

\begin{proof} [Proof of Lemma \ref{lem: B11}]

By adding and subtracting terms and by linearity of
$m(w,\rho)$ in $\rho,$%
\begin{align}
&  \hat{\theta}_{1}-\theta_{n}=\frac{1}{n_{1}}\sum_{i\in I_{1}}\left\{
m(W_{i},\hat{\rho}_{2})+\hat{\alpha}_{1}(X_{i})[Y_{i}-\hat{\rho}_{2}%
(X_{i})]-\theta_{n}\right\} \label{expansion}\\
&  =\hat{M}_{1}^{\top}\hat{\gamma}_{2}+\hat{\pi}_{1}^{\top}(\hat{\mu}%
_{1}-\hat{\Sigma}_{1}\hat{\gamma}_{2})-\theta_{n}=\hat{\mu}_{1}^{\top}%
\pi-\theta_{n}+\hat{\mu}_{1}^{\top}(\hat{\pi}_{1}-\pi)+\hat{M}_{1}^{\top
}\hat{\gamma}_{2}-\hat{\pi}_{1}^{\top}\hat{\Sigma}_{1}\hat{\gamma}%
_{2}\nonumber\\
&  =\hat{\mu}_{1}^{\top}\pi-\theta_{n}+(\hat{U}+\hat{\Sigma}_{1}%
\gamma)^{\top}(\hat{\pi}_{1}-\pi)+\hat{M}_{1}^{\top}\hat{\gamma}_{2}%
-\hat{\pi}_{1}^{\top}\hat{\Sigma}_{1}\hat{\gamma}_{2}\nonumber\\
&  =\hat{\mu}_{1}^{\top}\pi-\theta_{n}+\hat{U}^{\top}(\hat{\pi}_{1}%
-\pi)+\gamma^{\top}\hat{\Sigma}_{1}(\hat{\pi}_{1}-\pi)+\hat{M}_{1}^{\top
}\hat{\gamma}_{2}-\pi^{\top}\hat{\Sigma}_{1}\hat{\gamma}_{2}-(\hat{\pi}%
_{1}-\pi)^{\top}\hat{\Sigma}_{1}\hat{\gamma}_{2}\nonumber\\
&  =\hat{\mu}_{1}^{\top}\pi-\theta_{n}+\hat{U}^{\top}(\hat{\pi}_{1}%
-\pi)+(\gamma-\hat{\gamma}_{2})^{\top}\hat{\Sigma}_{1}(\hat{\pi}_{1}%
-\pi)+\hat{R}^{\top}\hat{\gamma}_{2}\nonumber\\
&  =\hat{M}_{1}^{\top}\gamma+\pi^{\top}(\hat{\mu}_{1}-\hat{\Sigma}%
_{1}\gamma)-\theta_{n}+\hat{R}^{\top}(\hat{\gamma}_{2}-\gamma)+\hat
{U}^{\top}(\hat{\pi}_{1}-\pi)+(\gamma-\hat{\gamma}_{2})^{\top}\hat{\Sigma
}_{1}(\hat{\pi}_{1}-\pi)\nonumber\\
&  =\frac{1}{n}\sum_{i=1}^{n}\varUpsilon_{n}(W_{i})+T_{2}+T_{3}+T_{1}\text{,}%
\nonumber
\end{align}
where
$$
T_{1}    =(\gamma-\hat{\gamma}_{2})^{\top}\hat{\Sigma}_{1}(\hat{\pi}%
_{1}-\pi),\text{ }T_{2}=\hat{R}^{\top}(\hat{\gamma}_{2}-\gamma),\text{
}T_{3}=\hat{U}^{\top}(\hat{\pi}_{1}-\pi).
$$

Since $\hat{M}_{1}=\E_{n,1} m(W_{i},\psi)$ and entries of $m(W_{i},\psi)$ are sub-Gaussian, the exponential-tail of sub-Gaussian variables and union bound yield $P(\Vert\hat{M}_{1}-M\Vert_{\infty
} > C_{0}\varepsilon_n)=o_{P}(1)$ for some constant $C_{0}>0$.
By Assumption  \ref{assp 9}(i) and Lemmas \ref{lem: B4} and \ref{lem: B7}
(applied to $\hat{\gamma},$ $\gamma_{\ast},$ $\gamma$ in place of $\hat{\pi},$
$\pi_{\ast}$, and $\pi$) and the triangle inequality, $\left\Vert \hat{\gamma
}_{2}-\gamma\right\Vert _{1}=O_{P}(\varepsilon_{n}%
^{(2\xi_1-1)/(2\xi_1+1)}).$ Then by   Lemma \ref{lem: B8} and  H\"older's inequality, 
$$
\left\vert T_{1}\right\vert    \leq\left\Vert \hat{\Sigma}_{1}(\hat{\pi}%
_{1}-\pi)\right\Vert _{\infty}\left\Vert \hat{\gamma}_{2}-\gamma\right\Vert
_{1}=O_{P}(\varepsilon_n)O_{P}(\varepsilon_{n}^{(2\xi_1-1)/(2\xi_1+1)})
  =o_{P}(n^{-1/2}),
$$ 
where the last step follows by $n\gg (\log p)^{(2\xi_1)/(\xi_1-1/2)} $. Also, by Lemma \ref{lem: B5},%
\begin{align}
\left\Vert \hat{R}\right\Vert _{\infty}  &  \leq\left\Vert \hat{M}%
_{1}-M\right\Vert _{\infty}+\left\Vert (\Sigma-\hat{\Sigma}_{1})\pi\right\Vert
_{\infty}+\left\Vert M-\Sigma\pi\right\Vert _{\infty}\label{rhat rate}\\
&  =O_{P}(\varepsilon_{n})+O_{P}(\varepsilon_{n})+0=O_{P}(\varepsilon
_{n}).\nonumber
\end{align}

Therefore it follows that
$$
\left\vert T_{2}\right\vert    \leq\left\Vert \hat{R}\right\Vert _{\infty
}\left\Vert \hat{\gamma}_2-\gamma\right\Vert _{1}=O_{P}(\varepsilon_{n}%
)O_{P}(\varepsilon_{n}^{(2\xi_1-1)/(2\xi_1+1)})=o_{P}(n^{-1/2}),
$$
where the last step follows by $n\gg (\log p)^{(2\xi_1)/(\xi_1-1/2)} $.

Next, note that for $\varepsilon_{i}=Y_{i}-\rho_{0}(X_{i})$ and $d_{i}%
=\rho_{0}(X_{i})-\rho_{n}(X_{i}),$%
\[
T_{3}=T_{31}+T_{32},\text{ }T_{31}=\frac{1}{n_{1}}\sum_{i\in I_{1}}%
\varepsilon_{i}\{\psi(X_{i})^{\top}(\hat{\pi}_{1}-\pi)\},\text{ }T_{32}%
=\frac{1}{n_{1}}\sum_{i\in I_{1}}d_{i}\{\psi(X_{i})^{\top}(\hat{\pi}_{1}%
-\pi)\}.
\]
For $\tilde{X}=\{X_{i}|i\in I_{1}\}$ we have $\mathbb{E}[\varepsilon_{i}|\tilde{X}]=0$,
so that%
\[
\mathbb{E}[T_{31}|\tilde{X}]=\frac{1}{n_{1}}\sum_{i\in I_{1}}\mathbb{E}[\varepsilon_{i}%
|\tilde{X}]\{\psi(X_{i})^{\top}(\hat{\pi}_{1}-\pi)\}=0.
\]
Also, by Lemma \ref{lem: B9} and the assumption of bounded $Var(Y\mid X)$, we have
\[
Var(T_{31}|\tilde{X})=\frac{1}{n_{1}^{2}}\sum_{i\in I_{1}}Var(Y_{i}%
|X_{i})\{\psi(X_{i})^{\top}(\hat{\pi}_{1}-\pi)\}^{2}\leq\frac{C}{n_{1}}%
(\hat{\pi}_{1}-\pi)^{\top}\hat{\Sigma}_{1}(\hat{\pi}_{1}-\pi)=o_{P}%
(1/n_{1}),
\]
so that $T_{31}=o_{P}(n^{-1/2}).$ Furthermore by the Cauchy-Schwartz and Markov's inequalities and Lemma \ref{lem: B9},
$$
\left\vert T_{32}\right\vert    \leq\sqrt{\frac{1}{n_{1}}\sum_{i\in I_{1}}d_{i}^{2}}\sqrt{(\hat{\pi}_{1}-\pi)^{\top}\hat{\Sigma}_{1}(\hat{\pi}%
_{1}-\pi)} =O_{P}(\left\Vert \rho_{0}-\rho_{n}\right\Vert _{2})\cdot o_P(1)=o_{P}(n^{-1/2}),
$$
where the last step follows by $\|\rho_0-\rho_n\|_2 \lesssim p^{-\xi_1}$ (due to $\rho_{0}\in \mathcal{M}_{C,\xi_{1}}$) and $p\gg n^{1/(2\xi_1)}$. 

The conclusion then follows by the triangle inequality. 
\end{proof}

We also give a result corresponding to Lemma \ref{lem: B11} for the case where $\xi
_{1}\xi_{2}>1/4.$


\begin{lemma}\label{lem: B12} 
Let Assumptions   \ref{assp 6} and  \ref{assp 9}(ii) hold. Then
\[
\hat{\theta}_{1}=\theta_{n}+\frac{1}{n_{1}}\sum_{i\in I_{1}}\varUpsilon_{n}
(W_{i})+o_{P}(n^{-1/2})=\theta_{n}+O_{P}(n^{-1/2}).
\]
\end{lemma}

\begin{proof} [Proof of Lemma \ref{lem: B12}]

 We consider the same remainder decomposition as the proof of Lemma \ref{lem: B11} with

$T_{1}=(\gamma-\hat{\gamma}_{2})^{\top}\hat{\Sigma}_{1}(\hat{\pi}_{1}-\pi),$
$T_{2}=\hat{R}^{\top}(\hat{\gamma}_{2}-\gamma),$ $T_{3}=\hat{U}^{\top
}(\hat{\pi}_{1}-\pi).$ Note that%

\begin{align*}
T_{1}  &  =\frac{1}{n_{1}}\sum_{i\in I_{1}}[\hat{\rho}_{2}(X_{i})-\rho
_{n}(X_{i})][\hat{\alpha}_{1}(X_{i})-\alpha_{n}(X_{i})],\\
T_{2}  &  =\frac{1}{n_{1}}\sum_{i\in I_{1}}[m(W_{i},\hat{\rho}_{2}-\rho
_{n})+\alpha_{n}(X_{i})\{\hat{\rho}_{2}(X_{i})-\rho_{n}(X_{i})\}],\\
T_{3}  &  =\frac{1}{n_{1}}\sum_{i\in I_{1}}[\hat{\alpha}(X_{i})-\alpha
_{n}(X_{i})]\{Y_{i}-\rho_{n}(X_{i})\}.
\end{align*}
We consider first $T_{1}$. By Lemmas \ref{lem: B2} and \ref{lem: B7} and the Markov inequality
\begin{align*}
&  \frac{1}{n_{1}}\sum_{i\in I_{1}}[\hat{\alpha}_{1}(X_{i})-\alpha_{n}%
(X_{i})]^{2}\\
&  =(\hat{\pi}_{1}-\pi)^{\top}\hat{\Sigma}_{1}(\hat{\pi}_{1}-\pi)\leq
2(\hat{\pi}_{1}-\pi_*)^{\top}\hat{\Sigma}_{1}(\hat{\pi}_{1}-\pi_*)+2(\pi-\pi_{\ast})^{\top}\hat{\Sigma
}_{1}({\pi}-\pi_{\ast})\\
&  =O_{P}(\varepsilon_{n}^{4\xi_{2}/(2\xi_{2}%
+1)})+O_{P}(\mathbb{E}[(\pi-\pi_{\ast})^{\top}\hat{\Sigma}_{1}({\pi}-\pi_{\ast
})])\\
&  =O_{P}(\varepsilon_{n}^{4\xi_{2}/(2\xi_{2}%
+1)})+O_{P}((\pi-\pi_{\ast})^{\top}\Sigma({\pi}-\pi_{\ast}%
))=O_{P}(\varepsilon_{n}^{4\xi_{2}/(2\xi_{2}+1)}).
\end{align*}
Similarly, 
\begin{align*}
&\mathbb{E}[\frac{1}{n_{1}}\sum_{i\in I_{1}}[\hat{\gamma}_{2}'\psi(X_{i})-\gamma'
\psi(X_{i})]^{2}|\{W_{i}\}_{i\in I_{2}}]  \\
& =(\hat{\gamma}_{2}-\gamma)^{\top}\Sigma (\hat{\gamma}_{2}-\gamma)\leq 2(\hat{\gamma}_{2}-\gamma_*)^{\top}\Sigma(\hat{\gamma}_{2}-\gamma_*)
+2(\gamma-\gamma_{\ast})^{\top}\Sigma(\gamma-\gamma_{\ast})\\
&  \lesssim \left\Vert \hat{\gamma}_{2}-\gamma_*\right\Vert _{2}^{2}+\varepsilon_{n}^{4\xi
_{1}/(2\xi_{1}+1)}=O_{P}(\varepsilon_{n}^{4\xi
_{1}/(2\xi_{1}+1)}).
\end{align*}
It then follows by the conditional Markov inequality that%
\[
\frac{1}{n_{1}}\sum_{i\in I_{1}}[\hat{\rho}_{2}(X_{i})-\rho_n(X_{i})]^{2}=\frac{1}{n_{1}}\sum_{i\in I_{1}}[\hat{\gamma}_{2}'\psi(X_{i})-\gamma'\psi(X_{i})]^{2}=O_{P}(\varepsilon_{n}^{4\xi_{1}/(2\xi
_{1}+1)}).
\]
Then by the Cauchy-Schwartz inequality,%
\begin{align*}
\left\vert T_{1}\right\vert  &  \leq\{\frac{1}{n_{1}}\sum_{i\in I_{1}}%
[\hat{\rho}_{2}(X_{i})-\rho_{n}(X_{i})]^{2}\}^{1/2}\{\frac{1}{n_{1}}%
\sum_{i\in I_{1}}[\hat{\alpha}_{1}(X_{i})-\alpha_{n}(X_{i})]^{2}\}^{1/2}\\
&  =O_{P}(\varepsilon_{n}^{2\xi_{1}/(2\xi_{1}%
+1)+2\xi_{2}/(2\xi_{2}+1)})\overset{\text{(i)}}{=}o_{P}(1/\sqrt{n}),
\end{align*}
where (i) follows by the assumption of $n^{\xi_1\xi_2-1/4}\gg (\log p)^{2\xi_1\xi_2+(\xi_1+\xi_2)/2}$.

Also

\begin{align*}
T_{2}  &  =T_{21}+T_{22},\text{ }T_{21}=\frac{1}{n_{1}}\sum_{i\in I_{1}%
}[m(W_{i},\hat{\rho}_{2}-\rho_{n})+\alpha_0(X_{i})\{\hat{\rho}_{2}%
(X_{i})-\rho_{n}(X_{i})\}],\\
T_{22}  &  =\frac{1}{n_{1}}\sum_{i\in I_{1}}[\{\alpha_{n}(X_{i})-\alpha_0(X_{i})\}\{\hat{\rho}_{2}(X_{i})-\rho_{n}(X_{i})\}.
\end{align*}
Note that  for $i\in I_1$,
\[
\mathbb{E}[m(W_{i},\hat{\rho}_{2}-\rho_{n})+\alpha_0(X)\{\hat{\rho}_{2}(X)-\rho
_{n}(X)\}|\{W_{i}\}_{i\in I_{2}}]=0.
\]
Let $\Delta_{2,\gamma}=\hat{\gamma}_2-\gamma$. Then $\hat{\rho}_{2}(X_{i})-\rho_{n}(X_{i})=\psi(X_{i})^{\top} \Delta_{2,\gamma}$. Notice that $\|\Delta_{\gamma,2}\|_{2}\leq\|\hat{\gamma}_{2}-\gamma_{*}\|_{2}+\|\gamma_{*}-\gamma\|_{2}=O_{P}(\varepsilon_{n}^{2\xi_{1}/(2\xi_{1}+1)})=o_{P}(1)$.  By the assumption of bounded eigenvalues of  $\mathbb{E}[\alpha_0(X)^{2} \psi(X)\psi(X)^{\top}]$ and $\mathbb{E} [m(W,\psi)m(W,\psi)^{\top}]$,  we have that for $i\in I_1$,
$$
\mathbb{E}[|m(W_{i},\hat{\rho}_{2}-\rho_{n})|^{2}|\{W_{i}\}_{i\in I_{2}}] 
  =\Delta_{2,\gamma}'\mathbb{E}[m(W,\psi)m(W,\psi)^{\top}]\Delta_{2,\gamma} \lesssim\|\Delta_{2,\gamma}\|_{2}^{2}=o_{P}(1)
$$
and 
$$
\mathbb{E}[\alpha_{0}(X_{i})^{2}\{\hat{\rho}_{2}(X_{i})-\rho_{n}(X_{i})\}^{2}|\{W_{i}\}_{i\in I_{2}}]=\Delta_{2,\gamma}'\mathbb{E}[\alpha_{0}(X)^2\psi(X)\psi(X)^{\top}]\Delta_{2,\gamma}
\lesssim\|\Delta_{2,\gamma}\|_{2}^{2}=o_{P}(1).
$$

Therefore it follows by independent observations that $\mathbb{E}[T_{21}^{2}%
|\{W_{i}\}_{i\in I_{2}}]=o_{P}(n^{-1}),$ so by the conditional Markov
inequality $T_{21}=o_{P}(n^{-1/2}).$ Similarly to the result for $T_{1}$, we observe that $\alpha_0\in \Mcal_{C,\xi_2} $ implies that $\mathbb{E}[(\alpha_0(X)-\alpha_n(X))^2]\lesssim p^{-2\xi_2} =o(n^{-1})$ due to $p\gg n^{1/(2\xi_2)}$. Then  $\frac{1}{n_{1}}\sum_{i\in I_{1}}%
[\alpha_0(X_{i})-\alpha_{n}(X_{i})]^{2}=o_P(n^{-1/2})$ and thus 
\begin{align*}
\left\vert T_{22}\right\vert  &  \leq\{\frac{1}{n_{1}}\sum_{i\in I_{1}}%
[\alpha_0(X_{i})-\alpha_{n}(X_{i})]^{2}\}^{1/2}\{\frac{1}{n_{1}}\sum_{i\in
I_{1}}[\hat{\rho}_{2}(X_{i})-\rho_{n}(X_{i})]^{2}\}^{1/2}\\
&  =o_{P}(n^{-1/2})O_{P}%
(\varepsilon_{n}^{2\xi_{1}/(2\xi_{1}+1)})=o_{P}(n^{-1/2}).
\end{align*}

Then by the triangle inequality $T_{2}=o_{P}(n^{-1/2})$. It also follows
similarly to the proof of Lemma \ref{lem: B11} that $T_3=o_P(n^{-1/2})$. The conclusion then follows by the triangle inequality.
\end{proof}

\begin{proof}[Proof of Theorem \ref{thm4}]

By Lemma \ref{lem: B11} or Lemma \ref{lem: B12} (and the same results for
$I_{2}$), we have 
$$
\hat{\theta}   =\frac{n_{1}}{n}\hat{\theta}_{1}+\frac{n_{2}}{n}\hat{\theta
}_{2}=\sum_{\ell=1}^{2}\frac{n_{\ell}}{n}\{\theta_{n}+\frac{1}{n_{\ell}}%
\sum_{i\in I_{\ell}}\varUpsilon_{n}(W_{i})+o_{P}(n^{-1/2})\}\\
  =\theta_{n}+\frac{1}{n}\sum_{i=1}^{n}\varUpsilon_{n}(W_{i})+o_{P}(n^{-1/2}%
).
$$

Notice that $\theta_{n}=\rho'\Sigma\gamma=\mathbb{E}\alpha_{n}(X)\rho_{n}(X)=\mathbb{E}\alpha_{0}(X)\rho_{n}(X)$
since $\E \psi(X)(\alpha_{n}(X)-\alpha_{0}(X))=0$. By $\rho_{0}\in\mathcal{M}_{C,\xi_{1}}$,
we have that $\mathbb{E}(\rho_{0}(X)-\rho_{n}(X))^{2}\lesssim p^{-2\xi_{1}}$
and thus 
\begin{align*}
\left|\theta_{0}-\theta_{n}\right| & =\left|\mathbb{E}\alpha_{0}(X)\rho_{0}(X)-\mathbb{E}\alpha_{0}(X)\rho_{n}(X)\right|\\
 & =\left|\mathbb{E}\alpha_{0}(X)(\rho_{0}(X)-\rho_{n}(X))\right|\\
 & \leq\sqrt{\mathbb{E}\alpha_{0}(X)^{2}}\cdot\sqrt{\mathbb{E}(\rho_{0}(X)-\rho_{n}(X))^{2}}\lesssim p^{-\xi_{1}}\overset{\text{(i)}}{=}o(n^{-1/2}),
\end{align*}
where (i) follows by $p\gg n^{1/(2\xi_{1})}$. Hence, 
\[
\hat{\theta}-\theta_{0}=n^{-1}\sum_{i=1}^{n}\varUpsilon_{n}(W_{i})+o_{P}(n^{-1/2}).
\]

Define $\tilde{\varUpsilon}_{n}(w)=m(w,\rho_{n})-\theta_{n}+\alpha_{n}(x)[y-\rho_{0}(x)]$.
Then by $p\gg n^{1/(2\xi_{1})}$, 
\begin{align*}
\mathbb{E}\left|n^{-1}\sum_{i=1}^{n}[\varUpsilon_{n}(W_{i})-\tilde{\varUpsilon}_{n}(W_{i})]\right| & \leq \E |\varUpsilon_{n}(W_{i})-\tilde{\varUpsilon}_{n}(W_{i})|\\
 & =\E |\alpha_{n}(X)[\rho_{n}(X)-\rho_{0}(X)]|\\
 & \leq\sqrt{\mathbb{E}\alpha_{n}(X)^{2}}\cdot\sqrt{\mathbb{E}(\rho_{0}(X)-\rho_{n}(X))^{2}}\\
 & \overset{\text{(i)}}{\lesssim}p^{-\xi_{1}}=o(n^{-1/2}),
\end{align*}
where (i) follows by $\sqrt{\mathbb{E}\alpha_{n}(X)^{2}}\leq\sqrt{\mathbb{E}\alpha_{0}(X)^{2}}+\sqrt{\mathbb{E}(\alpha_{n}(X)-\alpha_{0}(X))^{2}}=O(1)+O(p^{-\xi_{2}})$
(due to $\alpha_{0}\in\mathcal{M}_{C,\xi_{2}}$) and $\mathbb{E}(\rho_{0}(X)-\rho_{n}(X))^{2}\lesssim p^{-2\xi_{1}}$.
Therefore, 
\[
\hat{\theta}-\theta_{0}=n^{-1}\sum_{i=1}^{n}\tilde{\varUpsilon}_{n}(W_{i})+o_{P}(n^{-1/2}).
\]

Recall that  $\theta_{n}=\mathbb{E}\alpha_{0}(X)\rho_{n}(X)$ (shown above). Thus, 
\[
\mathbb{E}\tilde{\varUpsilon}_{n}(W)=\E m(W,\rho_{n})-\theta_{n}=\mathbb{E}\alpha_{0}(X)\rho_{n}(X)-\theta_{n}=0.
\]

We notice that 
\begin{align*}
\mathbb{E}\tilde{\varUpsilon}_{n}(W)^{2} & =\mathbb{E}[m(W,\rho_{n})-\theta_{n}]^{2}+\mathbb{E}\alpha_{n}(X)^{2}[Y-\rho_{0}(X)]^{2}\\
 & =\E m(W,\rho_{n})^{2}-\theta_{n}^{2}+\mathbb{E}\alpha_{n}(X)^{2}[Y-\rho_{0}(X)]^{2}\\
 & \leq \E m(W,\rho_{n})^{2}+\mathbb{E}\alpha_{n}(X)^{2}[Y-\rho_{0}(X)]^{2}\\
 & \overset{\text{(i)}}{=}\E m(W,\rho_{n})^{2}+O(1)\mathbb{E}\alpha_{n}(X)^{2}=\E m(W,\rho_{n})^{2}+O(1),
\end{align*}
where (i) follows by the assumption that $\E |Y-\rho_{0}(X)|^{2+c}$
is bounded. Notice that $\mathbb{E}\rho_{0}(X)^{2}=\mathbb{E}\rho_{n}(X)^{2}+\mathbb{E}(\rho_{0}(X)-\rho_{n}(X))^{2}$
and $\mathbb{E}\rho_{0}(X)^{2}\leq EY^{2}$ is bounded. Thus, $\mathbb{E}\rho_{n}(X)^{2}=O(1)$.
Since $\rho_{n}(X)=\psi(X)^{\top}\gamma$, we have that $\gamma^{\top}\Sigma\gamma=\mathbb{E}\rho_{n}(X)^{2}=O(1)$.
By the bounded eigenvalue assumption of $\Sigma^{-1}$, we have $\|\gamma\|_{2}=O(1)$.
Therefore, by Assumption \ref{assp 6},
\[
\E m(W,\rho_{n})^{2}=\gamma^{\top}\mathbb{E}[m(W,\psi)m(W,\psi)^{\top}]\gamma\lesssim\|\gamma\|_{2}^{2}=O(1).
\]

The above displays imply $\mathbb{E}\tilde{\varUpsilon}_{n}(W)^{2}=O(1)$. By $\mathbb{E}\tilde{\psi}(W)=0$,
we have $n^{-1}\sum_{i=1}^{n}\tilde{\varUpsilon}_{n}(W_{i})=O_{P}(n^{-1/2})$.
Therefore, $\hat{\theta}-\theta_{0}=O_{P}(n^{-1/2})$. 
\end{proof}

\begin{proof}[Proof of Theorem \ref{thm: asym norm gen}]
\jelenax{By Assumption \ref{assu: asy normal}, 
\begin{align*}
\E |\tilde{\varUpsilon}_{n}(W_{i})|^{2+c} & \lesssim \E |m(W,\rho_{n})|^{2+c}+\E |\alpha_{0}(X)|^{2+c}|Y-\rho_{0}(X)|^{2+c}\\
 & \overset{\text{(i)}}{\lesssim}\E |m(W,\rho_{n})|^{2+c}+\E |\alpha_{0}(X)|^{2+c}\lesssim1,
\end{align*}
where (i) follows by Assumption \ref{assu: asy normal}. Without loss of generality, we
can assume that $c\in(0,1)$; if $c\geq1$, we can reduce $c$ and
Assumption \ref{assu: asy normal} still holds. By Theorem \ref{thm4},
$\mathbb{E}\tilde{\varUpsilon}_{n}(W)=0$. Notice that $\mathbb{E}\tilde{\varUpsilon}_{n}(W)^{2}\geq \mathbb{E}\alpha_{0}(X)^{2}(Y-\rho_{0}(X))^{2}$,
which is assumed to be bounded below by a positive constant. By a
generalized Berry-Esseen bound (Theorem 6.3 in Chapter 7 of \cite{gut2013probability}),
it follows that 
\[
\sup_{x\in\mathbb{R}}\left|P\left(\frac{n^{-1/2}\sum_{i=1}^{n}\tilde{\varUpsilon}_{n}(W_{i})}{\sqrt{\mathbb{E}\tilde{\varUpsilon}_{n}(W)^{2}}}\leq x\right)-\Phi(x)\right|\lesssim\frac{n\E |\tilde{\varUpsilon}_{n}(W)|^{2+c}}{n^{1+c/2}}\lesssim n^{-c/2},
\]
where $\Phi(\cdot)$ is the cdf of the standard normal distribution
$N(0,1)$. The above display implies that 
\[
\frac{n^{-1/2}\sum_{i=1}^{n}\tilde{\varUpsilon}_{n}(W_{i})}{\sqrt{\mathbb{E}\tilde{\varUpsilon}_{n}(W)^{2}}}\overset{d}{\rightarrow}N(0,1).
\]
Notice that $\mathbb{E}\tilde{\varUpsilon}_{n}(W)^{2}=\mathbb{E}[m(W_{i},\rho_{n})-\theta_{n}]^{2}+\mathbb{E}\alpha_{n}(X_{i})^{2}[Y_{i}-\rho_{0}(X_{i})]^{2}$.
Define $\bar{V}=(n_{1}/n)\bar{V}_{1}+(n_{2}/n)\bar{V}_{2}$, where
$\bar{V}_{\ell}=\frac{1}{n_{\ell}}\sum_{i\in I_{\ell}}[(m(W_{i},\rho_{n})-\theta_{n})^{2}+\alpha_{n}(X_{i})^{2}[Y_{i}-\rho_{0}(X_{i})]^{2}]$
for $\ell\in\{1,2\}$. Since $\E |\tilde{\varUpsilon}_{n}(W_{i})|^{2+c}\lesssim1$,
the law of large numbers (e.g., Corollary 8.2 in Chapter 3 of \cite{gut2013probability})
implies that $\bar{V}=\mathbb{E}\tilde{\varUpsilon}_{n}(W)^{2}+o_{P}(1)$. To show
that $\hat{V}=n^{-1}\sum_{i=1}^{n}\tilde{\varUpsilon}_{n}(W_{i})^{2}+o_{P}(1)$,
we need to show $\hat{V}_{\ell}-\bar{V}_{\ell}=o_{P}(1)$ for $\ell\in\{1,2\}$.
We do so for $\ell=1$; the argument for $\ell=2$ is analogous. Notice
that 
\begin{align}
\hat{V}_{1}-\bar{V}_{1} & =\frac{1}{n_{1}}\sum_{i\in I_{1}}\left[\left(m(W_{i},\hat{\rho}_{2})-\hat{\theta}\right)^{2}-\left(m(W_{i},\rho_{n})-\theta_{n}\right)^{2}\right]\nonumber \\
 & \quad+\frac{1}{n_{1}}\sum_{i\in I_{1}}\left[\hat{\alpha}_{2}(X_{i})^{2}[Y_{i}-\hat{\rho}_{2}(X_{i})]^{2}-\alpha_{n}(X_{i})^{2}[Y_{i}-\rho_{0}(X_{i})]^{2}\right]\nonumber \\
 & =T_{1}+T_{2}+T_{3}+T_{4}+T_{5},\label{eq: thm asy norm eq 6}
\end{align}
 where 
\begin{align*}
T_{1} & =\frac{1}{n_{1}}\sum_{i\in I_{1}}\left[\left(m(W_{i},\hat{\rho}_{2})-\hat{\theta}\right)^{2}-\left(m(W_{i},\rho_{n})-\theta_{n}\right)^{2}\right]\\
T_{2} & =\frac{1}{n_{1}}\sum_{i\in I_{1}}\left(\hat{\alpha}_{2}(X_{i})^{2}-\alpha_{n}(X_{i})^{2}\right)[Y_{i}-\rho_{0}(X_{i})]^{2}\\
T_{3} & =\frac{2}{n_{1}}\sum_{i\in I_{1}}\hat{\alpha}_{2}(X_{i})^{2}(Y_{i}-\rho_{0}(X_{i}))(\rho_{0}(X_{i})-\rho_{n}(X_{i}))\\
T_{4} & =\frac{2}{n_{1}}\sum_{i\in I_{1}}\hat{\alpha}_{2}(X_{i})^{2}(Y_{i}-\rho_{0}(X_{i}))(\rho_{n}(X_{i})-\hat{\rho}_{2}(X_{i}))\\
T_{5} & =\frac{1}{n_{1}}\sum_{i\in I_{1}}\hat{\alpha}_{2}(X_{i})^{2}(\rho_{0}(X_{i})-\hat{\rho}_{2}(X_{i}))^{2}.
\end{align*}}

\jelenax{For $\ell\in\{1,2\}$, define $\delta_{\alpha,\ell}(X)=\hat{\alpha}_{\ell}(X)-\alpha_{n}(X)$
and $\delta_{\rho,\ell}(X)=\hat{\rho}_{\ell}(X)-\rho_{n}(X)$. We
bound $T_{1}$,...,$T_{5}$ in five steps. }

\jelenax{\textbf{Step 1:} show $T_{1}=o_{P}(1)$.}

\jelenax{Let $\Delta_{\gamma,2}=\hat{\gamma}_{2}-\gamma$. Since $\Delta_{\gamma,2}$
is independent of $\{W_{i}\}_{i\in I_{1}}$, by $\delta_{\rho,2}(x)=\psi(x)^{\top}(\hat{\gamma}_{2}-\gamma)$
and Assumption \ref{assp 6}, it follows that 
$$
\mathbb{E}\left[\frac{1}{n_{1}}\sum_{i\in I_{1}}m(W_{i},\delta_{\rho,2})^{2}\mid\{W_{i}\}_{i\in I_{2}}\right]=\Delta_{\gamma,2}^{\top}\mathbb{E}\left[m(W,\psi)(W,\psi)^{\top}\right]\Delta_{\gamma,2}
\lesssim\|\Delta_{\gamma,2}\|_{2}^{2}\overset{\text{(i)}}{=}o_{P}(1),
$$
where (i) follows by $\|\Delta_{\gamma,2}\|_{2}=o_{P}(1)$
as argued in the proof of Lemma \ref{lem: B12}. Thus, $\frac{1}{n_{1}}\sum_{i\in I_{1}}m(W_{i},\delta_{\rho,2})^{2}=o_{P}(1)$.
Since $\hat{\theta}-\theta_{n}=o_{P}(1)$ (shown in the proof of Theorem
\ref{thm4}), we have that $\frac{1}{n_{1}}\sum_{i\in I_{1}}[m(W_{i},\delta_{\rho,2})+\theta_{n}-\hat{\theta}]^{2}=o_{P}(1)$.
Therefore, 
\begin{multline*}
T_{1}=\frac{1}{n_{1}}\sum_{i\in I_{1}}\left(m(W_{i},\delta_{\rho,2})+\theta_{n}-\hat{\theta}\right)^{2}\\
+\frac{2}{n_{1}}\sum_{i\in I_{1}}\left(m(W_{i},\rho_{n})-\theta_{n}\right)\left(m(W_{i},\delta_{\rho,2})+\theta_{n}-\hat{\theta}\right)=o_{P}(1).
\end{multline*}}

\jelenax{\textbf{Step 2:} show $T_{2}=o_{P}(1)$.}

\jelenax{Since $\mathbb{E}([Y_{i}-\rho_{0}(X_{i})]^{2}\mid X_{i})$ is bounded (due to  the assumption of bounded $Var(Y\mid X)$), we have 
\begin{align*}
 & \mathbb{E}\left[\frac{1}{n_{1}}\sum_{i\in I_{1}}\left|\left(\hat{\alpha}_{2}(X_{i})^{2}-\alpha_{n}(X_{i})^{2}\right)[Y_{i}-\rho_{0}(X_{i})]^{2}\right|\mid\{W_{i}\}_{i\in I_{2}}\right]\\
 & \lesssim\frac{1}{n_{1}}\sum_{i\in I_{1}}\left|\hat{\alpha}_{2}(X_{i})^{2}-\alpha_{n}(X_{i})^{2}\right|\\
 & =\frac{1}{n_{1}}\sum_{i\in I_{1}}\left|(\hat{\alpha}_{2}(X_{i})-\alpha_{n}(X_{i}))(\hat{\alpha}_{2}(X_{i})+\alpha_{n}(X_{i}))\right|\\
 & \leq\sqrt{\frac{1}{n_{1}}\sum_{i\in I_{1}}[\hat{\alpha}_{2}(X_{i})-\alpha_{n}(X_{i})]^{2}}\cdot\sqrt{\frac{1}{n_{1}}\sum_{i\in I_{1}}[\hat{\alpha}_{2}(X_{i})+\alpha_{n}(X_{i})]^{2}}.
\end{align*}}

\jelenax{Let $\Delta_{\pi,2}=\hat{\pi}_{2}-\pi$ and note $\hat{\alpha}_{2}(X)-\alpha_{n}(X)=\psi(X)^{\top}\Delta_{\pi,2}$.
Notice that 
\begin{equation}
\mathbb{E}\left(\frac{1}{n_{1}}\sum_{i\in I_{1}}[\hat{\alpha}_{2}(X_{i})-\alpha_{n}(X_{i})]^{2}\mid\{W_{i}\}_{i\in I_{2}}\right)\lesssim\Delta_{\pi,2}^{\top}\Sigma\Delta_{\pi,2}\overset{\text{(i)}}{=}o_{P}(1),\label{eq: thm asy norm eq 9}
\end{equation}
where (i) follows by $\|\Delta_{\pi,2}\|_{2}\leq\|\hat{\pi}_{2}-\pi_{*}\|_{2}+\|\pi_{*}-\pi\|_{2}=O_{P}(\varepsilon_{n}^{2\xi_{2}/(2\xi_{2}+1)})=o_{P}(1)$;
Lemma \ref{lem: B7} gives $\|\hat{\pi}_{2}-\pi_{*}\|_{2}=O_{P}(\varepsilon_{n}^{2\xi_{2}/(2\xi_{2}+1)})$
and the proof of Lemma \ref{lem: B2} gives $(\pi_{*}-\pi)^{\top}\Sigma(\pi_{*}-\pi)=O(\varepsilon_{n}^{2\xi_{2}/(2\xi_{2}+1)})$,
which implies $\|\pi_{*}-\pi\|_{2}=O(\varepsilon_{n}^{2\xi_{2}/(2\xi_{2}+1)})$.
Therefore, the above two displays imply 
\[
|T_{2}|\leq\frac{1}{n_{1}}\sum_{i\in I_{1}}\left|\left(\hat{\alpha}_{2}(X_{i})^{2}-\alpha_{n}(X_{i})^{2}\right)[Y_{i}-\rho_{0}(X_{i})]^{2}\right|=o_{P}(1)\sqrt{\frac{1}{n_{1}}\sum_{i\in I_{1}}[\hat{\alpha}_{2}(X_{i})+\alpha_{n}(X_{i})]^{2}}.
\]}

\jelenax{Since $[\hat{\alpha}_{2}(x)+\alpha_{n}(x)]^{2}=[\hat{\alpha}_{2}(x)-\alpha_{n}(x)+2\alpha_{n}(x)]^{2}\leq2[\hat{\alpha}_{2}(x)-\alpha_{n}(x)]^{2}+8\alpha_{n}(x)^{2}$
and $\alpha_{n}(x)=[\alpha_{n}(x)-\alpha_{0}(x)+\alpha_{0}(x)]^{2}\leq2[\alpha_{n}(x)-\alpha_{0}(x)]^{2}+2\alpha_{0}(x)^{2}$,
we have 
\begin{eqnarray*}
&&\frac{1}{n_{1}}\sum_{i\in I_{1}}[\hat{\alpha}_{2}(X_{i})+\alpha_{n}(X_{i})]^{2}
\\
&\leq&\frac{2}{n_{1}}\sum_{i\in I_{1}}[\hat{\alpha}_{2}(X_{i})-\alpha_{n}(X_{i})]^{2}+\frac{16}{n_{1}}\sum_{i\in I_{1}}[\alpha_{n}(X_{i})-\alpha_{0}(X_{i})]^{2} 
+\frac{16}{n_{1}}\sum_{i\in I_{1}}\alpha_{0}(X_{i})^{2}\overset{\text{(i)}}{=}O_{P}(1),
\end{eqnarray*}
where (i) follows by $\frac{1}{n_{1}}\sum_{i\in I_{1}}[\hat{\alpha}_{2}(X_{i})-\alpha_{n}(X_{i})]^{2}=o_{P}(1)$,
$\mathbb{E}\frac{1}{n_{1}}\sum_{i\in I_{1}}[\alpha_{n}(X_{i})-\alpha_{0}(X_{i})]^{2}=\mathbb{E}[\alpha_{n}(X)-\alpha_{0}(X)]^{2}\lesssim p^{-2\xi_{2}}=o(1)$
and $\mathbb{E}\frac{1}{n_{1}}\sum_{i\in I_{1}}\alpha_{0}(X_{i})^{2}=\mathbb{E}\alpha_{0}(X)^{2}=O(1)$.
The above two displays imply $T_{2}=o_{P}(1)$. }

\jelenax{\textbf{Step 3:} show $T_{3}=o_{P}(1)$.}

\jelenax{Let $v_{i}=Y_{i}-\mathbb{E}(Y_{i}\mid X_{i})=Y_{i}-\rho_{0}(X_{i})$. Notice
that 
\begin{align}
\mathbb{E}\left(T_{3}^{2}\mid\{X_{i}\}_{i\in I_{1}},\{W_{i}\}_{i\in I_{2}}\right) & =\frac{4}{n_{1}^{2}}\sum_{i\in I_{1}}\hat{\alpha}_{2}(X_{i})^{4}(\rho_{0}(X_{i})-\rho_{n}(X_{i}))^{2}\mathbb{E}\left[v_{i}^{2}\mid\{X_{i}\}_{i\in I_{1}}\right]\nonumber \\
 & \overset{\text{(i)}}{\lesssim}\frac{4}{n_{1}^{2}}\sum_{i\in I_{1}}\hat{\alpha}_{2}(X_{i})^{4}(\rho_{0}(X_{i})-\rho_{n}(X_{i}))^{2},\label{eq: thm asy norm eq 11}
\end{align}
where (i) follows by Assumption \ref{assp 6}. Notice that 
\[
\mathbb{E}\left(n_{1}^{-1}\sum_{i\in I_{1}}\hat{\alpha}_{2}(X_{i})^{2}\mid\{W_{i}\}_{i\in I_{2}}\right)=\hat{\pi}_{2}^{\top}\Sigma\hat{\pi}_{2}\overset{\text{(i)}}{=}O_{P}(1),
\]
where (i) follows by $\|\hat{\pi}_{2}-\pi\|_{2}=\|\Delta_{\pi,2}\|_{2}=o_{P}(1)$
(shown in Step 2) and $\|\pi\|_{2}=O(1)$. Thus, $n_{1}^{-1}\sum_{i\in I_{1}}\hat{\alpha}_{2}(X_{i})^{2}=O_{P}(1)$,
which means that $\max_{i\in I_{1}}\hat{\alpha}_{2}(X_{i})^{2}\leq\sum_{i\in I_{1}}\hat{\alpha}_{2}(X_{i})^{2}=O_{P}(n)$.
Hence, 
\begin{align}
 & \mathbb{E}\left(T_{3}^{2}\mid\{X_{i}\}_{i\in I_{1}},\{W_{i}\}_{i\in I_{2}}\right)\nonumber \\
 & \lesssim\frac{4}{n_{1}^{2}}\left(\sum_{i\in I_{1}}(\rho_{0}(X_{i})-\rho_{n}(X_{i}))^{2}\right)\cdot\left(\max_{i\in I_{1}}\hat{\alpha}_{2}(X_{i})^{4}\right)\nonumber \\
 & =\frac{4}{n_{1}^{2}}\left(\sum_{i\in I_{1}}(\rho_{0}(X_{i})-\rho_{n}(X_{i}))^{2}\right)\cdot O_{P}(n^{2})\nonumber \\
 & =O_{P}(1)\cdot\sum_{i\in I_{1}}(\rho_{0}(X_{i})-\rho_{n}(X_{i}))^{2}\overset{\text{(i)}}{=}O_{P}(np^{-2\xi_{1}})\overset{\text{(ii)}}{=}o_{P}(1),\label{eq: thm asy nor 12}
\end{align}
where (i) follows by $\mathbb{E}\sum_{i\in I_{1}}(\rho_{0}(X_{i})-\rho_{n}(X_{i}))^{2}=n_{1}\mathbb{E}(\rho_{0}(X)-\rho_{n}(X))^{2}\lesssim np^{-2\xi_{1}}$
and (ii) follows by $p\gg n^{1/(2\xi_{1})}$. Therefore, $T_{3}=o_{P}(1)$. }

\jelenax{\textbf{Step 4:} show $T_{4}=o_{P}(1)$.}

\jelenax{Similar to (\ref{eq: thm asy norm eq 11}), we notice that
\begin{equation}
\mathbb{E}\left(T_{4}^{2}\mid\{X_{i}\}_{i\in I_{1}},\{W_{i}\}_{i\in I_{2}}\right)\lesssim\frac{4}{n_{1}^{2}}\sum_{i\in I_{1}}\hat{\alpha}_{2}(X_{i})^{4}(\rho_{n}(X_{i})-\hat{\rho}_{2}(X_{i}))^{2}.\label{eq: thm asy nor eq 13}
\end{equation}}

\jelenax{Suppose that Assumption \ref{assp 9}(i) holds. Notice that $\rho_{n}(x)-\hat{\rho}_{2}(x)=\psi(x)^{\top}(\gamma-\hat{\gamma}_{2})$.
As argued in the proof of Lemma \ref{lem: B11}, $\|\hat{\gamma}_{2}-\gamma\|_{1}=O_{P}(\varepsilon_{n}^{(2\xi_{1}-1)/(2\xi_{1}+1)})$.
Since $\xi_{1}>1/2$, we have $\|\hat{\gamma}_{2}-\gamma\|_{1}=o_{P}(1)$.
Due to the boundedness of $\|\psi(x)\|_{\infty}$, we have that $\max_{i\in I_{1}}|\rho_{n}(X_{i})-\hat{\rho}_{2}(X_{i})|\leq\|\hat{\gamma}_{2}-\gamma\|_{1}\cdot\max_{i\in I_{1}}\|\psi(X_{i})\|_{\infty}=o_{P}(1)$.
Thus, by (\ref{eq: thm asy nor eq 13}), we have 
\[
\mathbb{E}\left(T_{4}^{2}\mid\{X_{i}\}_{i\in I_{1}},\{W_{i}\}_{i\in I_{2}}\right)=o_{P}(1)\cdot\frac{4}{n_{1}^{2}}\sum_{i\in I_{1}}\hat{\alpha}_{2}(X_{i})^{4}\overset{\text{(i)}}{\lesssim}o_{P}(1)\cdot\frac{4}{n_{1}^{2}}\cdot O_{P}(n^{2})=o_{P}(1),
\]
where (i) follows by $\sum_{i\in I_{1}}\hat{\alpha}_{2}(X_{i})^{4}\leq\left(\sum_{i\in I_{1}}\hat{\alpha}_{2}(X_{i})^{2}\right)$
and $\sum_{i\in I_{1}}\hat{\alpha}_{2}(X_{i})^{2}=O_{P}(n)$ (shown
in Step 3). Thus, $T_{4}=o_{P}(1)$ under Assumption \ref{assp 9}. }

\jelenax{Suppose that Assumption \ref{assp 9}(ii) holds.  Similar to (\ref{eq: thm asy nor eq 13}),
we have 
\begin{eqnarray}
T_{4}^{2}&=&O_{P}(1)\cdot\frac{4}{n_{1}^{2}}\sum_{i\in I_{1}}\hat{\alpha}_{2}(X_{i})^{4}(\rho_{n}(X_{i})-\hat{\rho}_{2}(X_{i}))^{2}
\\
&\overset{\text{(i)}}{=}&O_{P}(1)\cdot\frac{1}{n_{1}}\sum_{i\in I_{1}}\hat{\alpha}_{2}(X_{i})^{2}(\rho_{n}(X_{i})-\hat{\rho}_{2}(X_{i}))^{2},\label{eq: thm asy norm eq 15} \nonumber 
\end{eqnarray}
where (i) follows by $\max_{i\in I_{1}}\hat{\alpha}_{2}(X_{i})^{2}=O_{P}(n)$
(shown in Step 3). }

\jelenax{Moreover, by 
$$\hat{\rho}_{2}(x)-\rho_{n}(x)=\psi(x)^{\top}\Delta_{\gamma,2},$$
we have $\frac{1}{n_{1}}\sum_{i\in I_{1}}(\rho_{n}(X_{i})-\hat{\rho}_{2}(X_{i}))^{2}=\Delta_{\gamma,2}^{\top}\hat{\Sigma}_{2}\Delta_{\gamma_{2}}$.
Since 
$$\mathbb{E}[\Delta_{\gamma,2}^{\top}\hat{\Sigma}_{2}\Delta_{\gamma_{2}}\mid\{W_{i}\}_{i\in I_{2}}]=\Delta_{\gamma,2}^{\top}\Sigma\Delta_{\gamma,2}\lesssim\|\Delta_{\gamma,2}\|_{2}^{2}$$
and $\|\Delta_{\gamma,2}\|_{2}=O_{P}(\varepsilon_{n}^{2\xi_{1}/(2\xi_{1}+1)})$
(shown in Step 1), we have that $\frac{1}{n_{1}}\sum_{i\in I_{1}}(\rho_{n}(X_{i})-\hat{\rho}_{2}(X_{i}))^{2}=O_{P}(\varepsilon_{n}^{4\xi_{1}/(2\xi_{1}+1)})$.
Similarly, $\frac{1}{n_{1}}\sum_{i\in I_{1}}[\hat{\alpha}_{2}(X_{i})-\alpha_{n}(X_{i})]^{2}=O_{P}(\|\Delta_{\pi,2}\|_{2}^{2})$.
By $\|\Delta_{\pi,2}\|_{2}=O_{P}(\varepsilon_{n}^{2\xi_{2}/(2\xi_{2}+1)})$
(shown in Step 2), we have $\frac{1}{n_{1}}\sum_{i\in I_{1}}[\hat{\alpha}_{2}(X_{i})-\alpha_{n}(X_{i})]^{2}=O_{P}(\varepsilon_{n}^{4\xi_{2}/(2\xi_{2}+1)})$.
Hence, 
\begin{align}
 & \frac{1}{n_{1}}\sum_{i\in I_{1}}(\hat{\alpha}_{2}(X_{i})-\alpha_{n}(X_{i}))^{2}(\rho_{n}(X_{i})-\hat{\rho}_{2}(X_{i}))^{2}\nonumber \\
 & \leq\frac{1}{n_{1}}\left(\sum_{i\in I_{1}}|\hat{\alpha}_{2}(X_{i})-\alpha_{n}(X_{i})|\cdot|\rho_{n}(X_{i})-\hat{\rho}_{2}(X_{i})|\right)^{2}\nonumber \\
 & \overset{\text{(i)}}{\leq}\frac{1}{n_{1}}\left(\sum_{i\in I_{1}}|\hat{\alpha}_{2}(X_{i})-\alpha_{n}(X_{i})|^{2}\right)\cdot\left(\sum_{i\in I_{1}}|\rho_{n}(X_{i})-\hat{\rho}_{2}(X_{i})|^{2}\right)\nonumber \\
 & =O(n^{-1})\cdot O_{P}\left(n^{2}\varepsilon_{n}^{4\xi_{1}/(2\xi_{1}+1)+4\xi_{2}/(2\xi_{2}+1)}\right)\overset{\text{(ii)}}{=}o_{P}(1),\label{eq: thm asy norm eq 16}
\end{align}
where (i) follows by Cauchy-Schwarz inequality and (ii) follows by
$n^{\xi_{1}\xi_{2}-1/4}\gg(\log p)^{2\xi_{1}\xi_{2}+(\xi_{1}+\xi_{2})/2}$
and $\xi_{1}\xi_{2}>1/4$ (due to Assumption \ref{assp 9}(ii)). We also notice that
\begin{align}
 & \frac{1}{n_{1}}\sum_{i\in I_{1}}(\alpha_{n}(X_{i})-\alpha_{0}(X))^{2}(\rho_{n}(X_{i})-\hat{\rho}_{2}(X_{i}))^{2}\nonumber \\
 & \leq\left(\max_{i\in I_{1}}(\rho_{n}(X_{i})-\hat{\rho}_{2}(X_{i}))^{2}\right)\cdot\frac{1}{n_{1}}\sum_{i\in I_{1}}(\alpha_{n}(X_{i})-\alpha_{0}(X))^{2}\nonumber \\
 & \overset{\text{(i)}}{\leq}O_{P}(n\varepsilon_{n}^{4\xi_{1}/(2\xi_{1}+1)})\cdot O_{P}(p^{-2\xi_{2}})\leq O_{P}(n)\cdot O_{P}(p^{-2\xi_{2}})\overset{\text{(ii)}}{=}o_{P}(1),\label{eq: thm asy norm eq 17}
\end{align}
where (i) follows by $\max_{i\in I_{1}}(\rho_{n}(X_{i})-\hat{\rho}_{2}(X_{i}))^{2}\leq\sum_{i\in I_{1}}(\rho_{n}(X_{i})-\hat{\rho}_{2}(X_{i}))^{2}=O_{P}(n\varepsilon_{n}^{4\xi_{1}/(2\xi_{1}+1)})$
(shown above) and $\mathbb{E}\frac{1}{n_{1}}\sum_{i\in I_{1}}(\alpha_{n}(X_{i})-\alpha_{0}(X))^{2}=\mathbb{E}(\alpha_{n}(X)-\alpha_{0}(X))^{2}=O(p^{-2\xi_{2}})$
(due to $\alpha_{0}\in\mathcal{M}_{C,\xi_{2}}$) and (ii) follows
by $p\gg n^{1/(2\xi_{2})}$ (due to Assumption \ref{assp 9}(ii)). }

\jelenax{Since $\hat{\rho}_{2}(x)-\rho_{n}(x)=\psi(x)^{\top}\Delta_{\gamma,2}$ and
$\Delta_{\gamma,2}$ is independent of $\{X_{i}\}_{i\in I_{1}}$,
we have 
\begin{align*}
 & \mathbb{E}\left[\frac{1}{n_{1}}\sum_{i\in I_{1}}\alpha_{0}(X)^{2}(\rho_{n}(X_{i})-\hat{\rho}_{2}(X_{i}))^{2}\mid\{W_{i}\}_{i\in I_{2}}\right]\\
 & =\Delta_{\gamma,2}\mathbb{E}[\alpha_{0}(X)^{2}\psi(X)\psi(X)^{\top}]\Delta_{\gamma,2}\overset{\text{(i)}}{\lesssim}\|\Delta_{\gamma,2}\|_{2}^{2}=o_{P}(1),
\end{align*}
where (i) follows by the boundedness of the eigenvalues of $\mathbb{E}[\alpha_{0}(X)^{2}\psi(X)\psi(X)^{\top}]$
(due to Assumption \ref{assp 6}). Thus, $\frac{1}{n_{1}}\sum_{i\in I_{1}}\alpha_{0}(X)^{2}(\rho_{n}(X_{i})-\hat{\rho}_{2}(X_{i}))^{2}=o_{P}(1)$.
This, (\ref{eq: thm asy norm eq 16}) and (\ref{eq: thm asy norm eq 17}),
displays imply 
\begin{equation}
\frac{1}{n_{1}}\sum_{i\in I_{1}}\hat{\alpha}_{2}(X_{i})^{2}(\rho_{n}(X_{i})-\hat{\rho}_{2}(X_{i}))^{2}=o_{P}(1).\label{eq: thm asy norm eq 17.5}
\end{equation}}

\jelenax{By (\ref{eq: thm asy norm eq 15}), $T_{4}^{2}=o_{P}(1)$ under Assumption \ref{assp 9}(ii).}

\jelenax{\textbf{Step 5:} show $T_{5}=o_{P}(1)$.}

\jelenax{Notice that 
\begin{align*}
 & \frac{1}{n_{1}}\sum_{i\in I_{1}}\hat{\alpha}_{2}(X_{i})^{2}(\rho_{0}(X_{i})-\rho_{n}(X_{i}))^{2}\\
 & \leq\left(\max_{i\in I_{1}}\hat{\alpha}_{2}(X_{i})^{2}\right)\cdot\left(\frac{1}{n_{1}}\sum_{i\in I_{1}}(\rho_{0}(X_{i})-\rho_{n}(X_{i}))^{2}\right)\\
 & \overset{\text{(i)}}{\leq}O_{P}(n)\cdot\left(\frac{1}{n_{1}}\sum_{i\in I_{1}}(\rho_{0}(X_{i})-\rho_{n}(X_{i}))^{2}\right)\overset{\text{(ii)}}{=}O_{P}(n)\cdot O_{P}(p^{-2\xi_{1}})\overset{\text{(iii)}}{=}o_{P}(1),
\end{align*}
where (i) follows by $\max_{i\in I_{1}}\hat{\alpha}_{2}(X_{i})^{2}=O_{P}(n)$
(shown in Step 3), (ii) follows by $\mathbb{E}\frac{1}{n_{1}}\sum_{i\in I_{1}}(\rho_{0}(X_{i})-\rho_{n}(X_{i}))^{2}=\mathbb{E}(\rho_{0}(X)-\rho_{n}(X))^{2}=O(p^{-2\xi_{1}})$
(due to $\rho_{0}\in\mathcal{M}_{C,\xi_{1}}$) and (iii) follows by
$p\gg n^{1/(2\xi_{1})}$. Therefore, by the definition of $T_{5}$,
it remains to verify that 
\begin{equation}
\frac{1}{n_{1}}\sum_{i\in I_{1}}\hat{\alpha}_{2}(X_{i})^{2}(\rho_{n}(X_{i})-\hat{\rho}_{2}(X_{i}))^{2}=o_{P}(1).\label{eq: thm asy norm eq 18}
\end{equation}}

\jelenax{In Step 4 (after (\ref{eq: thm asy nor eq 13})),
we already proved that $\max_{i\in I_{1}}|\rho_{n}(X_{i})-\hat{\rho}_{2}(X_{i})|=o_{P}(1)$ under Assumption \ref{assp 9}(i). Since $\frac{1}{n_{1}}\sum_{i\in I_{1}}\hat{\alpha}_{2}(X_{i})^{2}=O_{P}(1)$
(shown in Step 3), we obtain (\ref{eq: thm asy norm eq 18}) under Assumption \ref{assp 9}(i).}

\jelenax{We proved (\ref{eq: thm asy norm eq 17.5})
in Step 4 under Assumption \ref{assp 9}(ii). This is the same as (\ref{eq: thm asy norm eq 18}). The
proof is complete. }
\end{proof}

\section{Some General Lemmas}


\begin{lemma}\label{lem: C1}
For any $a\in
\mathbb{R}^{p}$ such that $\left\Vert a-b_{s}\right\Vert _{2}\leq Cs^{-r}%
$ for any $s\geq0$, where $C,r>0$\ are constants and
$b_{s}=\arg\min_{\left\Vert v\right\Vert _{0}\leq s}\left\Vert a-v\right\Vert
_{2}.$ If $r>1/2$ and $s\geq2$, then $\left\Vert
a-b_{s}\right\Vert _{1}\leq Ds^{1/2-r},$ where $D>0$ is a
constant depending only on $C$ and $r.$
\end{lemma}

\begin{proof}[Proof of Lemma \ref{lem: C1}]

Without loss of generality, we assume that $|a_{1}|\geq
|a_{2}|\geq\cdots\geq|a_{p}|\geq0$. Then clearly, $a-b_{k}=(0,0,\ldots
,0,a_{k+1},a_{k+2},\ldots,a_{p})^\top$ for $k\geq0$. By assumption, we have
that for any $k\geq0$,
\begin{equation}
\sum_{j=k+1}^{p}a_{j}^{2}\leq C^{2}k^{-2r}.
\label{eq: L1 bnd sparse apprx error 4}%
\end{equation}
Let $g\in\mathbb{N}$ be defined as $2^{g}<p/s\leq2^{g+1}$. With a slight abuse
of notation, we extend $a$ to be a $2^{g+1}s$-dimensional vector with
$a_{j}=0$ for $j>p$. Then we have that%
\begin{eqnarray*}
\sum_{j=s+1}^{p}a_{j}&=&\sum_{m=0}^{g}\sum_{j=2^{m}s+1}^{2^{m+1}s}a_{j}\leq
\sum_{m=0}^{g}\sqrt{2^{m}s\sum_{j=2^{m}s+1}^{2^{m+1}s}a_{j}^{2}}\\
&\leq&\sum
_{m=0}^{g}\sqrt{2^{m}s\sum_{j=2^{m}s+1}^{p}a_{j}^{2}}\\
&\overset{\text{(i)}}{\leq}& \sum_{m=0}^{g}\sqrt{2^{m}sC^{2}(2^{m}s)^{-2r}%
}=Cs^{1/2-r}\sum_{m=0}^{g}\left(  2^{1/2-r}\right)  ^{m}
\\
&<&Cs^{1/2-r}\sum
_{m=0}^{\infty}\left(  2^{1/2-r}\right)  ^{m}\overset{\text{(ii)}}{=}C\frac
{1}{1-2^{1/2-r}}s^{1/2-r}\text{,}%
\end{eqnarray*}
where (i) follows by (\ref{eq: L1 bnd sparse apprx error 4}) applied to $k=2^{m}s$ and (ii) follows by the fact
that $2^{1/2-r}<1$ (since $r>1/2$). 
\end{proof}

Let $X_{i}=(X_{i,1},...,X_{i,p})^{\top}\in{\mathbb{R}}^{p}$ and $X_{i,0}$ be
a scalar, with all random variables allowed to depend on $n.$


\begin{lemma}\label{lem: C2} 
\jelenax{Let $\{(X_{i},X_{i,0})\}_{i=1}^{n}$
be independent observations. Assume that there exist constants $\kappa,C_{1},C_{2}>0$
such that the sub-Gaussian norm of $X_{i,j}$ is bounded by $C_{1}$
for any $j\in\{1,...,p\}$ and $\\E \left[\,|X_{i,0}|^{2+\kappa}\,\right]\leq C_{2}$.
Suppose that $p\lesssim\exp(n^{\zeta})$ for some $\zeta<\kappa/(2\kappa+4)$.
Let $D_{i}=X_{i}X_{i,0}-\mathbb{E}[X_{i}X_{i,0}]$ and $\bar{D}=\sum_{i=1}^{n}D_{i}/n.$
Then, there exists a constant $\bar{C}>0$ (depending only on $C_{1}$,
$C_{2}$, $\zeta$, and $\kappa$) such that
\[
P\left(\|\bar{D}\|_{\infty}>\bar{C}\sqrt{\log(p)/n}\right)=o(1).
\]}
\end{lemma}
\begin{proof}[Proof of Lemma \ref{lem: C2}]

\jelenax{We prove this result using symmetrization. Let $\varepsilon_{1},...,\varepsilon_{n}$
be i.i.d Rademacher random variables independent of $X_{i}$ for for
all observations, i.e., $P(\varepsilon_{i}=1)=P(\varepsilon_{i}=-1)=1/2$.
Define the symmetrized quantity $W_{*,j}=\sum_{i=1}^{n}D_{i,j}\varepsilon_{i}$.}

\jelenax{Since the sub-Gaussian norm of $X_{i,j}$ is bounded by a constant
uniformly in $j$, $\max_{1\leq j\leq p}\mathbb{E}|X_{i,j}|^{2+4/\kappa}$
is also bounded by a constant. Then
\begin{multline*}
\max_{1\leq j\leq p}\mathbb{E}[D_{i,j}^{2}]\leq\max_{1\leq j\leq p}\mathbb{E}[X_{i,j}^{2}X_{i,0}^{2}]\leq\max_{1\leq j\leq p}\left(\mathbb{E}|X_{i,j}|^{2\cdot(2\kappa^{-1}+1)}\right)^{\kappa/(2+\kappa)}\cdot\left(\mathbb{E}|X_{i,0}|^{2\cdot(1+\kappa/2)}\right)^{2/(2+\kappa)}\\
=\max_{1\leq j\leq p}\left(\mathbb{E}|X_{i,j}|^{4\kappa^{-1}+2}\right)^{\kappa/(2+\kappa)}\cdot\left(\mathbb{E}|X_{i,0}|^{2+\kappa/}\right)^{2/(2+\kappa)}\leq C.
\end{multline*}}

\jelenax{Since $\varepsilon_{i}$ is a Rademacher variable, by Hoeffding's
lemma, for any $t\in{\mathbb{R}}$, $\mathbb{E}\operatorname{exp}(t\varepsilon_{i})\leq\operatorname{exp}(t^{2}/2)$.
Since $\{\varepsilon_{i}\}_{i=1}^{n}$ is independent of $X$ we have
\begin{multline*}
\mathbb{E}\left[\operatorname{exp}(tW_{\ast,j})|X\right]=\mathbb{E}\left[\prod_{i=1}^{n}\operatorname{exp}[tD_{i,j}\varepsilon_{i}]\mid X\right]=\prod_{i=1}^{n}\mathbb{E}\left[\operatorname{exp}[tD_{i,j}\varepsilon_{i}]\mid X\right]\\
\leq\prod_{i=1}^{n}\operatorname{exp}\left(t^{2}D_{i,j}^{2}/2\right)=\operatorname{exp}\left(t^{2}\sum_{i=1}^{n}D_{i,j}^{2}/2\right).
\end{multline*}
Similarly, apply the same argument to $-W_{\ast,j}$ to obtain $\mathbb{E}\left[\operatorname{exp}(-tW_{\ast,j})\mid X\right]\leq\operatorname{exp}\left(t^{2}\sum_{i=1}^{n}D_{i,j}^{2}/2\right).$
Since $\operatorname{exp}(t|W_{\ast,j}|)\leq\operatorname{exp}(tW_{\ast,j})+\operatorname{exp}(-tW_{\ast,j})$,
we have 
\[
\mathbb{E}\left[\operatorname{exp}(t|W_{\ast,j}|)\mid X\right]\leq2\operatorname{exp}\left(t^{2}\sum_{i=1}^{n}D_{i,j}^{2}/2\right).
\]}

\jelenax{Next let $z>0$ be a non-random quantity to be chosen later and $\Vert W_{\ast}\Vert_{\infty}=\max_{1\leq j\leq p}|W_{\ast,j}|$.
By Lemma 2.3.7 of \citet{van1996weak}, we have 
\begin{equation}
\left(1-\beta_{n}(z)\right)P\left(\max_{1\leq j\leq p}\left\vert \sum_{i=1}^{n}D_{i,j}\right\vert >z\right)\leq2P\left(\Vert W_{\ast}\Vert_{\infty}>z/4\right),\label{VVW}
\end{equation}
where 
\[
\beta_{n}(z)=1-4z^{-2}n\max_{1\leq j\leq p}\mathbb{E}[D_{i,j}^{2}]\geq1-Cz^{-2}n.
\]
The rest of the proof proceeds in two steps. We first bound $\max_{1\leq j\leq p}\sum_{i=1}^{n}D_{i,j}^{2}$
and then derive the final result.}

\jelenax{\textbf{Step 1:} bound $\max_{1\leq j\leq p}\sum_{i=1}^{n}D_{i,j}^{2}$.}

\jelenax{Notice that $D_{i,j}^{2}\leq2X_{i,j}^{2}X_{i,0}^{2}+2(\mathbb{E}[X_{i,j}X_{i,0}])^{2}$.
Since $|\mathbb{E}[X_{i,j}X_{i,0}]|\leq\sqrt{\mathbb{E}[X_{i,j}^{2}X_{i,0}^{2}]}$,
which has been shown to be bounded, it remains to bound $\max_{1\leq j\leq p}\sum_{i=1}^{n}X_{i,j}^{2}X_{i,0}^{2}$.
Notice that 
\begin{multline}
\max_{1\leq j\leq p}n^{-1}\sum_{i=1}^{n}X_{i,j}^{2}X_{i,0}^{2}\leq\max_{1\leq j\leq p}\left(n^{-1}\sum_{i=1}^{n}(X_{i,j}^{2})^{(\kappa+4)/\kappa}\right)^{\frac{\kappa}{4+\kappa}}\left(n^{-1}\sum_{i=1}^{n}(X_{i,0}^{2})^{1+\kappa/4}\right)^{\frac{4}{4+\kappa}}\\
=\left(\max_{1\leq j\leq p}n^{-1}\sum_{i=1}^{n}|X_{i,j}|^{2+8/\kappa}\right)^{\frac{\kappa}{4+\kappa}}\left(n^{-1}\sum_{i=1}^{n}|X_{i,0}|^{2+\kappa/2}\right)^{\frac{4}{4+\kappa}}.\label{eq: lem F2 eq 1}
\end{multline}}

\jelenax{Let 
$$S_{n}=\sum_{i=1}^{n}[|X_{i,0}|^{2+\kappa/2}-\mathbb{E}|X_{i,0}|^{2+\kappa/2}].$$
Let $\alpha=(2+\kappa)/(2+\kappa/2)=2(\kappa+2)/(\kappa+4)>1$. Notice
that $\mathbb{E}\left||X_{i,0}|^{2+\kappa/2}-\mathbb{E}|X_{i,0}|^{2+\kappa/2}\right|^{\alpha}$
is bounded since $\mathbb{E}|X_{i,0}|^{(2+\kappa/2)\alpha}=\mathbb{E}|X_{i,0}|^{2+\kappa}$
is bounded. By Rosenthal's inequality (Theorem 9.1 in Chapter 3.9
of \citet{gut2013probability}), 
$$\mathbb{E}|S_{n}|^{\alpha}\leq C_{1}(n+n^{\alpha/2})$$
for some constant $C_{1}>0$. We choose $C_{1}$ large enough such
that $\mathbb{E}|X_{i,0}|^{2+\kappa/2}\leq C_{1}$. Since $\alpha>1$,
it follows that 
\begin{equation}
P\left(n^{-1}\sum_{i=1}^{n}|X_{i,0}|^{2+\kappa/2}>2C_{1}\right)\leq P\left(n^{-1}S_{n}>C_{1}\right)\leq\frac{\mathbb{E}|S_{n}|^{\alpha}}{(nC_{1})^{\alpha}}\leq\frac{C_{1}(n+n^{\alpha/2})}{(nC_{1})^{\alpha}}=o(1).\label{eq: lem F2 eq 2}
\end{equation}}

\jelenax{On the other hand, let $\beta=2+8/\kappa$ and choose $C_{2}>0$ to
be a bound on the sub-Gaussian norm of $X_{i,j}$. Then for any $t\geq1$,
\begin{equation}
\mathbb{E}|X_{i,j}|^{t}\leq C_{3}^{t}t^{t/2},\label{eq: lem F2 eq 3}
\end{equation}
where $C_{3}>0$ is a constant depending only on $C_{2}$. Let $Z_{i,j}=|X_{i,j}|^{\beta}-\mathbb{E}|X_{i,j}|^{\beta}$.
Then for $m\geq1$, 
\begin{multline*}
\mathbb{E}|Z_{i,j}|^{m}\overset{\text{(i)}}{\leq}2^{m-1}\left[\mathbb{E}|X_{i,j}|^{\beta m}+\left(\mathbb{E}|X_{i,j}|^{\beta}\right)^{m}\right]\\
\leq2^{m-1}\left[\mathbb{E}|X_{i,j}|^{\beta m}+\mathbb{E}|X_{i,j}|^{\beta m}\right]=2^{m}\mathbb{E}|X_{i,j}|^{\beta m}\overset{\text{(ii)}}{\leq}2^{m}C_{3}^{\beta m}(\beta m)^{\beta m/2},
\end{multline*}
where (i) follows by the elementary inequality that $|a+b|^{m}\leq2^{m-1}(|a|^{m}+|b|^{m})$
for any $a,b\in\RR$ and $m\geq1$ and (ii) follows by (\ref{eq: lem F2 eq 3}).}

\jelenax{We apply Theorem 2 of \citet{nagaev1978some} (with $c=t=r$ in their
theorem). Then for any $r\geq2$, 
\begin{multline*}
\mathbb{E}\left|\sum_{i=1}^{n}Z_{i,j}\right|^{r}\leq r^{r}\left(\sum_{i=1}^{n}\mathbb{E}|Z_{i,j}|^{r}\right)+r^{r/2+1}\exp(r)\left(\sum_{i=1}^{n}\mathbb{E}Z_{i,j}^{2}\right)^{r/2}\int_{0}^{1}x^{r/2-1}(1-x)^{r/2-1}dx\\
\leq r^{r}\left(\sum_{i=1}^{n}\mathbb{E}|Z_{i,j}|^{r}\right)+r^{r/2+1}\exp(r)\left(\sum_{i=1}^{n}\mathbb{E}Z_{i,j}^{2}\right)^{r/2}.
\end{multline*}
}
\jelenax{
We will choose $r$ later. By (\ref{eq: lem F2 eq 3}), we also have
\[
\sum_{i=1}^{n}\mathbb{E}Z_{i,j}^{2}\leq n2^{2}C_{3}^{2\beta}(2\beta)^{\beta}=C_{4}^{2}n
\]
and 
\[
\sum_{i=1}^{n}\mathbb{E}|Z_{i,j}|^{r}\leq2^{r}C_{3}^{\beta r}(\beta r)^{\beta r/2}n=C_{5}r^{\beta r/2}n,
\]
where $C_{4}=\sqrt{2^{2}C_{3}^{2\beta}(2\beta)^{\beta}}$ and $C_{5}=2C_{3}^{\beta}\beta^{\beta/2}$.
The above three displays imply 
\[
\mathbb{E}\left|\sum_{i=1}^{n}Z_{i,j}\right|^{r}\leq C_{5}r^{(\beta/2+1)r}n+r^{r/2+1}\exp(r)C_{4}^{r}n^{r/2}.
\]}

\jelenax{Since $p\lesssim\exp(n^{\zeta})$, we have $\log p\leq C_{6}+n^{\zeta}$
for some $C_{6}>0$. Therefore, by the union bound, we have 
\begin{align*}
P\left(n^{-1}\max_{1\leq j\leq p}\left|\sum_{i=1}^{n}Z_{i,j}\right|\geq1\right) & =P\left(\max_{1\leq j\leq p}\left|\sum_{i=1}^{n}Z_{i,j}\right|^{r}\geq n^{r}\right)\\
 & \leq\frac{\mathbb{E}\left(\max_{1\leq j\leq p}\left|\sum_{i=1}^{n}Z_{i,j}\right|^{r}\right)}{n^{r}}\\
 & \leq n^{-r}\sum_{j=1}^{p}\mathbb{E}\left|\sum_{i=1}^{n}Z_{i,j}\right|^{r}\\
 & \leq n^{-r}p\left[C_{5}r^{(\beta/2+1)r}n+r^{r/2+1}\exp(r)C_{4}^{r}n^{r/2}\right]\\
 & =C_{5}\exp\left\{ \log p+(\beta/2+1)r\log r-(r-1)\log n\right\} \\
 & \quad+\exp\left\{ \log p+(r/2+1)\log r+r(1+\log C_{4})-(r/2)\log n\right\} \\
 & \leq C_{5}\exp\left\{ C_{6}+n^{\zeta}+(\beta/2+1)r\log r-(r-1)\log n\right\} \\
 & \quad+\exp\left\{ C_{6}+n^{\zeta}+(r/2+1)\log r+r(1+\log C_{4})-(r/2)\log n\right\} 
\end{align*}}

\jelenax{Let $\theta=\beta/2+1=(2\kappa+4)/\kappa$. By the assumption of $\zeta<\kappa/(2\kappa+4)$,
we have $\zeta<1/\theta$. Now choose $\sigma\in(\zeta,1/\theta)$
and set $r=n^{\sigma}$. Notice that $\sigma<1$. Then 
\[
n^{\zeta}+(\beta/2+1)r\log r-(r-1)\log n=n^{\zeta}+\sigma\theta n^{\sigma}\log n-(n^{\sigma}-1)\log n\rightarrow-\infty
\]
and 
\begin{multline*}
n^{\zeta}+(r/2+1)\log r+r(1+\log C_{4})-(r/2)\log n\\
=n^{\zeta}+\sigma(n^{\sigma}/2+1)\log n+n^{\sigma}(1+\log C_{4})-(n^{\sigma}/2)\log n\rightarrow-\infty.
\end{multline*}}

\jelenax{The above three displays imply $P\left(n^{-1}\max_{1\leq j\leq p}\left|\sum_{i=1}^{n}Z_{i,j}\right|\geq1\right)=o(1)$.
By (\ref{eq: lem F2 eq 3}), we have $\max_{1\leq j\leq p}n^{-1}\sum_{i=1}^{n}\mathbb{E}|X_{i,j}|^{\beta}\leq C_{3}$.
Then the above implies 
\begin{multline*}
P\left(n^{-1}\max_{1\leq j\leq p}\sum_{i=1}^{n}|X_{i,j}|^{\beta}\geq1+C_{3}\right)\\
\leq P\left(n^{-1}\max_{1\leq j\leq p}\left|\sum_{i=1}^{n}Z_{i,j}\right|+\max_{1\leq j\leq p}n^{-1}\sum_{i=1}^{n}\mathbb{E}|X_{i,j}|^{\beta}\geq1+C_{3}\right)=o(1).
\end{multline*}}

\jelenax{We now combine the above with (\ref{eq: lem F2 eq 1}) and (\ref{eq: lem F2 eq 2}),
obtaining 
\[
P\left(\max_{1\leq j\leq p}n^{-1}\sum_{i=1}^{n}D_{i,j}^{2}>2\tilde{M}\right)=P\left(\max_{1\leq j\leq p}n^{-1}\sum_{i=1}^{n}X_{i,j}^{2}X_{i,0}^{2}>2\tilde{M}\right)=o(1)
\]
for some constant $\tilde{M}>0$. Therefore, $P({\mathcal{A})\rightarrow}0$,
where ${\mathcal{A}}=\left\{ \max_{1\leq j\leq p}n^{-1}\sum_{i=1}^{n}D_{i,j}^{2}>2\tilde{M}\right\} $.}

\jelenax{\textbf{Step 2:} derive the final result.}

\jelenax{On the event ${\mathcal{A}}^{c}$ we have that 
\begin{multline}
\mathbb{E}\left[\operatorname{exp}\left(t\Vert W_{\ast}\Vert_{\infty}\right)\mid X\right]=\mathbb{E}\left[\operatorname{exp}\left(t\max_{1\leq j\leq p}\left\vert W_{\ast,j}\right\vert \right)\mid X\right]\leq\sum_{j=1}^{p}\mathbb{E}\left(\operatorname{exp}(t|W_{\ast,j}|)\mid X\right)\\
\leq2\sum_{j=1}^{p}\operatorname{exp}\left(t^{2}\sum_{i=1}^{n}D_{i,j}^{2}/2\right)\leq2p\operatorname{exp}\left(\tilde{M}t^{2}n\right),\label{eq:maineq5}
\end{multline}
Let $t>0$ be a non-random quantity to be chosen. By the Markov inequality
we have 
\begin{align*}
P\left(\Vert W_{\ast}\Vert_{\infty}>z/4\mid X\right){\mathbf{1}}_{{\mathcal{A}}^{c}} & \leq P\left(\operatorname{exp}(t\Vert W_{\ast}\Vert_{\infty})>\operatorname{exp}(tz/4)\mid X\right){\mathbf{1}}_{{\mathcal{A}}^{c}}\\
 & \leq\operatorname{exp}(-tz/4)\mathbb{E}\left[\operatorname{exp}(t\Vert W_{\ast}\Vert_{\infty})\mid X\right]\cdot{\mathbf{1}}_{{\mathcal{A}}^{c}}\\
 & \leq\operatorname{exp}(-tz/4)\cdot2p\operatorname{exp}\left(\tilde{M}t^{2}n\right)=\operatorname{exp}\left(-\frac{1}{4}tz+\log(2p)+\tilde{M}_{n}t^{2}n\right).
\end{align*}
Now choose $t=z/[8n\tilde{M}]$ to obtain 
\[
P\left(\Vert W_{\ast}\Vert_{\infty}>z/4\mid X\right){\mathbf{1}}_{{\mathcal{A}}^{c}}\leq\operatorname{exp}\left(-z^{2}/(nK)+\log(2p)\right),
\]
where $K=64\tilde{M}$. Therefore we have
\begin{align*}
P\left(\Vert W_{\ast}\Vert_{\infty}>z/4\mid X\right) & =P\left(\Vert W_{\ast}\Vert_{\infty}>z/4\mid X\right)({\mathbf{1}}_{{\mathcal{A}}^{c}}+{\mathbf{1}}_{{\mathcal{A}}})\\
 & \leq\operatorname{exp}\left(-z^{2}/(nK)+\log(2p)\right)+{\mathbf{1}}_{{\mathcal{A}}}.
\end{align*}
Then by iterated expectations and $\mathbb{E}[{\mathbf{1}}_{{\mathcal{A}}}]=P({\mathcal{A})\rightarrow}0,$
\[
P\left(\Vert W_{\ast}\Vert_{\infty}>z/4\right)\leq\operatorname{exp}\left(-z^{2}/(nK)+\log(2p)\right)+o(1).
\]}

\jelenax{Finally choose $z=2\sqrt{Kn\log(2p)}$ and note that 
\[
\beta_{n}(z):=1-Cz^{-2}n=1-\frac{C}{4K\log(2p)}\rightarrow1.
\]
Define $B_{n}:=\sqrt{n/\log(p)}\left\Vert \bar{D}\right\Vert _{\infty}.$
Then by (\ref{VVW}) we have for large enough $n$ 
\begin{align*}
P(B_{n} & \geq2\sqrt{K\frac{\log(2p)}{\log(p)}})=P(\left\Vert n\bar{D}\right\Vert _{\infty}>z)\leq CP\left(\Vert W_{\ast}\Vert_{\infty}>\frac{z}{4}\right)\\
 & \leq\operatorname{exp}\left(-z^{2}/(nK_{n})+\log(2p)\right)+o(1)=\operatorname{exp}\left(-4\log(2p)+\log(2p)\right)+o(1)\rightarrow0.
\end{align*}
Note that $\log(2p)/\log(p)=1+2/\log(p)$ is bounded by 2 for large
$p$. We have that $P(B_{n}>2\sqrt{2K})=o(1)$. }
\end{proof}
\begin{lemma}\label{lem: C3} 
\jelenax{Let $\{(X_{i},X_{i,0})\}_{i=1}^{n}$
be independent observations. Assume that there exist a constant $\kappa>0$
such that $X_{i,j}$ is bounded by $\kappa$ for any $(i,j)\in\{1,...,n\}\times\{1,...,p\}$.
Suppose that $\mathbb{E}X_{i,0}^{2}=o(1/\log(pn))$. Let 
$$D_{i}=X_{i}X_{i,0}-\mathbb{E}[X_{i}X_{i,0}]$$
and $\bar{D}=\sum_{i=1}^{n}D_{i}/n.$ Then, there exists a constant
$\bar{C}>0$ (depending only on $\kappa$) such that
\[
P\left(\|\bar{D}\|_{\infty}>\bar{C}\sqrt{\log(p)/n}\right)=o(1).
\]}
\end{lemma}
\begin{proof}[Proof of Lemma \ref{lem: C3}]
\jelenax{We follow the same argument and same notations as the proof of Lemma
\ref{lem: C2} and obtain
\[
\mathbb{E}\left[\operatorname{exp}(t|W_{\ast,j}|)\mid X\right]\leq2\operatorname{exp}\left(t^{2}\sum_{i=1}^{n}D_{i,j}^{2}/2\right)
\]
and 
\[
\left(1-\beta_{n}(z)\right)P\left(\max_{1\leq j\leq p}\left\vert \sum_{i=1}^{n}D_{i,j}\right\vert >z\right)\leq2P\left(\Vert W_{\ast}\Vert_{\infty}>z/4\right),
\]
where 
\[
\beta_{n}(z)=1-4z^{-2}n\max_{1\leq j\leq p}\mathbb{E}[D_{i,j}^{2}].
\]}

\jelenax{By the sub-Gaussian assumption and the union bound, $\max_{1\leq i\leq n,1\leq j\leq p}|X_{i,j}|=O_{P}(\sqrt{\log(pn)})$.
Thus, 
\[
\max_{1\leq j\leq p}n^{-1}\sum_{i=1}^{n}X_{i,j}^{2}X_{i,0}^{2}\leq\left(\max_{1\leq i\leq n,1\leq j\leq p}X_{i,j}^{2}\right)\cdot n^{-1}\sum_{i=1}^{n}X_{i,0}^{2}=O_{P}(\log(pn))\cdot o_{P}(1/\log(pn))=o_{P}(1).
\]
This completes Step 1 in the proof of Lemma \ref{lem: C2}. We also
have 
\[
\max_{1\leq j\leq p}\mathbb{E}[D_{i,j}^{2}]\leq\max_{1\leq j\leq p}\mathbb{E}[X_{i,j}^{2}X_{i,0}^{2}]\lesssim\mathbb{E}X_{i,0}^{2}=o(1).
\]
The rest of the argument follows the same in proof of Lemma \ref{lem: C2}.
}
\end{proof}

\end{appendix}

\end{document}